\documentclass[12pt]{report}

\usepackage{UOPhDthesis}
\usepackage{amsthm} 
\usepackage{amssymb}
\usepackage{amsmath}
\usepackage{epsfig}
\usepackage{makeidx}
\usepackage{rotating}
\usepackage{wasysym} 
\usepackage{pifont}

\newtheoremstyle{karendef}{\topsep}{\topsep}%
     {}%         Body font
     {}%         Indent amount (empty = no indent, \parindent =      
     {\bfseries}% Thm head font
     {.}%        Punctuation after thm head
     {.5em}%{.5em}     Space after thm head (\newline = linebreak)
     {\thmname{#1} \thmnumber{#2} (\thmnote{{\em #3}})}%         Thm head
\newtheoremstyle{marnidef}{\topsep}{\topsep}%
     {}%         Body font
     {}%         Indent amount (empty = no indent, \parindent =      
     {\bfseries}% Thm head font
     {.  }%        Punctuation after thm head
     {.5em}%     Space after thm head (\newline = linebreak)
     {\thmname{#1} \thmnumber{#2}}%         Thm head

\newcommand{\orphan}{\enlargethispage*{40pt}}
\newcommand{\halforphan}{\enlargethispage*{20pt}}

\newcommand{\remove}[1]{}

\theoremstyle{plain}
\newtheorem{thm}{Theorem}[section]
\newtheorem{lemma}[thm]{Lemma}%[section]
\newtheorem{cor}[thm]{Corollary}%[section]
\newtheorem{prop}[thm]{Proposition}%[section]

\theoremstyle{karendef}
\newtheorem{defn}[thm]{Definition}

\theoremstyle{marnidef}
\newtheorem{conj}[thm]{Conjecture}%[section]
\newtheorem{quest}[thm]{Question}

\newenvironment{exmp}{
\refstepcounter{thm}
\begin{quotation}\small\noindent{\bf Example \thethm .}}{\end{quotation}}

\newcommand{\fix}{\mathrm{fix}}
\renewcommand{\qed}{\hfill \ding{80}}

\makeindex %  \index{wordforinindex}
\makeglossary \glossary{wordinglossary}

\begin{document} 

{
\title{Covering arrays on graphs: qualitative independence graphs and extremal set partition
theory}
\author{Karen Meagher}
 \date{\today}

 % \copyrightfalse % do not produce a separate copyright page        % otherwise use \copyrighttrue
 %\figurespagefalse % do not produce a separate figures page
   %   \tablespagefalse  % do not produce a separate tables page
  \beforepreface 
  \signaturepage    

  \prefacesection{Abstract}{
     \thispagestyle{empty}
     
There has been a good deal of research on covering arrays over the
last 20 years. Most of this work has focused on constructions,
applications and generalizations of covering arrays. The main focus of
this thesis is a generalization of covering arrays, covering arrays on
graphs. The original motivation for this generalization was to improve
applications of covering arrays to testing systems and networks, but
this extension also gives us new ways to study covering arrays.

Two vectors $v,w$ in $\mathbb{Z}_k^n$ are {\em qualitatively
independent} if for all ordered pairs $(a,b) \in \mathbb{Z}_k \times
\mathbb{Z}_k$ there is a position $i$ in the vectors where $(a,b) =
(v_i,w_i)$. A {\em covering array} is an array with the property that
any pair of rows are qualitatively independent.  A {\em covering array
on a graph} is an array with a row for each vertex of the graph with
the property that any two rows which correspond to adjacent vertices
are qualitatively independent. A covering array on the complete graph
is a covering array. A covering array is {\em optimal} if it has the
minimum number of columns among covering arrays with the same number of rows.

The addition of a graph structure to covering arrays makes it possible
to use methods from graph theory to study these designs. In this
thesis, we define a family of graphs called the {\em qualitative
independence graphs}. A graph has a covering array, with given
parameters, if and only if there is a homomorphism from the graph to a
particular qualitative independence graph. Cliques in qualitative
independence graphs relate to covering arrays and independent sets
are connected to intersecting partition systems.

It is known that the exact size of an optimal binary covering array
can be determined using Sperner's Theorem and the Erd\H{o}s-Ko-Rado
Theorem. In this thesis, we find good bounds on the size of an optimal
binary covering array on a graph. In addition, we determine both the
chromatic number and a core of the binary qualitative independence
graphs. Since the rows of general covering arrays correspond to set
partitions, we give extensions of Sperner's Theorem and the
Erd\H{o}s-Ko-Rado Theorem to set-partition systems. These results are
part of a general framework to study extremal partition systems.

The core of the binary qualitative independence graphs can be
generalized to a subgraph of a general qualitative independence graph called
the {\em uniform qualitative independence graph}. Cliques in the
uniform qualitative independence graphs relate to {\em balanced} covering
arrays.  Using these graphs, we find bounds on the size of a
balanced covering array. We give the spectra for several of these
graphs and conjecture that they are graphs in an association scheme.

We also give a new construction for covering arrays which yields many
new upper bounds on the size of optimal covering arrays.

	\enlargethispage{40pt}
}

  \prefacesection{Acknowledgements}

   I would like to gratefully acknowledge some of the many people and
   organizations that helped me throughout this degree and in writing
   this thesis.

   First my two supervisors, Lucia Moura and Brett Stevens. It has
   been a pleasure to work with them and I feel very lucky to have had
   two such excellent supervisors.  I have learned an incredible
   amount from Lucia. She has been unbelievably patient, supportive
   and encouraging! This thesis would still be a scattered pile of
   rough notes if it was not for her guidance and diligence. It was
   Brett who convinced me to do this degree in the first place and
   started me on this research project;

   The examiners Charles Colbourn, Mateja \v{S}ajna, Steven Wang and
   Alan Williams who gave many valuable suggestions and raised several
   interesting questions for future work.

   The results in Chapter 6 of this thesis were greatly improved by
   valuable discussions with Chris Godsil and Mike Newman;

   The National Science and Engineering Research Council who provided
   generous financial support;

   My husband Paul Dechene, who not only moved to Ottawa for me, took
   over all the household duties, made me sushi and cookies,
   proof-read several chapters of this thesis, supported me through
   the more difficult days, but also helped me bribe Lucia with his
   delicious brownies;

   My family and Paul's family for being supportive through degree
   after degree;

   Marni Mishna, my math sister, who read several of the chapters and
   has always been an inspiration. Duff McLeod who helped edit a
   chapter, and hasn't been the same since.  The other math students I
   studied with at the University of Ottawa and Carleton University,
   especially Megan Dewar (who made the people who got us confused
   think I am a nicer person than I am), Mr. Paul Elliott-Magwood (whose
   enlightened conversation made coming into the office fun), Geetha
   Patil (who gave me the mysterious advice: ``a thesis is like a
   bollywood movie'') and Sebastian Raaphorst (whose gossip made
   coming into the office fun).
\begin{flushright}
\qed Karen Meagher \\
Ottawa, Canada, June 2005
\end{flushright}
} 
%    \prefacesection{Dedication}  I dedicate this work to my father and mother, for example.
  \afterpreface   % to print the Table of Contents,
\addtolength{\abovedisplayskip}{4pt}
\addtolength{\belowdisplayskip}{4pt}
\addtolength{\jot}{2pt}

\chapter{Introduction}\label{intro}
\thispagestyle{empty}

%currently 2 primary areas of research
%mathematics community - focus on building smaller designs of higher interaction strength

%software testing community - greedy search algorithms to build arrays
%in a more flexible environment.

Covering arrays, also known as qualitatively independent families and
surjective arrays, have been the focus of much research for the last
twenty years. They are natural generalizations of orthogonal arrays and
Sperner systems.  Covering arrays are computationally difficult to
find and useful in multiple applications, particularly for designing
test suites that test interactions of parameters in a
system~\cite{cohen:97, cohen:96, DM, williams:96, williams:00b}.

%definition of cas
For positive integers $n,r,k$, a covering array $CA(n,r,k)$ is an $r
\times n$ array with entries from $\mathbb{Z}_k$ and the property that
in any pair of rows, each of the $k^2$ ordered pairs from
$\mathbb{Z}_k \times \mathbb{Z}_k$ appear in at least one column.
This property is called {\em qualitative independence}, and thus, the
set of rows in a covering array is a qualitatively independent
family. The parameter $n$ is the {\em size} of the covering array and
$\mathbb{Z}_k$ is its {\em alphabet}. For most applications, it
is desirable to have covering arrays with smallest possible size.

Much of the work on covering arrays focuses on developing
constructions for them~\cite{Colbourn3, Colbourn2, Colbourn1,
brett:97b}, finding bounds on their minimum size~\cite{brett:97a} and
more recently, on improving heuristic searches for small covering
arrays~\cite{chateauneuf:99, Colbourn5, Nur, john:thesis}.  As a
result of this work, there are many bounds on the size of covering
arrays with specific parameters. For the special case of binary
covering arrays (covering arrays with $k=2$), the exact sizes of the
smallest array are known~\cite{katona:73, kleitman:73}.  For a fixed
$k$, the asymptotic growth, as $n$ increases, of the largest $r$ such
that a $CA(n,r,k)$ exists is known~\cite{gargano:92, MR94e:05024,
gargano:94}.

The first result in this thesis is a new construction that produces
some of the smallest known covering arrays for specific
parameters. This construction, called {\em the group construction of
covering arrays}, is presented in Section~\ref{sec:groupconst} and
published in~\cite{groupconstCA} and is also used in~\cite{karen2}.
Appendix~\ref{appendix} contains two tables of new upper bounds on
the minimal size of covering arrays.

To improve applications, generalizations of orthogonal arrays and of
covering arrays have been considered~\cite{constructiontestsuites,
MR2071905, MR2012427}. In this thesis, we examine covering arrays with
a graph structure.  Prior to this thesis, there have been few studies
of covering arrays on graphs.  Seroussi and Bshouty~\cite{Seroussi:01}
prove that finding the smallest binary covering array on a graph is an
NP-hard problem and Stevens gave some basic results on covering arrays
on graphs in his Ph.D. thesis~\cite{brett:thesis}.

For positive integers $n$ and $k$, a {\em covering array on a graph
$G$}, denoted $CA(n,G,k)$, is a $|V(G)| \times n$ array, whose entries
are from $\mathbb{Z}_k$, with the property that each row corresponds
to a vertex in $G$ and rows that correspond to adjacent vertices in
$G$ are qualitatively independent. A covering array on a complete
graph is a standard covering array.

Adding a graph structure to covering arrays adds much more than
just an improvement for applications of covering arrays to designing
test suites. It gives us a new way to study covering arrays and
provides interesting results and research problems in graph theory and
extremal set theory.

The goal of this thesis is to develop a new framework to study covering
arrays and covering arrays on graphs. The first way this is done is by
rephrasing a design theory problem as a graph theory problem. This is
achieved with a family of graphs called the {\em qualitative
independence graphs}.  The second way this is done is by generalizing
the results from extremal set theory that are used to find the exact
size of an optimal binary covering array.  To this end, we give
extensions of Sperner's Theorem and the Erd\H{o}s-Ko-Rado Theorem to
partition systems and we lay the foundation for a theory of extremal
partition systems.

The new results in this thesis are mostly contained in
Chapters~\ref{CAG}--\ref{higherparts}. Next, we give an overview
of the major results from each of these four chapters.

\subsubsection{The Qualitative Independence Graphs}

We start with covering arrays on graphs with a particular focus on
binary covering arrays.  We give an upper bound for the size of
a covering array on a graph based on the chromatic number of the graph
and a lower bound based on the size of the maximum clique of the
graph.  The upper bound is of particular interest because it gives a
method to construct covering arrays on graphs from standard covering
arrays and it raises a question that partially motivates this work:
can we find a graph for which this bound is not tight?

In Section~\ref{sec:QIgraphs}, we define a family of graphs, called
the {\em qualitative independence graphs}, denoted $QI(n,k)$, for
positive integers $n$ and $k$. It has become evident that the
qualitative independence graphs are an important family of
graphs. First, for infinitely many of these graphs, the upper bound on
the size of covering array on the graph from the chromatic number is
not tight.  Second, this family of graphs gives a good
characterization of covering arrays for all graphs, namely that for a
graph $G$ and positive integers $k$ and $n$, a $CA(n,G,k)$ exists if
and only if there is a graph homomorphism $G$ to $QI(n,k)$.
Finally, for all positive integers $k$ and $n$, a clique of size $r$
in the graph $QI(n,k)$ corresponds to a $CA(n,r,k)$. This rephrases a
design theory question to a graph theory question.  In this thesis, we
use techniques from graph theory, particularly algebraic graph theory,
to analyze the graphs $QI(n,k)$.

We give formulae for both the size of the maximum clique and the
chromatic number of the graphs $QI(n,2)$. These two formulae are used
to establish good bounds on the minimal size of a binary covering
array on a given graph. Also, using these formulae, we determine the
core of $QI(n,2)$. The structure of these cores implies that for a
graph $G$, if there exists a $CA(n,G,2)$, then there exists a
$CA(n,G,2)$ in which each row has exactly $\lfloor \frac{n}{2}
\rfloor$ ones. 

Finally, in Section~\ref{sec:QIk2k}, for a positive integer $k$, we
find upper bounds on the size of the maximum clique and the chromatic
number of $QI(k^2,k)$ and we prove that this graph is a
$(k!)^{k-1}$-regular graph.

\subsubsection{Eigenvalues and Association Schemes}

For positive integers $n$ and $k$, where $k$ divides $n$, it is clear
how to generalize the core we found for the binary qualitative
independence graphs to a subgraph of $QI(n,k)$.  These subgraphs are called
the {\em uniform qualitative independence graphs} and denoted
$UQI(n,k)$. It is not at all clear if $UQI(n,k)$ is a
core for $QI(n,k)$, nevertheless, the uniform qualitative independence
graphs are an intriguing family of graphs.

In Section~\ref{sec:boundsVT}, we give bounds on the size of the
maximum cliques in a uniform qualitative independence graph. These
bounds are not new, but the proof given here is new and uses
properties of the uniform qualitative independence graphs.

The qualitative independence graph $QI(9,3)$, which is also a uniform
qualitative independence graph, has previously appeared in the
literature in a different context. Mathon and Rosa~\cite{MR86e:05068}
prove that $QI(9,3)$ is part of an association scheme and give its
eigenvalues and their multiplicities.

We show in Section~\ref{sec:betterep} that the method used to find the
eigenvalues in~\cite{MR86e:05068} can be generalized to the uniform
qualitative independence graphs.  In particular, in
Chapter~\ref{chp:ass}, we give an equitable partition on the vertices
of the uniform qualitative independence graphs that can be used to
find their eigenvalues. The eigenvalues and multiplicities for the
graphs $UQI(3c,3)$ for $c=3,4,5,6$ and $QI(16,4)$ are stated.

It is interesting to ask if other qualitative independence graphs are
also classes in an association scheme and in particular, if the scheme
used in~\cite{MR86e:05068} can be generalized. In
Section~\ref{sec:otherschemes}, two sets of graphs are given that are
generalizations of the scheme from~\cite{MR86e:05068}. We conjecture
that these sets of graphs form association schemes and we give their
modified matrix of eigenvalues.

\subsubsection{Extremal Partition Theory}

Since each row of a binary covering array corresponds to a set, two
results from extremal set theory, Sperner's Theorem and the
Erd\H{o}s-Ko-Rado Theorem, can be used to find the exact size of the
smallest binary covering array, $CA(n,r,2)$ for all values of $r$.
The famous Erd\H{o}s-Ko-Rado Theorem~\cite{MR25:3839} is concerned
with the maximal cardinality of intersecting set systems as well as
with the structure of such maximal systems. The equally well-known
Sperner's Theorem is concerned with the cardinality and structure of
the largest system of incomparable sets. 

Since the rows of a non-binary covering array correspond to
partitions, it is desirable to have versions of Sperner's Theorem
and the Erd\H{o}s-Ko-Rado Theorem for partition systems.  In general,
this motivates the study of {\em extremal partition systems}. There is
an extensive theory of extremal finite set systems but there are
almost no similar results for systems of partitions. One goal of this
thesis is to lay down the foundation for the study of extremal
partition systems.

In Section~\ref{sec:spernerpartitions}, we give a new generalization
of Sperner's Theorem for partition systems.  For this generalization we give
an exact result for a class of partition systems, a bound
and an asymptotic result for all partition systems. Our extension of
Sperner's Theorem is different from other Sperner-type theorems for
systems of collections of sets that have been
examined~\cite{gargano:92, MR94e:05024, gargano:94, MR34:5690}.

In Section~\ref{intersectingpartitions} of this thesis, an extension
of the Erd\H{o}s-Ko-Rado Theorem to intersecting set partitions is
proven. We consider two different notions of intersecting partitions.
First, we consider two partitions to be intersecting if they both
contain a common class. With this type of intersection, we prove a
version of the Erd\H{o}s-Ko-Rado Theorem for partition systems. In
particular, we prove that the largest intersecting partition system is
a trivially intersecting system.  Another notion of intersecting
partitions is {\em partial intersection}. Two partitions are partially
$t$-intersecting if they have classes with intersection of cardinality
at least $t$. We conjecture a version of the Erd\H{o}s-Ko-Rado Theorem
for partially intersecting partitions, that is, a trivially partially
intersecting partition system is the largest partially intersecting
partition system. This conjecture is proven in certain restricted
cases in Section~\ref{PIP}. Different extensions of the
Erd\H{o}s-Ko-Rado Theorem to systems of intersecting partitions have
also been considered by P. L. Erd{\H{o}}s, Sz{\'e}kely and
L{\'a}szl{\'o}~\cite{MR2001f:05148}.

\subsubsection{General extremal partition systems}

Higher order extremal problems are extremal problems in which the
elements in the system are disjoint families of subsets of a finite
set called {\em clouds}, rather than sets.  Ahlswede, Cai and
Zhang~\cite{MR1739309} present a framework for several different types
of higher order extremal problems.
\remove{ Further results are given
in~\cite{intersectingsystems, MR1349259, MR1368831}. }

In Chapter~\ref{higherparts}, we examine higher order extremal
problems in which the elements of the system are partitions, rather
than clouds. We use a framework, similar to the one defined by
Ahlswede, Cai and Zhang, to present a variety of different types of
extremal problems for partition systems.  For many of these problems
we give an exact solution and for others we give bounds on the
cardinality of the maximum partition system.  These results include
our extension of Sperner's Theorem to partition systems and our
Erd\H{o}s-Ko-Rado Theorem for partition systems.  This chapter
includes many open problems and provides a unified framework for
problems considered in previous chapters.

\subsubsection{Overview of the Document}

%covering arrays
In Chapter~\ref{designs}, we introduce covering arrays, define three
closely related designs (orthogonal Latin squares, transversal designs
and orthogonal arrays) and describe three constructions for covering
arrays.  In Chapter~\ref{EST}, we show three common types of set
systems and state results from extremal set theory that are related to
covering arrays (Sperner's Theorem, Bollob\'as's Theorem and the
Erd\H{o}s-Ko-Rado Theorem), and we give an application of Sperner's
Theorem and the Erd\H{o}s-Ko-Rado Theorem to covering arrays.  In
Chapter~\ref{graphtheory}, we review some basic graph theory and
algebraic graph theory, including spectral theory of graphs and
association schemes. We use these to study the qualitative
independence graphs.  In Chapter~\ref{CAG}, we introduce the object
which motivates much of this work: covering arrays on graphs.  We also
introduce the qualitative independence
graphs. This chapter focuses on the binary qualitative independence
graphs.  We define the uniform qualitative independence graphs in
Chapter~\ref{chp:ass}. Using an equitable partition on the vertices of
this graph it is possible to find the spectrum of these graphs --- the
spectrum of several instances of these graphs is given in this chapter.
Then, we consider whether these graphs are part of an association
scheme.  In Chapter~\ref{partitionsystems}, extensions of Sperner's
Theorem and Erd\H{o}s-Ko-Rado Theorem are given for partition
systems. In Chapter~\ref{higherparts}, a general framework for
extremal partition systems is described. In Chapter~\ref{Questions},
we list open questions and conjectures that arise from this
study. Finally, in Appendix~\ref{appendix} we present a table of the new
best bounds on the size of covering arrays.  These new bounds come
from our new construction given in Chapter~\ref{designs}.

%%%%%%%%%%%%%%%%%%%%%%%%%%%%%%%%%%%%%%%%%%%%%%%%%%%%%%%

\chapter{Covering Arrays and Related Designs}\label{designs}
\thispagestyle{empty}

Covering arrays, which are also known as qualitatively independent
families and surjective arrays, are a generalization of the well-known
and well-studied orthogonal arrays~\cite{MR2000h:05042}. They are a
mathematically rich design with many applications.

Orthogonal arrays were first developed by Rao in the 1940s for
designing statistical experiments. Experiments based on orthogonal
arrays are particularly effective; such experiments allow for
independent estimation of each factor in the system and the effect of
different pairs of factors can be compared with the same accuracy (for
more information, see Section 11.4 of~\cite{MR2000h:05042}). But, for
many parameters, it is not possible to build an orthogonal
array. Covering arrays are a relaxation of orthogonal arrays (of index
1) that still provide good test design. Orthogonal arrays are closely
related to two other designs: Latin squares and transversal
designs. In this chapter, we give an overview of covering arrays and
related designs. The only new result in this chapter is a group
construction for covering arrays given in
Section~\ref{sec:groupconst}.

\section{Mutually Orthogonal Latin Squares}\label{sec:MOLS}

Latin squares are a very old and well-studied design. It is believed
that these were first studied by Euler around 1782. Latin squares have
also appeared in works of art, math puzzles and even game shows.
\enlargethispage{40pt}

\begin{defn}[Latin Square]\index{Latin square}
Let $n$ be a positive integer. A {\em Latin square of order $n$} is an
$n \times n$ array on $n$ symbols with the property that each symbol
occurs exactly once in each row and each column.
\end{defn}

Unless otherwise stated, in this thesis we use the set
$\{0,1,\dots, n-1\}$ as the symbols in a Latin square of order $n$.  A
Latin square is in {\em standard form} if the symbols in the first row
occur in their natural order.\index{standard form}

\begin{exmp}
Below are two Latin squares.  The first square, $L$, is a Latin square
of order 4 in standard form.  The second Latin square, $M$, is of order
5 and is also in standard form. We remark that it is the addition
table for $\mathbb{Z}_5$.

\begin{center}
\begin{tabular}{cc}
$L = \left[
\begin{tabular}{cccc}
0 & 1 & 2 & 3 \\
3 & 2 & 1 & 0 \\
1 & 0 & 3 & 2 \\
2 & 3 & 0 & 1 
\end{tabular}
\right] $
& 
$M = \left[
\begin{tabular}{ccccc}
0 & 1 & 2 & 3 & 4 \\
1 & 2 & 3 & 4 & 0 \\
2 & 3 & 4 & 0 & 1 \\
3 & 4 & 0 & 1 & 2 \\
4 & 0 & 1 & 2 & 3
\end{tabular}\right]$
\end{tabular}
\end{center}
\end{exmp}

For every $n$, there exists a Latin square of order $n$; the addition table
for $\mathbb{Z}_n$ is always a Latin square of order $n$.  Any
permutation of the rows or the columns of a Latin square is still a
Latin square. Also, any relabelling of the symbols in a Latin square
is still a Latin square. Thus any Latin square can be put in standard
form by permuting the columns or relabelling the symbols.

\begin{defn}[Orthogonal Latin Squares]
Let $n$ be a positive integer and let $A = [a_{i,j}]$ and
$B=[b_{i,j}]$ be Latin squares of order $n$.  Consider the $n \times n$
array with entries $[ (a_{i,j},b_{i,j} )]$; if, in this array, each of
the $n^2$ ordered pairs of symbols occurs exactly once, then the squares $A$
and $B$ are {\em orthogonal}\index{orthogonal}.
\end{defn}

If $A$ and $B$ are a pair of orthogonal Latin squares, then any
relabelling of the symbols in either $A$ or $B$, produces another pair
of Latin squares which are orthogonal. So for any pair of orthogonal
Latin squares, both can be put in standard form.

\begin{exmp}
Below are two orthogonal Latin squares of order 4. Note that both
Latin squares are in standard form.

\begin{center}
\begin{tabular}{cc}
 $ L=\left[
\begin{tabular}{cccc}
0 & 1 & 2 & 3 \\
3 & 2 & 1 & 0 \\
1 & 0 & 3 & 2 \\
2 & 3 & 0 & 1 
\end{tabular}
\right]$, 
&$ L' = \left[
\begin{tabular}{cccc}
0 & 1 & 2 & 3 \\
2 & 3 & 0 & 1 \\
3 & 2 & 1 & 0 \\
1 & 0 & 3 & 2 
\end{tabular} 
\right]$
\end{tabular}
\end{center}

\begin{center}
\[
[ (L_{i,j},L'_{i,j} )] =
\left[ 
 \begin{tabular}{cccc}
(0,0) & (1,1) & (2,2) & (3,3) \\
(3,2) & (2,3) & (1,0) & (0,1) \\
(1,3) & (0,2) & (3,1) & (2,0) \\
(2,1) & (3,0) & (0,3) & (1,2) 
\end{tabular} 
\right]
\]
\end{center}

\end{exmp}

A set of Latin squares with the property than any two distinct squares
from the set are orthogonal is called a set of \index{MOLS} {\em
mutually orthogonal Latin squares} (MOLS).  Let $N(n)$
denote the cardinality of the largest set of MOLS of order $n$. For
all $n$,
\begin{eqnarray}\label{eq:maxMOLS}
1 \leq N(n)\leq n-1.
\end{eqnarray}

To derive this upper bound, consider a set of $N(n)$ Latin squares in
standard form.  For each of the Latin squares, the $(1,1)$ entry is 0
and the $(1,i+1)$ entry is $i$.  Consider the symbol which occurs in
the $(2,1)$ position of each Latin square in this set.  This symbol
must be one of $n-1$ letters from the set $\{1,\dots, n-1\}$.  Assume
two Latin squares $A$ and $B$ from the set have the entry $i \in
\{1,\dots, n-1\}$ in the $(2,1)$ position.  Then $(A_{1,i+1}, B_{1,
i+1}) = (i,i) = (A_{2,1}, B_{2,1})$, so the Latin squares $A$ and $B$
are not be orthogonal.  Thus each symbol in $\{1, \dots ,n-1\}$
corresponds to at most one Latin square in the set and there can be at
most $n-1$ MOLS of order $n$.

A set of $n-1$ MOLS of order $n$ is called a {\em complete set of
MOLS}\index{complete set of MOLS}. When $n$ is a prime power, that is,
$n=p^a$, where $p$ is a prime number and $a$ is a positive integer,
there exists a complete set of MOLS. We do not give a proof of this
theorem here; rather in Section~\ref{sec:finitefieldconstruction}, we
give a construction of an equivalent design.

\begin{thm}[see \cite{street:87}]\label{MOLSprimepower}
For a prime power $n$, $N(n) = n-1$.
\end{thm}

\begin{exmp}\label{exampleMOLS4}
The following is a complete set of MOLS of order 4:
\begin{center}
\begin{tabular}{ccc}
$A_1$ =
\begin{tabular}{|cccc|}
\hline
 0 & 1 & 2 & 3 \\
 1 & 0 & 3 & 2 \\
 2 & 3 & 0 & 1 \\
 3 & 2 & 1 & 0 \\
\hline
\end{tabular}
&
$A_2$ =
\begin{tabular}{|cccc|}
\hline
 0 & 1 & 2 & 3 \\
 2 & 3 & 0 & 1 \\
 3 & 2 & 1 & 0 \\
 1 & 0 & 3 & 2\\
\hline
\end{tabular}
&
$A_3$ =
\begin{tabular}{|cccc|}
\hline
 0 & 1 & 2 & 3 \\
 3 & 2 & 1 & 0 \\
 1 & 0 & 3 & 2 \\
 2 & 3 & 0 & 1 \\
\hline
\end{tabular}
\end{tabular}
\end{center}
\end{exmp}

It is possible to construct orthogonal Latin squares from smaller orthogonal Latin squares.
This construction provides a lower bound on $N(n)$.

\begin{thm}[\cite{street:87}]
For all positive integers $m,n$, we have 
\[
N(mn) \geq \min \{ N(m), N(n)\}.
\]
\end{thm}
\begin{proof}
Let $l = \min \{ N(m), N(n) \}$.  Let $M_1, M_2, \dots, M_{l}$ be the first
$l$ orthogonal Latin squares of order $m$, and let $N_1, N_2,
\dots, N_{l}$ be the first $l$ orthogonal Latin squares of order $n$.
For $0\leq r \leq m-1$, let $N_i^r$ be the Latin square obtained from
$N_i$ by adding $rn$ to each entry of $N_i$.  In each Latin square
$M_i$, for $i=1,\dots, l$, replace each occurrence of the symbol $r$
by the square $N_i^r$.  This set of $l$ squares will be a set of
mutually orthogonal Latin squares of order $mn$.
\end{proof}

For $n=6$ the previous theorem states that $N(6) \geq \min\{N(2), N(3)\}
= 1$.  Since $N(n) \geq 1$ for all $n$, this is not a useful result.
It is non-trivial to conclude that it is impossible to build a pair of
MOLS of order 6 (indeed, $N(6)=1$, see Section 6.4~\cite{street:87}).
Euler conjectured that $N(n)=1$ for all $n \equiv 2\pmod{4}$. This
conjecture was proven false by Bose, Shrikhande and
Parker~\cite{MR0122729}.  The methods they used were generalized (in
the following article in the same issue of the journal!) by Chowla, Erd\H{o}s and
Straus~\cite{MR0122730} to show $\lim_{n \rightarrow \infty} N(n) =
\infty$.  Currently, it is known that $N(2) = 1$, $N(3) = 2$, $N(6) =
1$ and $N(10)\geq 2$, and for all other values of $n$, $N(n) \geq
3$~\cite{street:87}.

\section{Transversal Designs}\label{sec:td}

Transversal designs are a generalization of mutually orthogonal Latin
squares.  Often, results for Latin squares are more conveniently
rephrased in terms of transversal designs (for example, the proof
of the non-existence of a pair of orthogonal Latin squares of order
6 in~\cite{street:87}).

\begin{defn}[Transversal Design] \index{transversal design}
Let $r, \lambda, k$ be positive integers.
A {\em transversal design}, $T[r,\lambda;k]$, is a triple $(X, \mathcal{G}, \mathcal{B})$ 
with the following properties:
\begin{enumerate}
\item $X$ is a set of $rk$ symbols, called {\em varieties}.\index{varieties}

%\item $\mathcal{G}$ is a set of $r$ $k$-subsets of $X$ (called groups), 
%$\mathcal{G} = \{G_1, G_2, \dots, G_r\}$, such that
%the groups $G_i$ partition the set $X$.

\item $\mathcal{G}$ is a partition of the set $X$, 
$\mathcal{G} = \{G_1, G_2, \dots, G_r\}$, where each $G_i$ is a
$k$-subset of $X$. The sets $G_i$ are called {\em groups}.

\item $\mathcal{B}$ is a collection of $r$-sets
of $X$, called {\em blocks}, with the property that each block
intersects each group $G_i$, for $i=1,\dots,r$, in exactly one
variety.

\item Any pair of varieties from different groups occur in the 
same number of blocks.  This number is called the {\em index} and is denoted by $\lambda$.
\end{enumerate}

\end{defn}

\begin{exmp}\label{exmp:td}
The following is a transversal design, $T[4,1;3]$.

\begin{center}
$X = \{0, \; a,\;\alpha,\;\heartsuit,\; 1,\; b, \;\beta, \;\clubsuit,\;
 2,\; c,\; \gamma,\; \diamondsuit\}$ 
\vspace{.5cm}

\begin{tabular}{cc}
$\mathcal{G} = \left\{ 
\begin{array}{cccc}
G_1 = & \{0 & 1 & 2\} \\
G_2 = & \{a & b & c\} \\
G_3 = & \{\alpha & \beta & \gamma\} \\
G_4 = & \{\heartsuit & \clubsuit & \diamondsuit\}
\end{array}
\right\}
$ 
& 
$\mathcal{B} = \left\{ 
\begin{array}{ccccc}
B_1 = & \{0 & a & \alpha & \heartsuit \} \\
B_2 = & \{0 & b & \beta  & \clubsuit \}  \\
B_3 = & \{0 & c & \gamma & \diamondsuit \} \\
B_4 = & \{1 & a & \gamma & \clubsuit \}  \\
B_5 = & \{1 & b & \alpha & \diamondsuit \} \\
B_6 = & \{1 & c & \beta  & \heartsuit \} \\
B_7 = & \{2 & a & \beta  & \diamondsuit \} \\
B_8 = & \{2 & b & \gamma & \heartsuit \} \\
B_9 = & \{2 & c & \alpha & \clubsuit \}  
\end{array}
\right\}$
\end{tabular}
\end{center}
\end{exmp}

\remove{
\begin{exmp}\label{exmp:td}
The following is a transversal design, $T[4,1;3]$.
\begin{center}
\begin{tabular}{ccc}
\begin{tabular}{ccccc}
$X =$&$\{$0, & a, & $\alpha$, & $\heartsuit$,\\
           & 1, & b, & $\beta$, & $\clubsuit$,\\
           & 2, & c, & $\gamma$, & $\diamondsuit \}$
\end{tabular}
& 
\begin{tabular}{cccc}
$G_1$: & 0 & 1 & 2 \\
$G_2$: & a & b & c \\
$G_3$: & $\alpha$ & $\beta$ & $\gamma$ \\
$G_4$: & $\heartsuit$ & $\clubsuit$ & $\diamondsuit$ 
\end{tabular}
& 
\begin{tabular}{ccccc}
$B_1$: & 0 & a & $\alpha$ & $\heartsuit$  \\
$B_2$: & 0 & b & $\beta$  & $\clubsuit$   \\
$B_3$: & 0 & c & $\gamma$ & $\diamondsuit$\\
$B_4$: & 1 & a & $\gamma$ & $\clubsuit$   \\
$B_5$: & 1 & b & $\alpha$ & $\diamondsuit$\\
$B_6$: & 1 & c & $\beta$  & $\heartsuit$  \\
$B_7$: & 2 & a & $\beta$  & $\diamondsuit$\\
$B_8$: & 2 & b & $\gamma$ & $\heartsuit$  \\
$B_9$: & 2 & c & $\alpha$ & $\clubsuit$   
\end{tabular}
\end{tabular}
\end{center}
\end{exmp}
}

\begin{thm}[\cite{street:87}]
Let $r,k$ be positive integer with $r\geq 3$.
A transversal design $T[r,1;k]$ exists if and only if there are
$r-2$ MOLS of order $k$.
\end{thm}
\begin{proof}
Let $(X, \mathcal{G}, \mathcal{B})$ be a $T[r,1;k]$.
For each block $B \in \mathcal{B}$ in the transversal design, use the first pair of
symbols in the block to index the position in a Latin square. Then use
the remaining $r-2$ elements in the block $B$ as the entry in the
given position for each of the $r-2$ Latin squares.  Similarly, a transversal
design can be constructed from a set of MOLS.
\end{proof}

\begin{exmp}
For the transversal design in Example~\ref{exmp:td} the corresponding pair of MOLS of order 3 are:

\begin{center}
\begin{tabular}{ccc}

\begin{tabular}{lccc}
  & 0 & 1 & 2 \\ \cline{2-4}
a & \multicolumn{1}{|c}{$\alpha$} & \multicolumn{1}{c}{$\gamma$} & \multicolumn{1}{c|}{$\beta$}  \\
b & \multicolumn{1}{|c}{$\beta$}  & \multicolumn{1}{c}{$\alpha$} & \multicolumn{1}{c|}{$\gamma$} \\
c & \multicolumn{1}{|c}{$\gamma$} & \multicolumn{1}{c}{$\beta$}  & \multicolumn{1}{c|}{$\alpha$} \\ \cline{2-4}
\end{tabular}
& \phantom{and} &
\begin{tabular}{lccc}
    & 0 & 1 & 2 \\ \cline{2-4}
a  & \multicolumn{1}{|c}{$\heartsuit$}   & \multicolumn{1}{c}{$\clubsuit$}    & \multicolumn{1}{c|}{$\diamondsuit$} \\
b  & \multicolumn{1}{|c}{$\clubsuit$}    & \multicolumn{1}{c}{$\diamondsuit$} & \multicolumn{1}{c|}{$\heartsuit$}   \\
c  & \multicolumn{1}{|c}{$\diamondsuit$} & \multicolumn{1}{c}{$\heartsuit$}   &  \multicolumn{1}{c|}{$\clubsuit$}   \\ \cline{2-4}
\end{tabular}
\end{tabular}
\end{center}

\end{exmp}

Two well-known generalizations of transversal designs are {\em transversal
covers} and {\em point-balanced transversal covers}.  A transversal
cover \index{transversal cover}is a triple $(X,
\mathcal{G}, \mathcal{B})$ with the same properties as a transversal
design $T[r,1;k]$ except that the property that any pair of varieties from
different groups occur in {\em exactly one} block is relaxed to any
pair of varieties from different groups occur in {\em at least one} block.
Transversal covers are denoted by $TC[r,k]$.
The fewest blocks possible in a  $TC[r,k]$ is denoted by $tc[r,k]$.

A point-balanced transversal cover\index{point-balanced
transversal cover} is a transversal cover $TC[r,k]$ with the additional property
that every variety appears in the same number of blocks.
Point-balanced transversal covers are denoted by $PBTC[r,k]$ and
the fewest blocks possible in a $PBTC[r,k]$ is denoted by
$pbtc[r,k]$.

Stevens, Moura and Mendelsohn~\cite{brett:97b} prove several bounds
for transversal covers and point-balanced transversal covers. Four of
these bounds are stated here. These bounds are not stated in full
generality; there are refinements that can be made which are not
included here.

\begin{lemma}[\cite{brett:97b}]\label{tcbounds}
Let $r,k$ be positive integers, then
\begin{enumerate}
\item  $\frac{pbtc(r,k)}{k} + k(k-1) \leq tc(r,k) \leq pbtc(r,k)$; \label{tcbound1}
\item $tc(r,k) \geq \left\lceil \frac{k \log_2 r}{2} \right\rceil$; \label{tcbound2}
\item  $tc(r,k) \geq k^2+2 \mbox{ for all } r \geq k+2 \mbox{ and } k \geq 3$; \label{tcbound3}
\item { if there exists a $PBTC(r,k)$ with $b$ blocks, then
      \[
       r \leq \left\lfloor 
                 \frac{  {b \choose \frac{b}{k}-(k-2)}}{k { \frac{b}{k} \choose k-2}  } 
            \right\rfloor.
      \] \label{tcbound4} 
}
\end{enumerate}
\end{lemma}

\section{Orthogonal Arrays}\label{OA}

This section introduces a third design, orthogonal arrays. The results
in this section (and many more which are not included here) can be
found in Hedayat, Sloane and Stufken's text, {\em Orthogonal
Arrays}~\cite{MR2000h:05042}.

\begin{defn}[Orthogonal Array]\index{orthogonal array}Let $n,r,k,t$ be positive integers with $t\leq r$.
An {\em orthogonal array}, $OA(n,r,k,t)$, of index $\lambda$, with
strength $t$ and alphabet size $k$, is an $r \times n$ array with
entries from $\{ 0,1, \dots, k-1 \}$ and the property that any $t
\times n$ subarray has all $k^t$ possible $t$-tuples occurring as
columns exactly $\lambda$ times.
\end{defn}

The parameter $\lambda$ is not included in the notation since in
any $OA(n,r,k,t)$, $\lambda = n/k^t$.

The rows of an orthogonal array are also called {\em factors}.  Like MOLS and transversal
designs, constructing general orthogonal arrays is difficult. Often the
best results for orthogonal arrays are bounds on the sizes of either
$n$ or $r$ when all other parameters are fixed.

\remove{
The following are a few straightforward facts about orthogonal arrays
that are stated without proof.
\begin{itemize}
\item An orthogonal array of strength $t$ with index $\lambda$ is also an orthogonal array
      of strength $t'$ with index $\lambda k^{t-t'}$ for all $t'< t$.
\item If an $OA(n,r,k,t)$ exists then a  $OA(n,r-1,k,t)$ exists.
\item Permuting symbols in any row of an orthogonal array gives an orthogonal 
      array with the same parameters.
\item Permuting the columns or the rows of an orthogonal array gives an orthogonal 
      array with the same parameters.
\item From an $OA(n,r,k,t)$ it is possible to construct an $OA(n/k,r-1,k,t-1)$:
      Select all the columns from $OA(n,r,k,t)$ that begin with a
      given letter and remove the first row from these columns.
\end{itemize}
}

A strength-2 orthogonal array, $OA(n,r,k,2)$, is equivalent to a
transversal design $T[r,\lambda;k]$, where $\lambda = n/k^2$.  For
each row $r_i$ of an $OA(n,r,k,2)$ with $i \in \{0,\dots, r-1\}$
rewrite each letter $a$ in the row as $a_i$. Then, the $kr$ varieties
in the transversal design are $a_i$ for $a \in \{0,1, \dots, k-1\}$
and $i \in \{0,1, \dots, r-1\}$.  Define the $r$ groups of the
transversal design by $G_i = \{a_i \; :\; a = 0,1, \dots, k-1\}$ for
$i=0,\dots, r-1$. Finally, the blocks in the transversal design are the
columns in the covering array.  Similarly, it is possible to construct
an $OA(n,r,k,2)$ from a $T(r,\lambda;k)$ where $n = k^2 \lambda$.

\begin{exmp}
The following $OA(9,4,3,2)$ is equivalent to the transversal design
$TD[4,1;3]$ from Example~\ref{exmp:td} in Section~\ref{sec:td}.
\begin{center}
$\left[
\begin{tabular}{ccc ccc ccc}
0 & 0 & 0 & 1 & 1 & 1 & 2 & 2 & 2 \\
0 & 1 & 2 & 0 & 1 & 2 & 0 & 1 & 2 \\
0 & 1 & 2 & 2 & 0 & 1 & 1 & 2 & 0 \\
0 & 1 & 2 & 1 & 2 & 0 & 2 & 0 & 1 
\end{tabular}
\right]$
\end{center}

\end{exmp}

This equivalence between transversal designs and strength-2 orthogonal
arrays implies that the existence of an $OA(k^{2},r,k,2)$ is
equivalent to the existence of $r-2$ MOLS of order $k$.  The
equivalence between orthogonal arrays, transversal designs and MOLS
are summarized in the theorem below.

\begin{thm}[see \cite{stinson:02}]\label{thm:equivdesigns}
Let $r,k$ be positive integers with $k\geq 2$ and $r \geq 3$.  Then
the existence of any one of the following designs implies the
existence of the other two designs:
\begin{enumerate}
\item{an $OA(k^2,r,k,2)$,}
\item{a $T[r,1;k]$,}
\item{$r-2$ $MOLS(k)$.}
\end{enumerate}
\end{thm}

Results on the largest set of MOLS can be translated to results for
orthogonal arrays (for a table of largest known MOLS see Section II.2
of~\cite{colbourn:96}). For example, Inequality~(\ref{eq:maxMOLS}) from
Section~\ref{sec:MOLS} implies the following corollary.

\begin{cor}\label{OAnonprimepower}
For $k$ a positive integer, if an
$OA(k^2,r,k,2)$ exists, then $r \leq k+1$. 
\end{cor}

Further, Theorem~\ref{MOLSprimepower} can be translated to the
following corollary.  Again, the proof of this is delayed until
Section~\ref{sec:finitefieldconstruction} where a construction for an
equivalent design is given.

\begin{cor}\label{OAprimepower}
For $k$ a prime power, there exists an $OA(k^{2},k+1,k,2)$. 
\end{cor}

The requirement that $k$ be a prime power can not be dropped from
Corollary~\ref{OAprimepower}. For example, since $N(6) = 1$, it is not
possible to build an $OA(36,4,6,2)$ so no $OA(36,r,6,2)$ for
$r\geq 4$ exists.

\section{Covering Arrays}\label{CA}

For fixed positive integers $n,k$ and $t$, and for some values $r$, it is
impossible to build an orthogonal array $OA(n,k,r,t)$.  For example,
for an $OA(9,r,3,2)$, Inequality~(\ref{tcbound4}) from
Lemma~\ref{tcbounds} shows that
\[
r \leq  \frac{{9 \choose \frac{9}{3}-1}}{3 { 3 \choose 1}} = 4.
\]
Thus, it is impossible to build an orthogonal array with $k=3,t=2$ and
$\lambda=1$ with more than four rows.  It may be possible to build an
orthogonal array $OA(9\lambda, 5,3,2)$ for some $\lambda$
sufficiently large, but for many applications increasing the $\lambda$
produces an array with too many columns.  Instead, the condition on
the array that each $t$-tuple occurs {\em exactly} once can
be relaxed to be that each $t$-tuple occurs {\em at least} once.  These arrays
are called {\em covering arrays}, as all $t$-tuples are covered in
the array.

\begin{defn}[Covering Array]\label{cak}\index{covering array}
Let $n,r,k,t$ be positive integers with $t \leq r$.  A {\em covering
array}, $t$-$CA(n,r,k)$, with strength $t$ and alphabet size $k$ is an $r
\times n$ array with entries from $\{ 0,1, \dots, k-1\}$ and the property
that any $t \times n$ subarray has all $k^t$ possible $t$-tuples
occurring at least once.
\end{defn}

\begin{exmp}\label{caexample}\label{exmp:optimalCA}
The following array is a 2-$CA(11,5,3)$
\begin{center}
$\left[
\begin{tabular}{ccccccccccc}
0 & 0 & 2 & 1 & 1 & 1 & 0 & 1 & 2 & 2 & 2 \\
0 & 1 & 0 & 2 & 1 & 1 & 2 & 0 & 1 & 2 & 2 \\
0 & 1 & 1 & 0 & 2 & 1 & 2 & 2 & 0 & 1 & 2 \\
0 & 1 & 1 & 1 & 0 & 2 & 2 & 2 & 2 & 0 & 1 \\
0 & 2 & 1 & 1 & 1 & 0 & 1 & 2 & 2 & 2 & 0 
\end{tabular}
\right].$
\end{center}
It is possible to construct an $OA(18,5,3)$ which has every pair
occurring exactly twice, but this would require 7 more columns than the
above covering array.
\end{exmp}

The number of columns, $n$, in a $t\mbox{-}CA(n,r,k)$ is the {\em size}\index{size of a
covering array} of the covering array.  For many applications it is
best to use the covering array with the smallest size.  Thus, the
smallest possible size of a covering array is denoted by $t\mbox{-}CAN(r,k)$, that is
\[
t\mbox{-}CAN(r,k)=\min_{l \in \mathbb{N}} \{l: \exists \;t\mbox{-}CA(l,r,k) \}. 
\]
A covering array $t\mbox{-}CA(n,r,k)$ with $n = t\mbox{-}CAN(r,k)$ is said to be
\index{optimal covering array} {\em optimal}.

It is often useful to consider the maximum number of rows
possible in a covering array with $n$ columns on a given alphabet.
This is denoted by $t\mbox{-}N(n,k)$, that is
\[
t\mbox{-}N(n,k)=\max_{r \in \mathbb{N}} \{r: \exists \;t\mbox{-}CA(n,r,k) \}. 
\]

Throughout this thesis only strength-2 covering arrays are
considered. Thus, the $t$ will be dropped from the notation, so
$CA(n,r,k)$, $CAN(r,k)$ and $N(n,k)$ will be used to denote 
2-$CA(n,r,k)$, 2-$CAN(r,k)$ and 2-$N(n,k)$.

Just as strength-2 orthogonal arrays correspond to transversal
designs, covering arrays correspond to transversal covers. The bounds
from Lemma~\ref{tcbounds} apply to covering arrays. In particular,
from Inequality~(\ref{tcbound3}) of Lemma~\ref{tcbounds}, $CAN(k+2,k,2)
\geq k^2+2$, so the covering array in Example~\ref{exmp:optimalCA} is
an optimal $CA(11,5,3)$.

There is a different way to consider the rows of a covering array
which we also use throughout this thesis.

\begin{defn}[Qualitatively Independent Vectors]\label{defn:qi1}\index{qualitatively independent vectors}
Let $k,n$ be positive integers.  Two vectors $u,v \in \mathbb{Z}_k^n$
are {\em qualitatively independent} if for each one of the possible
$k^2$ ordered pairs $(a,b) \in \mathbb{Z}_k \times \mathbb{Z}_k$, there
is an index $i$ so $(u_i,v_i) = (a,b)$. A set of vectors is
qualitatively independent if any two distinct vectors in the set are
qualitatively independent.
\end{defn}

The set of rows in a covering array $CA(n,r,k)$ is a set of $r$
pairwise qualitatively independent vectors from $\mathbb{Z}_k^n$.

In Section~\ref{sec:qisets}, qualitatively independent {\em sets} are
considered; these are defined in Definition~\ref{defn:qi2}.
Qualitatively independent sets are equivalent to qualitatively
independent binary vectors --- that is, qualitatively independent vectors
from $\mathbb{Z}_2^n$ (see Section~\ref{sec:qisets}).  Qualitative
independence can also be extended to partitions
(Definition~\ref{defn:qi3}), and, for all positive integers $k$ and
$n$, qualitatively independent $k$-partitions are equivalent to
qualitatively independent vectors in $\mathbb{Z}_k^n$ (see
Section~\ref{sec:qipartitions}).

The most commonly cited application of covering arrays is for testing
systems (usually software or networks).  Assume that the system to be
tested has $r$ parameters and that each parameter can take $k$
different values. Each row in the array corresponds to a parameter or
a variable in the system. The letters in a $k$-alphabet, $\{0,1,
\dots, k-1\}$, correspond to the different values the parameters can
be assigned. Each column of the covering array describes the values
for the parameters in a test run for the system.  The test suite
described by all the columns in a covering array will, for every pair
of parameters, test all $k^{2}$ possible values of those two
parameters. This means, any two parameters will be tested completely
against one another.  There are several references that
support the practical effectiveness of test suites that guarantee this
pairwise coverage~\cite{cohen:97, cohen:96, DM, williams:96, williams:00b}.

\subsection{Constructions for Covering Arrays}\label{CAconstructions}

In light of this application, there has been much research to find
constructions of covering arrays with the fewest possible
columns~\cite{chateauneuf:99,sloane:93,brett:97a}.
Colbourn~\cite{Colbourn4} gives a comprehensive survey of results for
covering arrays. Determining the exact value of $CAN(r,k)$, for a
given $r$ and $k$, is, in general, a difficult problem. A construction
for covering arrays gives an upper bound on $CAN(r,k)$.  In this
section, three constructions will be given: the finite field
construction, the block-size recursive construction and the group
construction.  This last construction is new and improves many of the
upper bounds on the size of covering arrays.

\subsubsection{Finite Field Construction}\label{sec:finitefieldconstruction}

The construction given in this section is the very well-known {\em
finite field construction for orthogonal arrays}\index{finite field
construction}. It is essentially this construction that is used to
construct a complete set of MOLS of order $k$, where $k$ is a prime
power (Theorem~\ref{MOLSprimepower}).  This construction builds the
$OA(k^2,k+1,k,2)$ from Corollary~\ref{OAprimepower} which is also a $CA(k^2,k+1,k)$.

\begin{lemma}\label{mols}
Let $k$ be a prime power, then $CAN(k+1,k) = k^2$.
\end{lemma}
\begin{proof}
We will prove this by constructing a $CA(k^2,k+1,k)$, call this $C$.
In this proof, we will index the rows and columns of $C$ starting from 0.

Let $GF[k]$ be the finite field of order $k$.  Fix an ordering of the elements of
$GF[k]$, $\{ f_0, f_1, f_{2}, \dots, f_{k-1} \}$ with $f_0 =0$ and $f_1=1$. 

Let the first row of $C$ be each element in the field repeated $k$
times in the fixed order, so the entry in column $c$ of the first row is
$f_l$, where $l=\lfloor c/k \rfloor$.  For $i = 0, \dots ,k-1$, set the
entry in row $i+1$ and column $c$ of the covering array to be $f_if_l
+ f_j$ where $l = \lfloor c/k \rfloor$ and $j \equiv c \pmod{k}$.

Any row of $C$ is qualitatively independent from any other row of $C$.
First, it is not hard to see that the first row is qualitatively
independent from all other rows. Next, assume rows $i_0+1$ and $i_1+1$
(for any distinct $i_0, i_1 \in \{0,1, \dots , k-1\}$) are not
qualitatively independent. Then for columns $r$ and $s$, a pair is
repeated between two distinct rows $i_0+1$ and $i_1+1$. In particular,
$( C_{i_0+1,r}, C_{i_1+1,r}) = (C_{i_0+1,s}, C_{i_1+1,s})$.  Thus,
\[ f_{i_0}f_{l_r} + f_{j_r} = f_{i_0}f_{l_s} + f_{j_s},\] and 
\[ f_{i_1}f_{l_r} + f_{j_r} = f_{i_1}f_{l_s} + f_{j_s},\]
where $l_r = \lfloor r/k \rfloor$ and $j_r \equiv r \pmod{k}$, and 
$l_s = \lfloor s/k \rfloor$ and $j_s \equiv s \pmod{k}$.

As $GF[k]$ is a field, $f_{l_r} = f_{l_s}$, or $l_r = l_s$.  From
this it also follows that $ f_{j_r} = f_{j_s}$ and $j_r=j_s$.  These
together imply that $f_{i_0} = f_{i_1}$ which contradicts $i_0+1 \neq
i_1+1$.

\end{proof}

Inequality~(\ref{eq:maxMOLS}) from Section~\ref{sec:MOLS} and
Lemma~\ref{mols} can be restated in terms of $N(n,k)$.

\begin{lemma}\label{omegag2g}
For any positive integer $k$, $N(k^2,k) \leq k+1$, and for $k$ a prime power, $N(k^2,k)
= k+1$.
\end{lemma}

\begin{exmp}\label{ffcaexample}
The following array is a $CA(16,5,4)$ from the finite field
construction on the finite field with four elements, labelled
$\{0,1,2,3\}$.  This covering array corresponds to the complete set of
MOLS of order four in Example~\ref{exampleMOLS4} in Section~\ref{sec:MOLS}.
\begin{center}
$\left[
\begin{tabular}{cccc cccc cccc cccc}
0 & 0 & 0 & 0 & 1 & 1 & 1 & 1 & 2 & 2 & 2 & 2 & 3 & 3 & 3 & 3 \\
0 & 1 & 2 & 3 & 0 & 1 & 2 & 3 & 0 & 1 & 2 & 3 & 0 & 1 & 2 & 3 \\
0 & 1 & 2 & 3 & 1 & 0 & 3 & 2 & 2 & 3 & 0 & 1 & 3 & 2 & 1 & 0 \\
0 & 1 & 2 & 3 & 2 & 3 & 0 & 1 & 3 & 2 & 1 & 0 & 1 & 0 & 3 & 2 \\
0 & 1 & 2 & 3 & 3 & 2 & 1 & 0 & 1 & 0 & 3 & 2 & 2 & 3 & 0 & 1 \\
\end{tabular}
\right]$
\end{center}

\end{exmp}

Note that for each row, other than the first row, the first four entries
are (in order) $0,1,2,3$.  This leads to the concept of {\em disjoint
columns}\index{disjoint columns}.  In a covering array, two columns are
disjoint if, for each row, the two columns have different
entries.  An example of two columns that are disjoint is a
column of all 0s and a column of all 1s.  A covering array has $m$
disjoint columns if there is a set of at least $m$ columns that are
pairwise disjoint. 

The term ``disjoint columns'' comes from the correspondence between
covering arrays and transversal covers --- in particular, each column
of a covering array corresponds to a block in a transversal cover.  If
two columns in the covering array are disjoint, then the two blocks in
the transversal cover corresponding to these columns are disjoint
sets.

\begin{cor}\label{lem:disjoint}
For any prime power $k$, there exists a covering array
$CA(k^2,k,k)$ with $k$ disjoint columns.
\end{cor}
\begin{proof}
Let $C$ be the covering array $CA(k^2,k+1,k)$ built by the finite
field construction. Let $C'$ be the covering array formed by removing
the first row of $C$.  By the finite field construction, for columns
$j=0,\dots, k-1$ the entry on row $i$ of $C'$ is $f_i 0 + f_j =
f_j$. Thus the first $k$ columns of $C'$ are disjoint.
\end{proof}

\begin{exmp}\label{ffcaexampledisjoint}
The following array is a $CA(16,4,4)$ with four disjoint columns.
\begin{center}
$\left[
\begin{tabular}{|cccc| cccc cccc cccc} \cline{1-4}
0 & 1 & 2 & 3 & 0 & 1 & 2 & 3 & 0 & 1 & 2 & 3 & 0 & 1 & 2 & 3 \\ 
0 & 1 & 2 & 3 & 1 & 0 & 3 & 2 & 2 & 3 & 0 & 1 & 3 & 2 & 1 & 0 \\
0 & 1 & 2 & 3 & 2 & 3 & 0 & 1 & 3 & 2 & 1 & 0 & 1 & 0 & 3 & 2 \\
0 & 1 & 2 & 3 & 3 & 2 & 1 & 0 & 1 & 0 & 3 & 2 & 2 & 3 & 0 & 1 \\ \cline{1-4}
\end{tabular}
\right]$
\end{center}

\end{exmp}

\subsubsection{Block-size Recursive Construction}\label{sec:blockrecursive}

The next construction is known as \index{block-size recursive
construction} the {\em block-size recursive construction}; it appears
in~\cite{poljak:83} and~\cite{brett:97a}. This construction uses two
covering arrays with the same alphabet, $A$ a $CA(n,r,k)$ and $B$ a
$CA(m,s,k)$, to build a $CA(m+n,rs,k)$. Let $a_i$ for $i=0,1, \dots,
r-1$ denote the rows of $A$, and let $b_j$ for $j =0,1, \dots, s-1$
denote the rows of $B$. So,

\begin{center}
\begin{tabular}{ccc}
$A = 
\overbrace{
\left[ 
\begin{array}{c}
{a_0} \\
{a_1} \\
\vdots \\
{a_{r-1}}
\end{array}
\right]
}^{n \textrm{ columns}}$
& and &
$B =
\overbrace{
\left[ 
\begin{array}{c}
{b_0} \\
{b_1} \\
\vdots \\
{b_{s-1}}
\end{array}
\right]
}^{m \textrm{ columns}}$.
\end{tabular}
\end{center}

Construct a $CA(n+m,rs,k)$ as follows.
The first $s$ rows of the $CA(n+m,rs,k)$ are row $a_0$ of $A$ concatenated
with row $b_j$ of $B$ for $j=0, \dots , s-1$. The next $s$ rows of $CA(n+m,rs,k)$ are
row $a_1$ of $A$ concatenated with row $b_j$ of $B$ for $j=0,\dots, s-1$.  In general, row $t$
of $CA(n+m,rs,k)$ is row $a_i$, where $i= \lfloor t/s \rfloor$,
concatenated with row $b_j$, where $j \equiv t \pmod{s}$. Thus,

\[
CA(n+m,rs,k) = 
\overbrace{
\left[ 
\begin{array}{c|c}
{a_0} & {b_0}  \\
{a_0} & {b_1}  \\
\vdots & \vdots \\
{a_0} & {b_{s-1}}  \\
\hline 
{a_1} & {b_0}  \\
{a_1} & {b_1}  \\
\vdots & \vdots\\
{a_1} & {b_{s-1}}  \\
\hline 
\vdots & \vdots\\
\hline
{a_{r-1}} & {b_0}  \\
{a_{r-1}} & {b_1}  \\
\vdots & \vdots \\
{a_{r-1}} & {b_{s-1}}  
\end{array}
\right]
}^{n+m \textrm{ columns}}.
\]

This is indeed a covering array.  Any
two distinct rows $t_0$ and $t_1$ of this $CA(m+n,rs,k)$ are of the
form $a_{i_0}b_{j_0}$ and $a_{i_1}b_{j_1}$. If $i_0 \neq i_1$ then all
possible pairs in the $k$-alphabet occur in the first $n$ columns
between $a_{i_0}$ and $a_{i_1}$. If $i_0 = i_1$, then $j_0 \neq j_1$ and all
possible pairs in the $k$-alphabet occur in the last $m$ columns
between $b_{i_0}$ and $b_{i_1}$.

\begin{thm}[\cite{poljak:83},\cite{brett:97a}]\label{thm:construction}
If there exists a $CA(m,s,k)$ and a $CA(n,r,k)$, then there
exists a $CA(m + n, sr,k)$.
\end{thm}

Let $C$ be a $CA(m + n, sr,k)$ built by the block-size recursive
construction.  For any two rows of $C$, for all $i \in
\{0,\dots,k-1\}$, the pairs $(i,i) \in \mathbb{Z}_k \times
\mathbb{Z}_k$ must occur in the first $n$ columns of $C$. This is
clear because each pair of rows has either the first $n$ columns the
same or the first $n$ columns are a pair of distinct rows from a
covering array. Either way, any pair $(i,i)$ is covered.  Similarly,
for any two rows of $C$, the pairs $(i,i)$, for all $i \in
\{0,\dots,k-1\}$, also occur in the last $m$ columns.

It is not necessary to cover these pairs twice; if we could remove
some of these pairs we could improve the block-size recursive
construction.  In particular, if one of the covering arrays in the
construction has a set of disjoint columns, the letters in each row of
this array can be relabelled so that the disjoint columns are constant
columns; then these columns can be removed from the final covering
array.  For example, if the block-size recursive construction is used
with the covering array $CA(k^2,k+1,k)$ from the finite field
construction (Lemma~\ref{mols}) and the $CA(k^2,k,k)$ with $k$
disjoint columns (Corollary~\ref{lem:disjoint}), all pairs $(i,i) \in
\mathbb{Z}_k \times \mathbb{Z}_k$ are covered in the first $k^2$
columns, so the first $k$ columns of $CA(k^2,k,k)$ do not need to be
included in the final covering array.

\begin{exmp}
Let $k$ be a prime power. Denote the elements of the finite field of
order $k$ by $\{0,1,...,k-1\}$.  The following covering array is
constructed with the block-size recursive construction using the
covering array $CA(k^2,k+1,k)$ from the finite field construction,
Lemma~\ref{mols}, and the $CA(k^2,k,k)$ with $k$ disjoint columns from
Corollary~\ref{lem:disjoint}. Note that columns $k^2+1$ through
$k^2+k$ can be removed so this is a $CA(k^2-k,k+1,k)$,
\orphan
\[
\scriptsize{
\left[ 
\begin{tabular}{c@{}c@{}c@{}c  c@{}c@{}c@{}c c  c@{}c@{}c@{}c  | c@{}c@{}c@{}c c@{}c@{}c@{}c  c  c@{}c@{}c@{}c}
0&0& \dots & 0  &  1 &1& \dots& 1   & \dots & k-1& k-1& \dots & k-1  \, &
                 0&1& \dots & k-1\hspace{-1.2cm}-------------&0&1 &\dots  & k-1 & \dots & 0 & 1 & \dots & k-1 \\
0&0& \dots & 0  &  1 &1& \dots& 1   & \dots & k-1& k-1& \dots & k-1  &  
                 0&1& \dots & k-1\hspace{-1.2cm}-------------&1&2 &\dots  & 0   & \dots &  k-1 & 0 & \dots & k-2  \\
\vdots &&&&&&&&&&&&&&&&&&&&&&&&& \\
0&0& \dots & 0  &  1 &1& \dots& 1   & \dots & k-1& k-1 &\dots & k-1  & 
                 0&1& \dots & k-1\hspace{-1.2cm}-------------& &&  & && &\dots  &\\
\hline 
0&1& \dots &k-1 &  0 &1& \dots& k-1 & \dots & 0  & 1   & \dots & k-1 &
                 0&1& \dots & k-1\hspace{-1.2cm}------------&0&1& \dots & k-1 & \dots & 0 & 1  &\dots & k-1\\
0&1& \dots &k-1 &  0 &1& \dots& k-1 & \dots & 0  & 1   &\dots & k-1  & 
                 0&1& \dots & k-1\hspace{-1.2cm}-------------&1&2& \dots &  0  & \dots & k-1 & 0  &\dots & k-2 \\
\vdots  &&&&&&&&&&&&&&&&&&&&&&&&& \\
0&1& \dots &k-1 &  0 &1& \dots &k-1 & \dots & 0  &1    & \dots & k-1 & 
                 0&1& \dots & k-1\hspace{-1.2cm}-------------& &&  & && &\dots  & \\
\hline 
0&1& \dots &k-1 &    &  &\dots &    &       &    &     &       &     & 
                 0&1& \dots & k-1\hspace{-1.2cm}-------------&0&1& \dots &k-1 & \dots & 0 & 1 & \dots & k-1 \\
0&1& \dots &k-1 &    &  &\dots &    &       &    &     &       &     & 
                 0&1& \dots & k-1\hspace{-1.2cm}-------------&1&2& \dots & 0  & \dots &  k-1  & 0  &\dots  &k-2  \\
\vdots &&&&&&&&&&&&&&&&&&&&&&&&&  \\
0&1& \dots &k-1 &    &  &\dots &    &       &    &      &      &     & 
                 0&1& \dots & k-1\hspace{-1.2cm}-------------&&&  & && &\dots  &  \\
\end{tabular}
\right].
}
\]
\end{exmp}

The result of this construction is that it is possible to build small 
covering arrays with a large number of rows when the alphabet is a prime power.
\begin{lemma}[\cite{poljak:83}, \cite{brett:97a}]\label{blockrec}
For a prime power $k$, there exists a $CA(2k^2-k, k(k+1), k)$.
Equivalently, for any prime power $k$ and any integer $r \leq k(k+1)$,
\[
CAN(r,k) \leq 2k^2-k.
\]
\end{lemma}

Further, the block-size recursive method applied to the covering array
$CA(2k^2-k, k(k+1),k)$, from Lemma~\ref{blockrec} and the covering
array $CA(k^2, k,k)$ with $k$ disjoint columns produces a covering
array $CA(3k^2-2k, k^2(k+1),k)$. Applying the block-size recursive
construction with this new covering array and the array $CA(k^2, k,k)$
with $k$ disjoint columns $i$ times produces a $CA(k^2+i(k^2-k),
k^i(k+1),k)$.  Moreover, in each row of this covering array, each
letter occurs exactly $k + i(k-1)$ times.

A covering array $CA(n,r,k)$ with the property that each letter
occurs exactly $n/k$ times in every row is a {\em balanced covering
array}\index{balanced covering array}. 

%Balanced covering arrays will be denoted by $BCA(n,r,k)$.

\begin{lemma}\label{fullblockrec}
For $k$ a prime power and any positive integer $i$,
there exists a balanced covering array $CA(k^2+i(k^2-k), k^i(k+1),k)$.
Thus,
\[
CAN( k^i(k+1),k) \leq k^2+i(k^2-k).
\]
\end{lemma}

Two extensions of this construction are given by Colbourn, Martirosyan,
Mullen, Shasha, Sherwood and Yucus~\cite{Colbourn3}.

\subsubsection{Group Construction of Covering Arrays}\label{sec:groupconst}

For a prime power $k$ we know from Lemma~\ref{mols} that there exists
$CA(k^2, k+1, k)$ and from Lemma~\ref{blockrec} that there exists a
$CA(2k^2-k, k(k+1), k)$.  We give here a new construction that can be
used to build covering arrays $CA(n,r,k)$, where $k^2 \leq n \leq
2k^2-k$ and $k+1 \leq r \leq k(k+1)$. This construction appeared
in~\cite{karen2, groupconstCA}.

This construction involves selecting a subgroup of the symmetric group
on $k$ elements, $G < Sym(k)$, and finding a {\em starter vector}, $v \in
\mathbb{Z}_k^r$, (which depends on the group $G$).  The vector is used to
form a circulant matrix $M$. The group acting on the matrix $M$
produces several matrices which are concatenated to form a covering
array.  Often it will be necessary to add a small matrix, $C$,
to complete the covering conditions.

This group construction extends the method of Chateauneuf, Colbourn
and Kreher~\cite{chateauneuf:99}. This method has several parameters
that depend on each other; the group $G$, the vector $v$ and the small
array $C$. We first choose a specific group $G$ which will determine
$C$ and then, knowing $G$, we find $v$ by a computer search. To
motivate the choice for $G$, an example of how this method works is
given.

\begin{exmp}\label{ex:cyclic}
Let 
\[
G = \{ e, (12) \} < Sym(3)  \textrm{ and } v= (0,1,1,1,2) \in \mathbb{Z}_3^5.
\]

Build the following circulant matrix taking $v$ as the first column,
\[
M = 
\left[
\begin{array}{ccccc}
0 & 2 & 1 & 1 & 1 \\
1 & 0 & 2 & 1 & 1 \\
1 & 1 & 0 & 2 & 1 \\
1 & 1 & 1 & 0 & 2 \\
2 & 1 & 1 & 1 & 0 \\
\end{array}
\right].
\]

The elements of $G$ acting on the matrix $M$ produce
\[
\begin{array}{cc}
M_e = 
\left[
\begin{array}{ccccc}
0 & 2 & 1 & 1 & 1 \\
1 & 0 & 2 & 1 & 1 \\
1 & 1 & 0 & 2 & 1 \\
1 & 1 & 1 & 0 & 2 \\
2 & 1 & 1 & 1 & 0 \\
\end{array}
\right]
& \textrm{ and~~}
M_{(12)} = 
\left[
\begin{array}{ccccc}
0 & 1 & 2 & 2 & 2 \\
2 & 0 & 1 & 2 & 2 \\
2 & 2 & 0 & 1 & 2 \\
2 & 2 & 2 & 0 & 1 \\
1 & 2 & 2 & 2 & 0 \\
\end{array}
\right].
\end{array}
\]
The vector $ C = \left[ 0\;0\;0\;0\;0 \right]^\top$ is needed to
ensure the coverage of all pairs.

From this, a $CA(11,5,3)$ is built by concatenating the matrices $C$,
$M_e$ and $M_{(1,2)}$
\begin{displaymath}
\left[
\begin{array}{c|ccccc|ccccc}
0 & 0 & 2 & 1 & 1 & 1 & 0 & 1 & 2 & 2 & 2 \\
0 & 1 & 0 & 2 & 1 & 1 & 2 & 0 & 1 & 2 & 2 \\
0 & 1 & 1 & 0 & 2 & 1 & 2 & 2 & 0 & 1 & 2 \\
0 & 1 & 1 & 1 & 0 & 2 & 2 & 2 & 2 & 0 & 1 \\
0 & 2 & 1 & 1 & 1 & 0 & 1 & 2 & 2 & 2 & 0 \\
\end{array}
\right].
\end{displaymath}

\end{exmp}

Consider the group action of $G$ on the pairs from $\mathbb{Z}_k$, let
$\mathcal{O}$ be the set of orbits from this action.  For the group
action of $G$ on the matrix $M$ to produce a covering array, it is
necessary that for any orbit $O \in \mathcal{O}$, and for any two
rows $r_0,r_1$ of $M$, there must be at least one column $c$ of $M$
such that $( M_{r_0,c},M_{r_1,c}) \in O$.

To determine which vectors $v=(v_0,v_1, \dots,v_{r-1})$ are starter vectors, consider the sets $d_{i}$,
for $i = 1,\dots,r-1$,
\begin{equation}
d_{i} = \{ ( v_{j}, v_{j+i} )  |\;  j = 0, 1, \dots, r-1 \},
\end{equation}
where the subscripts are taken modulo $r$.

For $v$ to be a starter vector,\index{starter vector} each set $d_{i}$
(for $i = 1, \dots, r-1$) must contain at least one element from each
orbit of the group action of $G$ on the ordered pairs from $\mathbb{Z}_k$.

\remove{\begin{lemma}
If a starter vector $v$ exists in $\mathbb{Z}_k^r$ with respect to a
group $G$ then there exists a $CA( r|G|,r,k)$.
\end{lemma}}

To produce small covering arrays, the group $G$ should be small and
have few orbits.  These two properties of the group action are
connected by Burnside's Theorem (also known as Cauchy-Frobenius
Theorem).  This theorem states that if a group $G$ acting on a set $X$
has $N$ orbits, then
\begin{equation}
N = \frac{1}{|G|} \sum_{g \in G} \fix(g)   \nonumber
\end{equation}
where $\fix(g)$ is the number of elements $ x \in X$ for which the 
action of $g \in G$ fixes $x$ (this means $g(x)=x$). 

For this construction, it is best to use a group that has $ N |G| =
\sum_{g \in G} \fix(g)$ small. The group used in both~\cite{karen2}
and~\cite{groupconstCA} is the $(k-1)$-element group $G= \langle (1\;2\;
\dots \; k-1) \rangle < Sym(k)$. This group was chosen because each
element, other than the identity, fixes exactly one pair from
$\mathbb{Z}_k$, the pair $(0,0)$.

Consider the action of the subgroup $G= \langle (1\;2\; \dots \; k-1)
\rangle < Sym(k)$ on the pairs from $\mathbb{Z}_k$.  There are $k+2$
orbits generated by this group action:
\begin{enumerate}
\item $\{(0,0)\}$, 
\item $\{ (0,i) \; : \; i = 1,\dots , k-1\}$, 
\item $\{ (i,0) \; : \; i = 1,\dots , k-1\}$,
\item $\left\{ (a,b) \; : \; a,b \in \mathbb{Z}_k \backslash \{0\} \mbox{ and } a-b \equiv j \pmod {k-1}\right\}$ 
                   for $j = 0,\dots, k-2$.
\end{enumerate}

Consider a  vector $v=(0,v_1,v_2, \dots ,v_{r-1})$ from $\mathbb{Z}_k^r$
where only the first entry of the vector is 0.
Then consider the sets of pairs
\[
d_{i} = \{ ( v_{j}, v_{j+i} )  |\;  1 \leq j \leq  r-1  \mbox{ and } i+j \neq 0 \}
\mbox{ for } 1 \leq i \leq r-1,
\]
(where the subscripts are taken modulo $r$). Let $v$ be a vector in
$\mathbb{Z}^r_k$ with the property that the set of differences $v_j -
v_{j+i}$ for all $(v_{j}, v_{j+i}) \in d_{i}$ covers
$\mathbb{Z}_{k-1}$ for every $i = 1, \dots, r-1$. Then the set 
\[
d^\ast_{i} = \{ ( v_{j}, v_{j+i} )  |\;  j = 0, \dots, r-1 \},
                 \quad i = 1, \dots, r-1
\]
intersects all the orbits except $\{(0,0)\}$.

Set $M$ to be the $r \times r$ circulant matrix generated from the
vector $v$.  Then, for any pair of rows in $M$, at least one element
from each of the orbits of $G$ (except $\{(0,0)\}$) occurs in some column.
For each $g \in G$, let $M_g$ be the matrix formed by the action of
$g$ on the elements of $M$. As any two rows in $M$ have a
representative from each orbit of $G$, for any two rows all pairs from
$\mathbb{Z}_k$ (except $\{(0,0)\}$) will occur in $M_g$ for some $g \in
G$.  Finally let $C$ be the $r \times 1$ matrix with all entries equal
to 0.  Then the array formed by concatenating $C$ and $M_g$, for all $g
\in G$, is a covering array.

Starter vectors exist for many values of $r$ and $k$.
Table~\ref{results1:tab} is an exhaustive list of all values of $r$
and $k$ with $k \leq 10$ with $k+1 \leq r\leq 2(k+1)$ for which a
starter vector exists for the given group $G$.  These were found by an
exhaustive search for all values with $k = 3, \dots, 10$ and all $r$
with $k+1 \leq r \leq 2(k+1)$; if values in this range are not listed
in Table~\ref{results1:tab} then a starter vector does not exist for
those values.

\begin{table}[ht]
\begin{displaymath}
\begin{tabular}{|l|r|}\hline
$k$ & $r$  \\
\hline
3 & 5,8 \\
4 & 5--10 \\
5 & 7--12 \\
6 & 9--14 \\
7 & 10--16 \\
8 & 9, 11--18 \\
9 & 13--20 \\
10 & 15--22\\\hline
\end{tabular}
\end{displaymath}
\caption{Exhaustive list of values for which a starter vector exists with $k\leq 10$ and $k+1 \leq r \leq 2(k+1)$.
}% \label{results1:tab}}
 \label{results1:tab}
\end{table}

\begin{table}[ht]
\begin{displaymath}
\begin{tabular}{|l|r|}\hline
$k$ & $r$  \\
\hline
11 & 17 \\
12 & 18--19 \\
13 & 21--27 \\
14 & 23--29 \\
15 & 27--31 \\
16 & 29--33 \\
17 & 33--35 \\
18 & 35--37\\\hline
\end{tabular}
\end{displaymath}
\caption{A list of values (not exhaustive) for which a starter vector exists with $11 \leq k \leq 18$ and $k+1 \leq r \leq 2(k+1)$.
 \label{results2:tab}}
\end{table}

Many of these starter vectors give an upper bound on the size of a
covering array which is better than the previous best known bound. For
example, when $k=6$ and $r=9$, a $CA(46, 9, 6)$ can be constructed
with this method. Previously, the smallest known covering array with
$r=9$ and $k=6$ was a $CA(48,9,6)$.  A table of the new upper bounds
obtained with this construction is given in Appendix~\ref{appendix}.

For values of $k\geq 11$ the exhaustive search was found to be
inefficient. Replacing the exhaustive search, a simple hill-climbing
algorithm has improved bounds for covering arrays with $11 \leq k
\leq 20$ (these bounds are in~\cite{karen2}).  These new upper bounds are listed in
also Appendix~\ref{appendix}. Table~\ref{results2:tab} gives an
non-exhaustive list of values for which a starter vector exists with
$11 \leq k \leq 18$ and $k+1 \leq r \leq 2(k+1)$.

The covering arrays that we construct with the group construction using our group $G$
have the property that the $k$ columns that have 0 in the first row
are mutually disjoint, except for the first row. Removing the first
row from the covering arrays, $CA(n,r,k)$, that we constructed with the group
construction method produces a $CA(n,r-1,k)$ with $k$ disjoint
columns.  This method not only produces small covering arrays, it also
produces arrays that can be used in the block-size recursive
construction.

For example, if $k=4$, from Lemma~\ref{fullblockrec}, it is possible
to construct a $CA(16+12i, 5(4^i), 4)$, for all positive integers $i$.
From the group construction, there exist covering arrays with four
disjoint columns with the following parameters: $CA(19,5,4)$,
$CA(22,6,4)$, $CA(25,7,4)$, $CA(28,8,4)$ and $CA(31,9,4)$.  Using
these covering arrays with the block-size recursive construction and
the $CA(16,5,4)$, from Example~\ref{ffcaexample},
Section~\ref{sec:finitefieldconstruction}, it is possible to construct
the following covering arrays for all positive integers $i$:
$CA(16+22i,5(6^i),4)$, $CA(16+25i,5(7^i),4)$, $CA(16+28i,5(8^i),4)$
and $CA(16+31i,5(9^i),4)$.

\remove{This provides many new bounds, but the asymptotic growth of $N(n,k)$ 
relative to $n$ is the largest for $CA(16+12i, 5(4^i), 4)$. This means
these new bounds are useful, but do not give new information about the asymptotic
growth of balanced covering arrays.
}

%%%%%%%%%%%%%%%%%%%%%%%%%%%%%%%%%%%%%%%%%%%%%%%%%%%%%%%%%%%%%%%%%%
\subsection{Asymptotic Results}\label{asymptotics}

With the results we have seen so far, it is possible to show one bound
on the asymptotic growth of $N(n,k)$. From Lemma~\ref{fullblockrec},
for $k$ a prime power and $i$ a positive integer, there exists a
covering array $CA(k^2+i(k^2-k), k^i(k+1),k)$.  For any integer $i$,
set $n = k^2+i(k^2-k)$, then,
\[
N(n,k) \geq k^i(k+1) > k^{i+1} = k^{(n-k)/(k^2-k)}.
\]

To motivate the next result, a short digression on Shannon theory is
needed. Shannon theory was developed by C. Shannon~\cite{MR0089131}
and has resulted in many interesting and difficult combinatorial
problems. One of the goals of this theory was to develop methods to
identify encoded messages that may have some of their bits changed. To
do this, we need to be able to identify when two messages are ``really
different'', and not just the same message with some bits changed.

For a graph $G$ with vertex set $V$, two
vectors $v,w \in V^n$ are considered ``really different'' if for some
$i \in \{1,\dots,n\}$, $(v_i,w_i) \in E(G)$, where $E(G)$ is the edge
set of $G$.  Denote by $N(G,n)$ the size of the largest set of vectors
from $V^n$ that are pairwise ``really different.''  Then, the {\em
zero error capacity}\index{zero error capacity} of a graph $G$ is
defined to be
\[
C(G) = \limsup_{n \rightarrow \infty} \frac{\log_2 N(G,n)}{n}.
\]

We consider the same asymptotic growth of $N(n,k)$, that is,
\[
\limsup_{n \rightarrow \infty} \frac{\log_2 N(n,k)}{n}.
\] When $k$ is a prime power,
\begin{eqnarray*}
\limsup_{n \rightarrow \infty} \frac{\log_2 N(n,k)}{n}
 &\geq& \limsup_{n \rightarrow \infty} \frac{\log_2( k^{(n-k)/(k^2-k)})}{n} \\
 &=& \limsup_{n \rightarrow \infty} \frac{ n \log_2 k-k\log_2 k  }{n (k^2-k) }\\
& = & \frac{\log_2 k}{k^2-k}.
\end{eqnarray*}
Poljak and Tuza~\cite{poljak:89} show that for a prime power $k$, we have the following bounds:
\begin{eqnarray}\label{eq:poljak}
 \frac{\log_2 k}{k^2-k} \leq \limsup_{n \rightarrow \infty} \frac{\log_2 N(n,k)}{n} \leq \frac{2}{k}.
\end{eqnarray}
In fact, the upper bound holds for all $k$ (this is shown in Section~\ref{asympbollobas}).

Gargano, K\"orner and Vaccaro show that 
\begin{eqnarray}\label{eq:asymptotic}
\limsup_{n \rightarrow \infty} \frac{\log_2 N(n,k)}{n} = \frac{2}{k}. 
\end{eqnarray}

The proof of Equation~(\ref{eq:asymptotic}) is very complex,
non-constructive and is contained in the three
papers~\cite{gargano:92, MR94e:05024, gargano:94}.

%%%%%%%%%%%%%%%%%%%%%%%%%%%%%%%%%%%%
% Extremal Set Theory

\chapter{Extremal Set Theory}\label{EST}

One of the central problems in extremal set theory is the problem of
finding a system of sets with the largest cardinality given some
restriction on the sets of the system. Usually the restriction is of
the form that any two distinct sets in the system satisfy some
constraint.

In this chapter, two such problems are discussed.  The problem
considered in Section~\ref{sec:sperner} is to find the maximum
cardinality of a set system, over a finite ground set, satisfying the
constraint that any two distinct sets in the system are incomparable
(one is not contained in the other). The major result for this type of
restriction is Sperner's Theorem. The problem that we describe in
Section~\ref{sec:EKR} is to find the set system with the largest
cardinality satisfying the constraint that any two distinct sets in
the system are intersecting. The main result for this type of
problem is the Erd\H{o}s-Ko-Rado Theorem.

Some design theory problems can be rephrased as extremal set theory
problems. For example, the minimal size of a binary covering array can
be found by using Sperner's Theorem and the Erd\H{o}s-Ko-Rado Theorem.
General covering arrays are not related to set
systems, instead they correspond to systems of partitions. In
Chapter~\ref{partitionsystems}, a theory of extremal partitions is
introduced, and several of the results presented in this chapter are
extended to include partition systems; such extensions are also known
as {\em higher order} extremal problems.

%%%%%%%%%%%%%%%%%%%%%%%%%%%%%%%%%%%%%%%%%%%%%%%%%%%%%%%%%%%%%%%%%%%%%%%%%%%%%%
\section{Set Systems}\label{setsystems}

For $n$ a positive integer, let $X = \{1,2,\dots, n\}$ be an
$n$-set\index{$n$-set}.  The {\em power set} of $X$, denoted $P(X)$,
is the collection of all subsets of $X$. A {\em set system}\index{set
system} on an $n$-set is a collection of sets from $P(X)$.  For a
positive integer $k \leq n$, a {\em $k$-set}\index{$k$-set} is a set
$A \in P(X)$ with $|A|=k$. The collection of all $k$-sets of an
$n$-set is denoted by ${[n]\choose k}$. A {\em $k$-uniform set
system}\index{$k$-uniform set system} on an $n$-set is a collection of
sets from ${[n]\choose k}$. It is possible to arrange $P(X)$ in a
poset ordered by inclusion. In this poset, for every positive
integer $k\leq n$, the sets ${[n]\choose k}$ are the {\em level
sets}\index{level sets}.  In this section, three types of set systems
are defined: $t$-designs, set partitions and hypergraphs.

\subsection{$t$-Designs}\label{sec:blockdesigns}

In this section, we introduce $t$-designs and give some basic
results. For more on $t$-designs, see~\cite{beth:85,
stinson:02, street:87}.

\begin{defn}[$t$-$(n,k,\lambda)$ Design]\index{$t$-$(n,k,\lambda)$ design}
Let $t,n,k$ and $\lambda$ be positive integers such that $2\leq t \leq
k < n$.  A {\em $t$-$(n,k,\lambda)$ design}\index{$t$-$(n,k,\lambda)$
design} is a $k$-uniform set system $\mathcal{A}$ on an $n$-set with
the property that any $t$-set of the $n$-set is contained in exactly
$\lambda$ sets in $\mathcal{A}$.
\end{defn}

A $t$-$(n,k,\lambda)$ design with $t=2$ is also known as a
$(n,k,\lambda)$-BIBD\index{ $(n,k,\lambda)$-BIBD} (BIBD stands for
{\em balanced incomplete block design}\index{balanced incomplete block
design}). The sets in a design are also called the {\em blocks} of the
design.

\begin{exmp}
The most famous $t$-design is the ubiquitous Fano plane. The Fano
plane\index{Fano plane} is a 2-$(7,3,1)$ design (and thus a $(7,3,1)$-BIBD). The Fano
plane is represented in the picture below. The vertices are the
elements of the $7$-set and each of the blocks is represented by a line in the plane.

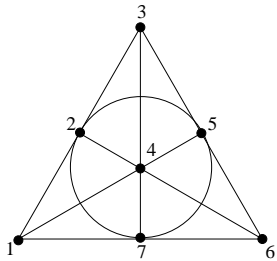
\centerline{ \psfig{figure=fano.eps,height=3.5cm}}\hfill

That is, the set system $\mathcal{A}$ below is a 2-$(7,3,1)$ design,
\[
\mathcal{A} = \{ \{123\},\{145\},\{167\},\{246\},\{257\},\{347\},\{356\}\}.
\]
\end{exmp}

\begin{thm}[see \cite{stinson:02}]
In a $t$-$(n,k,\lambda)$ design there are exactly
\[
b = \frac{\lambda {n \choose t}}{{k \choose t}}
\]
blocks.
Moreover, every element of the $n$-set is contained in exactly
\[
r = \frac{kb}{n} = \frac{ \lambda {n-1 \choose t-1} }{{k-1 \choose t-1}}
\]
sets.
\end{thm}

Clearly, if a $t$-$(n,k,\lambda)$ design exists, then $b$ and $r$ in
the above theorem must be integers. This gives a necessary condition
for the existence of a $t$-$(n,k,\lambda)$ design.

\remove{For the special case of a
$(n,3,\lambda)$-BIBD, if $b$ and $r$ (as defined above) are integer then a
$(n,3,\lambda)$-BIBD exists.  A $(n,3,1)$-BIBD is also called a {\em
Steiner triple system}\index{Steiner triple system} and is denoted by
$STS(n)$.}

\begin{defn}[Resolvable Design]\index{resolvable $t$-$(n,k,\lambda)$ design}\label{defn:resolvable}
Suppose $\mathcal{A}$ is a $t$-$(n,k,\lambda)$ design.  A {\em
parallel class}\index{parallel class in a design} in $\mathcal{A}$ is a collection
of pairwise disjoint sets from $\mathcal{A}$ whose union is the entire
$n$-set. A partition of $\mathcal{A}$ into $r = \frac{\lambda {n-1
\choose t-1}}{{k-1 \choose t-1}}$ parallel classes is called a {\em
resolution}\index{resolution}. A $t$-$(n,k,\lambda)$ design is {\em
resolvable} if a resolution exists.
\end{defn}

For many values of $t, n,k, \lambda$, it is not possible to construct a
$t$-$(n,k,\lambda)$ design. On the other hand, it is always possible to construct a
design called a {\em packing}.

\begin{defn}[$t$-$(n,k, \lambda)$ packing]\index{$t$-$(n,k, \lambda)$ packing}\label{defn:packing}
Let $t,n,k$ and $\lambda$ be positive integers such that $2\leq t \leq
k < n$. A {\em $t$-$(n,k,\lambda)$ packing} is a $k$-uniform set
system $\mathcal{A}$ on an $n$-set with the property that any $t$-set
of the $n$-set is contained in at most $\lambda$ sets in
$\mathcal{A}$.
\end{defn}

The definition of a {\em resolvable packing} follows.

\begin{defn}[Resolvable Packing]\index{resolvable $t$-$(n,k,\lambda)$ packing}\label{defn:resolvablepacking}
Suppose $\mathcal{A}$ is a $t$-$(n,k,\lambda)$ packing.  A {\em
parallel class}\index{parallel class in a packing} in $\mathcal{A}$ is a collection
of pairwise disjoint sets from $\mathcal{A}$ whose union is the entire
$n$-set. 
A partition of $\mathcal{A}$ into parallel classes is called a {\em
resolution}\index{resolution}. A $t$-$(n,k,\lambda)$ packing is {\em
resolvable} if a resolution exists.
\end{defn}

\subsection{Set Partitions}\label{setpartitions}

A {\em set partition}\index{set partition} of an $n$-set is a set of
disjoint non-empty subsets (called classes) of the $n$-set whose union
is the $n$-set.  Throughout this thesis, set partitions will be
referred to simply as {\em partitions}\index{partition}.  A partition
$P$ is called a {\em $k$-partition} if it contains $k$ classes, that
is $P=\{P_1,P_2, \dots , P_k\}$, or, equivalently, $|P|=k$.  For
positive integers $k,n$, let $\mathcal{P}^n_k$ denote the set of all
$k$-partitions of an $n$-set. Let $S(n,k)=|\mathcal{P}^n_k|$; the
values $S(n,k)$ are the {\em Stirling numbers of the second type}.
\index{Stirling numbers of the second type}

A partition $P \in \mathcal{P}^n_k$ is {\em uniform}
\index{uniform partition} if every class $P_i \in P$, $i=1,\dots ,k$, has the same
cardinality, that is, $|P_i|=n/k$ for all $P_i\in P$.  For positive
integers $n,c,k$ with $n=ck$, denote by $\mathcal{U}^n_k$ the set of
all uniform $k$-partitions in $\mathcal{P}^n_k$. Analogous to the Stirling
numbers of the second type, denote $U(n,k)=|\mathcal{U}^n_k|$, that is, for $n=ck$,
\[U(n,k)= \frac{1}{k!} {n \choose c} {n-c \choose c}
\cdots {n-(k-1)c \choose c}.\]

If $k$ does not divide $n$, it is not possible for a partition in
$\mathcal{P}^n_k$ to be uniform, in this case {\em almost-uniform
partitions} are considered. For positive integers $n,k,c$ with $n=ck+r$ where $0\leq r < k$, a
partition $P \in \mathcal{P}^n_k$ is {\em almost-uniform}
\index{almost-uniform partition} if every class of $P$ has cardinality
$c$ or $c+1$. Note that in an almost-uniform partition, there are $r$
classes of cardinality $c+1$ and $k-r$ classes of cardinality $c$.
Denote by $\mathcal{AU}^n_k$ the set of all almost-uniform partitions
in $\mathcal{P}^n_k$.  Denote $AU(n,k)=|\mathcal{AU}^n_k|$. For
$n=ck+r$,
\begin{eqnarray*}
AU(n,k) &=& \frac{1}{r!(k-r)!} {n \choose c} {n-c \choose c} \cdots {n-(k-r-1)c \choose c}  \\
         &&  {n-(k-r)c \choose c+1}{n-(k-r)c-(c+1) \choose c+1} \cdots {c+1 \choose c+1}.
\end{eqnarray*}

If $k$ divides $n$, then $\mathcal{AU}^n_k=\mathcal{U}^n_k$ and $U(n,k) = AU(n,k)$.

\subsection{Hypergraphs}

Hypergraphs are generalizations of graphs, these will be used in
Section~\ref{sec:spernerpartitions}. A {\em
hypergraph}\index{hypergraph} on $n$ vertices is a collection of
subsets of an $n$-set. The elements of the $n$-set are the {\em
vertices}\index{vertices of a hypergraph} of the hypergraph and the
subsets are the {\em edges} of the hypergraph\index{edges of a
hypergraph}. A hypergraph is {\em uniform}\index{uniform hypergraph}
if all the edges in the hypergraph have the same cardinality.  A
uniform hypergraph with all edges of cardinality two is a graph.  For
positive integers $n$ and $c$ with $c \leq n$, the {\em complete
$c$-uniform hypergraph}\index{complete uniform hypergraph}, denoted
$K_{n}^{(c)}$, has vertex set $\{1,2,\dots,n\}$ and edge set all
$c$-sets from the $n$-set, that is, ${[n] \choose c}$.

The {\em degree} of a vertex in a hypergraph is the number of edges
that contain the vertex. If every vertex in a hypergraph has degree $d$, then the
hypergraph is {\em $d$-regular}\index{$d$-regular hypergraph}.

A {\em 1-factor}\index{1-factor of a hypergraph} of a hypergraph on an
$n$-set is a set of disjoint edges whose union is the entire
$n$-set. A {\em 1-factorization}\index{1-factorization} of a
hypergraph is a partition of its edges into a collection
of pairwise disjoint 1-factors. The {\em size of a
1-factorization}\index{size of a 1-factorization} is the number of
1-factors in the 1-factorization. For positive integers $c,n$, if $c$
divides $n$, then a 1-factorization of $K_n^{(c)}$ is possible.

\begin{thm}[1-factorization of hypergraphs~\cite{MR88g:05001}]\label{thm:1factorization}
Let $c,n$ be positive integers such that $c$ divides $n$, then the
hypergraph $K_n^{ (c) }$ has a 1-factorization.
\end{thm}

Let $c,k,n$ be positive integers such that $n=ck$, the size of a
1-factorization of $K_{n}^{(c)}$ is $\frac{1}{k}{n
\choose c} = {n-1 \choose c-1}$. 

Assume for positive integers $k,c,n$ that $n=ck$. A 1-factor of the
complete uniform hypergraph $K_n^{(c)}$ is equivalent to a uniform
$k$-partition of an $n$-set. Further, for any $t \leq c$, a
1-factorization of the complete hypergraph $K_{n}^{(c)}$ corresponds
to a resolvable $t$-$(n,c,\lambda)$ design with $\lambda = {n-t
\choose c-t}$.  In particular, it is a resolvable $(n,c,\lambda)$-BIBD, where
$\lambda = {n-2 \choose c-2}$.

\section{Sperner Theory}\label{sec:sperner}

An important class of problems in extremal set theory deals with the
maximum cardinality of a set system with some restriction on the sets
in the system.  The first restriction that we consider is that any two
distinct sets from the system must be {\em incomparable}.  Two subsets
$A$ and $B$ of an $n$-set are {\em comparable} if $A
\subseteq B$ or $B\subseteq A$.  If $A$ and $B$ are not comparable
then they are {\em incomparable}.

\begin{defn}[Sperner Set System]\index{Sperner set system}\label{defn:sperner}
Let $n$ be a positive integer.  A {\em Sperner set system}
$\mathcal{A}$ on an $n$-set is a set system on an $n$-set with the
\orphan
property that any two distinct sets in $\mathcal{A}$ are incomparable.
\end{defn}

In 1928, Sperner proved that a Sperner set system cannot have more
sets than are in the largest level set~\cite{sperner:28}. We give the
proof of this here since a similar proof is used in Lemma~\ref{thm:kuniform}
and Theorem~\ref{thm:QIchi}.

\begin{thm}[Sperner's Theorem~\cite{sperner:28}] \label{thm:sperner}
Let $n$ be a positive integer. If $\mathcal{A}$ is a Sperner set
system on an $n$-set, then $|\mathcal{A}| \leq {n \choose \lfloor
\frac{n}{2}\rfloor}$. 
\end{thm}
\begin{proof}
A chain in a poset is a collection of sets in the poset with the
property that any two distinct sets in the chain are ordered in the
poset.  We will show that the poset formed by $P(\{1,2,\dots,n\})$
ordered by inclusion can be
decomposed into ${n \choose \lfloor
\frac{n}{2}\rfloor}$ disjoint chains. Any Sperner system can intersect
such a chain in at most one set, and thus, has cardinality no more than ${n \choose
\lfloor \frac{n}{2}\rfloor}$.

Let $r$ be a positive integer with $r < \lfloor \frac{n}{2}\rfloor$.
We will show that every $r$-set can be matched to an $(r+1)$-set such
that the $r$-set is a subset of the $(r+1)$-set.  Construct a
bipartite graph as follows.  For every $r$-set of the $n$-set there is
a corresponding vertex in the first part of the graph.  For each
$(r+1)$-set in the $n$-set there is a corresponding vertex in the
second part of the graph. Two vertices are adjacent in this bipartite
graph if and only if one of the corresponding sets is contained in the
other set. All the vertices in the first part of the graph have degree
$n-r$ and the vertices in the second part all have degree $r+1$.  Let
$S$ be a set of vertices from the first part of the graph. Let $N(S)$
be the set of vertices in the second part of the graph adjacent to any
vertex in $S$. Since every vertex in the first part of the graph has
degree $n-r$ and every vertex in the second part of the graph has
degree $r+1$, counting the number of edges between $S$ and $N(S)$ we
have that $|S|(n-r) \leq |N(S)|(r+1)$. Since $\frac{r+1}{n-r} \leq 1$,
we have that $|S| \leq |N(S)|$ and by Hall's Theorem, there is a
one-to-one matching from the first part of this graph to the second
part.  Thus, there is a one-to-one matching from ${ [n] \choose r}$ to
${ [n] \choose r+1}$ for $r \leq \lfloor \frac{n}{2}\rfloor$.

Similarly, for any positive integer $r$ with $r > \lfloor
\frac{n}{2}\rfloor$, there is a one-to-one matching from ${ [n]
\choose r}$ to ${ [n] \choose r-1}$.
Finally, by the same argument, if $n$ is odd, there is one-to-one matching from ${ [n]
\choose \frac{n-1}{2}}$ to ${ [n] \choose \frac{n+1}{2}}$.

Two sets are in the same chain if they are matched in one of these
matchings. Then these matching define ${n \choose \lfloor
\frac{n}{2}\rfloor}$ disjoint chains which partition the poset, since
\halforphan
each chain has exactly one set in the set ${ [n] \choose \lfloor
\frac{n}{2} \rfloor}$.
\end{proof}
 
Moreover, $|\mathcal{A}| = {n \choose \lfloor
\frac{n}{2}\rfloor}$ if and only if $\mathcal{A} = { [n] \choose k}$
where $k= \lfloor \frac{n}{2} \rfloor$ or $k=  \lceil \frac{n}{2} \rceil$.

A sharper inequality is the {\em LYM Inequality} named after
Lubell~\cite{MR33:2558}, Meshalkin~\cite{MR27:52} and
Yamamoto~\cite{MR16:668c}, who each independently established the
result. 

\begin{thm}[LYM Inequality]\label{theorem:lym}
Let $n$ be a positive integer and $\mathcal{A}$ be a Sperner set
system on an $n$-set. Then
\[
\sum_{A \in \mathcal{A}} {n \choose |A|}^{-1} \leq 1.
\]
\end{thm}
\remove{
\begin{proof}
Let $\mathcal{A}$ be a Sperner set system.  Let $\sigma = (a_1,a_2,\dots
a_n)$ be a permutation of the $n$-set, a $k$-set $A$ is {\em
contained} in $\sigma$ if $A= \{a_1,a_2, \dots , a_k\}$. A $k$-set 
$A \in \mathcal{A}$ is contained in exactly $k!(n-k)!$ permutations. Since
$\mathcal{A}$ is a Sperner set system, each permutation can contain at most
one set from $\mathcal{A}$. There are $n!$ permutations of an $n$-set,
\[
\sum_{A \in \mathcal{A}} |A|!(n-|A|)! \leq n!
\]
and the result follows.
\end{proof}
}

The LYM Inequality implies Sperner's Theorem. Since for all $k \leq
n$, ${n \choose k} \leq {n \choose \lceil \frac{n}{2} \rceil}$, we have that
for a Sperner set system $\mathcal{A}$ on an $n$-set
\[
|\mathcal{A}|{n \choose \lceil \frac{n}{2} \rceil}^{-1} 
        \leq \sum_{A \in \mathcal{A}} {n \choose |A|}^{-1} \leq 1.
\]

The LYM Inequality can be rearranged into another form. Let $p_i$ denote the
number of subsets in $\mathcal{A}$ of size $i$, then
\begin{eqnarray}\label{equationLYM}
\sum_{i=1}^n \frac{p_i}{ {n \choose i} } \leq 1.
\end{eqnarray}

\subsection{Bollob\'as's Theorem}

An important theorem on set systems is Bollob\'as's Theorem. This theorem
implies Sperner's Theorem and the LYM Inequality. Further, it is an
example of a {\em higher order} Sperner Theorem, meaning that it gives
a Sperner-like result for a system of families of sets rather than systems of sets.

\begin{thm}[Bollob\'as's Theorem~\cite{MR32:1133}]\label{bollobas}
Let $A_1,A_2,\dots ,A_m$ and $B_1,B_2,\dots , B_m$ be
sequences of sets such that $A_i \cap B_j = \emptyset$ if and only if $i =j$.
Then
\[
\sum_{i=1}^m {a_i +b_i \choose a_i}^{-1} \leq 1,
\]
where $a_i = |A_i|$ and $b_i = |B_i|$.
\end{thm}

In~\cite{MR34:5690} it is noted that the natural extension of this
theorem to three families of sets is not true.  In
Section~\ref{sec:spernerpartitions}, we give a different higher order
version of Sperner's Theorem.

\label{asympbollobas} 

Recall that for given positive integers $n$ and $k$, $N(n,k)$ is the
largest $r$ such that a $CA(n,r,k)$ exists. Bollob\'as's Theorem was
used by Stevens, Moura and Mendelsohn~\cite{brett:97b} and by Poljak
and Tuza~\cite{poljak:89} (Equation~(\ref{eq:poljak}),
Section~\ref{asymptotics}) to find an upper bound on $N(n,k)$.

Let $C$ be a $CA(n,r,k)$. For each $i$, $1 \leq i \leq r$, let $k_i^1$
be a letter which occurs least often in row $i$ of $C$ and $k_i^2$ be
a letter which occurs second least often in row $i$ of $C$, where ties
are broken arbitrarily. Define two subsets $A^1_i$ and $A^2_i$ by $j
\in A^1_i$ if and only if $k_i^1$ occurs in column $j$ and
$j \in A^2_i$ if and only if $k_i^2$ occurs in column $j$.

Define sets $A_i$ and $B_i$ for $i=1,\dots, 2r$ by $A_i = A^1_i$ for
$1\leq i \leq r$ and $A_i = A^2_{i-r}$ for $r+1\leq i \leq 2r$ and
$B_i = A^2_i$ for $1\leq i \leq r$ and $B_i = A^1_{i-r}$ for $r+1\leq
i \leq 2r$.

For all $i\in \{1,\dots,2r\}$, $|A_i| \leq \lfloor n/k \rfloor$,
$|A_i| + |B_i| \leq \lfloor 2n/k \rfloor$ and $A_i \cap B_i = \emptyset$.
For all distinct $i,j \in \{ 1,\dots,2r\}$, $A_i \cap B_j \neq
\emptyset$.  By Bollob\'as's Theorem, we have $2r \leq {\lfloor{2n/k}\rfloor
\choose \lfloor{n/k}\rfloor}$, and using Stirling's formula, $n! \sim
\sqrt{2 \pi n} \left( \frac{n}{e}\right)^n$, we obtain
\[
N(n,k) \leq \frac{1}{2}{ \lfloor{2n/k}\rfloor \choose \lfloor{n/k}\rfloor}
           = O \left(\frac{4^{n/k}}{n^{1/2}}\right).
\]

Then
\[
\limsup_{n \rightarrow \infty} \frac{\log_2 N(n,k)}{n}
      \leq \limsup_{n \rightarrow \infty} \frac{1}{n} \log_2 \left( \frac{4^{n/k}}{n^{1/2}}\right)
   =\limsup_{n \rightarrow \infty} \frac{1}{n}\left(\frac{2n}{k}-\frac{1}{2}\log_2{n}\right)
   = \frac{2}{k}.
\]

As stated in Section~\ref{asymptotics}, Gargano, K\"orner and
Vaccaro~\cite{gargano:92,MR94e:05024,gargano:94} show that for all $k \geq 2$, 
\[
\limsup_{n \rightarrow \infty} \frac{\log_2 N(n,k)}{n} = \frac{2}{k}.
\]

%%%%%%%%%%%%%%%%%%%%%%%%%%%%%%%%%%%%%%%%%%%%%%%%%%%%%%%%%%%%%%%%%%%%%%%%%%%%%%%%
\section{The Erd\H{o}s-Ko-Rado Theorem} \label{sec:EKR}

The Erd\H{o}s-Ko-Rado Theorem is a key result in extremal set
theory. It is a theorem that has many proofs, applications and
extensions. We state this theorem without a proof, a particularly nice
proof can be found in~\cite{MR46:3316}. There are many generalizations
of this theorem; two such results will be stated in
Section~\ref{sec:EKRgeneralizations}. In Section~\ref{sect:EKRapps}, we
show how the Erd\H{o}s-Ko-Rado Theorem, with Sperner's Theorem, can be
used to find the exact value of $CAN(r,2)$ for all $r$.

The Erd\H{o}s-Ko-Rado Theorem has been extended to combinatorial
objects other than sets. Rands in 1982~\cite{MR84i:05024} gives a
similar result for intersecting blocks of $t$-$(v,k,\lambda)$
designs. In Section~\ref{intersectingpartitions}, we give an extension
of the Erd\H{o}s-Ko-Rado Theorem to uniform partition systems.

\subsection{Intersecting Set Systems}\label{sec:ISS}

For fixed positive integers $t,k,n$, let $I(t,k,n)$ \index{$I(t,k,n)$} denote the
collection of all set systems $\mathcal{A}$ on an $n$-set with the
following properties: for all $A \in \mathcal{A}$, $|A| \leq k$; and
for all distinct $A,B \in \mathcal{A}$, $A \not \subset B$; for
all $A,B \in \mathcal{A}$, $|A \cap B| \geq t$.

The set systems in $I(t,k,n)$ are also known as {\em $t$-intersecting
set systems}\index{$t$-intersecting set system} and if $t=1$, these
are also called {\em intersecting set systems}\index{intersecting set
system}. The requirement that for all distinct sets $A,B$ in a system
that $A \not\subset B$ forces all set systems in $I(t,k,n)$ to be
Sperner set systems. If $2k-t \geq n $, then any two
$k$-sets from the $n$-set have at least $t$ elements in common.

\begin{defn}[$k$-Uniform $t$-Intersecting Set System]\index{$k$-uniform $t$-intersecting set system}
For positive integers $n,k,t$ a {\em $k$-uniform $t$-intersecting set
system} is a $k$-uniform set system, $\mathcal{A}$, on an $n$-set with the property
that for all distinct $A,B \in \mathcal{A}$, we have $A \not\subset B$; and
for all $A,B \in \mathcal{A}$, we have $|A \cap B| \geq t$.
\end{defn}

For $n$ sufficiently large, if there exists a $t$-intersecting set
system on an $n$-set, then there exists a $k$-uniform $t$-intersecting
set system that has at least the same cardinality. In particular, it
is possible to replace each set of size less than $k$ in a
$t$-intersecting set system on an $n$-set by a set of size $k$, which
is $t$-intersecting with all the sets in the system.  The proof of
this is not included since it is similar to the proof of Sperner's
Theorem.

\begin{lemma}[\cite{MR25:3839}]\label{thm:kuniform}
If $2k \leq n$, then for any $\mathcal{A} \in I(t,k,n)$ there exists a
$k$-uniform $t$-intersecting set system $\mathcal{A'} \in I(t,k,n)$
with $|\mathcal{A}| \leq |\mathcal{A'}|$.
\end{lemma}
\remove{
\begin{proof}
Let the size of the smallest set in $\mathcal{A}$ be $l$.  We will
prove this by induction on $k-l$. If $k-l=0$, then $\mathcal{A} \in
I(t,k,n)$ is a $k$-uniform set system and we are done.  

Assume $k-l >0$. Let $\mathcal{A}^l$ denote all the sets in
$\mathcal{A}$ of size $l$. Define
\[
\partial (\mathcal{A}^l) = \{ B\in {[n] \choose l+1} \; : \; A \subset B \mbox{ for some } A \in \mathcal{A}^l  \}.
\]
Then $|\partial (\mathcal{A}^l)| \geq |\mathcal{A}^l|(n-l)(l+1)^{-1}$, to
see, this consider the matching on the bipartite graph from the proof
of Sperner's Theorem, Theorem~\ref{thm:sperner}. 

Set $\mathcal{A'} = \mathcal{A} \backslash \mathcal{A}^l \cup \partial (\mathcal{A^l})$.
Then $\mathcal{A'} \in I(t,k,n)$ and
\[
|\mathcal{A'}| \geq |\mathcal{A}| - |\mathcal{A}^l| + |\mathcal{A}^l|\frac{n-l}{l+1}.
\]
Since $l < k \leq n/2$, $|\mathcal{A'}| \geq |\mathcal{A}|$.
\end{proof}
}

For positive integers $n,k,t$ with $t\leq k \leq n$, a set system on
an $n$-set is a {\em $k$-uniform trivially $t$-intersecting set
system}\index{$k$-uniform trivially $t$-intersecting set system} if it is equal,
up to a permutation on $\{1,\dots,n\}$, to
\[
\mathcal{A} =\left\{A \in {[n]\choose k} :   \{1,2,\dots, t\}\subseteq A\right\}. 
\]
The cardinality of a $k$-uniform trivially $t$-intersecting set system
is ${n-t \choose k-t}$.  If $t=1$ then a $k$-uniform trivially
$t$-intersecting set system is simply called a {\em $k$-uniform
trivially intersecting set system}.

%%%%%%%%%%%%%%%%%%%%%%%%%%%%%%%%%%%%% Proof of EKR
\subsection{The Erd\H{o}s-Ko-Rado Theorem}\label{sect:ERKproof}

The Erd\H{o}s-Ko-Rado Theorem proves that the set system in $I(t,k,n)$
with the largest cardinality, provided that $n$ is sufficiently large,
is a $k$-uniform trivially $t$-intersecting set system.
\enlargethispage{20pt}

\begin{thm}[Erd\H{o}s, Ko and Rado~\cite{MR25:3839}]\label{thm:ekr}
Let $k$ and $t$ be positive integers, with $0< t < k$.  There exists a
function $f(t,k)$ such that if $n$ is a positive integer with $n >
f(t,k)$, then for any $\mathcal{A} \in I(t,k,n)$
\[
|\mathcal{A}| \leq {n-t \choose k-t}.
\]
Moreover, equality holds if and only if $\mathcal{A}$ is a $k$-uniform trivially
$t$-intersecting set system.
\end{thm}

Erd\H{o}s, Ko and Rado, in 1961, proved that $f(t,k) \leq t+(k-t){k
\choose t}^3$ and stated this bound was certainly not the best possible.
In 1984, Wilson~\cite{MR86f:05007} improved this to the exact bound
$f(t,k) = (t+1)(k-t+1)$. It is interesting to note that Wilson's proof
is very algebraic; it uses the eigenvalues of the Johnson scheme (see
Example~\ref{exmp:johsonscheme}, Section~\ref{sec:introassscheme}) and
ratio bounds (see Section~\ref{sec:ratiobounds}).

\subsection{Generalizations of the Erd\H{o}s-Ko-Rado Theorem}\label{sec:EKRgeneralizations}

In this section, we state two extensions of the Erd\H{o}s-Ko-Rado
Theorem.  The first, by Hilton and Milner~\cite{MR36:2510}, gives the
largest non-trivial intersecting set system for $n \geq 2k$. For the proof of
Proposition~\ref{prop:HMpartitions}, we use a partition system similar
to the non-trivial set system given by Hilton and Milner. The second
extension of the Erd\H{o}s-Ko-Rado Theorem is by Ahlswede and
Khachatrian~\cite{MR97m:05251}.  In 1997, they gave a complete
description of all the maximal $k$-uniform $t$-intersecting set systems for all
values of $n,k$ and $t$.  Conjecture~\ref{conj:completeintersection}
in Section~\ref{sec:completetheorem} is a generalization of this
result by Ahlswede and Khachatrian to partition systems.

\begin{thm}[Hilton and Milner~\cite{MR36:2510}]
Let $k$ and $n$ be positive integers with $2 \leq k \leq
\frac{n}{2}$. Let $\mathcal{A} \in I(1,k,n)$ such that 
$\cap_{A \in \mathcal{A}} A = \emptyset$. 
Then
\begin{eqnarray*}\label{eq:HM}
|\mathcal{A}| \leq {n-1 \choose k-1} - {n-k-1 \choose k-1 }+1.
\end{eqnarray*}
\end{thm}

It is not difficult to construct an intersecting $k$-set system
$\mathcal{A}$ with cardinality $ {n-1 \choose k-1} - {n-k-1 \choose
k-1 }+1$ and $\cap_{A \in \mathcal{A}} A = \emptyset$. 
For $2 \leq k \leq \frac{n}{2}$, define a set system $\mathcal{A}'$ by
$A \in \mathcal{A}'$ if $1 \in A$ and $A \cap \{2, 3,\dots, k+1\} \neq
\emptyset$.  Define $\mathcal{A} = \mathcal{A}'\cup \{ \{2,3, \dots ,
k+1\}\}$. Then $|\mathcal{A}| = {n-1 \choose k-1} - {n-k-1 \choose k-1
}+1$ and $\cap_{A \in \mathcal{A}} A = \emptyset$.

\begin{thm}[Ahlswede and Khachatrian~\cite{MR97m:05251}]\label{thm:ahlswede}
Let $t,k,n$ be positive integers with $1 \leq t \leq k\leq n$ and $r$
be an integer with $r \leq k-t$. Define
\[
\mathcal{A}_r = \left\{A \in {[n] \choose k}  : |A \cap \{1,2 ,\dots, t+2r\} | \geq t+r \right\}.
\]
If
\[
(k-t+1) \left(2 + \frac{t-1}{r+1}\right) < n < (k-t+1)\left(2 + \frac{t-1}{r}\right),
\]
then $\mathcal{A}_r$ is the unique (up to a permutation on
$\{1,\dots,n\}$) system in $I(t,k,n)$ with maximal cardinality.  By
convention, $\frac{t-1}{r} = \infty$ for $r=0$.

If 
\[
n = (k-t+1)\left(2 + \frac{t-1}{r+1}\right),
\]
then $|\mathcal{A}_r|=|\mathcal{A}_{r+1}|$ and these two systems are
the unique (up to a permutation on $\{1,\dots,n\}$) systems with maximal
cardinality in $I(t,k,n)$.
\end{thm}

Note that $\mathcal{A}_0$ is a trivially $t$-intersecting set system.
Thus, when $r=0$, Theorem~\ref{thm:ahlswede} gives the Erd\H{o}s-Ko-Rado
Theorem with the exact lower bound for $n$.

\section{Application of the Erd\H{o}s-Ko-Rado Theorem}\label{sect:EKRapps}

Sperner's Theorem and the Erd\H{o}s-Ko-Rado Theorem can be used to
find the exact value of $CAN(r,2)$ for all $r$. In 1973, both
Katona~\cite{katona:73} and Kleitman and Spencer~\cite{kleitman:73}
gave similar proofs of this. The proof given here is closer to the
proof given by Kleitman and Spencer.

\subsection{Qualitatively Independent Subsets}\label{sec:qisets}

\begin{defn}[Qualitatively Independent Subsets]\label{defn:qi2}
\index{qualitatively independent subsets}
Two subsets $A$ and $B$ of an $n$-set are {\em qualitatively independent subsets} if 
\[
A \cap B \neq \emptyset, \quad
A \cap \overline{B} \neq \emptyset, \quad 
\overline{A} \cap B \neq \emptyset, \quad 
\overline{A} \cap \overline{B} \neq \emptyset. 
\]
\end{defn}

The term ``qualitative independence'' comes from probability
theory: if two sets $A$ and $B$ are qualitatively independent, then the
information that $x \in A$ gives no information about
whether or not $x \in B$.

Definition~\ref{defn:qi2} for qualitatively independent sets is
equivalent to Definition~\ref{defn:qi1} for binary vectors.  To see
this, let $A$ and $B$ be qualitatively independent subsets of an
$n$-set. Define the vector corresponding to $A$ to be $u \in
\mathbb{Z}_2^n$, with $u =(u_1,u_2,\dots, u_n)$, $u_k =1$ if $k \in A$
and $u_k = 0$ otherwise. Similarly let $v \in \mathbb{Z}_2^n$ be the
vector corresponding to the set $B$.  Since $A$ and $B$ are
qualitatively independent, $A \cap B \neq \emptyset$.  So there exists
some $i \in A \cap B$, and so $(u_i,v_i) = (1,1)$.  Further, since $A
\cap \overline{B} \neq \emptyset$, there exists some $j \in A \cap
\overline{B}$, and so $(u_j,v_j) = (1,0)$. Similarly, for $k \in
\overline{A} \cap B$ we get $(u_k,v_k) = (0,1)$ and for $l \in
\overline{A} \cap \overline{B}$ we get $(u_l,v_l) = (0,0)$.  Thus the
vectors $u$ and $v$ are qualitatively independent.

Conversely, if $u,v \in \mathbb{Z}_2^n$ are qualitatively independent,
then the sets $A$ and $B$, defined by $i\in A$ if and only if $u_i=1$
and $i\in B$ if and only if $v_i=1$, are also qualitatively
independent.

\begin{thm}[\cite{katona:73, kleitman:73}]\label{thm:katona73}
If $\mathcal{A} =\{A_1,A_2 ,\dots, A_k\}$ is a qualitatively
independent set system of an $n$-set, then
\[
|\mathcal{A}| \leq {n-1 \choose \lfloor n/2 \rfloor -1}.
\]
\orphan
Further, this bound is attained by a $\lfloor n/2 \rfloor$-uniform trivially 1-intersecting set system.
\end{thm}

\begin{proof}
First, assume $n$ is even. Define
\[
\mathcal{A}^\ast = \{A_i, \overline{A_i} \; : \; A_i \in \mathcal{A} \}.
\]
The system $\mathcal{A}^\ast$ is a Sperner set system so, 
$|\mathcal{A}^\ast| \leq {n \choose n/2 }$. Hence 
$|\mathcal{A}| \leq 1/2{n \choose n/2 } = {n-1 \choose n/2 -1}$. 
This bound is attained by the set system
\[
\mathcal{A} = \{A \in {[n] \choose \frac{n}{2}} : 1 \in A\}.
\]

Next, assume $n$ is odd.  If $A_i \in \mathcal{A}$ and $|A_i| \geq n/2$
then replace $A_i$ with $\overline{A_i}$. This does not affect the
pairwise qualitative independence of $\mathcal{A}$. So, we can
assume that each $A_i \in
\mathcal{A}$ has $|A_i| \leq \lfloor n/2 \rfloor$.

By the definition of qualitative independence, $\mathcal{A}$ is a
1-intersecting set system and by the Erd\H{o}s-Ko-Rado Theorem,
$|\mathcal{A}| \leq {n-1 \choose \lfloor n/2 \rfloor -1}$.  This bound
is attained by the set system
\[
\mathcal{A} = \{A \in {[n] \choose \frac{n-1}{2}} : 1 \in A\}.
\]
\end{proof}

From Theorem~\ref{thm:katona73}, the exact size of the optimal
binary covering array with $r$ rows can be found for all $r$.

\begin{thm}[\cite{kleitman:73}]\label{thm:binaryCAN}
Let $r$ be a positive integer, then
\[
CAN(r,2) = \min \left\{n : {n-1 \choose \lfloor n/2 \rfloor-1 } \geq r \right\}.
\]
\end{thm}

%%%%%%%%%%%%%%%%%%%%%%%%%%%%%%%%%%%%%%%%%%%%%%%%%%%%%%%%%%%%
\chapter{Graph Theory}\label{graphtheory}

As stated in Section~\ref{CA}, covering arrays can be used to design
test suites for systems and networks. To improve efficiency in such
applications, several variations on covering arrays have been
considered, such as mixed-level covering arrays and variable-strength
covering arrays~\cite{ constructiontestsuites, MR2071905,MR2012427}.
The generalization considered in this thesis adds a graph structure to
the covering array, obtaining a {\em covering array on a graph}.  These
are defined in Chapter~\ref{CAG}, where we also define a family of
graphs called the {\em qualitative independence graphs}. These graphs
are used to find bounds both on covering arrays on graphs and
covering arrays. For this reason, we need to introduce graph concepts that
will be used throughout this thesis.

The first half of this chapter is a review of basic graph theory. The
second half is a review of algebraic graph theory, which includes spectral
theory of graphs and association schemes. These will be used in
Chapter~\ref{chp:ass} in the study of the qualitative independence
graphs. Association schemes at first look seem to be purely
algebraic objects, but they are very closely related to designs. In
fact, they were first defined by Bose and Shimamoto in
1952~\cite{MR0048772} to generalize block designs. For more on the
history of block designs and association schemes see Chapter 13 of
Bailey's {\em Association Schemes: Designed Experiments, Algebra and
Combinatorics}~\cite{MR2047311}.

\section{Basic Graph Theory}

Throughout this chapter, $G$ will denote a simple graph, unless
otherwise stated. The vertex set of $G$ will be denoted by $V(G)$, and
the edge set will be denoted by $E(G)$.

\subsection{Graph Homomorphism}

This section introduces graph homomorphism; for more information
on graph homomorphism see~\cite{MR2002f:05002,hahn:97, MR2089014}.

\begin{defn}[Graph Homomorphism]\index{graph homomorphism}\index{homomorphism}
Let $G$ and $H$ be graphs.  A mapping $\phi$ from $V(G)$ to $V(H)$ is
a {\em graph homomorphism} from $G$ to $H$ (or simply a {\em
homomorphism}) if for all $v,w\in V(G)$, the vertices $\phi(v)$ and
$\phi(w)$ are adjacent in $H$ whenever $v$ and $w$ are adjacent in
$G$.
\end{defn}

For graphs $G$ and $H$, if there exists a homomorphism from $G$ to
$H$, then we write $G \rightarrow H$.

For graphs $G$ and $H$, a map $\phi$ from $V(G)$ to $V(H)$ is a {\em graph
isomorphism} (or simply an {\em isomorphism}) if $\phi$ is
a bijection such that $x,y \in V(G)$ are adjacent in $G$ if and only if
$\phi(x)$ and $\phi(y)$ are adjacent in $H$.
\index{graph isomorphism}\index{isomorphism} If there
exists an isomorphism between two graphs, then we say the graphs are
{\em isomorphic}\index{isomorphic graphs}.

A homomorphism from a graph $G$ to itself is a {\em graph
endomorphism}\index{graph endomorphism} (or an {\em endomorphism}).  An
isomorphism from a graph $G$ to itself is a {\em graph automorphism}
\index{graph automorphism} (or an {\em automorphism}). The {\em
automorphism group} for a graph $G$ is the group of all automorphisms
of $G$; it is denoted $Aut(G)$.\index{automorphism group}

A {\em fibre} of a homomorphism $\phi:G\rightarrow H$ is the preimage
$\phi^{-1}(w)$ of some vertex $w \in V(H)$, that is,
\[
\phi^{-1}(w) = \{v \in V(G) \; : \; \phi(v)= w \}.
\]

\subsection{Colourings, Cliques and Independent Sets}\label{sec:cliquesandchromatic}

The {\em complete graph on $n$ vertices}\index{complete graph}, $K_n$,
is the graph with $n$ vertices and with an edge between any two distinct
vertices.

A {\em proper colouring}\index{proper colouring} of $G$ with $n$
colours is a map from $V(G)$ to a set of $n$ colours such that no
two adjacent vertices are assigned the same colour.  A proper
colouring of a graph $G$ with $n$ colours is equivalent to a
homomorphism from $G$ to $K_n$.  The {\em chromatic number}
\index{chromatic number} of a graph $G$, denoted $\chi(G)$, is the smallest $n$ such that
$G \rightarrow K_n$.

A {\em clique}\index{clique} in a graph $G$ is a set of vertices from
$V(G)$ in which any two distinct vertices are adjacent in $G$.  The
size of a maximum clique in $G$ is denoted by $\omega(G)$.  If $G$ has a
clique of size $n$, then there is a homomorphism $K_n \rightarrow G$,
and the size of a maximum clique in $G$ is the largest $n$ for which
$K_n \rightarrow G$. An {\em $n$-clique} is a clique of size
$n$.\index{$n$-clique}

For graphs $G$ and $H$, if there is a homomorphism $G \rightarrow H$,
then $\chi(G) \leq \chi(H)$ and $\omega(G) \leq \omega(H)$.  Also, for
all graphs $G$, $\omega(G) \leq \chi(G)$.

An {\em independent set} \index{independent set} in a graph $G$ is a
set of vertices from $V(G)$ in which no two vertices are adjacent in
$G$.  The size of a largest independent set in a graph $G$ is
denoted by $\alpha(G)$.

The vertices which are assigned the same colour in a proper colouring
form an independent set. In fact, a proper colouring on a graph $G$
partitions the vertices of $G$ into independent sets called {\em
colour classes}\index{colour classes}. A proper colouring corresponds
to a binary function on the independent sets of a graph: each
independent set that is a colour class in the proper colouring is
assigned a value of 1 and all other independent sets are assigned a
value of 0 by the function. Further, each vertex is in exactly one independent set
which has an assigned value of 1. The fractional relaxation of this
binary function is used in Section~\ref{knesergraphs} to define a
fractional colouring.

\subsection{Vertex-Transitive Graphs}

The {\em degree}\index{degree} of a vertex $v \in V(G)$ is the number
of vertices in $G$ which are adjacent to $v$. If every vertex in $G$
has the same degree, then we say that $G$ is {\em regular}\index{regular}.
Specifically, if every vertex in $G$ has degree $k$, we say $G$ is {\em
$k$-regular}\index{$k$-regular}.

A graph is {\em vertex transitive} \index{vertex transitive} if its
automorphism group acts transitively on the set of vertices. This means that for
any two distinct vertices, there is an automorphism on the graph that maps one
vertex to the other. If a graph is vertex transitive, then it must also be regular.
\enlargethispage{20pt}

There are many other types of transitivity. One example is {\em
arc transitivity}\index{arc transitive}. An {\em arc}\index{arc} in a
graph is an ordered pair of adjacent vertices. A graph is
arc transitive if its automorphism group acts transitively on the arcs
of the graph.  This means for any two arcs in the graph
there is an automorphism that maps one arc to the other.

\subsection{Cores of a Graph}\label{subsec:cores}

\begin{defn}[Core Graph]\index{core graph}
A graph $G$ is a {\em core} if any endomorphism on $G$ is an
automorphism.
\end{defn}

If a graph $G$ is a core, then there is no homomorphism from $G$ to a
proper subgraph of $G$. Indeed, for any positive integer $n$, the complete
graph $K_n$ is a core.

\begin{defn}[Core of a Graph]\index{core of a graph}
A {\em core of a graph} $G$ is a subgraph
$G^\bullet$ of $G$ such that $G^\bullet$ is a core and there is a
homomorphism $G \rightarrow G^\bullet$.
\end{defn}

If $G^\bullet_1$ and $G^\bullet_2$ are both cores of a graph $G$, then
$G^\bullet_1$ and $G^\bullet_2$ are isomorphic (see Lemma 6.2.2,~\cite{MR2002f:05002}). Further, a core of a
graph preserves many important properties of the graph. For example,
\[
\chi(G) = \chi(G^\bullet) \textrm{ and } \omega(G) = \omega(G^\bullet).
\]

Any automorphism on $G$ induces an automorphism on $G^\bullet$,
so a core of a graph also preserves vertex transitivity.

\begin{thm}[\cite{hahn:97}]\label{thm:corevt}
If $G$ is a vertex-transitive graph, then $G^\bullet$ is also
vertex transitive.
\end{thm}

If $G$ is a vertex-transitive graph, then for any homomorphism
from $G$ to $G^\bullet$ each fibre of the homomorphism has the same
cardinality. This means that the preimage of every vertex in $G^\bullet$ has
the same cardinality. From this we can conclude the following theorem.

\begin{thm}[\cite{hahn:97}]\label{thm:fibres}
If $G$ is a vertex-transitive graph, then $|V(G^\bullet)|$ divides $|V(G)|$.
\end{thm}

\subsection{Kneser Graphs}\label{knesergraphs}

In this section, a family of vertex-transitive graphs, the Kneser
graphs, is defined. The vertex set of a Kneser graph is the set of
all $r$-subsets of an $n$-set, so this graph is closely related to set
systems.  In Section~\ref{CAG}, we introduce a similar family of
graphs, called the {\em qualitative independence graphs}, whose vertex
set is a set of partitions rather than sets.

\begin{defn}[Kneser Graphs]\index{Kneser graphs}
For positive integers $r,n$ with $r\leq n$, the {\em Kneser graph}
$K_{n:r}$ is the graph whose vertex set is ${[n] \choose r}$ and
$r$-subsets are adjacent if and only if they are disjoint.
\end{defn}

Kneser graphs are related to a generalization of proper colouring called {\em
fractional colouring}.  From Section~\ref{sec:cliquesandchromatic}, a
proper colouring of a graph is usually considered to be a homomorphism
from $G$ to a complete graph, but it can also be considered a
binary function on the independent sets of the graph. A fractional
colouring is also a function on the independent sets of a graph, but
the requirement that the function be binary is relaxed to that the
function be non-negative.

For a graph $G$ with $v \in V(G)$, define $\mathcal{I}(G,v)$ to be the
set of all independent sets in $G$ that contain the vertex $v$, and
define $\mathcal{I}(G)$ to be the set of all independent sets in $G$.

\begin{defn}[Fractional Colouring]\index{fractional colouring}
A {\em fractional colouring} of a graph $G$ is a non-negative function
$f$ on the independent sets of $G$ with the property that for any
vertex $v \in V(G)$,
\[
\sum_{S \in \mathcal{I}(G,v)} f(S) \geq 1.
\]
\end{defn}

\begin{defn}[Fractional Chromatic Number] \index{fractional chromatic number}
Let $G$ be a graph and $f$ a fractional colouring on $G$.  The {\em
weight}\index{weight of a fractional colouring} of $f$ is the sum of
the values of $f$ over all independent sets in $G$, that is, $\sum_{S \in
\mathcal{I}(G)}f(S)$.  The {\em fractional chromatic number} of $G$,
denoted $\chi^\ast(G)$, is the minimum weight over all fractional
colourings of $G$.
\end{defn}

The connection between fractional chromatic number and Kneser graphs
is that, like the chromatic number, the fractional chromatic number can
be determined by a homomorphism. While the chromatic number is determined
by homomorphisms to complete graphs, the fractional chromatic number is
determined by homomorphisms to Kneser graphs.

\begin{thm}[see Theorem 7.4.5,~\cite{MR2002f:05002}]\label{thm:fractionalchromatic}
For any graph $G$,
\[
\chi^\ast(G) = \min \left\{ \frac{n}{r} : G \rightarrow K_{n:r} \right\}.
\]
\end{thm}
 
For any graph $G$, the fractional chromatic number satisfies the following bounds:
\begin{eqnarray}\label{eqn:omegafracchi}
\omega(G) \leq \chi^\ast(G) \leq \chi(G).
\end{eqnarray}
The fractional chromatic number, like the chromatic number and maximum clique
size, is preserved by cores, that is, for any graph $G$,
\begin{equation}\label{eqn:fraccores}
\chi^\ast(G) = \chi^\ast(G^\bullet).
\end{equation}

For the special case when $G$ is a vertex-transitive graph, there is
another formula for $\chi^\ast(G)$. This formula, is derived from the
fact that there exists a fractional colouring with the minimum weight
that assigns the same value to each maximum independent set in $G$.
 
\begin{cor}[see Chapter 7,~\cite{MR2002f:05002}]\label{cor:vtbound}
If $G$ is a vertex-transitive graph, then
\begin{eqnarray*}
\chi^\ast(G) = \frac{|V(G)|}{\alpha(G)}.
\end{eqnarray*} 
\end{cor}

For a vertex-transitive graph $G$, this can be used with
Inequality~(\ref{eqn:omegafracchi}) as a bound on the size of a maximum
clique in $G$
\begin{eqnarray}\label{eq:vt2}
\omega(G) \leq \chi^\ast(G) = \frac{|V(G)|}{\alpha(G)}.
\end{eqnarray}

The next theorem states that for $n > 2r$ the Kneser graph $K_{n:r}$
is a core. The proof of this is included for two reasons: first, it is
a nice application of the Erd\H{o}s-Ko-Rado Theorem; and second, a
similar, but more complicated proof is used to prove
Lemma~\ref{lem:qi93core}, but we do not include it in this thesis.
Before we can prove that a Kneser graph is a core, we need a short lemma.

\begin{lemma}[Theorem 7.5.4,~\cite{MR2002f:05002}]\label{lem:preimageindy}
Let $G$ and $H$ be vertex-transitive graphs with the same fractional chromatic number
with a homomorphism $\phi : G \rightarrow H$.  Then, if $S$ is
a maximum independent set in $H$, the preimage of $S$,
$\phi^{-1}(S) = \{ v \in V(G) \; : \; \phi(v) \in S\}$, is a
maximum independent set in $G$.
\end{lemma}
\remove{
\begin{proof}
{}From the definition of a homomorphism, $\phi^{-1}(S)$
is an independent set.  All that is needed is to show that $|\phi^{-1}(S)|=\alpha(G)$.

Since $G$ is vertex transitive, the size of the fibres of $\phi$ is the
same for every vertex in $G^\bullet$. Let $a$ be the size of the
fibres, so $|V(G)| = a |V(G^\bullet)|$ and $|\phi^{-1}(S)| = a|S| = a \alpha(G^\bullet)$.

From Equation~(\ref{eqn:fraccores}), $\chi^\ast(G) =\chi^\ast
(G^\bullet)$, and since both $G$ and $G^\bullet$ are vertex transitive,
\[
\frac{|V(G)|}{\alpha(G)} = \frac{|V(G^\bullet)|}{\alpha(G^\bullet)}.
\]
Thus $|\phi^{-1}(S)| = a\alpha(G^\bullet) = \frac{a |V(G^\bullet)| \alpha(G)}{ |V(G)|} =\alpha(G)$.
\end{proof}
}

\begin{thm}[Theorem 7.9.1,~\cite{MR2002f:05002}]\label{knesercore}
For positive integers $n,r$ with $n > 2r$, $K_{n:r}$ is a core.
\end{thm}
\begin{proof}
Let $S$ be a maximum independent set in $K_{n:r}$. The vertices of
$K_{n:r}$ are $r$-subsets of an $n$-set, and $S$ is an $r$-uniform set
system.  As $S$ is an independent set, for any $A,B \in S$, $A \cap B
\neq \emptyset$.  Thus, $S$ is a maximum $r$-uniform 1-intersecting
set system of an $n$-set.  As $n>2r$, by the Erd\H{o}s-Ko-Rado
Theorem, $S$ must be a trivially 1-intersecting set system of
cardinality ${n-1 \choose r-1}$, and all $r$-subsets in $S$ contain a
common element.

For each $i=\{1,\dots, n\}$, let $S_i \subset V(K_{n:r})$ be the system of
all $r$-subsets that contain the element $i$. The sets $S_i$ are all
the maximum independent sets in $K_{n:r}$.

Let $\phi$ be an endomorphism $\phi:K_{n:r}\rightarrow K_{n:r}$.  By
Lemma~\ref{lem:preimageindy}, the preimage $\phi^{-1}(S_i)$ is also a
maximum independent set. This means there is a $j_i \in \{1,\dots,n\}$
such that $\phi^{-1}(S_i) = S_{j_i}$.

The subset $B = \{1,2,\dots, r\}$ is the unique vertex in
the intersection of the set systems $S_1 \cap S_2 \cap \dots \cap S_r$.
Then 
\[
\phi^{-1}(B) = \phi^{-1}(S_1) \cap \phi^{-1}(S_2) \cap \dots \cap \phi^{-1}(S_r) \\
= S_{j_1} \cap S_{j_2} \cap \dots \cap S_{j_r}.
\]  

There is at least one $r$-subset $A \in V(K_{n:r})$ with $\{j_1,j_2,\dots,
j_r\} \subseteq A$ (there are more subsets if the elements of
$\{j_1,j_2,\dots, j_r\}$ are not distinct). Thus, for every $B \in V(K_{n:r})$,
$\phi^{-1}(B) \neq \emptyset$. So every endomorphism $\phi$ is an onto map
and thus, is an automorphism.
\end{proof}

\section{Algebraic Graph Theory}

In this section, some concepts from algebraic graph theory are
introduced. This theory is used in Chapter~\ref{chp:ass} to study the
qualitative independence graphs. This section follows Godsil and
Royle's {\em Algebraic Graph Theory}~\cite{MR2002f:05002} and Godsil's
{\em Algebraic Combinatorics}~\cite{MR94e:05002}.

%%%%%%%%%%%%%%%%%%%%%%%%%%%%%%%%%%%%%%%%%%  spectral theory
\subsection{Spectral Theory of Graphs}

Let $G$ be a simple graph on $n$ vertices labelled by $1,2,\dots,
n$.  The {\em adjacency matrix} \index{adjacency matrix of a graph} of
$G$ is the $n \times n$ matrix with a 1 in the $i,j$ position if
vertices $i$ and $j$ are adjacent in $G$, and 0 if vertices $i$ and
$j$ are not adjacent in $G$. The adjacency matrix of a graph $G$ is
denoted by $A(G)$.  The {\em characteristic polynomial}
\index{characteristic polynomial of a graph} of a graph $G$ is the characteristic
polynomial of the matrix $A(G)$. The eigenvalues and eigenvectors of
the graph $G$ are the eigenvalues and eigenvectors of the matrix
$A(G)$.  The {\em spectrum}\index{spectrum} for a graph is the set of
all the eigenvalues with their multiplicities.

For any undirected graph $G$ on $n$ vertices, $A(G)$ is a real
symmetric matrix.  Thus, the eigenvalues of $G$ are real numbers and
the eigenvectors for $G$ span $\mathbb{R}^n$. If two graphs are
isomorphic, they have the same spectrum. But, two graphs can have the
same spectrum and not be isomorphic.

The eigenvalues of a graph $G$ can give information about $G$. In
particular, we will consider cases where eigenvalues can be used to
find upper bounds on $\alpha(G)$ and $\omega(G)$. In general, finding
the eigenvalues of a graph is difficult. However, for special classes
of graphs we know more, as stated in the following theorems.

\begin{thm}\label{evaluescomplete}
The eigenvalues of the complete graph $K_n$ are $n-1$ and $-1$ with
multiplicities $1$ and $n-1$, respectively.
\end{thm}
 
\begin{thm}[see Theorem 2.4.2,~\cite{MR94e:05002}]\label{regularevalues}
If $G$ is a $k$-regular graph, then $k$ is an
eigenvalue of $G$. Moreover, $k$ is the largest eigenvalue of $G$.
\end{thm}

In a $k$-regular graph $G$, the multiplicity of the
eigenvalue $k$ is equal to the number of components in $G$.

\begin{thm}[see Lemma 8.5.1,~\cite{MR2002f:05002}]\label{thm:evaluescomp}
If $G$ is a $k$-regular graph on $n$ vertices with eigenvalues
$k,\lambda_2, \dots,\lambda_n$, then the eigenvalues of the complement $\overline{G}$
of $G$ are $n-k-1, -1-\lambda_2, \dots, -1-\lambda_n$. Moreover, the
graphs $G$ and $\overline{G}$ have the same eigenvectors.
\end{thm}

\subsection{Equitable partitions}\label{sec:ep}

In this section, we introduce a method that reduces the amount of calculation
needed to find the eigenvalues of a graph.

For a graph $G$, a partition $\pi$ of $V(G)$ with $r$ classes
$C_1,C_2,\dots,C_r$ is an {\em equitable partition}\index{equitable
partition} if the number of vertices in $C_j$ that are adjacent to $v
\in C_i$ is a constant $b_{ij}$, independent of $v$. For a graph $G$
with an equitable partition $\pi$, the {\em quotient graph of $G$ over
$\pi$} \index{quotient graph} is the directed multi-graph whose
vertices are the $r$ classes $C_i$ with $b_{ij}$ arcs from the
$i^{th}$ class to the $j^{th}$ class. This graph is denoted by
$G/\pi$. The adjacency matrix of $G/\pi$ is given by
\[
(A(G/\pi))_{i,j} = b_{ij}.
\]

There are many equitable partitions for a graph $G$.  For
any subgroup $X \leq Aut(G)$, the set of orbits formed by the group
action of $X$ on the vertices of $G$ is an equitable partition.  To
see this, consider two vertices $v,v'$ which are in the same orbit
under the group action of $X$. In this case, there is some $x \in X$ such that
$x(v) = v'$. For any $w$ adjacent to $v$, the vertex $x(w)$ is adjacent
to $v'$ and is in the same orbit as $w$.  Thus the number of vertices
adjacent to $v$ in any orbit is the same as the number of vertices in
that orbit adjacent to $v'$.

\begin{thm}[see Theorem 9.3.3, \cite{MR2002f:05002}]\label{thm:epevalues}
If $\pi$ is an equitable partition of a graph $G$, then the
characteristic polynomial of $G/\pi$ divides the characteristic
polynomial of $G$. 
\end{thm}

If $\lambda$ is an eigenvalue for $G/\pi$ with multiplicity $m$,
then $\lambda$ is an eigenvalue for $G$ with multiplicity at least
$m$. In some special circumstances, the eigenvalues of $G/\pi$ are exactly
the eigenvalues of $G$.
 
\begin{thm}[see Theorem 9.4.1, \cite{MR2002f:05002}]\label{singlecell}
Let $G$ be a vertex-transitive graph and $\pi$ be a partition on
$V(G)$ generated by the orbits of some subgroup $X \leq Aut(G)$. If
$\pi$ has a class of size 1 (a singleton class\index{singleton
class}), then every eigenvalue of $G$ is an eigenvalue of $G/\pi$.
\end{thm}
%\todo{learn a proof of this}
This can be used to find all the eigenvalues of the Kneser graphs.

\begin{exmp}\label{exmp:kneserevalues}
Let $r,n$ be positive integers with $r\leq n$ and let $K_{n:r}$ be a
Kneser graph.  Fix an arbitrary $r$-subset $A \in
V(K_{n,r})$ of the $n$-set. For $i$ $=0,\dots,r$, define $C_i$ to be the collection of
$r$-subsets of the $n$-set (the vertices in $K_{n,r}$) that have
intersection of size $r-i$ with the subset $A$.  This is an equitable
partition; it corresponds to the orbit partition of the group $\fix(A)
= \{\sigma \in Sym(n) \; : \; \sigma(A) = A\}$. The class $C_{0}$
contains exactly one element, $A$.  It is possible to build the
adjacency matrix $A(K_{n:r}/\pi)$ for the quotient graph.  Let $B$ be
an $r$-subset that intersects with $A$ in exactly $r-i$ elements, so
$B\in C_i$.  The $i,j$ entry of $A( K_{n:r}/\pi)$ is the number of
$r$-subsets in $C_j$ that are disjoint from $B$.  This is the number of
$r$-subsets that meet $A$ in exactly $r-j$ positions and are disjoint
from $B$.  The number of such subsets is
\[
{i \choose r-j}{n-r-i \choose j}.
\]
Since this depends only on $i$ and $j$, the partition is equitable.

From Theorem~\ref{singlecell} the eigenvalues of $A(K_{n:r}/\pi)$ are
exactly the eigenvalues of the Kneser graph $K_{n:r}$. Using some
computations with binomial coefficients (Section 9.4,
\cite{MR2002f:05002}), it is possible to show that the eigenvalues of
$A(K_{n:r}/\pi)$, and hence $K_{n:r}$, are
\begin{equation}\label{kneserevalues}
(-1)^i{ n-r-i \choose r-i } \;\mbox{  for  }\; i = 0,1,\dots,r.
\end{equation}

\end{exmp}

\subsection{Association Schemes}\label{sec:introassscheme}

%\todo{connect to designs}
%history.
\remove{In this section, we give two definitions for association schemes. The
first used graphs and the second used matrices. These definitions are
equivalent and we will use which ever form is more convenient and we need them.}

\begin{defn}[Association Scheme]\index{association scheme}\label{defn:assscheme}
An {\em association scheme} on a set $X$ is a set of $d$ graphs
$G_1,G_2,\dots, G_d$, all with vertex set $X$, which has the following
properties:
\begin{enumerate}
\item{distinct elements $x,y \in X$ are adjacent in exactly one graph $G_i$;}
\item{for all $x,y \in X$ and $1 \leq i,j \leq d$, the number of elements $z$, 
with $x$ and $z$ adjacent in $G_i$ and $y$ and $z$ adjacent in $G_j$
depends only on $i,j,k$ where $G_k$ is the graph in which $x$ and $y$
adjacent. This number is denoted by $p_{i,j}^k$.}
\end{enumerate}
\end{defn}

Note that the second requirement implies that each of the graphs in
the scheme is regular. The graphs in the scheme are also called the
{\em classes}\index{classes of an association scheme} of the scheme.
For any set $X$, the complete graph on $|X|$ vertices is a scheme with just one
class. This scheme is called the {\em trivial scheme}\index{trivial scheme}.

Association schemes can also be defined in terms of matrices.  In this
case, the classes of the scheme are represented by 01-matrices.  The
matrix $J$ denotes the matrix with all entries equal to 1.

\begin{defn}[Association Scheme (definition 2)]\label{defn:assscheme2}
Let $n$ be a positive integer. An {\em association scheme}, with $d$
classes is a set of $n \times n$ 01-matrices $\{A_0, A_1, \dots,
A_d\}$, with the following properties:
\begin{enumerate}
\item{$A_0 = I$;}
\item{$\sum_{i=0}^{d}A_i = J$;}
\item{$A_i^\top = A_i$, for all $i = 1,\dots, d$;}\label{symmetric}
\item{for all $i,j \in \{0,\dots,d\}$, $A_iA_j$ is a linear combination of $\{A_0,\dots, A_d \}$;}
\label{closed}
\item{$A_iA_j = A_jA_i$, for all $i,j \in \{0,\dots,d\}$.}\label{commute}
\end{enumerate}
\end{defn}

These two definitions of association schemes are equivalent.  Assume
the graphs $G_i$, for $i=1,\dots,d$, form an association scheme on a set
$X$ with cardinality $n$.  Set $A_i = A(G_i)$ for $i = 1,\dots, d$
and $A_0 = I$. The first three properties of matrices in an
association scheme hold since each $G_i$ is a graph and any pair
of distinct elements $x,y \in X$ are adjacent in exactly one graph
$G_i$.

Define $G_0$ to be the graph with vertex set $X$ and a loop
at each vertex.  Let $x,y \in X$ be adjacent in the graph $G_k$.
Then, for any $i,j \in \{0,\dots,d\}$, the $(x,y)$ entry of the matrix $A_iA_j$ is
\begin{eqnarray*}
(A_iA_j)_{x,y} &=& \sum_{z \in X} (A_i)_{x,z}(A_j)_{z,y}  \\
            &=& |\{z :\mbox{ $x$ and $z$ adjacent in $G_i$, $y$ and $z$ adjacent in $G_j$}\}| \\
            &=& p_{i,j}^k.
\end{eqnarray*}
Thus $A_iA_j = \sum_{k=0}^{d} p_{i,j}^k A_k$, and the final two
conditions for an association scheme hold.

Conversely, if the 01-matrices $A_i$, for $i = 0,\dots, d$, form an
association scheme, the set of $d$ graphs that have $A_i$, for
$i=1,\dots, d$, as their adjacency matrices form an association scheme.

The trivial scheme in matrix form is the set of matrices $\{I, J-I\}$. Note
that there are two matrices in this scheme, but the scheme only has one
class.

In some of the literature (\cite{MR882540,MR0384310}), association
schemes as defined here are referred to as {\em symmetric association
schemes}. In these references, a general association scheme is the same
as in Definition~\ref{defn:assscheme2}, but with Condition~\ref{symmetric}
replaced with the condition that $A_i^T \in \mathcal{A}$. In most
recent papers, association schemes are assumed to be symmetric and
schemes with Condition~\ref{symmetric} replaced with the condition
that $A_i^T \in \mathcal{A}$ are call {\em asymmetric association
schemes}\index{asymmetric association scheme}.
This thesis will follow this trend and use ``association
scheme'' to mean a symmetric association scheme.

\begin{exmp}{\bf(Johnson Scheme). }\label{exmp:johsonscheme}
Let $i,r,n$ be positive integers with $i \leq r \leq n$.  Define the
{\em generalized Johnson graph}\index{generalized Johnson graph}
$J(n,r,r-i)$ to be the graph whose vertex set is the set of all
$r$-subsets of an $n$-set and two $r$-sets are adjacent if and only if
the subsets have exactly $i$ elements in common. The graph $J(n,r,r)$
is isomorphic to the Kneser graph $K_{n:r}$ and the graph $J(n,r,0)$
has a loop at each vertex and no other edges. If $n<2r-i$ the graph
$J(n,r,r-i)$ is the empty graph.

The $r$ graphs $J(n,r,i)$, for $i = 1,\dots, r$, form an association
scheme called the {\em Johnson scheme}\index{Johnson scheme}.

If the $r$-subsets $A,B$ are adjacent in the graph $J(n,r,r-i)$, then
$|A\cap B| = i$. For any $r$-subset $C$ that is adjacent to $A$ in the
graph $J(n,r,r-j)$ and adjacent to $B$ in the graph $J(n,r,r-k)$, we have
\[
|A \cap C| = j \mbox{ and } |B \cap C| = k.
\]
The total number of such $r$-subsets $C$ is
\[
 \sum_{ \stackrel{x+y = j}{x+z = k} }  {n-(2r-i) \choose r-(x+y+z)}
    {i \choose x}{r-i \choose y}{r-i \choose z}.
\]
As the number of subsets $C$ depends only on $i,j$ and $k$, the Johnson scheme
is an association scheme.

\end{exmp}

Let $\mathcal{A}$ be an association scheme with $d+1$ matrices of order
$n$. From Condition~\ref{closed} in
Definition~\ref{defn:assscheme2}, the span (over $\mathbb{R}$) of the
matrices in $\mathcal{A}$ is a $(d+1)$-dimensional algebra.  This
algebra is called the {\em Bose-Mesner Algebra}\index{Bose-Menser
algebra}.  The matrices in the association scheme are a basis for this
algebra, and hence a convenient way to consider the algebra.  Often it
is more useful to use the matrix definition of association schemes
than the graph definition.

Any matrix in an association scheme is a real symmetric matrix; thus,
it is diagonalizable and all the eigenvalues are real numbers.
Moreover, from Condition~\ref{commute} in
Definition~\ref{defn:assscheme2}, any two matrices in an association
scheme commute; thus, for any $i,j \in \{1,\dots, d\}$, any eigenspace
of the matrix $A_i$ is {\em $A_j$-invariant} for any matrix $A_j$ in
the scheme. That is, if $E_\lambda$ is an eigenspace of the matrix
$A_i$, then for any $v \in E_\lambda$ we have that $A_j v \in
E_\lambda$ for any matrix $A_j$ in the scheme.

It is possible to find $d+1$ subspaces $U_0,U_1,\dots, U_d$ of
$\mathbb{R}^n$ that partition $\mathbb{R}^n$ with the property that
for every $i,j \in \{0,\dots, d\}$, the space $U_j$ is contained in an
eigenspace of $A_i$.  For each $i,j \in \{0,\dots, d\}$, denote the
eigenvalue of $A_i$ on the subspace $U_j$ by $\lambda_i(j)$. Then the
{\em matrix of eigenvalues}\index{matrix of eigenvalues for an
association scheme} for an association scheme is defined to be the matrix $(P)_{j,i} =
\lambda_i(j)$. In particular, the $i^{th}$ column of the matrix of
eigenvalues contains the list of eigenvalues of the matrix $A_i$ in the
scheme.

As the matrix $A_0$ is the identity matrix, the first column of the
matrix of eigenvalues has each entry equal to one. Often, this first
column is replaced with a column containing the dimension of the
subspace $U_i$; these values give the multiplicities of the eigenvalues. The
matrix of eigenvalues with the first column replaced by the
multiplicities of the eigenvalues is called the {\em modified matrix
of eigenvalues}\index{modified matrix of eigenvalues} for the
association scheme.

%%%%%%%%%%%%%%%%%%%%%%%%%%%%%%%%%%%%%%%%%%%%%%%%%%%%%%%%%%%%%%%%%%

\subsection{Ratio Bounds}\label{sec:ratiobounds}

There are two bounds using eigenvalues called the {\em ratio bounds}.
The first ratio bound is an upper bound on the size of the maximum
independent set in a graph and the second ratio bound is an upper
bound on the size of the maximum clique. The ratio bound on
independent sets was established for graphs in an association scheme
by Delsarte~\cite{MR0384310}. The proof given here is from Section 9.6
of~\cite{MR2002f:05002}. For an interesting examination of several
results from the ratio bounds, see Newman's
Ph.D. thesis~\cite{mikesthesis}.

\begin{lemma}[Ratio Bound for $\alpha(G)$, see Section 9.6,~\cite{MR2002f:05002}]\label{lem:smallevbound}
\index{ratio bound for $\alpha(G)$} 
Let $G$ be a $d$-regular graph on $n$ vertices and let $\tau$ be the
least eigenvalue of $G$. Then
\[
\alpha(G) \leq \frac{n}{1-\frac{d}{\tau}}.
\] 
\end{lemma}
\remove{
\begin{proof}
Let $\alpha = \alpha(G)$.  Let $S$ be a maximum independent set in
$G$.  Define $\pi$ to the be the partition of the vertices of $G$ with
two classes, $S$ and $V(G) \backslash S$.  There are $\alpha$ vertices
in $S$ and $n-\alpha$ vertices in $V(G)
\backslash S$.  Each vertex in $S$ is adjacent to $d$ vertices in
$V(G) \backslash S$.  Thus, each vertex in $V(G) \backslash S$ is adjacent
to $\frac{\alpha n}{n-\alpha}$ vertices in $S$ and $d -
\frac{\alpha n}{n-\alpha }$ vertices in $V(G) \backslash S$. The quotient
matrix of $A(G/\pi)$ is
\[ 
\left(\begin{matrix} 0 & d \\ \frac{\alpha n}{n-\alpha } 
         & d - \frac{\alpha n}{n-\alpha} \end{matrix} \right).
\]
The eigenvalues for $A(G/\pi)$ are $d$ and
$\frac{-d\alpha}{n-\alpha}$.  From Lemma~\ref{lem:evaluepartition},
the eigenvalues of $A(G/\pi)$ interlace the eigenvalues of $G$, the
smallest eigenvalue of $G$ is smaller than the smallest eigenvalue of
$A\pi$.  Thus $\tau \leq \frac{-d\alpha}{n-\alpha}$ and the result
follows.
\end{proof}
}

Any independent set in $G$ with cardinality
$\frac{n}{1-\frac{d}{\tau}}$ is called {\em
ratio-tight}\index{ratio-tight}.  If an independent set is ratio-tight,
then more is known about the independent set.

\begin{thm}[\cite{mikesthesis}]\label{lem:completeratiobound}
Let $G$ be a $d$-regular graph on $n$ vertices with the least eigenvalue
$\tau$. Let $A(G)$ be the adjacency matrix of the graph $G$. Let $S$
be an independent set in $G$ of size $s$ and ${\bf z}$ the characteristic vector of $S$.  If
\[
s = \frac{n}{1-\frac{d}{\tau}},
\] 
then
\[
A(G) \left( {\bf z}-\frac{s}{n} {\bf 1} \right) = \tau \left( {\bf z} - \frac{s}{n} {\bf 1}\right).
\]
Moreover, $s = \frac{n}{1-\frac{d}{\tau}}$ holds if and only
if the partition of $V(G)$, $\{S, V(G) \backslash S\}$, is
an equitable partition.
\end{thm}

A proof of this theorem can be found in~\cite{mikesthesis}, where it is
attributed to Godsil. This proof uses the fact that the matrix $A-\tau
I$ is a positive semi-definite matrix.

From this theorem, the characteristic vector of a ratio-tight independent
set is a linear combination of the all-ones vector, {\bf 1}, and an
eigenvector corresponding to the smallest eigenvalue.
Also, for an independent set $S$ of size $s$, if the partition $\{S, V(G) \backslash S\}$
of $V(G)$ is equitable, then the smallest
eigenvalue of $G$ is $\tau = \frac{d}{1-\frac{n}{s}}$.

There is another ratio bound on the size of the maximum clique in a
graph. The set of graphs where the ratio bound on the clique size is
valid is more restrictive than for the ratio bound on the size of the
maximum independent set.

\begin{thm}[Ratio Bound for $\omega(G)$, see~\cite{Godsil:2003}]\label{thm:tauomegabound}
\index{ratio bound for $\omega(G)$} 
Assume $G$ is a $d$-regular graph 
which is either arc transitive or a single graph in an association scheme.
Let $\tau$ be the least eigenvalue of $G$, then
\[
\omega(G) \leq 1-\frac{d}{\tau}.
\]
\end{thm}

In Chapter~\ref{chp:ass}, we use the ratio bounds to find bounds on
the size of the maximum independent sets and the maximum cliques for a
family of {\em qualitative independence graphs} and several {\em uniform
qualitative independence graphs}. These graphs are defined in the next
chapter.

%%%%%%%%%%%%%%%%%%%%%%%%%%%%%%%%%%%%%%%%%%%%%%%%%%%%%%%%%%%%%%%%%%ay??

\remove{
In testing applications, each row in the array represents a particular
component of the system being tested. It may be an input variable, a
network node, a subroutine or a hardware component.  Each column in
the array corresponds to a test on the system. The goal is to produce
an array with the fewest number of columns, hence tests. Strength two
covering arrays test all pairwise interactions. This requires far
fewer tests than complete testing but in practice provides good test
coverage \cite{cohen:96,DM}.

}

%%%%%%%%%%%%%%%%%%%%%%%%%%%%%%%%%%%%%%%%%%%%%%%%%%%%%%%%%%%%%%%%%%
\chapter[Covering Arrays on Graphs]{Covering Arrays on Graphs and Qualitative Independence Graphs}
\label{CAG}

\setcounter{thm}{0}

In this chapter, we extend the definition of a covering array to
include a graph structure.  Recall from Definition~\ref{cak}, for
positive integers $n,r,k$, a covering array, $CA(n,r,k)$, is an $r
\times n$ array with entries from $\{0,1, \dots, k-1\}$ with the
property that any two rows in the arrays are qualitatively
independent.  One application of covering arrays is to design test
suites for systems or networks (see Section~\ref{CA}). In such an
application, each row of the array corresponds to a parameter in the
system, each column corresponds to a test run, and the entries from
$\{0,1, \dots, k-1\}$ correspond to the values the parameters are
assigned in the test run.  Such a test suite completely tests any pair
of parameters in the system against one another.

If, in a system, it is not necessary to test all pairs of parameters
against each other, a graph can be used to describe which pairs of
parameters need to be tested. In particular, each vertex in such a graph
corresponds to a parameter in the system and vertices are adjacent if
and only if the corresponding pairs of parameters need to be tested
against each other. Adding such a graph structure to covering arrays
makes it possible to build test suites that are more efficient for a
specific system and provides a way to use the internal structure of
the system to optimize covering arrays.

\begin{defn}[Covering Array on a Graph] \label{caG}\index{covering array on a graph}
Let $G$ be a graph and $n$ and $k$ be positive integers.  A {\em covering
array on the graph G}, $CA(n,G,k)$, is a $|V(G)| \times n$ array with
entries from $\{0,1,\dots, k-1\}$ whose rows correspond to the vertices of
$G$ and any two rows corresponding to vertices
which are adjacent in $G$ are qualitatively independent.
\end{defn}

The size of the smallest covering array on a graph $G$ with alphabet $k$
will be denoted by $CAN(G,k)$, that is,
\[
CAN(G,k)=\min_{l \in \mathbb{N}} \{l: \exists \; CA(l,G,k) \}.
\]
A $CA(n,G,k)$ with $n = CAN(G,k)$ is an \index{optimal covering array
on a graph}{\em optimal} covering array on $G$.

\begin{exmp}
An optimal covering array on a graph $G$ with $CAN(G,2) = 4$.
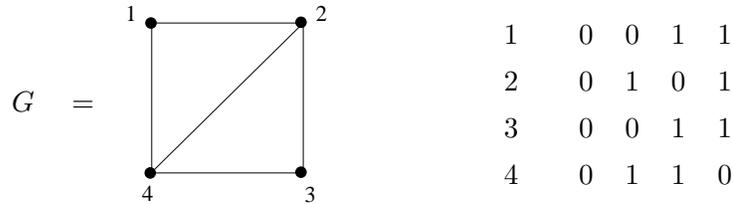
\begin{center}
\begin{tabular}{cc@{\extracolsep{2cm}} c}
\raisebox{1.2cm}{$G\quad =$} &
\psfig{figure=new_square.eps,height=2.75cm}
&
\raisebox{1.2cm}{
\begin{tabular}{lcccc}
1\quad \quad& 0 & 0 & 1 & 1 \\
2 & 0 & 1 & 0 & 1 \\
3 & 0 & 0 & 1 & 1 \\
4 & 0 & 1 & 1 & 0 \\
\end{tabular}
} 
\end{tabular}
\end{center}

\end{exmp}

Covering arrays on graphs are extensions of standard covering arrays;
in particular, for $K_k$ the complete graph on $k$ vertices, 
$CAN(k,g) = CAN(K_k,g)$. 

Binary covering arrays on graphs have been studied by Seroussi and
Bshouty, who proved that determining the existence of an optimal binary
covering array on a graph is an NP-complete
problem~\cite{Seroussi:01}.  General covering arrays on graphs are
introduced in the conclusion of Stevens's thesis
\cite{brett:thesis}.

In Section~\ref{sec:homomor}, we show that for all graphs $G$ and all positive integers $k$,
\begin{equation*}
CAN(K_{\omega(G)},k)    \leq CAN(G,k) \leq  CAN(K_{\chi(G)},k).\label{clique:col}
\end{equation*}
The upper bound is of particular interest because it gives a method
to construct covering arrays on graphs (see
Section~\ref{sec:homomor}). This also raises a question that partially motivates
this work: can determining $CAN(G,k)$ be reduced to determining
$\chi(G)$ and $CAN( \chi(G) ,k)$? We show that it cannot.  In
particular, we look for graphs $G$ so that
\begin{equation}
CAN(G,k) <  CAN(K_{\chi(G)},k)\label{chi:bound}.
\end{equation}

In Section~\ref{sec:QIgraphs}, we define a family of graphs, called
the {\em qualitative independence graphs}.  We show that this family
of graphs gives a good characterization of covering arrays for all
graphs, namely that for a graph $G$ and positive integers $k$ and $n$,
a $CA(n,G,k)$ exists if and only if there is a graph homomorphism $G
\rightarrow QI(n,k)$.  Moreover, for all positive integers $k$ and
$n$, a clique of size $r$ in the graph $QI(n,k)$ corresponds to a
$CA(n,r,k)$.  This new family converts a problem in combinatorial
systems into a question about homomorphisms on
graphs~\cite{MR2002f:05002}.

In Section~\ref{sec:cores}, we show
\[
\omega(QI(n,2)) = {n-1 \choose \lfloor \frac{n}{2} \rfloor-1}
\mbox{   and   }
\chi(QI(n,2)) = \bigg{\lceil} \frac{1}{2} {n \choose \lfloor \frac{n}{2}
\rfloor} \bigg{\rceil}.
\]
Using this formula for chromatic number, we show that if $n$ is odd,
then the graph $QI(n,2)$ satisfies
Inequality~(\ref{chi:bound}). Additionally, we present necessary
conditions on a graph $G$ to satisfy Inequality~(\ref{chi:bound}) for
$k=2$. These are: if $CAN(K_{\chi(G)},2) = c$, then $c$ must be even
and
\[ 
\frac{1}{2}{c-2 \choose  \frac{c-2}{2}} 
< \chi(G)
\leq \bigg{\lceil} \frac{1}{2} {c-1 \choose  \frac{c}{2}-1}
\bigg{\rceil}. 
\]

In Section~\ref{sec:QIcores}, we give a graph that is a core of
$QI(n,2)$. The structure of these cores implies that if there exists a
$CA(n,r,2)$ or a $CA(n,G,2)$, then there exists a covering array with
the same parameters in which each row has exactly $\lfloor
\frac{n}{2} \rfloor$ ones. 

Finally, in Section~\ref{sec:QIk2k}, for a positive integer $k$, we
find upper bounds on $\omega( QI(k^2,k) )$ and $\chi( QI(k^2,k) )$ and
prove that $QI(k^2,k)$ is a $(k!)^{k-1}$-regular graph.

Most of the results in this chapter are published in~\cite{brett:ca5}.

\section{Bounds from Homomorphisms}\label{sec:homomor}

From Section~\ref{sec:cliquesandchromatic}, for graphs $G$ and $H$, if
there is a graph homomorphism $G \rightarrow H$, then $\chi(G) \leq
\chi(H)$ and $\omega(G) \leq \omega(H)$.  We can obtain a similar
bound on the optimal size of a covering array.

\begin{thm}\label{thm:hombound}
Let $k$ be a positive integer and $G$ and $H$ be graphs.  If there is a
graph homomorphism $\phi:G \rightarrow H$, then
\[
  CAN(G,k) \leq   CAN(H,k).
\]
\end{thm}
\begin{proof}
For $n$ a positive integer, assume that there exists a covering array
$CA(n,H,k)$. This covering array will be used to construct a covering
array $CA(n,G,k)$. For each $i \in \{1,\dots,|V(G)|\}$, row $i$ of
$CA(n,G,k)$ corresponds to a vertex $v_i \in V(G)$.  Set row $i$ of
$CA(n,G,k)$ to be identical to the row corresponding to the vertex
$\phi(v_i) \in V(H)$ in $CA(n,H,k)$. Since $\phi$ is a homomorphism,
any pair of adjacent vertices in $G$ are mapped to adjacent vertices
in $H$.  Thus, any pair of rows in $CA(n,G,k)$, which correspond to
vertices that are adjacent in $G$, will be qualitatively independent,
since the rows are qualitatively independent in $CA(n,H,k)$.
\end{proof}

For any graph $G$, there are homomorphisms between the following complete graphs
\[
K_{\omega(G)} \rightarrow G \rightarrow K_{\chi(G)}.
\]
These homomorphisms can be used to find bounds on $CAN(G,k)$.

\begin{cor}\label{cor:completebounds}
For all positive integers $k$ and all graphs $G$, 
\[
CAN(K_{\omega(G)},k)    \leq CAN(G,k) \leq  CAN(K_{\chi(G)},k).
\]
\end{cor}

\section{Qualitative Independence Graphs}\label{sec:QIgraphs}

In Section~\ref{sec:defnQIgraphs}, we define the {\em qualitative
independence graphs}. The vertex set for these graphs is a set of
partitions, so before defining the qualitative independence graphs, we
need to extend the definition of qualitative independence to
partitions.

We will use the notation for set partitions from
Section~\ref{setpartitions}. In particular, a $k$-partition of an
$n$-set is a set of $k$ disjoint non-empty classes whose union is the
$n$-set and $\mathcal{P}^n_k$ denotes the set of all $k$-partitions of
an $n$-set.

\subsection{Qualitatively Independent Partitions}\label{sec:qipartitions}

We have defined qualitative independence for vectors
(Definition~\ref{defn:qi1}) and for sets
(Definition~\ref{defn:qi2}). There is a third (also equivalent)
definition of qualitative independence for partitions.
   
\begin{defn}[Qualitatively Independent Partitions]
\label{defn:qi3}\index{qualitatively independent partitions}
Let $n$ and $k$ be positive integers with $n \geq k^2$.
Let $A,B \in \mathcal{P}^n_k$ be two $k$-partitions of an $n$-set.
Assume $A =\{A_1, A_2 ,\dots, A_k\}$
and $B = \{B_1, B_2 ,\dots, B_k\}$.
The partitions $A$ and $B$ are {\em qualitatively independent} if
\[
A_i \cap B_j \neq \emptyset \quad \textrm{for all $i$ and $j$.}
\]
\end{defn}

If $k$-partitions $A = \{ A_1, \dots, A_k\}$ and $B = \{B_1, \dots,
B_k\}$ are qualitatively independent, then for each $i \in \{
1,\dots,k\}$, $|A_i| \geq k$ and $|B_i| \geq k$.

For all positive integers $k$ and $n$, this definition for qualitatively
independent $k$-partitions of an $n$-set is equivalent to
Definition~\ref{defn:qi1} for qualitatively independent vectors in
$\mathbb{Z}_k^n$.  Assume vectors $u, v \in\mathbb{Z}_k^n$, with
$u = (u_1,u_2,\dots ,u_n)$ and $v=(v_1,v_2,\dots,v_n)$, are
qualitatively independent.  For the vector $u$, define a $k$-partition $P =
\{P_1,P_2,\dots, P_k\}$ by $a \in P_i$ if and only if $u_a = i$. Similarly
let $Q$ be the partition corresponding to the vector $v$.  Since $u$
and $v$ are qualitatively independent, for all ordered pairs $(i,j)
\in \mathbb{Z}_k \times \mathbb{Z}_k$ there is an index $a$ such that
$(u_a,v_a) = (i,j)$.  By the definition of $P$ and $Q$, $a \in P_i$
and $a \in Q_j$. So, for all $i,j \in \{1,\dots, n\}$, $P_i \cap Q_j
\neq \emptyset$. Conversely, if $P$ and $Q$ are qualitatively
independent $k$-partitions, then the vectors $u,v \in \mathbb{Z}_k^n$
corresponding to $P$ and $Q$ are qualitatively independent.

\remove{a partition $P =\{ P_1,P_2,\dots,P_k \} \in \mathcal{P}^n_k$ can be
used to define a vector $v \in\mathbb{Z}_k^n$ by $v_a = i$ if and only
if $a \in P_i$. Assume partitions $P,Q \in \mathcal{P}^n_k$ are qualitatively
independent. Let $u,v \in \mathbb{Z}_k^n$ be the vectors corresponding
to $P$ and $Q$, respectively.  Since $P$ and $Q$ are qualitatively
independent for any $i,j \in \{1, \dots, n\}$, there exists some
element $a \in P_i \cap Q_j$. By definition of $u,v$, $u_a = i$ and
$v_a =j$. So for any $i,j \in \{1,\dots, n\}$ there exists a $a$, such
that $(u_a,v_a)=(i,j)$.}

\subsection{Definition of Qualitative Independence Graphs}\label{sec:defnQIgraphs}

\begin{defn}[Qualitative Independence Graph] \index{qualitative independence graph}\label{defn:qigraphs}
Let $n$ and $k$ be positive integers with $n \geq k^2$. Define the
{\em qualitative independence graph} $QI(n,k)$ to be the graph whose
vertex set is the set of all $k$-partitions of an $n$-set with the
property that every class of the partition has size at least $k$.
Vertices are adjacent if and only if the corresponding partitions are
qualitatively independent.
\end{defn}

\begin{exmp}
Consider the graph $QI(4,2)$.  There are only three partitions in
$\mathcal{P}^4_2$ with the property that every class has size at least two:
\[
1\;2\;\mid\;3\;4, \quad 1\;3\;\mid\;2\;4, \quad \;1\;4\;\mid\;2\;3. 
\]
These partitions are all pairwise qualitatively independent, so the
graph $QI(4,2)$ is isomorphic to the complete graph on three vertices,
$K_3$.

\begin{figure}
\centerline{\psfig{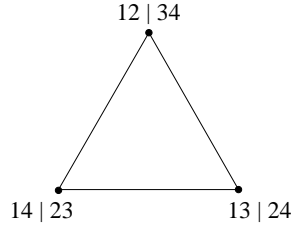}}
\caption{The graph $QI(4,2)$.}
\end{figure}

\end{exmp}

It would be possible to define the vertex set of $QI(n,k)$ to be the
set of all $k$-partitions of an $n$-set, $\mathcal{P}^n_k$.  But then,
any $k$-partition with a class of size smaller than $k$ would be an
isolated vertex in $QI(n,k)$.  Since the vertex set of $QI(n,k)$
excludes these partitions, the qualitative independence graphs have no
isolated vertices, in fact, they are connected and have small
diameter.

\begin{lemma}\label{QIdiameter2}
Let $k,n$ be positive integers with $n\geq k^2$. The graph
$QI(n,k)$ has diameter 2.
\end{lemma}
\begin{proof}
Given any two partitions, $P,Q\in V(QI(n,k))$, it is possible to construct a
partition $R \in V(QI(n,k))$ which is qualitatively independent from $P$ and $Q$.

Define a bipartite multi-graph $H$ as follows: each part of $H$ has $k$
vertices. The $k$ vertices in the first part of $H$ correspond to the
classes $P_i$ of $P$, and the $k$ vertices in the second part of $H$
correspond to the classes $Q_i$ of $Q$. 

For every $a \in \{1,\dots,n\}$, there exists a unique class $P_{i_a} \in
P$ with $a \in P_{i_a}$ and there exists a unique class $Q_{i_a} \in Q$ with
$a \in Q_{i_a}$. For each such $a$, add an edge to $H$ between the
vertices corresponding to $P_{i_a}$ and $Q_{i_a}$. Label this edge by $a$. The
multi-graph $H$ has exactly $n$ edges, one for each $a \in
\{1,\dots,n\}$. Further, each class $P_i \in P$ and each class $Q_i
\in Q$ have cardinality at least $k$. Thus, the degree of each
vertex in $H$ is at least $k$. By a dual of K\"{o}nig's Theorem~\cite{MR80d:05023},
there are $k$ edge-disjoint 1-factors in $H$.

Use these $k$ 1-factors to build the $k$-partition $R$. Each class
of $R$ will correspond to a 1-factor of $H$. In particular, $a \in
R_i$ for $i=1,\dots,k$ if and only if the edge labelled by $a$ is in
the 1-factor that corresponds to the class $R_i \in R$. Place the remaining elements
from $\{1,\dots ,n\}$ in classes of $R$ arbitrarily. Each 1-factor has
exactly $k$ edges, so the classes of $R$ will have size at least $k$.
Thus, $R$ is a $k$-partition of an $n$-set with each class of size at least
$k$, so $R \in V(QI(n,k))$.

Finally, we need to show that the partitions $P$ and $R$ are
qualitatively independent. For any $j =1,\dots, k$, the class $R_j$
corresponds to a 1-factor of $H$, this means that for any $i \in \{1,
\dots, k\}$, there is an edge $a$ in the 1-factor incident to $P_i$.
By the definition of $R$, $a \in R_j$. Since the edge labelled $a$ is
incident with the vertex corresponding to $P_i$, we also know that $a
\in P_i$. This means for any $i,j \in \{1, \dots, n\}$, $a \in P_i
\cap R_j$ so $P_i \cap R_j \neq \emptyset$. Thus $P$ and $R$ are
qualitatively independent. Similarly, $Q$ and $R$ are qualitatively
independent.
\end{proof}

A clique in $QI(n,k)$ is equivalent to a collection of pairwise
qualitatively independent $k$-partitions of an $n$-set. Thus, an
$r$-clique in $QI(n,k)$ corresponds to a $CA(n,r,k)$. An independent
set in $QI(n,k)$ is a collection of partitions in which no two are
qualitatively independent.

For any positive integers $k,n$, it is possible to construct a $CA(n,
QI(n,k), k)$.  Each row of a $CA(n,QI(n,k),k)$ corresponds to a
partition in $V(QI(n,k))$, which in turn corresponds to a vector in
$\mathbb{Z}_k^n$. Set each row in $CA(n, QI(n,k), k)$ to be the
corresponding vector. Two partitions are adjacent in $QI(n,k)$ if and
only if they are qualitatively independent and the vectors
corresponding to the partitions are qualitatively independent if and
only if the partitions are. Thus, two rows in $CA(n, QI(n,k), k)$ are
qualitatively independent if and only if the partitions corresponding
to the rows are adjacent in $QI(n,k)$.

\begin{lemma}\label{lem:nbound}
For positive integers $n,k$ with $n \geq k^2$, we have 
\[
CAN(QI(n,k) , k) \leq n.
\]
\end{lemma}

We conjecture that a $CA(n,QI(n,k),k)$ is an optimal covering array.

\begin{conj}\label{conj:canQIn}
For positive integers $n,k$ with $n \geq k^2$, we have 
\[
CAN(QI(n,k) , k) = n.
\]
\end{conj}

\remove{something about QI(10,3)?}

\subsection{Homomorphisms and Qualitative Independence Graphs}\label{sec:gencolouring}

Recall from Section~\ref{sec:cliquesandchromatic}, that a proper
colouring with $n$ colours of a graph $G$ is equivalent to a
homomorphism from $G$ to the complete graph $K_n$. We give a similar
characterization for covering arrays on graphs.

\begin{thm} \label{crysthom:lem}
For a graph $G$ and positive integers $k$ and $n$, a $CA(n,G,k)$
exists if and only if there exists a graph homomorphism $G
\rightarrow QI(n,k)$.
\end{thm}
\begin{proof}
Assume that there exists a $CA(n,G,k)$, call this $C$.  For a vertex
$v \in V(G)$, let the vector $C_v \in \mathbb{Z}_k^n$ be the row in
$C$ corresponding to $v$.  Consider a mapping $\phi:V(G) \rightarrow
V(QI(n,k))$ which takes a vertex $v \in V(G)$ to the partition $P_v
\in V(QI(n,k))$ which corresponds to the vector $C_v$.

The map $\phi$ is a homomorphism. To see this, let $v,w \in
V(G)$ be adjacent vertices. Since $C$ is a covering array, the vectors
$C_v$ and $C_w$ are qualitatively independent, thus the corresponding
partitions $P_v$ and $P_w$ are qualitatively independent and adjacent
in $QI(n,k)$.

Conversely, assume there is a graph homomorphism $\phi: G \rightarrow
QI(n,k)$.  For each $v \in V(G)$, $\phi(v)$ is a $k$-partition and 
has a corresponding vector $C_v \in \mathbb{Z}_k^n$.  Build a
$CA(n,G,k)$ by using the vector $C_v$ as the row corresponding to $v
\in V(G)$.  If the vertices $v,w \in V(G)$ are adjacent in
$G$, then the vertices corresponding to partitions $\phi(v),\phi(w)$
are adjacent in $QI(n,k)$. This means the partitions $\phi(v)$,
$\phi(w)$ are qualitatively independent, thus the corresponding
vectors $C_v,C_w$ are qualitatively independent.
\end{proof}

\begin{exmp}\label{exmp:qi42}
For a graph $G$, there exists a $CA(4,G,2)$ if and only if there is a
homomorphism from $G$ to $QI(4,2)$. Since $QI(4,2)$ is isomorphic to
$K_3$, there exists a $CA(4,G,2)$ if and only if $G$
is 3-colourable.  This gives a characterization of the graphs for
which a covering array of size 4 on a 2-alphabet exists (they are the 3-colourable graphs).

Further, this shows that determining if $CAN(G,2) = 3$ is as hard as
determining whether a graph is 3-colourable.  This is the approach
used by Seroussi and Bshouty~\cite{Seroussi:01} to prove that finding
$CAN(G,k)$ for a given $G$ and $k$ is an NP-hard problem.
\end{exmp}

Similar to the definition of chromatic number for a graph $G$,
\[
\chi(G) = \min \{ n : G \rightarrow K_n \},
\]
and the characterization of fractional chromatic number
(Theorem~\ref{thm:fractionalchromatic}),
\[
\chi^\ast(G) = \min \{ n/r : G \rightarrow K_{n:r} \},
\]
$CAN(G,k)$ can be characterized by a homomorphism.

\begin{cor}
For any graph $G$, and any positive integer $k$,
\[
CAN(G,k) = \min_{n\in\mathbb{N}} \{ n \; : \; G \rightarrow  QI(n,k) \}.
\]
\end{cor}

Knowing the chromatic number and the maximum clique size of $QI(n,k)$
will give information on which graphs have covering arrays of size $n$
on an alphabet of size $k$.

\begin{cor}
Let $G$ be a graph and $n,k$ be positive integers. If there exists a $CA(n,G,k)$, then
\[
\chi(G) \leq \chi(QI(n,k)) \mbox{ and } \omega(G) \leq \omega(QI(n,k)). 
\]
\end{cor}

With Theorem~\ref{crysthom:lem}, Conjecture~\ref{conj:canQIn} can be
rephrased in terms of homomorphisms.
\begin{conj}\label{conj:nohomo}
For all integers $n,k$, there is no homomorphism $QI(n,k) \rightarrow QI(n-1,k)$.
\end{conj}

It is in general difficult to prove that no homomorphism exists
between two graphs. One method to prove this would be to show either
$\chi(QI(n,k)) >  \chi(QI(n-1,k))$ or $\omega(QI(n,k)) >  \omega(QI(n-1,k))$.

With this motivation, we try to find bounds and exact values for
$\chi(QI(n,k))$ and $\omega(QI(n,k))$.

\section{Binary Qualitative Independence Graphs}\label{sec:binary}

From Example~\ref{exmp:qi42}, for a graph $G$, a $CA(4,G,2)$ exists if and
only if $G$ is 3-colourable. It is not true in general that the
existence of a $CA(n,G,k)$ can be characterized by the chromatic
number of $G$. To see this we consider the graph $QI(5,2)$.

\begin{exmp}\index{$QI(5,2)$}
There are 10 vertices in $\mathcal{P}^5_2$ that have each class of size at least 2:

\begin{tabular}{ccccc}
1\;2\; $\mid$ \;3\;4\;5, &  1\;3\; $\mid$ \;2\;4\;5, & 
1\;4\; $\mid$ \;2\;3\;5,  & 1\;5\; $\mid$  \;2\;3\;4, & 
1\;2\;3\; $\mid$ \;4\;5, \\
1\;2\;4\; $\mid$ \;3\;5, & 1\;2\;5\; $\mid$ \;3\;4, & 
1\;3\;4\; $\mid$ \;2\;5, & 1\;3\;5\; $\mid$ \;2\;4, & 
1\;4\;5\; $\mid$ \;2\;3.
\end{tabular}
\vspace{.25cm}

A representation of $QI(5,2)$ is given below. The vertices
1\;2\;5\;$\mid$\;3\;4 and 1\;3\;4\;$\mid$\;2\;5 are repeated (in
grey) to make the graph easier to read.

\begin{figure}
\centerline{\psfig{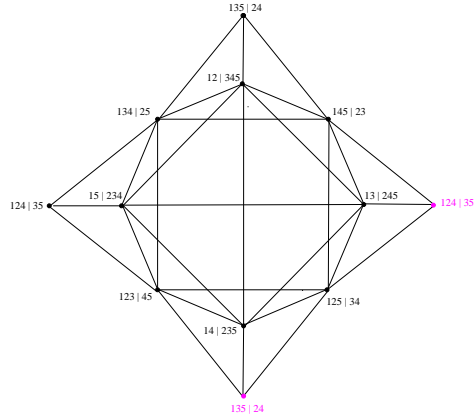}}
\caption{The graph $QI(5,2)$.}
\end{figure}

By inspection, $\omega(QI(5,2)) = 4$ and $\chi(QI(5,2))=5$. A 5-colouring
of $QI(5,2)$ is given below.
\begin{center}
\begin{tabular}{|ccc|} \hline
colour  & \multicolumn{2}{c|}{vertices}  \\
\hline
1 & 00111 & 00011 \\
2 & 01011 & 01010 \\
3 & 01101 & 00101 \\
4 & 01110 & 01100 \\
5 & 01001 & 00110 \\
\hline
\end{tabular}
\end{center}

From Theorem~\ref{crysthom:lem}, $QI(5,2)$ has a covering array of
size 5.  This graph does not have a covering array of size 4 since
otherwise there would exist a homomorphism $QI(5,2) \rightarrow QI(4,2)$,
contradicting that
$\chi(QI(5,2))>3$. Thus, $CAN( QI(5,2),2) = 5$.  From
Theorem~\ref{thm:binaryCAN}, $CAN(K_5,2) = CAN(5,2) =6$.

Therefore, for $QI(5,2)$, the inequality from
Corollary~\ref{cor:completebounds} holds strictly,
\[ 5 = CAN(QI(5,2),2) <  CAN(K_{\chi(QI(5,2))},2) =  6. \]

If for a graph $G$ a $CA(5,G,2)$ exists, then there is a homomorphism
from $G$ to $QI(5,2)$. Since $QI(5,2)$ is not isomorphic to a complete
graph, determining if a graph has a binary covering array of size 5 is
not equivalent to finding a proper colouring of the graph.

Finally, it is interesting to note that the complement of $QI(5,2)$ is
the Petersen graph, which is also known as the Kneser graph $K_{5:2}$.
\end{exmp}

\subsection{Formulae for $\omega(QI(n,2))$ and $\chi(QI(n,2))$}
\label{sec:cores}

In this section, the following formulae for maximum clique size
and chromatic number of $QI(n,2)$ are established:
\[
\omega(QI(n,2)) = {n-1 \choose \lfloor \frac{n}{2} \rfloor - 1 } 
\quad \textrm{ and } \quad
 \chi(QI(n,2))   =  \left\lceil \frac{1}{2} 
                     {n \choose  \lfloor \frac{n}{2} \rfloor } 
                   \right\rceil. 
\]

The vertices of $QI(n,2)$ are 2-partitions of an $n$-set with the
property that each class of the partition has at least 2 elements.
For any 2-partition $P= \{ P_1, P_2 \}$, it is clear that $P_2 =
\overline{P_1}$.  So, each 2-partition can be unambiguously described
by a subset of an $n$-set. Moreover, by choosing the smaller of $P_1$
and $P_2$, any 2-partition of an $n$-set can be described as an
$r$-subset of an $n$-set where $r \leq n/2$.

Let $P= \{ P_1, P_2\}$ and $Q=\{ Q_1, Q_2 \}$ be 2-partitions of an
$n$-set. The partitions $P$ and $Q$ are qualitatively independent if
and only if the sets $P_1$ and $Q_1$ are intersecting and incomparable
(Section~\ref{sect:EKRapps}). The exact value of the maximum clique
for $QI(n,2)$ is the maximum size of a 1-intersecting set system; this
is known from Sperner's Theorem and the Erd\H{o}s-Ko-Rado Theorem. In
fact, Theorem~\ref{thm:binaryCAN} can be restated in terms of a
maximum clique in $QI(n,2)$.
 
\begin{thm}[\cite{kleitman:73}]\label{thm:maxclique}
For integer $n\geq 4$,
\[
\omega(QI(n,2)) = \max_{r \in \mathbb{N}} \{r : \exists \;CA(n,K_{r},2)\} 
               =  {n-1 \choose \left \lfloor \frac{n}{2} \right \rfloor - 1}.
\]
\end{thm}

When $n$ is even, the set of uniform 2-partitions of an $n$-set is a
maximum clique.  When $n$ is odd the set of all almost-uniform
2-partitions of an $n$-set that have a common element in the smaller
class is a maximum clique.

Next, we determine the chromatic number of the graphs $QI(n,2)$.

\begin{thm} \label{thm:QIchi}
For all positive integers $n$,
\[ 
\chi ( QI(n,2) ) = \bigg\lceil \frac{1}{2} {n \choose \lfloor \frac{n}{2} \rfloor }\bigg\rceil.
\]
\end{thm}
\begin{proof} 
Consider the vertices in $QI(n,2)$ not as partitions, but as
subsets of an $n$-set of size no more than $\lfloor n/2 \rfloor$ (as
described above).

{}From the proof of Sperner's Theorem (Theorem~\ref{thm:sperner}), the poset of
subsets of an $n$-set ordered by inclusion can be decomposed into ${n
\choose \lfloor n/2 \rfloor}$ disjoint chains and each chain contains
exactly one set of size $\lfloor n/2 \rfloor$. Call these chains $C_i$,
where $i \in \{1,\dots, {n \choose \lfloor n/2 \rfloor} \}$.  For any
$i$, if the sets $A, B \in C_i$, then $A$ and $B$ are not qualitatively
independent.  In particular, any two vertices in $QI(n,k)$ which
correspond to sets that are in the same chains are not adjacent in
$QI(n,k)$.

It is possible to pair the ${n \choose \lfloor n/2 \rfloor}$ chains so
that any subset of size no more than $n/2$ in one chain is disjoint
from any subset of size no more than $n/2$ in the other chain. To see
this, consider two cases, first when $n$ is even and second when $n$ is
odd.

Assume $n$ is even. For each chain $C_i$, let $A_i$ be the set of size
$n/2$. Match the chains $C_i$ and $C_j$ where $A_i = \overline{A_j}$.

Assume $n$ is odd.  For each chain $C_i$, let $A_i$ be the set of size
$(n-1)/2$.  The sets $A_i$ are the vertices of the Kneser graph
$K_{n:\frac{n-1}{2}}$ (see Section~\ref{knesergraphs}). The graph
$K_{n:\frac{n-1}{2}}$ is vertex transitive so there exists a matching
that is perfect or is missing just one vertex (Section 3.5
of~\cite{MR2002f:05002}). So each set $A_i$ (except possibly one set)
is matched to another set of size $\frac{n-1}{2}$, call it $A_i'$.
The set $A_i' \subset \overline{A_i}$.  Match the chain $C_i$ which
contains $A_i$ with the chain $C_{i'}$ that contains the set $A_i'$.

Any two sets in a matched pair of chains have the property that either
one set contains the other or one set contains the complement of the
other. In either case, the partitions are not qualitatively independent.
All sets in the paired chains can be assigned the same colour in a
proper colouring of $QI(n,2)$. This produces a proper 
$\left \lceil \frac{1}{2} {n \choose \lfloor n/2 \rfloor } \right \rceil$-colouring on the
graph $QI(n,2)$.

To see that this is the smallest possible colouring of $QI(n,2)$,
consider the vertices of $QI(n,2)$ that correspond to $\lfloor n/2
\rfloor$-sets.  Two such vertices may be assigned the same colour if
and only if the subsets are disjoint. It is clear that there can not
be three mutually disjoint subsets of an $n$-set with size $\lfloor n/2
\rfloor$. So it is not possible to properly colour these vertices with
fewer than $\left\lceil \frac{1}{2} {n \choose \lfloor n/2 \rfloor } \right
\rceil$ colours.
\end{proof}

From the formulae for maximum clique size and chromatic number, for
all $n \geq 4$,
\[
\chi( QI(n-1,2) ) < \omega ( QI(n,2)) 
        \leq \chi(QI(n,2)) < \omega(QI(n+1,2)).
\]
Note that when $n$ is even, we have $ \omega ( QI(n,2)) = \chi(QI(n,2))$.
From this inequality we have the following result, which confirms
Conjecture~\ref{conj:canQIn} for $k=2$.

\begin{cor}\label{caGn22isn}
For all $n \geq 4$, we have $CAN(QI(n,2),2) = n$.

\end{cor}
\begin{proof}
Assume that there exists a $CA(n-1,QI(n,2),2)$, then by
Theorem~\ref{crysthom:lem} there is a homomorphism
\[
QI(n,2) \rightarrow QI(n-1,2),
\]
which contradicts $\chi( QI(n-1,2) ) < \chi(QI(n,2))$.

From Theorem~\ref{crysthom:lem}, there exists a $CA(n, QI(n,2),2)$, so 
\[
CAN(QI(n,2),2) = n.
\] 
\end{proof}

With the exact value of $\chi(QI(n,2))$ we can prove that the second
inequality in Corollary~\ref{cor:completebounds} is strict for
$QI(n,2)$ for all odd $n$.

\begin{cor}
For $n$ odd,
\[ 
CAN(QI(n,2),2)  <  CAN(K_{\chi(QI(n,2))},2). 
\]
\end{cor}
\begin{proof}
{}From Corollary~\ref{caGn22isn}, $CAN(QI(n,2),2) = n$.
{}From Theorem~\ref{thm:QIchi} and Theorem~\ref{thm:binaryCAN},
\begin{eqnarray*}
CAN(K_{\chi(QI(n,2))},2) 
&=& CAN\bigg( \bigg\lceil \frac{1}{2} {n \choose \frac{n-1}{2} } \bigg\rceil, 2\bigg)  \\
&=& \min \bigg\{ m : {m-1 \choose \lfloor m/2 \rfloor-1 } \geq \bigg\lceil \frac{1}{2} {n \choose \frac{n-1}{2} } \bigg\rceil \bigg\} \\
&=& n+1.
\end{eqnarray*}
\end{proof}

This gives an infinite family of graphs for which
Inequality~(\ref{chi:bound}) holds with strict inequality.  Further, we
can use this to get a lower bound on $CAN(G,2)$.

\begin{cor}
For any graph $G$  
\[
CAN(K_{ \chi(G) },2)-1 \leq CAN(G,2) \leq CAN(K_{ \chi(G) },2).
\]
Moreover, if $CAN(K_{ \chi(G) },2)$ is odd, then 
\[
CAN(G,2) = CAN(K_{ \chi(G) },2).
\]
\end{cor}

\begin{proof}
{}From Corollary~\ref{cor:completebounds}, $CAN(G, 2) \leq CAN(K_{\chi(G)}, 2)$. 

Assume that $m = CAN(K_{\chi(G)}, 2)$ is even and that $CAN(G,2) \leq
m - 2$.  Then there is a homomorphism from $G$ to $QI(m-2,2)$ and
\[
\chi(G) \leq \chi(QI(m-2,2)) = \frac{1}{2}{ m-2 \choose \frac{m-2}{2}}.
\]
Since 
\[
CAN(K_{ \chi(G) },2) = \min \bigg\{l : {l-1 \choose \lfloor \frac{l}{2} \rfloor -1} 
                          \geq \chi(G)\bigg\},
\]
$CAN(K_{ \chi(G) },2) \leq m-1$.
This is a contradiction with $m=CAN(K_{ \chi(G) },2)$.

Next, assume that $CAN(K_{\chi(G)}, 2)$ is odd and set $m +1 = CAN(K_{\chi(G)}, 2) $.
If $CAN(G, 2) \leq m$, then there is a homomorphism from $G$ to $QI(m,2)$.
{}From Theorem~\ref{thm:QIchi}, $\chi(QI(m,2)) = \frac{1}{2}{m \choose m/2}={m-1 \choose m/{2}}$,
and $\chi(G) \leq {m-1 \choose m/2}$.
By definition,
\[
CAN( K_{\chi(G)}, 2) = \min \bigg\{ l : {l-1 \choose \lfloor \frac{l}{2} \rfloor - 1} \geq  \chi(G)\bigg\}.
\]
Since $\chi(G) \leq {m-1 \choose {m}/{2}}$, the minimum occurs when 
$l \leq m$. This means that 
\[
CAN(K_{\chi(G)}, 2) \leq m <CAN(K_{\chi(G)},2),
\]
and this contradiction gives that $CAN(G, 2) = CAN(K_{\chi(G)}, 2)$.
\end{proof}

\remove{The result of these theorems are two necessary conditions on a graph
$G$ for when $CAN(G,2) < CAN(K_{\chi(G)},2)$.
These are, for $c=CAN(K_{\chi(G)},2)$: 
\begin{enumerate}
\item $c$ must be even;
\item $ {c-2 \choose  \frac{c}{2}}  < \chi(G)
\leq \bigg{\lceil} \frac{1}{2} {c-1 \choose  \frac{c}{2}}
\bigg{\rceil}. $
\end{enumerate}}

If $CAN(G,2) < CAN(K_{\chi(G)},2)$, then by Corollary~\ref{cor:bounds}
$c=CAN(K_{\chi(G)}, 2)$ is even and 
\[
CAN(G,2) = CAN(K_{\chi(G)},2)-1.
\] 
This means there is a homomorphism from
$G$ to $QI(c-1,2)$. Since $\chi(QI(c-1,2)) = \lceil \frac{1}{2} {c-1
\choose \frac{c-2}{2}}\rceil$, we have that $\chi(G) \leq \lceil \frac{1}{2} {c-1
\choose \frac{c-2}{2}}\rceil$.

If $\chi(G) \leq \frac{1}{2}{c-2 \choose \frac{c-2}{2}}$, then there
is a homomorphism from $G$ to $K_{\frac{1}{2}{c-2 \choose
\frac{c-2}{2}}}$.  For $c$ even, $K_{\frac{1}{2}{c-2 \choose
\frac{c-2}{2}}}$ is isomorphic to the graph $UQI(c-2,2)$. Thus, there
is a homomorphism from $G$ to $UQI(c-2,2)$, and in particular,
$CAN(G,2) \leq c-2 = CAN(K_{\chi(G)},2)-2$, contradicting
Corollary~\ref{cor:bounds}.

Thus we have two necessary
conditions on $G$ for
\[
CAN(G,2) < CAN(K_{\chi(G)},2).
\] 
These are: 
\begin{enumerate}
\item $CAN(K_{\chi(G)}, 2)=c$ must be even, and
\item $ \frac{1}{2}{c-2 \choose  \frac{c-2}{2}} 
< \chi(G)
\leq \bigg{\lceil} \frac{1}{2} {c-1 \choose  \frac{c}{2}-1}
\bigg{\rceil}. $
\end{enumerate}

\subsection{Cores of the Binary Qualitative Independence Graphs}\label{sec:QIcores}

{}From Section~\ref{subsec:cores}, a core of a graph $G$ is an induced
subgraph, denoted $G^\bullet$, with the property that there exists a
homomorphism from $G$ to $G^\bullet$ and every endomorphism on
$G^\bullet$ is an automorphism.  A core of a graph is useful since it
preserves any property of the graph defined by homomorphisms (for
example, chromatic number, maximum clique size, and fractional
chromatic number) and for any graph $H$, if there exists a
homomorphism from $H$ to $G$, then there also exists a homomorphism
from $H$ to $G^\bullet$.

In this section, we find a core for the graphs $QI(n,2)$. First we
need some notation.  Recall from Section~\ref{setpartitions}, that a
{\em uniform partition} is a partition in which all the classes have
the same size and an {\em almost-uniform partition} is a partition in
which the sizes of the classes differ by at most one. For $n$ even, let
$UQI(n,2)$ be the subgraph of $QI(n,2)$ induced by the vertices that
correspond to uniform partitions. For $n$ odd, let $AUQI(n,2)$ be the
subgraph of $QI(n,2)$ induced by the vertices that correspond to
almost-uniform partitions. These graphs are further explored in
Section~\ref{UQI}. 

%karen intro here

\begin{prop}\label{core1}
For $n$ an even positive integer, there exists a homomorphism $QI(n,2) \rightarrow UQI(n,2)$.
For $n$ an odd positive integer, there exists a homomorphism $QI(n,2) \rightarrow AUQI(n,2)$.
\end{prop}

\begin{proof}
Recall from Theorem~\ref{thm:QIchi} that the poset of subsets of an $n$-set can
be decomposed into  ${n \choose \lfloor \frac{n}{2} \rfloor}$ disjoint chains.

For every $P \in V(QI(n,2))$ with $P=\{P_1,P_2\}$, assume $|P_1| \leq
|P_2|$.  The chain in the poset of subsets of an $n$-set that contains
$P_1$ contains a unique $\lfloor \frac{n}{2} \rfloor$-set, call this
set $P'_1$. Let $P'=\{P'_1,\overline{P'_1}\}$.  Define a map $\phi$
from $QI(n,2)$ to $UQI(n,2)$ (or $AUQI(n,2)$ if $n$ is odd) by
$\phi(P)= P'$.

The map $\phi$ is a homomorphism. To see this, assume partitions
$P=\{P_1,P_2\}$ and $Q=\{Q_1,Q_2\}$ are qualitatively
independent. Then $P_1 \cap Q_1 \neq \emptyset$, and $P_1$ and $Q_1$
are in distinct chains. Let $P'_1$ and $Q'_1$ be the unique sets of
size $\lfloor \frac{n}{2} \rfloor$ in the chains that contain $P_1$
and $Q_1$, respectively. Then $P'_1$ and $Q'_1$ are distinct, and
since $P \subseteq P'_1$ and $Q\subseteq Q'_1$, $P_1 \cap Q_1 \subseteq P'_1
\cap Q'_1$. Thus $P'_1\cap Q'_1 \neq \emptyset$. Finally, since $P'_1$
and $Q'_1$ are distinct $\lfloor \frac{n}{2} \rfloor$-sets, they are incomparable. Thus the
partitions $P'=\{P'_1,P'_2\}$ and $Q'=\{Q'_1,Q'_2\}$ are
qualitatively independent.
\end{proof}

From Theorem~\ref{crysthom:lem}, for a graph $G$ and positive integers
$k,n$, a $CA(n,G,k)$ exists if and only if there is a homomorphism $G
\rightarrow QI(n,k)$. Further, from Proposition~\ref{core1},
a $CA(n,G,2)$ exists if and only if there is a homomorphism from $G$ to
$UQI(n,2)$, if $n$ is even, or to $AUQI(n,2)$, if $n$ is odd.

\begin{thm}\label{balancedCA}
For a graph $G$ and a positive integer $n$, if there exists a
$CA(n,G,2)$, then it is possible to find a covering array
$CA(n,G,2)$ in which the rows have exactly $\lfloor {n/2} \rfloor$
0s. Moreover, if $n$ is even, then it is possible to find such a covering
array with the rows all beginning with 0.
\end{thm}

\begin{proof}
From Proposition~\ref{crysthom:lem}, if there exists a $CA(n,G,2)$,
then there exists a homomorphism $\phi:G \rightarrow QI(n,2)$.

If $n$ is even, then there is an endomorphism $\psi:QI(n,2) \rightarrow
UQI(n,2)$.  For each vector $v \in V(G)$, replace the row of $CA(n,G,2)$
corresponding to $v$ by the vector corresponding to the uniform
2-partition $\psi(\phi(v))$. If this vector does not start with 0,
simply relabel the 0s by 1s and the 1s by 0s.

If $n$ is odd, then there is an endomorphism $\psi:QI(n,2) \rightarrow
AUQI(n,2)$.  For each vector $v \in V(G)$, replace the row of
$CA(n,G,2)$ corresponding to $v$ by the vector corresponding to the
almost-uniform partition $\psi(\phi(v))$.
\end{proof}

In the next two theorems, we give a core of $QI(n,2)$.

\begin{thm}\label{core2}
For $n$ even, $UQI(n,2)$ (which is isomorphic to $K_{\frac{1}{2} {n
\choose n/2 }}$) is a core of $QI(n,2)$.
\end{thm}
\begin{proof}
For $n$ even, any two distinct uniform 2-partitions of an $n$-set are
qualitatively independent. Thus, the graph $UQI(n,2)$ is isomorphic to
$K_{\frac{1}{2} {n \choose n/2 }}$. From Theorem~\ref{core1} there is
a homomorphism from $QI(n,2)$ to $UQI(n,2)$.  Since
$K_{\frac{1}{2} {n \choose n/2 }}$ is a core, the theorem holds.
\end{proof}

\begin{thm}\label{core3}
For $n$ odd, $AUQI(n,2)$ is a core of $QI(n,2)$.
\end{thm}

\begin{proof}
From Theorem~\ref{core1}, there is a homomorphism from $QI(n,2)$ to
$AUQI(n,2)$.  All that is needed is to show that $AUQI(n,2)$
is itself a core.

Let $AUQI^{\bullet}(n,2)$ be a core of $AUQI(n,2)$. The graph
$AUQI(n,2)$ is vertex transitive. Thus, a core $AUQI^{\bullet}(n,2)$
is also vertex transitive (Theorem~\ref{thm:corevt}) and from
Theorem~\ref{thm:fibres},
\[
|V(AUQI^{\bullet}(n,2))| \quad \mathrm{divides} \quad 
|V( AUQI(n,2))| = {n \choose \frac{n-1}{2}}.
\]

Since the graph $AUQI(n,k)$ is an induced subgraph of $QI(n,k)$ and
there is a homomorphism from $QI(n,k)$ to $AUQI(n,k)$, we have that
$\chi(AUQI(n,k)) =\chi(QI(n,k))$.
As $AUQI^{\bullet}(n,2)$ is a core, we also have that 
\[ \chi(AUQI^\bullet(n,k)) =\chi(QI(n,k)).\]
Thus,
\[
\chi(AUQI^{\bullet}(n,2))   = \chi(AUQI(n,2))
                   =  \bigg{\lceil} \frac{1}{2} 
                      {n \choose  \frac{n-1}{2} } 
                      \bigg{\rceil}.
\]

These facts give two possibilities for $|V(AUQI^{\bullet}(n,2))|$, either it is 
$\frac{1}{2}{n \choose  \frac{n-1}{2}}$ or ${n \choose \frac{n-1}{2}}$.

If $|V(AUQI^{\bullet}(n,2))| = \frac{1}{2}{n \choose \frac{n-1}{2}}$
then, from the chromatic number of $AUQI^{\bullet}(n,2)$, the graph
$AUQI^{\bullet}(n,2)$ would have to be the complete graph. This is not
the case since
\[
\omega(AUQI^{\bullet}(n,2)) = \omega(AUQI(n,2)) 
  ={n-1 \choose \frac{n-1}{2} -1} < \frac{1}{2}{n \choose \frac{n-1}{2} }.
\]
Thus $|V(AUQI^{\bullet}(n,2))| ={ n \choose (n-1)/2}$, and as it is an
induced subgraph, $AUQI^{\bullet}(n,2) = AUQI(n,2)$.  This means $AUQI(n,2)$ is
a core and in particular, it is a core of $QI(n,2)$.
\end{proof}

For $n$ odd, the graph $AUQI(n,2)$ is isomorphic to the complement of
the Kneser graph $K_{n:\frac{n-1}{2}}$.

\section{The graphs $QI(k^2,k)$}\label{sec:QIk2k}

For $k$ a positive integer, the family of graphs $QI(k^2,k)$ is a
special family. First, for each $k$, the graph $QI(k^2,k)$ is the
qualitative independence graph with the smallest number of vertices.
Further, a clique in $QI(k^2,k)$ is equivalent to an orthogonal array,
thus from Corollary~\ref{OAprimepower}, for $k$ a prime power, we have
$\omega( QI(k^2,k) ) = k+1$, and from Corollary~\ref{OAnonprimepower}
we have $\omega( QI(k^2,k) ) \leq k+1$ in general.  We study some properties
of these graphs.

\begin{lemma}
For any integer $k\geq 2$, the qualitative independence graph
$QI(k^2,k)$ is vertex transitive.
\end{lemma}
\begin{proof}
The vertices of $QI(k^2,k)$ correspond to uniform
$k$-partitions of $\{1,\dots,k^2\}$.
For any pair of uniform $k$-partitions $P,Q$, there is a 
permutation in $\sigma \in Sym(k^2)$ such that $\sigma(P) =Q$.
\end{proof}

\begin{thm}\label{biggchrom}
For any integer $k$, 
\[
\chi( QI(k^2,k) ) \leq {k+1 \choose 2}. 
\]
\end{thm}
\begin{proof}
Assign one colour to each of the ${k+1 \choose 2}$ pairs from the set
$\{1,2,\dots, k+1\}$.  By the pigeon-hole principle, each partition
$P\in V(QI(k^2,k))$ has at least two distinct elements from
$\{1,2,\dots, k+1\}$ in the same class $P_i \in P$. Assign $P$ any
colour that corresponds to a pair of distinct elements from $\{1,2,\dots,
k+1\}$ that occur in the same class of $P$.

This colouring is a proper colouring of $QI(k^2,k)$.
If partitions $P$ and $Q$ are assigned the same colour, then for some
pair $a,b \in\{ 1,2,\dots, k+1\}$ (with $a\neq b$) there are $i,j$ such that $a,b \in
P_i$ and $a,b \in Q_j$. This means that $P_i$ can not intersect all
$k-1$ sets $Q_l$, $l\neq j$, so $P$ and $Q$ are not qualitatively independent.
\end{proof}

\begin{thm}\label{thm:fracchrom}
For any positive integer $k$, 
\[
\chi^\ast( QI(k^2,k) ) \leq k+1. 
\]
\end{thm}
\begin{proof}
There is a homomorphism from $QI(k^2,k)$ to the
Kneser graph $K_{k^2-1: k-1}$. Each vertex in $QI(k^2,k)$ corresponds
to a $k$-partition of a $k^2$-set. For a partition $P$, let $P_1$ be
the class with the element $k^2 \in P_1$. Map the partition $P$ to the $(k-1)$-subset
$P_1 \backslash \{k^2\}$. 

If partitions $P$ and $Q$ are adjacent in $QI(k^2,k)$, then, for each
class $P_i \in P$ and $Q_j \in Q$, $|P_i \cap Q_j| =1$. In particular,
if $k^2 \in P_1 \in P$ and $k^2 \in Q_1 \in Q$, then 
$(P_1\backslash \{k^2\}) \cap (Q_1\backslash \{k^2\}) = \emptyset$. Thus, this map
defines a homomorphism from $QI(k^2,k)$ to $K_{k^2-1: k-1}$.

\end{proof}

Also by Inequality~(\ref{eq:vt2}), Section~\ref{knesergraphs}, Theorem~\ref{thm:fracchrom}
implies that $\omega( QI(k^2,k) ) \leq k+1$, which is another proof of
Lemma~\ref{omegag2g}.

\begin{lemma}\label{QIregular}
For all $k$, the graph $QI(k^2,k)$ is $(k!)^{k-1}$-regular.
\end{lemma}

\begin{proof}
Let $P$ be a vertex in  $QI(k^2,k)$, so $P$ is a $k$-uniform $k$-partition.

Any $k$-uniform $k$-partition that is qualitatively independent with
$P$ must have each of the $k$ elements in each of the $k$ classes of
$P$ in different classes. There are $k!$ ways to place the $k$
elements of a class of $P$ in different classes. Since there are $k$
classes in $P$, this gives $(k!)^k$ ways to arrange all the elements
of all the classes. Finally, we need to divide by $k!$ so as not to
over count the partitions with the same $k$ classes in different orders.
Thus there are
$\frac{1}{k!}(k!)^k$ partitions qualitatively independent with $P$.
Thus each vertex of $QI(k^2,k)$ has degree $(k!)^{k-1}$.
\end{proof}

From Lemma~\ref{QIdiameter2}, the graphs $QI(k^2,k)$ have diameter 2.
For each vertex $v \in V(QI(k^2,k))$ there are $(k!)^{k-1}$
vertices at distance 1 and 
\[
\frac{1}{k!}{k^2 \choose k}{k^2-k \choose k}\dots{k \choose k}-(k!)^{k-1} - 1
\]
vertices at distance 2.

%%%%%%%%%%%%%%%%%%%%%%%%%%%%%%%%%%%%%%%%%%%%%%%%%%%%

\chapter{Uniform Qualitative Independence Graphs}\label{chp:ass}

In the previous chapter, we saw that a core of the graph $QI(n,2)$,
for $n$ even, is the subgraph induced by vertices that are uniform
2-partitions of an $n$-set. Motivated by this, we consider the
subgraph of $QI(n,k)$ which is induced by the uniform $k$-partitions. This subgraph
is called the {\em uniform qualitative independence graph}.

These graphs are regular and vertex transitive.  A clique in the
uniform qualitative independence graph corresponds to a {\em balanced
covering array}; that is, a covering array with the property that each
letter occurs the same number of times in each row.

We start this chapter by proving bounds on the size of a maximum clique
in a uniform qualitative independence graph. These bounds are a
result of the fact that the uniform qualitative independence graphs
are vertex transitive. Some of these bounds are known, but the method 
we use to prove them is new.

Next, we consider the two ratio bounds for graphs from
Section~\ref{sec:ratiobounds}, which we restate here: 
For a $d$-regular graph $G$ with smallest eigenvalue $\tau$, we have
\begin{eqnarray}\label{ratio1}
\alpha(G) \leq \frac{|V(G)|}{1-\frac{d}{\tau}};
\end{eqnarray}
further, if $G$ is arc transitive or a single graph in an association scheme, then
\begin{eqnarray}\label{ratio2}
\omega(G) \leq 1-\frac{d}{\tau}.
\end{eqnarray}

To use these bounds for a graph $G$, it is necessary to know the
largest and smallest eigenvalue of $G$. For certain special cases of
the uniform qualitative independence graphs, it is possible to find
the eigenvalues, but in general they are difficult to find.  In this
chapter, we find an equitable partition on the vertex set of the
uniform qualitative independence graphs that reduces the calculations
for the eigenvalues. With this partition, we give the eigenvalues (and
their multiplicities) for several small uniform qualitative
independence graphs.

Eigenvalues and their multiplicities can give more information about
the graph than just these two ratio bounds.  For example, Godsil and
Newman~\cite{Godsil:Newman} use the multiplicity of the smallest
eigenvalue of $QI(9,3)$ and the properties of the maximum independent
sets to characterize all of the maximum independent sets in the
graph. This, in turn, is used to prove that the graph $QI(9,3)$ is a
core.  This result is restated with more detail in
Section~\ref{assscheme93}.

Mathon and Rosa~\cite{MR86e:05068} have shown that the graph $QI(9,3)$
is a graph in an association scheme on $\mathcal{U}^9_3$. We can
generalize this to $\mathcal{U}^{ck}_k$ for all positive integers $k$
and $c$, although it is not clear if this generalization produces an
association scheme. For the graphs $UQI(12,3)$ and $UQI(15,3)$, we
give the modified matrix of eigenvalues of the graphs in this
generalization and conjecture that they are part of an association
scheme.

%%%%%%%%%%%%%%%%%%%%%%%%%%%%%%%%%%%%%%%%%%%%%%%%%%%%%%%%%%%%%%%%%%%%%%%%%%%%%%%%%%%
\section{Uniform Qualitative Independence Graphs}\label{UQI}

Recall from Section~\ref{setpartitions} that, for positive integers $n,k,c$ with
$n=ck$, a {\em uniform $k$-partition of an $n$-set} is a partition of
an $n$-set into $k$ classes each of size $c$. 

\begin{defn}[Uniform Qualitative Independence Graph]\index{uniform qualitative independence graph}
For positive integers $n,k$, with the property that $k$ divides $n$
and $n \geq k^2$, the {\em uniform qualitative independence graph}
$UQI(n,k)$ is the graph whose vertex set is the set of all uniform
$k$-partitions of an $n$-set. Vertices are adjacent if and only if the
corresponding partitions are qualitatively independent.
\end{defn}

\orphan
The graphs $UQI(n,k)$ are called {\em partition
graphs}\index{partition graphs} in \cite{Godsil:Newman} and are
denoted $P(c^k)$, where $c=n/k$.

The uniform qualitative independence graphs are vertex transitive, and
from Section~\ref{setpartitions},
\[
|V( UQI(n,k) )| = U(n,k)= \frac{1}{k!} {n \choose c} {n-c \choose c}
\cdots {n-(k-1)c \choose c}.
\]
  
If $k$ does not divide $n$, it is not possible to have uniform
$k$-partitions of an $n$-set. In this case, {\em almost-uniform
partitions} are considered.  Again, recall from
Section~\ref{setpartitions} that, for positive integers $n,k,c$ with
$n=ck+r$ where $0\leq r < k$, an almost-uniform $k$-partition of an
$n$-set is a partition of an $n$-set into $k$ classes, each of
size $c$ or $c+1$. If $r=0$ in the above definition, then
an almost-uniform $k$-partition of an $n$-set is also a uniform
$k$-partition of an $n$-set.

\begin{defn}[Almost-Uniform Qualitative Independence Graph] 
\index{almost-uniform qualitative independence graph}
Let $n,k,c$ be positive integers such that $n=ck+r$ where $0\leq r <
k\leq c$.  The {\em almost-uniform qualitative independence graph}
$AUQI(n,k)$ is the graph whose vertex set is the set of all
almost-uniform $k$-partitions of an $n$-set. Vertices are adjacent if
and only if the corresponding partitions are qualitatively
independent.
\end{defn}

The almost-uniform qualitative independence graphs are vertex transitive.
From Section~\ref{setpartitions} the number of vertices in this graph is
\begin{small}
\begin{eqnarray*}
|V(AUQI(n,k))|= AU(n,k) &=& \frac{1}{r!(k-r)!} {n \choose c} {n-c \choose c} \cdots {n-(k-r-1)c \choose c} \\
               &&\; {r(c+1) \choose c+1}{(r-1)(c+1) \choose c+1} \cdots {c+1 \choose c+1}.
\end{eqnarray*}
\end{small}

If $k$ and $n$ are positive integers, and $k$ divides $n$, then
$AUQI(n,k) = UQI(k,n)$; in the special case that $n=k^2$, $AUQI(k^2,k)
= UQI(k^2,k) = QI(k^2,k)$.

\section{Bounds from Vertex Transitivity} \label{sec:boundsVT}

In this section, we give upper bounds on $\omega(UQI(n,k))$ and
$\omega(AUQI(n,k))$ derived from the fact that the uniform and
almost-uniform qualitative independence graphs are vertex
transitive. An upper bound on $\omega(UQI(n,k))$ (or
$\omega(AUQI(n,k))$) gives a lower bound on the size of a balanced (or
an almost-balanced covering array). To see this, assume $\omega(UQI(ck,k))
\leq r$. Then, there does not exist a balanced $CA(ck,r+1,k)$, which means
a balanced covering array with $r+1$ rows on a $k$-alphabet must have
more than $ck$ columns.

These results are not new, as versions of these bounds for transversal
covers were proven by Stevens, Moura and Mendelsohn~\cite{brett:97b}
(these are restated in Lemma~\ref{tcbounds}).  The proof
from~\cite{brett:97b} uses Bollob\'{a}s's Theorem,
Theorem~\ref{bollobas}.

Before proving these bounds, we need some notation. Let $n,k$ be
positive integers and set $c = \lfloor\frac{n}{k}\rfloor$.  For a set 
$A \subseteq \{1,\dots,n\}$, define the subset $S_A \subseteq \mathcal{AU}^n_k$ by
\[
S_A = \{P \in \mathcal{AU}^n_k : A \subseteq P_0, \mbox{ for some class } 
               P_0 \in P \mbox{ with } |P_0| = c \}.
\]
These sets will be seen again in Section~\ref{PIP}, where they are called
{\em trivially partially intersecting partition systems}.

\begin{thm}\label{thm:generalvtbound}
For positive integers $n,k,c$ with $n=ck+r$ where $0\leq r < k \leq c$, 
\[
\omega(AUQI(n,k)) \leq \frac{n!(k-2)!}{(k-r)(n-c+k-2)!c!}.
\]
\end{thm}

\begin{proof}
Let $A = \{1,2,\dots, c-(k-2)\}$.  Consider the set $S_A \subseteq
\mathcal{AU}^n_k$. Assume $P,Q \in S_A$, and let $P_0 \in P$ and $Q_0 \in
Q$ be the classes with $A \subseteq P_0$ and $A \subseteq Q_0$ and
$|P_0|=|Q_0|=c$. There are at most $c - (c-(k-2)) = k-2$ elements in
$P_0$ that are not also in $Q_0$.  Thus, the class $P_0$ can intersect at
most $k-2$ classes of $Q$ other than $Q_0$.  Since $Q$ has $k$
classes, $P_0$ can not intersect every class in $Q$.
This means the partitions $P$ and $Q$ are not qualitatively independent, and the set
$S_A$ is an independent set in $AUQI(n,k)$. Since
\[
|S_A| = {n-(c-(k-2)) \choose k-2} AU(n-c, k-1),
\]
this produces a lower bound on  $\alpha(AUQI(n,k))$,
\[
 {n-(c-(k-2)) \choose k-2} AU(n-c, k-1) \leq \alpha(AUQI(n,k)).
\]

Since the graph $AUQI(n,k)$ is vertex transitive, by
Inequality~(\ref{eq:vt2}) in Section~\ref{knesergraphs}

\begin{eqnarray*}
\omega(AUQI(n,k)) \leq \frac{|V(AUQI(n,k)|}{\alpha(AUQI(n,k)) } 
&\leq& \frac{AU(n, k)}{{n-(c-(k-2)) \choose k-2} AU(n-c, k-1)} \\
&\leq&  \frac{\frac{1}{ (k-r)!}{n \choose c}}  
            {\frac{1}{(k-1-r)!}{ n-(c-(k-2)) \choose k-2 }}  \\
&=& \frac{n!(k-2)!}{(k-r)(n-c+k-2)!c!}.
\end{eqnarray*}
\end{proof}

Next, we give several bounds for specific values of $n$. Starting with a
bound for $\omega(QI(k^2,k))$, which is a restatement of
Lemma~\ref{omegag2g}.

%%%%%%%%%%%%%%%%%%%%%%%%%%%%%%%%%%%%%%%%%%%%%%%%%%%%%%%%%%%%%%%%%%%%%%%%%%%%%%%%%%%%%%%%%%%%%%%%%%%
\begin{cor}\label{cor:vtboundk2}
For all positive integers $k$, we have $\omega(QI(k^2,k)) \leq k+1$.
\end{cor}

\begin{proof}
Use Theorem~\ref{thm:generalvtbound} with $n=k^2$, $c=k$ and $r=0$
along with the fact that $QI(k^2,k) = AUQI(k^2,k)$.
\remove{
\[
\omega(QI(k^2,k)) = \omega(AUQI(n,k)) \leq \frac{n!(k-2)!}{(k-r)(n-c+k-2)!c!}
          \remove{= \frac{k^2!(k-2)!}{(k)(k^2-2)!k!}}  =
          \frac{k^2(k^2-1)}{k^2(k-1)} = k+1.
\]
}
\end{proof}

%%%%%%%%%%%%%%%%%%%%%%%%%%%%%%%%%%%%%%%%%%%%%%%%%%%%%%%%%%%%%%%%%%%%%%%%%%%%%%%%%%%%%%%%%%%%%%%%%%%
% for K^2+1
This next bound is for $n=k^2+1$. A stronger version of this
bound is given in~\cite{brett:97b} for transversal covers and is
restated as Part~\ref{tcbound3} of Lemma~\ref{tcbounds}.

\begin{cor}\label{cor:vtboundk2+1}
For all positive integers $k \geq 6$, we have $\omega(QI(k^2+1,k)) \leq k+2$.
\end{cor}
\begin{proof}
The graph $QI(k^2+1,k) = AUQI(k^2+1,k)$. Using
Theorem~\ref{thm:generalvtbound} with $n=k^2+1$, $c=k$ and $r=1$,
\begin{eqnarray*}
\omega(QI(k^2+1,k))  &\leq&  \frac{n!(k-2)!}{(k-r)(n-c+k-2)!c!} \\
%         &=&  \frac{( k^2+1)!(k-2)!}{(k-1)(k^2-1)!k!} \\
    &=&  \frac{(k^2+1)k}{(k-1)(k-1)} \\
    &=& k+2 + \frac{4k-2}{k^2-2k+1}.
\end{eqnarray*}
For $k \geq 6$, we have $4k-2 < k^2-2k+1$ and the corollary holds.
\end{proof}

Since the rows of any covering array $CA(k^2+1,s,k)$ with $s>1$
correspond to almost-uniform $k$-partitions of a $(k^2+1)$-set, this
bound can be translated into a bound on $s$.  In particular,
Corollary~\ref{cor:vtboundk2+1} implies that if $k\geq 6$, then for
any $CA(k^2+1,s,k)$, $s \leq k+2$. Equivalently, for $k\geq 6$, if a
$CA(n,s,k)$ exists with $s \geq k+3$, then $n \geq k^2+2$.  The result
from~\cite{brett:97b} is: if $s \geq k+2$, then $n \geq k^2+2$.

%%%%%%%%%%%%%%%%%%%%%%%%%%%%%%%%%%%%%%%%%%%%%%%%%%%%%%%%%%%%%%%%%%%%%%%%%%%%%%%%%%%%%%%%%%%%%%%%%%%
%For K^2+2

\begin{cor}\label{cor:vtboundk2+2}
For $k \geq 13$, we have $\omega(AUQI(k^2+2,k)) \leq k+3$.
\end{cor}
\begin{proof}
Using Theorem~\ref{thm:generalvtbound} with $n=k^2+2$, $c=k$ and $r=2$,
\begin{eqnarray*}
\omega(AUQI(k^2+1,k))  &\leq&  \frac{n!(k-2)!}{(k-r)(n-c+k-2)!c!} \\
%        &=&  \frac{ (k^2+2)!(k-2)!}{(k-2)(k^2)!k!} \\
    &=&  \frac{(k^2+1)(k^2+2)}{(k-2)(k-1)k } \\
&=& k+3 +\frac{10k^2 - 6k+2}{k^3-3k^2+2k}.
\end{eqnarray*}
For $k \geq 13$, we have $10k^2 - 6k+2 < k^3-3k^2+2k$ and the corollary holds.
\end{proof}

Corollary~\ref{cor:vtboundk2+2} can not be directly translated to a
lower bound on the size of a covering array.  The rows of a covering
array $CA(n,s,k)$ correspond to a set of cardinality $s$ of
qualitatively independent $k$-partitions of an $n$-set, but when $n
\geq k^2+2$ these partitions are not necessarily almost-uniform
partitions. For example, such a partition could have one class of
cardinality $k+2$ and all other classes with cardinality $k$.

However, for positive integers $n,k,s$ such that $k$ divides $n$ and
$s>1$, if a $CA(n,s,k)$ is a balanced covering array, then the rows of
$CA(n,s,k)$ correspond to uniform $k$-partitions of an $n$-set. An upper
bound on $\omega(UQI(n,k))$ gives a lower bound on the size of a
balanced covering array. 

With Theorem~\ref{thm:generalvtbound}, we have such a bound on general
uniform qualitatively independence graphs. This bound is not new, but it
is, surprisingly, equivalent to the bound in Part~\ref{tcbound4}
of Lemma~\ref{tcbounds} for point-balanced transversal covers which
originally appeared in~\cite{brett:97b}.

\begin{cor}\label{vtbounduniform}
For positive integers $n,k,c$ with $n=ck$ and $k \leq c$, 
\[
\omega(UQI(n,k)) \leq \frac{(ck)! (k-2)!}{k(ck-c+k-2)!c!}.
\]
\end{cor}
\begin{proof}
This is Theorem~\ref{thm:generalvtbound} with $n=ck$, and $r=0$.
\end{proof}

This corollary can be restated in terms of covering arrays.  For positive
integers $n,k,c,s$ with $n=ck$ and $s>1$, if $CA(n,s,k)$ is a balanced
covering array, then $s \leq \frac{(ck)!  (k-2)!}{k(ck-c+k-2)!c!}$.
Alternately, if a $CA(n,s,k)$ exists with $s > \frac{(ck)!
(k-2)!}{k(ck-c+k-2)!c!}$, then the $CA(n,s,k)$ is not a balanced
covering array.

There are several questions concerning these bounds. For all
integers $n,k$, the graph $AUQI(n,k)$ is an induced subgraph of
$QI(n,k)$, and, in particular, 
\[
\omega(AUQI(n,k)) \leq \omega(QI(n,k)).
\]
The first question is: how big can $\omega(QI(n,k)) -\omega(AUQI(n,k))$ be? 

From Section~\ref{sec:QIcores}, the graph $AUQI(n,2)$ is a core of
$QI(n,2)$; in this case $\omega(QI(n,2)) = \omega(AUQI(n,2))$.  If
this is true in general --- that is, if $AUQI(n,k)$ is a core for
$QI(n,k)$ --- then $\omega(QI(n,k)) = \omega(AUQI(n,k))$ and an upper
bound on $\omega(AUQI(n,k))$ gives a lower bound on the size of a
covering array. Moreover, there would exist a balanced covering array
with same size as an optimal covering array.

\begin{quest}\label{quest:generalcores}
For positive integers $k$ and $n$, is the graph $AUQI(n,k)$ a core of
$QI(n,k)$?  For a positive integer $r$, if $n=CAN(r,k)$, 
does there exist a balanced $CA(n,r,k)$?
\end{quest}

Another question is: are the bounds given here close to the actual
size of the maximum cliques in $AUQI(n,k)$?  This question is much
easier to answer.  A clique in $QI(n,k)$ is a qualitatively
independent family of partitions so $\omega(QI(n,k))= N(n,k)$ (recall
from Section~\ref{CA} that $N(n,k)$ is the maximum number of rows $s$ in a
$CA(n,s,k)$).  Since $AUQI(n,k)$ is an induced subgraph of $QI(n,k)$,
$\omega(AUQI(n,k)) \leq N(n,k)$. The value of $N(n,k)$ is not known in
general, but from Equation~(\ref{eq:asymptotic}) in
Section~\ref{asymptotics},
\[
\limsup_{n \rightarrow \infty} \frac{\log_2 N(n,k)}{n} = \frac{2}{k}.
\]

Let $f(n,k)$ be the bound on $\omega(AUQI(n,k))$ from Lemma~\ref{thm:generalvtbound}.
Then
\[
\limsup_{n \rightarrow \infty} \frac{\log_e f(n,k)}{n} 
         = \log_e \left( \frac{k}{k-1}\right) + \frac{\log_e(k-1)}{k}.
\]
The details of this calculation are omitted, but a similar calculation
is shown in full detail in Section~\ref{SpernerAsymptotics}.  For $k
>2$, it is clear from the asymptotic growth of $f(n,k)$ that the bound
on $\omega(AUQI(n,k))$ from Theorem~\ref{thm:generalvtbound} is far
from the real value of $\omega(AUQI(n,k))$ for large $n$.

\section{Eigenvalues of Qualitative Independence Graphs}

In this section, we give the eigenvalues for the binary qualitative
independence graphs and we find the largest and smallest eigenvalues
for $QI(k^2,k)$.

\subsection{Eigenvalues for $UQI(n,2)$ and $AUQI(n,2)$}

For $n$ even, $UQI(n,2)$ is isomorphic to the complete graph on
$\frac{1}{2}{n \choose n/2}$ vertices.  By
Theorem~\ref{evaluescomplete}, the eigenvalues of the complete graph
on $\frac{1}{2}{n \choose n/2}$ vertices are $\frac{1}{2} {n \choose
n/2 }-1$ and -1. Further, the complete graph is regular and part of an
association scheme (the trivial scheme) so both ratio bounds
(where $d$ is the largest eigenvalue and $\tau$ is the smallest)
\[
\alpha(UQI(n,2)) \leq \frac{|V(UQI(n,2))|}{1-\frac{d}{\tau}} 
\mbox {\; and \;}
\omega(UQI(n,2)) \leq 1-\frac{d}{\tau}
\]
hold, and in fact, both hold with equality.

If $n$ is odd, it is an entirely different story.  The graph
$AUQI(n,2)$ is isomorphic to the complement of the Kneser graph
$K_{n:\frac{n-1}{2}}$ (Section~\ref{sec:QIcores}). As stated in
Example~\ref{exmp:kneserevalues} in Section~\ref{sec:ep}, the
eigenvalues of $K_{n:\frac{n-1}{2}}$ are
\begin{eqnarray*}
(-1)^i{\frac{n+1}{2}-i \choose \frac{n-1}{2}-i} = (-1)^i \left( \frac{n+1}{2} - i \right), &\quad& 
  \mbox{ for }i=0,1,\dots,\frac{n-1}{2}. 
\end{eqnarray*}
From Theorem~\ref{thm:evaluescomp}, the eigenvalues of $AUQI(n,2)$ are
\[
{n \choose \frac{n-1}{2}}- \frac{n+1}{2} - 1 \mbox{ and }
-1-(-1)^i\left(\frac{n+1}{2}-i \right), \mbox{ for } i=1,\dots,\frac{n-1}{2}.
\]
In particular, the largest eigenvalue of $AUQI(n,2)$ is 
${n \choose \frac{n-1}{2}}- \frac{n+1}{2}  - 1$
and the smallest is $\frac{1-n}{2}$ (when $i=2$).

Since $AUQI(n,2)$ is vertex transitive, the ratio bound for independent sets holds:
\[
\alpha(AUQI(n,2)) \leq \frac{|V(AUQI(n,2))|}{1-\frac{d}{\tau}}  \\
   = \frac{{n \choose \frac{n-1}{2}}}
           {1-\frac{{n \choose \frac{n-1}{2}}- \left( \frac{n+1}{2} \right) - 1}{-\frac{n-1}{2}}}
   = \frac{n-1}{2 - 4{ n \choose \frac{n-1}{2} }^{-1}}.
\]
It is not hard to see that $\alpha(AUQI(n,2)) = 2$ for all $n$. So, this
bound is tight only for $n=5$, and for large $n$ this bound is not
good at all.

Next, consider the ratio bound for maximum cliques,
Inequality~(\ref{ratio2}). From Theorem~\ref{thm:maxclique},
$\omega(AUQI(n,2)) = \frac{1}{2}{n \choose \frac{n-1}{2}}$.  Let $d$
be the largest eigenvalue of $AUQI(n,2)$ and $\tau$ the smallest, then
$1-\frac{d}{\tau} = \frac{2}{n-1}\left( {n \choose \frac{n-1}{2}} -2 \right)$. For
$n>5$, the ratio bound for maximum cliques does not hold. This means
that, for $n$ odd and $n > 5$, the graph $AUQI(n,2)$ is not a single
graph in an association scheme nor is it arc transitive
(Theorem~\ref{thm:tauomegabound}).

In conclusion, it seems that the ratio bounds are more appropriate for the
uniform qualitative independence graphs. Motivated by this, we will
focus on these graphs.

\subsection{Eigenvalues for $QI(k^2,k)$}

In this section, we consider the qualitative independence graphs
$QI(k^2,k)$. In this case, $QI(k^2,k)=UQI(k^2,k)$. Cliques
in $QI(k^2,k)$ correspond to orthogonal arrays. The size of the largest clique in
$QI(k^2,k)$ is at most $k+1$ (Corollary~\ref{cor:vtboundk2},
or, equivalently, Lemma~\ref{omegag2g}), and for $k$ a prime power
$\omega(QI(k^2,k)) = k+1$ (Corollary~\ref{OAprimepower} and Lemma~\ref{mols}). The
equivalence between orthogonal arrays and transversal designs and
MOLS, from Theorem~\ref{thm:equivdesigns}, makes this case particularly
interesting.

For all positive integers $k$, we find the largest and the
smallest eigenvalues for the graph $QI(k^2,k)$. These can be used with
the ratio bounds to find upper bounds on $\alpha(QI(k^2,k))$ and
$\omega(QI(k^2,k))$.

For the rest of this section, $k$ will be a positive integer and ${\bf
1}$ will denote the all ones vector of length $U(k^2,k)$.
Following the notation from Section~\ref{sec:boundsVT}, for any
positive integer $n$, and $a,b\in\{1,\dots,n\}$, $S_{ \{a,b\} }
\subseteq \mathcal{U}^n_k$ denotes the set of all $k$-partitions of an
$n$-set that have a class that contains both the elements $a$ and
$b$. For any distinct $a,b \in \{1,2,\dots, k^2\}$, the set $S_{
\{a,b\} }$ is an independent set in $QI(k^2,k)$. A proof of this is
contained in the proof of Theorem~\ref{thm:generalvtbound}.

\begin{lemma}\label{lem:evUQIlargest}
For positive integer $k$, the largest eigenvalue for $QI(k^2,k)$ is
$(k!)^{k-1}$ with multiplicity one.  The corresponding eigenvector is
$\bf{1}$.
\end{lemma}
\begin{proof}
From Lemma~\ref{QIregular} and Lemma~\ref{QIdiameter2}, $QI(k^2,k)$ is
$(k!)^{k-1}$-regular and connected. 
The result follows from Theorem~\ref{regularevalues}.
\end{proof}

\begin{thm}\label{lem:evUQIsmallest}
For any $k$, $\frac{-(k!)^{k-1}}{k}$ is the smallest eigenvalue for $QI(k^2,k)$.
Further, for any distinct $a,b \in \{1,2,\dots,k^2\}$, if $v$ is the
characteristic vector for $S_{ \{a,b\} }$, then $v-\frac{1}{k+1}{\bf 1}$ is
an eigenvector for $\frac{-(k!)^{k-1}}{k}$.
\end{thm}

\begin{proof}
Let $A$ be the adjacency matrix for $QI(k^2,k)$ and let $x =
v-\frac{1}{k+1}{\bf 1}$.  To show that $x$ is an eigenvector corresponding
to the eigenvalue $\frac{-(k!)^{k-1}}{k}$, all that is needed is to
show that $Ax=\frac{-(k!)^{k-1}}{k}x$.

For $i \in \{1,2, \dots, U(k^2,k)\}$, denote row $i$ of $A$ by
$A_i$.  Each row $A_i$ corresponds to a vertex in $QI(k^2,k)$, which
is a $k$-uniform $k$-partition; call this partition $P^i$.

Consider two cases, first when $P^i \in S_{ \{a,b\} }$ and second when
$P^i \not\in S_{ \{a,b\} }$. Since the degree of $P^i$ is $(k!)^{k-1}$
(Lemma~\ref{QIregular}), if $P^i \in S_{ \{a,b\} }$, then $P^i$ is
adjacent to $(k!)^{k-1}$ vertices in $QI(k^2,k)$, none of which are
in $S_{ \{a,b\} }$. Thus
\[
A_i \cdot x = (k!)^{k-1}\left( \frac{-1}{k+1} \right) 
            = \frac{-(k!)^{k-1}}{k} \left( 1 - \frac{1}{k+1} \right).
\]

If $P^i \not\in S_{ \{a,b\} }$ then $P^i$ is adjacent to $(k!)^{k-1}$
vertices in $QI(k^2,k)$.  By a counting argument, similar to, but
slightly more complicated than, the counting argument for the degree of
each vertex (see the proof of Lemma~\ref{QIregular}),
$(k-1)!(k!)^{k-2}$ of these are in $S_{ \{a,b\} }$ and
$(k-1)^2(k-2)!(k!)^{k-2}$ are not in $S_{ \{a,b\} }$.  Thus
\begin{eqnarray*}
A_i \cdot x &=& (k-1)^2(k-2)!(k!)^{k-2}\left(\frac{-1}{k+1}\right) +  (k-1)!(k!)^{k-2}\left(1-\frac{1}{k+1}\right) \\
 %           &=& (k!)^{k-2}(k-1)! \left( (k-1)\left(\frac{-1}{k+1}\right) + \left(1-\frac{1}{k+1}\right) \right) \\
            &=& \frac{-(k!)^{k-1}}{k}\left( \frac{-1}{k+1}\right).
\end{eqnarray*}
Thus, 
\[
Ax = A\left(v-\frac{1}{k+1}{\bf 1}\right) = \frac{-(k!)^{k-1}}{k} \left(v-\frac{1}{k+1}{\bf 1}\right) = \frac{-(k!)^{k-1}}{k}x.
\]

Finally, we need to show that $\frac{-(k!)^{k-1}}{k}$ is the smallest
eigenvalue for $QI(k^2,k)$.  From above, a partition $P \in S_{
\{a,b\} }$ is adjacent to no partitions in $S_{ \{a,b\} }$ and
$(k!)^{k-1}$ partitions not in $S_{ \{a,b\} }$; a partition $P
\not\in S_{ \{a,b\} }$ is adjacent to $(k-1)!(k!)^{k-2}$ partitions in
$S_{ \{a,b\} }$ and $(k-1)^2(k-2)!(k!)^{k-2}$ partitions not in $S_{
\{a,b\} }$. This means the partition 
$\{S_{ \{a,b\} }, V(QI(k^2,k)) \backslash S_{ \{a,b\} } \}$
is an equitable partition (see Section~\ref{sec:ep}). 
From Theorem~\ref{lem:completeratiobound},
\[
|S_{ \{a,b\} }| = \frac{|V(QI(k^2,k))| }{1-\frac{(k!)^{k-1} }{\tau}}
\]
where $\tau$ is the least eigenvalue.  Since the cardinality of $S_{ \{a,b\} }$ is ${k^2-2
\choose k-2}U(k^2-k,k-1)$ and $|V(QI(k^2,k))| = U(k^2,k)$, the least
eigenvalue $\tau = \frac{-(k!)^{k-1}}{k}$.
\end{proof}

It is possible to see that the partition $\{S_{ \{a,b\} }, V(QI(k^2,k))
\backslash S_{ \{a,b\} } \}$ is an equitable partition in another manner.
Recall from Section~\ref{sec:ep} that the orbits of a group action on
the vertices of a graph form an equitable partition.  For any vertex
$P \in V(QI(k^2,k))$ and any permutation $\sigma \in Sym(k^2)$ let
$\sigma(P)$ be the partition with $\sigma(a) \in (\sigma(P))_i$ if and
only if $a \in P_i$.  In this manner, any subgroup of $Sym(k^2)$ induces a subgroup of
$Aut(QI(k^2,k))$.  Let $H' = \{ \sigma \in Sym(k^2) : \sigma(\{a,b\}) = \{a,b\}
\}$.  Then $H'$ is a subgroup of $Sym(k^2)$ and $H'$ induces a subgroup
$H$ of $Aut(QI(k^2,k))$. The group action of $H$ on the
vertices of $QI(k^2,k)$ has two orbits: $S_{ \{a,b\} }$ and
$V(QI(k^2,k)) \backslash S_{ \{a,b\} }$.

From Lemma~\ref{lem:evUQIsmallest} and Lemma~\ref{lem:evUQIlargest},
we know the largest and the smallest eigenvalues of $QI(k^2,k)$, so
now we can use the ratio bound on the size of the largest independent
set in $QI(k^2,k)$.  This result will be seen again in
Section~\ref{PIP}.

\begin{lemma}\label{rationprimepower}
For all $k$,
\begin{eqnarray}
\alpha( QI(k^2,k) ) = {k^2-2 \choose k-2} U(k^2-k, k-1).
\end{eqnarray}
\end{lemma}
\begin{proof}
From the ratio bound for independent sets (Inequality~(\ref{ratio1})),
with the largest and smallest eigenvalues for $QI(k^2,k)$
(Theorem~\ref{lem:evUQIlargest} and Theorem~\ref{lem:evUQIsmallest}),
we have
\[
\alpha( QI(k^2,k) ) \leq \frac{U(k^2, k)}{1-\frac{(k!)^{k-1}}{\frac{-(k!)^{k-1}}{k}}} = {k^2-2 \choose k-2} U(k^2-k, k-1).
\]
For any distinct $i,j \in \{1,\dots,k^2\}$, the set $S_{ \{i,j\} } \subset \mathcal{U}^n_k$ is an
independent set of size ${k^2-2 \choose k-2} U(k^2-k, k-1)$.
\end{proof}

With the exact value of $\alpha( QI(k^2,k) )$ and
Corollary~\ref{cor:vtbound}, we also have the exact value of the
fractional chromatic number of $QI(k^2,k)$.

\begin{cor}
For any positive integer $k$, 
\[
\chi^\ast( QI(k^2,k) ) = k+1. 
\]
\end{cor}

Before we can apply the ratio bound for maximum cliques
(Inequality~(\ref{ratio2})), we need to prove that $QI(k^2,k)$ is
arc transitive.

\begin{thm}
For all positive integers $k$, the graph $QI(k^2,k)$ is arc transitive.
\end{thm}
\begin{proof}
To prove that $QI(k^2,k)$ is arc transitive, we build an automorphism
that takes an arbitrary arc to any other arc. Let $u$ and $v$ be arcs
in $QI(k^2,k)$. Then $u=(P,Q)$ and $v = (R,S)$, where $P,Q,R,S \in
\mathcal{U}^{k^2}_k$ and both $P$ and $Q$ are qualitatively independent, and $R$
and $S$ are qualitatively independent.  Let $P_i$ and $Q_j$, $i,j \in
\{1,2, \dots, k\}$, be the classes in the partitions $P$ and $Q$, and
similarly, $R_i$ and $S_j$ the classes of $R$ and $S$. Then for any
$i,j \in \{1,2, \dots, k\}$, $|P_i \cap Q_j| =1$ and $|R_i \cap S_j|
=1$.

Define a permutation $\phi \in Sym(k^2)$ by $\phi(P_i \cap Q_j)
= R_i \cap S_j$ for all $i,j \in \{1,2, \dots, k\}$.
The permutation $\phi$ induces an automorphism on
$QI(k^2,k)$ with the property that $\phi(P) = R$ and $\phi(Q) = S$.
\end{proof}

This proof can be extended to the graph $QI(k^2+1,k)$ but not to any
other qualitative independence graph.

\begin{thm}
For all positive integers $k$, the graph $QI(k^2+1,k)$ is
arc transitive.  
\end{thm}

\remove{
\begin{conj}
For $n,k$ positive integers with $n \geq k^2+2$, the
graph $AUQI(n,k)$ is not arc transitive.
\end{conj}
}

Since the graph $QI(k^2+1,k)$ is arc transitive, the ratio bound for
cliques can be used to bound $\omega(QI(k^2+1,k))$. This makes it
particularly interesting to try to find the largest and smallest eigenvalues of
$QI(k^2+1,k)$. 
\remove{Perhaps this could be done using the same method used
in the proofs of Theorem~\ref{lem:evUQIlargest} and Theorem~\ref{lem:evUQIsmallest}.}

\begin{quest}
What are the eigenvalues of the graph $QI(k^2+1,k)$?
\end{quest}

Since the graph $QI(k^2,k)$ is arc transitive, the ratio bound for
cliques, Inequality~(\ref{ratio2}), holds. The ratio bound for cliques
gives the now very well-known bound
\[
\omega(QI(k^2,k)) \leq 1 - \frac{(k!)^{k-1}}{\frac{-(k!)^{k-1}}{k}} = k+1.
\]

We have now found this bound in three different ways. The first proof
was Lemma~\ref{omegag2g}. This lemma used the fact that a clique of
size $r$ in $QI(k^2,k)$ corresponds to an orthogonal array
$OA(n,r,k,2)$ which in turn corresponds to a set of $(r-2)$-$MOLS$ of
order $k$. Since there can be at most $k-1$ MOLS of order $k$, we have the bound $r \leq
k+1$. This method can not be generalized to $AUQI(n,k)$, as cliques in this graph do
not correspond to sets of MOLS.

The second method used to find this bound was in
Corollary~\ref{cor:vtboundk2}. This proof used bounds from
vertex transitivity of the graph $QI(k^2,k)$. This method can be used
to bound the number of rows in a balanced covering array
(Corollary~\ref{vtbounduniform}). The comments on the asymptotic
growth of this bound indicate that this bound does not seem to be very
good.

The third method (Theorem~\ref{thm:tauomegabound}) uses the ratio
bound for cliques and the eigenvalues of the graph. This method has
two problems: first, finding the eigenvalues of graphs can be
difficult, and second, it is not clear if the ratio bound for cliques
holds for all uniform qualitative independence graphs. Indeed,
Theorem~\ref{thm:tauomegabound} holds for arc-transitive graphs and for
graphs that are a single graph in an association scheme.  As it is
clear that the first two methods cannot be extended, we turn our
attention to this third method. In the next section, we find the
eigenvalues for several qualitative independence graphs and also
consider whether or not they are graphs in association schemes.

\section{Equitable Partitions}

Throughout this section, we assume that $k, c$ and $n$ are positive
integers with $n=kc$.

For $k>2$, the eigenvalues of the graphs $QI(n,k)$ and $UQI(n,k)$ are not known
in general. One notable exception is the graph $QI(9,3)$. For this
graph, all eigenvalues and their multiplicities are known. These are
found using an equitable partition on the vertices of $QI(9,3)$. This
example is particularly interesting because the eigenvalues and their
multiplicities can be used to prove that $QI(9,3)$ is a core.

In this section, we give two equitable partitions of the vertices of
$UQI(n,k)$ that can be used to find some of the eigenvalues of
$UQI(n,k)$. The first is a simple partition motivated by the comments
following the proof of Theorem~\ref{lem:evUQIsmallest}. The second is a
natural extension of the partition used to find the eigenvalues of
$QI(9,3)$. This second equitable partition reduces the work to find
eigenvalues but it does not solve the problem completely. Computation
is still needed, and for large graphs the computation takes too long.  Hence only the eigenvalues for the graphs $QI(9,3)$,
$UQI(12,3)$, $UQI(15,3)$, $UQI(18,3)$ and $QI(16,4)$ are given.

\subsection{Eigenvalues for $QI(9,3)$}\label{assscheme93}

One of the qualitative independence graphs, $QI(9,3)$, has previously
appeared in the literature.  Mathon and Rosa~\cite{MR86e:05068} give
an association scheme which has $QI(9,3)$ as one of the graphs, and in
this paper they also state all the eigenvalues and their
multiplicities for $QI(9,3)$.  Mathon and Rosa do not focus on the graph
$QI(9,3)$, rather they focus on a different graph in the scheme that
is strongly regular.

For two partitions $P$ and $Q$ the {\em meet of $P$ and $Q$}
\index{meet of two partitions} is defined by
\[
P \wedge Q = |\{ (i,j) \; :\; P_i \cap Q_j \neq \emptyset \}|.
\]
If $P,Q \in \mathcal{U}^9_3$, the value of $P \wedge Q$ is one of 3,
5, 6, 7 or 9.  In addition, $P \wedge Q = 3$ if and only if $P=Q$; and $P \wedge Q = 9$
if and only if $P$ and $Q$ are qualitatively independent.

Define an association scheme with graphs $G_i$ for $i=1,\dots, 4$ as follows.
The vertex set of each graph is the set of all uniform 3-partitions of an 9-set,
$\mathcal{U}_3^9$. The graphs $G_i$ for $i =1,2,3,4$ have an edge between
vertices if and only if the meet of the corresponding partitions is 5, 6, 7
or 9 respectively.

The modified matrix of eigenvalues for this scheme,
Table~\ref{eq:evQI93}, was given by Mathon and Rosa
in~\cite{MR86e:05068}. In this table, the last four columns contain
all of the eigenvalues of the graphs $G_i$ for $i =1,\dots,4$.  The
first column contains the multiplicities of the eigenvalues (see
Section~\ref{sec:introassscheme}).

\begin{table}[h]
\[
\left(
\begin{array}{r|rrrr}
  1 & 27 &162 & 54 & 36 \\
 27 & 11 & -6 &  6 &-12 \\
 48 &  6 & -6 & -9 &  8 \\
120 & -3 & -6 &  6 &  2 \\
 84 & -3 & 12 & -6 & -4
\end{array}
\right)
\]
\caption{The modified matrix of eigenvalues for the association scheme on $\mathcal{U}^9_3$ \label{eq:evQI93}}
\end{table}

These eigenvalues can be found using equitable partitions.  This is
similar to how the eigenvalues of the Kneser graphs are found
(Example~\ref{exmp:kneserevalues}, Section~\ref{sec:ep}). We will
define an equitable partition $\pi$ on $\mathcal{U}^9_3$ (which is
the vertex set of the graphs $G_i$ for $i=1,\dots, 4$).  Fix a vertex
$P \in
\mathcal{U}^9_3$.  Partition the other vertices of $\mathcal{U}^9_3$
into 4 classes, $C_j$, $j= 1,\dots,4$, such that all vertices in a
class have the same meet with $P$.

In Theorem~\ref{meettable}, it will be shown that this vertex partition is
also formed by the orbits of a group acting on the partitions of
$\mathcal{U}^9_3$, so for all $i=1,\dots, 4$, this vertex partition is
equitable on $G_i$. It is possible to build the adjacency matrix of
the quotient graph $G_i/\pi$; indeed, for $j,k \in \{1,\dots,4\}$, the
$(j,k)$-entry of matrix $A(G_i / \pi)$ is the number of vertices in
the class $C_k$ which are adjacent to a single (but arbitrary) vertex
in class $C_j$ in graph $G_i$. The fact that these numbers can be
found and do not depend on the arbitrary vertex in class $C_j$ also
proves that the partition is equitable.

The adjacency matrices of the quotient graphs $G_i/\pi$, for $i=1,\dots,
4$, are listed below:

\begin{eqnarray*}
A(G_1 / \pi) = \left(
\begin{array}{ccccc}
0 & 1 & 0 & 0 & 0  \\
27 & 8 & 2 & 3 & 0 \\
0 & 12 & 15 & 18 & 18 \\
0 & 6 & 6 & 6 & 0 \\
0 & 0 & 4 & 0 & 9 
\end{array}
\right),
&
A(G_2 / \pi) = \left(
\begin{array}{ccccc}
0 & 0 & 1 & 0 & 0 \\
0 & 12 & 15 & 18 & 18 \\
162 & 90 & 96 & 90 & 90 \\
0 & 36 & 30 & 30 & 36 \\
0 & 24 & 20 & 24 & 18 
\end{array}
\right),
\\
A(G_3 / \pi) = \left(
\begin{array}{ccccc}
0 & 0 & 0 & 1 & 0 \\
0 & 6 & 6 & 6 & 0 \\
0 & 36 & 30 & 30 & 36 \\
54 & 12 & 10 & 9 & 12 \\
0 & 0 & 8 & 8 & 6
\end{array}
\right),
&
A(G_4 / \pi) = \left(
\begin{array}{ccccc}
0 & 0 & 0 & 0 & 1 \\
0 & 0 & 4 & 0 & 9 \\
0 & 24 & 20 & 24 & 18 \\
0 & 0 & 8 & 8 & 6 \\
36 & 12 & 4 & 4 & 2
\end{array}
\right).
\end{eqnarray*}

This last matrix is the adjacency matrix for $QI(9,3)/ \pi$.
Table~\ref{eq:evQI93} is the modified matrix of eigenvalues for this
scheme, so the first column contains the multiplicities of the eigenvalues
and the last four columns contain the eigenvalues of these four matrices.
How the multiplicities of the eigenvalues are found is discussed in
Section~\ref{sec:betterep}.

The eigenvalues for $QI(9,3)$ are $(36,-12,8,2,-4)$ with
corresponding multiplicities $(1,27,48,120,84)$.  The largest
eigenvalue, 36, is the degree of the vertices in $QI(9,3)$.  The
smallest eigenvalue is $\tau = -12 = -\frac{(3!)^{{3-1}}}{3}$, which 
is expected from Lemma~\ref{lem:evUQIsmallest}. 

\remove{ Further, from Theorem~\ref{lem:evUQI}
\[
\alpha(QI(9,3)) \leq \frac{280}{1-\frac{36}{-12}} = 70.
\]
This bound is exactly the result from Lemma~\ref{rationprimepower}. }

We have more information than just the eigenvalues. But first, we need
some notation. Let $v$ be the characteristic vector for the set $S_{
\{a,b\} }$, for distinct $a,b \in \{1,2,\dots,k^2\}$.  From
Theorem~\ref{lem:evUQIsmallest}, $v-\frac{1}{k+1}{\bf 1}$ is an
eigenvector for $QI(k^2,k)$ corresponding to the smallest eigenvalue,
$\frac{-(k!)^{k-1}}{k}$. Let $V_k \subset \mathbb{R}^{U(k^2,k)}$ be the vector space spanned by the
set of ${k^2 \choose 2}$ characteristic vectors of $S_{ \{a,b\} }$,
for all distinct $a,b \in \{1,2,\dots,k^2\}$.

Godsil and Newman~\cite{Godsil:Newman} show that the dimension of the
vector space $V_3$ is 27.  Using the fact that this is exactly the
dimension of the eigenspace corresponding to the smallest eigenvalue
of $QI(9,3)$, and properties of the maximum
independent sets in $QI(9,3)$, they prove that all independent sets in
$QI(9,3)$ are of the form $S_{ \{i,j\} }$ for distinct $i,j
\in \{1,2,\dots, 9\}$. With this result, they are also able to prove that
the graph $QI(9,3)$ is a core.  This proof is similar to the proof of
Lemma~\ref{knesercore}, which shows that the Kneser graphs are cores.
 
\begin{lemma}[\cite{Godsil:Newman}]\label{lem:qi93core}
The graph $QI(9,3)$ is a core.
\end{lemma}

Moreover, in~\cite{mikesthesis}, Newman shows that the dimension of the
vector space $V_{k}$ is ${k^2 \choose 2} - {k^2 \choose 1}$ for all
$k$.  If it was known that the dimension of the eigenspace corresponding
to the smallest eigenvector was ${k^2 \choose 2} - {k^2 \choose 1}$,
then we might be able to show that all the independent sets in
$QI(k^2,k)$ are sets $S_{\{i,j\} }$ for some distinct $i,j \in
\{1,2,\dots, k^2\}$. Further, it may be possible to use this to prove
that the graphs $QI(k^2,k)$ are cores in general.

It is not clear whether or not this can be done directly. One problem is that for
$k>3$ it is not clear that it is possible to construct an association
scheme where $QI(k^2,k)$ is a graph (this is discussed in
Section~\ref{sec:otherschemes}). Another problem is that we would need
to generalize the extra properties of the independent sets. But this
motivates trying to find the eigenvalues and the multiplicities of the
qualitative independence graphs. It also motivates searching for
association schemes (or possibly asymmetric association schemes) that
have $QI(k^2,k)$ as one of the graphs.  We end this section with
several conjectures on extending Lemma~\ref{lem:qi93core} to all
$QI(k^2,k)$.

\begin{conj}\label{conj:dimespace}
For all positive integers $k$, the eigenspace corresponding to the
smallest eigenvector of $QI(k^2,k)$, $ \frac{-(k!)^{k-1}}{k}$, has
dimension ${k^2 \choose 2} - {k^2 \choose 1}$.
\end{conj}

%\begin{conj}
%The set of characteristic vectors for $S_{ \{a,b\} }$ in $QI(n,k)$ for all
%distinct $a,b \in \{1,\dots,n\}$ spans a vector space of dimension ${n
%\choose 2}-{n \choose 1}$.
%\end{conj}

\begin{conj}
For all positive integers $k$, all maximum independent sets in
$QI(k^2,k)$ correspond to sets $S_{ \{a,b\} }$ for distinct $a,b \in
\{1,\dots,k^2\}$.
\end{conj}

Conjecture~\ref{conj:partialintersect} in Section~\ref{PIP} is a
generalization of this conjecture to {\em intersecting partition
systems}.

\begin{conj}
For any $k$, the graph $QI(k^2,k)$ is a core.
\end{conj}

\subsection{A Simple Equitable Partition}
 
From the proof of Theorem~\ref{lem:evUQIsmallest}, the partition $\pi
= \{S_{ \{1,2\} } , V(QI(k^2,k)) \backslash S_{ \{1,2\} }\}$ is an
equitable partition. This partition can not be used to find all
eigenvalues for $QI(k^2,k)$, but it can be used to find their largest
and smallest eigenvalues.  

\remove{In this section we extend this equitable partition for the vertices of
$QI(k^2,k)$ to all $UQI(ck,k)$.}

The adjacency matrix for the quotient graph $QI(k^2,k) / \pi$ is 
\[
\left( 
\begin{array}{cc}
 0 & k!^{k-1} \\
 \frac{k!^{k-1}}{k} &  k!^{k-1} - \frac{k!^{k-1}}{k} 
\end{array}
\right)
\]
which has eigenvalues $k!^{k-1}$ and $-\frac{k!^{k-1}}{k}$.
These are the largest and smallest eigenvalues of $QI(k^2,k)$.

It would be interesting to find a similar equitable partition on the
vertices of the graph $UQI(ck,k)$, where $c$ is a positive integer
with $k < c$. From the proof of Theorem~\ref{thm:generalvtbound}, for
a set $A \subset \{1,\dots,ck\}$ with $|A| = c-(k-2)$, $S_A$ is an
independent set in $QI(ck,k)$. But, the partition $\{S_A, V(UQI(ck,k))
\backslash S_A \}$ is not equitable. So instead, we consider the set $S_{
\{1,2\} }$, which is not an independent set in $UQI(ck,k)$ for $k < c$. But,
from the comments following the proof of
Lemma~\ref{lem:evUQIsmallest}, the vertex partition
\[
\{S_{ \{1,2\} }, V(QI(k^2,k)) \backslash S_{ \{1,2\} } \}
\]
is the orbit partition formed by the 
subgroup $H < Aut(QI(k^2,k))$ induced by the group $H'= \{\sigma
\in Sym(k^2) : \sigma(\{1,2\}) = \{1,2\} \}$.  Similarly, the partition 
\[
\{S_{ \{1,2\} }, V(UQI(ck,k)) \backslash S_{ \{1,2\} }\}
\]
is the orbit partition formed by
the subgroup $H < Aut(QI(ck,k))$ induced by the group $H'=
\{\sigma \in Sym(ck) : \sigma(\{1,2\}) = \{1,2\} \}$. From the comments in
Section~\ref{sec:ep}, this means $\{S_{ \{1,2\} }, V(UQI(ck,k)) \backslash
S_{ \{1,2\} }\}$ is an equitable partition.

For a given vertex $P \in S_{ \{1,2\} }$ in $UQI(ck,k)$, it is much harder to count
the number of vertices in $S_{ \{1,2\} }$ adjacent to $P$ and the
number of vertices not in $S_{ \{1,2\} }$ adjacent to $P$.  So, let
$a$ be the number of partitions with 1 and 2 in the same class
adjacent to a fixed partition also with 1 and 2 in the same class. Let
$b$ be the number of partitions with 1 and 2 in the same class and
adjacent to a fixed partition with $1$ and $2$ in different
classes. Let $d$ be the degree of the vertices in $UQI(ck,k)$.

Then the adjacency matrix for this quotient
graph is 
\[
\left(
\begin{array}{cc}
 a & d-a \\
 b & d-b
\end{array}
\right).
\]
The eigenvalues for this quotient graph are $a-b$ and $d$.  It is no
surprise that $d$ is an eigenvalue since the graph is $d$-regular.
Since the partition is equitable, $a-b$ is an eigenvalue of
$UQI(ck,k)$.

\begin{quest}
For integers $c,k$ with $k \leq c$, what is the exact value of $a$ and $b$ for $UQI(ck,k)$?
\end{quest}

For $QI(k^2,k)$ the value of $a-b$ is $-\frac{k!^{k-1}}{k}$. From
Lemma~\ref{lem:evUQIsmallest}, the eigenvalue $a-b$ is the smallest
eigenvalue for $QI(k^2,k)$. Is this true for all $UQI(ck,k)$?
\begin{quest}
For integers $c,k$ with $k \leq c$, is the eigenvalue $a-b$ from above
the smallest eigenvalue for $UQI(ck,k)$?
\end{quest} 

One obvious problem with this equitable partition is that it can only
be used to find at most two eigenvalues of the graph $UQI(ck,k)$. In
the next section, we give an equitable partition that can be used to
find all eigenvalues of $UQI(ck,k)$.

\subsection{A Better Equitable Partition}\label{sec:betterep}

Recall from Theorem~\ref{singlecell}, that if $G$ is a
vertex-transitive graph and $\pi$ is the orbit partition for some
subgroup of $Aut(G)$ with the property that $\pi$ has a singleton class
$\{u\}$, then every eigenvalue of $G$ is an eigenvalue of $G/\pi$.

A partition with these properties will be used to simplify the
calculations of the eigenvalues of $UQI(n,k)$.
\remove{ A subgroup of the
automorphism group of $UQI(n,k)$ is given and the orbits of this group
will produce an equitable partition, $\pi$, of $V(UQI(n,k))$ that has
a singleton class.}

For $P \in \mathcal{U}^n_k$, define a subgroup
\[
\fix'(P) = \{\sigma \in Sym(n) \; : \; \sigma(P) = P\},
\]  
where $\sigma(P)$ is defined by: $\sigma(a) \in (\sigma(P))_i$ if and only if $a \in P_i$.  The
group $\fix'(P)$ induces a subgroup of $Aut(UQI(n,k))$, which we call
$\fix(P)$.  The orbits of the group action of $\fix(P)$ on the vertices
of $UQI(n,k)$ contain the partition $P$ as a singleton class.  Thus
the orbit partition of $\fix(P)$ forms an equitable partition of the
vertex set of $UQI(n,k)$ with a singleton class.

\begin{cor}
Fix an arbitrary partition $P \in V(UQI(n,k))$, let $\pi$ be the partition of $V(UQI(n,k))$
induced by the group $\fix(P)$. Then every eigenvalue of $UQI(n,k)$ is
an eigenvalue of the quotient graph $UQI(n,k) / \pi$.
\end{cor}

It is easier to calculate the eigenvalues of the quotient graph
$UQI(n,k)/ \pi$ than $UQI(n,k)$. Since the group $\fix(P)$ is very
large and difficult to work with, we give a characterization for when
two partitions in $V(UQI(n,k))$ are in the same orbit under the group
action of $\fix(P)$.  For the graph $QI(9,3)$, this gives the same
partition of the vertices as the one defined using the meet of two
partitions in Section~\ref{assscheme93}.

For partitions $P,Q \in V(QI(n,k))$ define the {\em meet table of P
and Q}\index{meet table} to be the $k \times k$ array with the $i,j$
entry $|P_i \cap Q_j|$. Denote the meet table of $P$ and $Q$ by
$M_{P,Q}$. Partitions $P$ and $Q$ are qualitatively independent if and
only if the meet table of $P$ and $Q$ has all entries non-zero. 

Two meet tables are {\em isomorphic}\index{isomorphic meet tables} if
some permutation of the rows and columns of one array produces the
other array\remove{(note we do not include taking the transpose in our
definition of isomorphic)}.  Different orderings on the classes in the
partitions $P$ and $Q$ could produce different meet tables, but the
tables would be isomorphic.  If $P=Q$, then the meet table $M_{P,Q}$
is isomorphic to the array with diagonal entries $c = n/k$ and all
other entries are zero.

\begin{thm}\label{meettable}
Let $P,Q,R \in V(QI(n,k))$. The meet table for $P$
and $Q$ is isomorphic to the meet table for $P$ and $R$ if and
only if there exists $g \in \fix(P)$ such that $g(Q)=R$.
\end{thm}
\begin{proof}
For a partition $P$, let $P_i$ for $i\in \{0,\dots, k-1\}$ denote the
classes of $P$.
  
Assume the meet tables $M_{P,Q}$ and $M_{P,R}$ are isomorphic. Then for
some permutations $\sigma, \phi \in Sym(k)$
\[
 [M_{P,Q}]_{i,j} =  [M_{P,R}]_{\sigma(i),\phi(j)} \mbox{, for } i,j\in\{0,1,\dots, k-1\}.
\]
For all $i,j \in \{1,\dots,k\}$, we have $|P_i \cap Q_j| = |P_{\sigma(i)} \cap
R_{\phi(j)}|$.  For fixed $i,j$, set $m = |P_i \cap Q_j|$. Denote $P_i \cap Q_j =
\{a_1,a_2,\dots, a_m\}$ and $P_{\sigma(i)} \cap R_{\phi(j)} =
\{b_1,b_2,\dots, b_m\}$.  Define $g_{i,j}$ to be the mapping that
assigns $a_l$ to $b_l$ for $l = 1, \dots, m$. For distinct pairs
$(i_0,j_0)$ and $(i_1,j_1)$ the mappings $g_{i_0,j_0}$ and
$g_{i_1,j_1}$ have disjoint domains and map these sets to disjoint sets.
In particular, for $j_0 \neq j_1$, $g_{i,j_0}(P)$ and $g_{i,j_1}(P)$
are disjoint.

\remove{By definition, for all $i,j \in \{0,\dots, k-1\}$, $g_{i,j}(P_i)
\subseteq P_{\sigma(i)}$ and $g_{i,j}(Q_j) \subseteq R_{\phi(j)}$.}

\remove{Thus, $\cup_{j=0}^{k-1}(P_i \cap Q_j) = P_i$.}
\remove{\[ \left| \bigcup_{j = 0}^{k-1} g_{i,j}(P_i) \right| = c.\]}

Define $g_i =  \Pi_{j = 0}^{k-1}g_{i,j}$, then $g_i (P_i) = P_{\sigma(i)}$.
Define $g = \Pi_{i=0}^{k-1}g_i$, so $g(P_i) = P_{\sigma(i)}$ for all
$i \in \{0, \dots, k-1\}$.  
Similarly, for all $j \in \{ 0, \dots, k-1\}$, $g(Q_j) = R_{\phi(j)}$.
Thus, the permutation $g$ is in $\fix(P)$ and $g(Q) =R$.

Next, assume that there exists a $g \in \fix(P)$ such that $g(Q) = R$ and
show that $M_{P,Q}$ is isomorphic to $M_{P,R}$.  Define a permutation
on the rows $i = 0,\dots,k-1$ of $M_{P,Q}$ by $\sigma(i) = i'$ if
and only if $g(P_i) = P_{i'}$. Similarly, define a permutation $\phi$
on the columns $i=0,\dots, k-1$ of $M_{P,Q}$ by $\phi(j) = j'$ if and
only if $g(Q_j) = R_{j'}$.  Thus,
\[
 [M_{P,Q}]_{\sigma(i),\phi(j)} =  [M_{P,R}]_{i,j} \mbox{, for } i,j\in\{0,1,\dots, k-1\}
\]
and the meet tables are isomorphic.
\end{proof}

With this characterization of an equitable partition, we use a
computer program to build the quotient graph of several of the graphs
$UQI(ck,k)$.  Then, using {\texttt Maple} we can find the eigenvalues
of these graphs.  First, we fix a partition $P \in V(UQI(ck,k))$.  Let
$\pi$ be the partition of the vertices in $UQI(ck,k)$ formed by the
orbits of $\fix(P)$.  Next, we build the list of the classes in $\pi$.
To do this, we create a list of all non-isomorphic meet tables.  This
is done by going through all the partitions in $\mathcal{U}^n_k$ and,
for each one, building the meet table with $P$. For each meet table we
construct, we check if it is isomorphic to a meet table already in the
list. If not, we add the meet table to the list. The classes of $\pi$
(the non-isomorphic meet tables) are stored as a single partition
that has the particular meet table with $P$; we call this partition
the {\em representative of the orbit}. Next, the adjacency matrix for
$UQI(ck,k)/\pi$ is built. Again, we go through the list of all
partitions in $\mathcal{U}^{ck}_k$, and for each orbit we count the
number of partitions which are qualitatively independent with the
representative of the orbit.

Since this process requires going through the entire set of partitions
$\mathcal{U}^n_k$ twice, it can take a very long time to
run. Currently, the largest graph our program can complete in a
reasonable time (2 days) is $QI(16,4)$, which has more than 2.6
million vertices. Fortunately, the quotient graphs are small: the
graph $QI(16,4)/\pi$ has only 43 vertices.

\subsection{Multiplicities of the Eigenvalues}

Once the adjacency matrix for $UQI(ck,k)/ \pi$ is built, we also have
enough information to find the multiplicities of the eigenvalues.
First we need some notation. For a graph $G$, let $\phi(G,x)$ denote
the characteristic polynomial of the adjacency matrix $A(G)$, and let
$\phi'(G,x)$ denote the derivative of this polynomial.  Let $ev(G)$
denote the set of distinct eigenvalues of $G$, and for $\lambda \in
ev(G)$, let $m_\lambda$ be the multiplicity of $\lambda$.

\begin{cor}[Section 5.3,~\cite{MR94e:05002}]
Let $\pi = \{C_1,C_2,\dots, C_r\}$ be an equitable partition of the
vertices of a graph $G$ with the property that $C_1$ is a singleton
class. Then
\[
\frac{\phi'(G,x)}{\phi(G,x)} = \frac{ |V(G)| \phi( (G/\pi)\backslash C_1, x )}{\phi(G/\pi, x)}.
\]
\end{cor}
 
It is a straightforward application of the multiplication rule for
derivatives to see that
\[
\frac{\phi'(G,x)}{\phi(G,x)} = \sum_{\lambda \in ev(G) } \frac{m_\lambda}{x-\lambda}.
\]

With the adjacency matrix for $UQI(n,k)/\pi$, using {\texttt Maple}, it
is simple to find the partial fraction expansion of $\frac{|V(G)|
\phi( (G/\pi)\backslash C_1, x)}{\phi(G/\pi, x)}$ where $G =
UQI(ck,k)$, and $C_1 = \{P\}$. Then the numerators in the partial fraction
expansion are the multiplicities of the eigenvalues for $UQI(ck,k)$.

\subsection{Eigenvalues for Several Small Uniform Qualitative Independence Graphs}

\remove{        $QI(9,3)$, $UQI(12,3)$, $UQI(15,3)$, $UQI(18,3)$ and $QI(16,4)$    }

The equitable partition described in Section~\ref{sec:betterep} can be
used to reduce the calculation needed to find the eigenvalues of
$UQI(ck,k)$ and their multiplicities.  As stated in the comments
following the proof of Theorem~\ref{meettable}, the adjacency matrix
for $UQI(ck,k)/\pi$ can be built using a computer program, but this
can be very time consuming. For large $c$ and $k$, this method takes
too much time. To date, we have been able to find the eigenvalues for
the graphs $QI(9,3)$, $UQI(12,3)$, $UQI(15,3)$, $UQI(18,3)$ and
$QI(16,4)$. These are given, along with their multiplicities, in
Table~\ref{evaluestable}.

\begin{table}[h]
\begin{tabular}{|l|l|} \hline
Graph & Eigenvalues\\ \cline{2-2}
&  Corresponding multiplicities (in the same order as the eigenvalues) \\
\hline \hline
$UQI(9,3)$   & (-4, 2, 8, -12, 36) \\ \cline{2-2}
            & (84, 120, 48, 27, 1) \\
\hline
$UQI(12,3)$ & (0, 8, -12, 18, -27, 48, 108, -252, 1728) \\ \cline{2-2}
          & (275, 2673, 462, 616, 1408, 132, 154, 54, 1)\\
\hline
$UQI(15,3)$ & (4, 8, -10, -22, 29, 34, -76, 218, -226, 284, 1628, -5060, 62000) \\ \cline{2-2}
          &  \small{(1638, 21450, 910, 25025, 32032, 22113, 11583, 1925, 7007, 2002, 350, 90,1)}\\
\hline
$UQI(18,3)$ &\small{(8, 15, 18, -60, 60, -102, -120, 120, 368, 648, -655, -2115,}\\
          &\small{2370, -2115, 2370, 2460, -4140, 24900, -89550, 1876500, $954 \pm 18\sqrt{10209}$)}\\ \cline{2-2}
         & \footnotesize{(787644, 678912, 136136, 87516, 331500, 259896, 102102, 219912, 99144, 11934,}\\ 
         & \footnotesize{88128, 22848, 4641, 5508, 2244, 663, 135, 1, 9991)} \\ 
\hline
$UQI(16,4)$ 
&\footnotesize{(-72, $-56\pm 8\sqrt{193}$, $-96\pm 96\sqrt{37}$, $24\pm 24\sqrt{97}$, -96, 96, -288, 8, -144, 24,}\\
 &\footnotesize{192, 32, 1728, -64, -16, 432, 48, 1296,-48, -576, 128, -3456,}\\
 &\footnotesize{576, 13824, -1152, 144)} \\ \cline{2-2}
 &\footnotesize{(266240, 137280, 7280, 76440, 69888, 91520, 24960, 262080, 73920, 24024,}\\
 &\footnotesize{65520, 150150, 440, 51480, 753324, 20020, 420420, 1260, 23100, 10752, 60060, 104,}\\
 &\footnotesize{4070, 1, 1260, 32032)}\\
\hline
\end{tabular}
\caption{Eigenvalues and their multiplicities for small, uniform qualitative independence graphs \label{evaluestable}}
\end{table}

The largest and smallest eigenvalues for each of the
graphs in Table~\ref{evaluestable} are particularly interesting.
Since the graphs $UQI(ck,k)$ are vertex transitive, the ratio bound
for independent sets (Inequality~(\ref{ratio1})) gives an upper bound on the size of the maximum
independent set for each of the graphs in
Table~\ref{evaluestable}. Moreover, from the proof of
Theorem~\ref{thm:generalvtbound}, the set $S_A$ with $|A| = c-k+2$ is
an independent set of $UQI(ck,k)$, and the cardinality of $S_A$ is a
lower bound for $\alpha(UQI(ck,k))$.

\begin{cor} \hfill
\begin{enumerate}
\item{$\alpha(QI(9,3))=70$ }
\item{$ 315 \leq \alpha(UQI(12,3)) \leq  735 $}
\item{$1386 \leq \alpha(UQI(15,3)) \leq  9516$}
\item{$ 6006 \leq \alpha(UQI(18,3)) \leq 130215$}
\orphan
\item{$\alpha(QI(16,4))= 525525$}
\end{enumerate}
\end{cor}

\begin{quest}
What is the size of a largest independent set in each of $UQI(12,3)$,
$UQI(15,3)$ and $UQI(18,3)$?
\end{quest}

Considering Conjecture~\ref{conj:dimespace}, it is also interesting to
note the multiplicity of the smallest eigenvalue for each of the
graphs in Table~\ref{evaluestable}.  The smallest eigenvalue for
$UQI(16,4)$ is -3456 and it has multiplicity 104. This confirms
Conjecture~\ref{conj:dimespace} for $k=4$.  In fact, for each of
the graphs $QI(3c,3)$, $c=3,4,5,6$, the dimension of the eigenspace
corresponding to the smallest eigenvalue is ${ck \choose 2} - {ck
\choose 1}$. This suggests that Conjecture~\ref{conj:dimespace} can be
extended from $QI(k^2,k)$ to $UQI(ck,k)$.

\begin{conj}\label{conj:espacen}
For all positive integers $c,k$, the eigenspace corresponding to the
smallest eigenvalue of $UQI(ck,k)$ has dimension ${ck \choose 2} - {ck
\choose 1}$.
\end{conj}

\section{Other Schemes}\label{sec:otherschemes}

The association scheme from Section~\ref{assscheme93} that has
$QI(9,3)$ as one of its graphs can also be defined by meet tables.
There are only five possible values for the meet of two partitions in
$\mathcal{U}^9_3$; for each value, there is (up to isomorphism) exactly one
meet table. In this association scheme there are four graphs $G_i$,
one for each non-isomorphic meet table, except the table corresponding
to meet 3 (two partitions have meet 3 if and only if they are
identical). Two partitions $P,Q \in \mathcal{U}^9_3$ are adjacent in
the graph $G_i$ if the meet table $M_{P,Q}$ is isomorphic to the meet table
corresponding to $G_i$.

This can be generalized to $\mathcal{U}^n_k$ for positive integers $k$
and $n$, where $k$ divides $n$. Assume that there are $m$ pairwise
non-isomorphic meet tables for the partitions in
$\mathcal{U}^n_k$. For each meet table from this list, define the
graph $G_i$, for $i=0,\dots, m-1$ as follows: the vertex set of $G_i$
is the set of all uniform partitions $\mathcal{U}^n_k$, and two
partitions $P,Q
\in \mathcal{U}^n_k$ are adjacent in $G_i$ if the meet table $M_{P,Q}$ is
isomorphic to the meet table corresponding to $G_i$.

The big question here is: Does this, in general, define an association scheme for
$\mathcal{U}^n_k$?

The first problem is that for $P,Q \in \mathcal{U}^n_k$ it is possible
that the meet tables $M_{P,Q}$ and $M_{Q,P}$ are not isomorphic.  For
example, consider the partitions from $UQI(18,3)$
\begin{eqnarray*}
P &=& 1\; 2\; 3\; 4\; 5 \; 6\; \mid \; 7 \; 8\; 9\; 10\; 11\; 12\; \mid \; 13\; 14\; 15\; 16\; 17\; 18 \\
Q &=& 1\; 2\; 3\; 4\; 7 \;13\; \mid \; 5 \; 6\; 8\;  9\; 14\; 15\; \mid \; 10\; 11\; 12\; 16\; 17\; 18\; .
\end{eqnarray*}

The meet arrays for $P$ and $Q$  and for $Q$ and $P$ are
\begin{center}
\begin{tabular}{ccc} 
$M_{P,Q} = \left[
\begin{tabular} {ccc}
 4 & 2 & 0 \\
 1 & 2 & 3 \\
 1 & 2 & 3
\end{tabular} \right]$
& and &
$M_{Q,P} = \left[
\begin{tabular} {ccc}
 4 & 1 & 1 \\
 2 & 2 & 2 \\
 0 & 3 & 3
\end{tabular}
\right].$
\end{tabular}
\end{center}
These two tables are transposes of each other, but they are not
isomorphic.  This means that the graphs $G_i$, as defined above, may
actually be directed graphs. So this method does not in general define
an association scheme. But, it may define an {\em asymmetric association
scheme} (see Section~\ref{sec:introassscheme} for definition).

This example of meet tables $M_{P,Q}$ and $M_{Q,P}$ which are not
isomorphic is the smallest example for $UQI(3c,3)$. For $c=4,5$, for all
partitions $P, Q \in \mathcal{U}^{3c}_3$ the meet tables
$M_{P,Q}$ and $M_{Q,P}$ are isomorphic.  The two
cases $\mathcal{U}^{12}_3$ and $\mathcal{U}^{15}_3$ are considered in the next sections.

\subsection{Partitions in $\mathcal{U}^{12}_3$ }\label{sec:U12-3}

For $\mathcal{U}^{12}_3$ there are 9 non-isomorphic meet tables.  Let
$P=$ 1 2 3 4 $\mid$ 5 6 7 8 $\mid$ 9 10 11 12 be the fixed
partition. Each of the following nine partitions is a representative
for a class of non-isomorphic meet tables. Only the last
representative is qualitatively independent with $P$.

\[
\begin{tabular}{l@{\extracolsep{1cm}}l}
$1\; 2\; 3\; 4\; \mid \; 5\; 6\; 7\; 8\; \mid \; 9\; 10\; 11\; 12$
&
$1\; 2\; 3\; 4\; \mid \; 5\; 6\; 7\; 9\; \mid \; 8\; 10\; 11\; 12$\\
$1\; 2\; 3\; 4\; \mid \; 5\; 6\; 9\; 10\; \mid \; 7\; 8\; 11\; 12$
&
$1\; 2\; 3\; 5\; \mid \; 4\; 6\; 7\; 9\; \mid \; 8\; 10\; 11\; 12$\\
$1\; 2\; 3\; 5\; \mid \; 4\; 6\; 9\; 10\; \mid \; 7\; 8\; 11\; 12$
&
$1\; 2\; 3\; 5\; \mid \; 4\; 9\; 10\; 11\; \mid \; 6\; 7\; 8\; 12$\\
$1\; 2\; 5\; 6\; \mid \; 3\; 4\; 9\; 10\; \mid \; 7\; 8\; 11\; 12$
&
$1\; 2\; 5\; 9\; \mid \; 3\; 4\; 6\; 10\; \mid \; 7\; 8\; 11\; 12$\\
$1\; 2\; 5\; 9\; \mid \; 3\; 6\; 7\; 10\; \mid \; 4\; 8\; 11\; 12$ & \\
\end{tabular}
\]

For each of the nine non-isomorphic meet
tables, define a graph $G_i$ with vertex set $\mathcal{U}^{12}_3$. Two
partitions will have an edge in graph $G_i$ if and only if their meet
table is isomorphic to the meet table corresponding to the graph
$G_i$.  Using the equitable partition described in
Section~\ref{sec:betterep}, it is possible to calculate the
eigenvalues of these graphs. The table below is the modified matrix of
eigenvalues of these graphs. As only the last representative is
qualitatively independent with $P$, the final column of this table is
the set of eigenvalues of the graph $UQI(12,3)$.

\begin{table}[h]
\begin{eqnarray*}
\left(
\begin{tabular}{c|ccccccccc}
275  & 12 & 21 & -48 & -24 & -16 & 18 &  36 &    0 \\
2673 & -4 &  1 &  -8 &  24 &   0 & -2 & -20 &    8 \\
462  &-12 &  9 &  36 & -72 &   8 &  6 &  36 &  -12 \\
616  &  6 & -9 & -18 & -36 &  20 & 18 &   0 &   18 \\
1408 &  3 & -6 &   6 &   3 &  -7 & -9 &  36 &  -27 \\
132  &  6 & -9 &  72 & -36 & -40 & 48 & -90 &   48 \\
154  & 18 &  9 &  36 & -72 &   8 &-54 & -54 &  108 \\
54   & 26 & 21 &  92 & 144 &  40 & 18 & -90 & -252 \\
1    & 48 & 54 & 576 &1728 & 128 &216 &1296 & 1728 
\end{tabular}
\right)
\end{eqnarray*}
\caption{Modified matrix of eigenvalues for the graphs $G_i$ (with vertex set $\mathcal{U}^{12}_3$)}
\end{table}

If the set of graphs $\{ G_i \: : \: i=1,\dots, 9\}$ describes an
association scheme, then the ratio bound for cliques,
Inequality~(\ref{ratio2}), would hold for the graph $UQI(12,3)$. If this
bound held, then $\omega(UQI(12,3)) \leq 7$. This particular bound is
interesting because the exact value of $N(12,3)$ is not known; it is
known that $7 \leq N(12,3)$ and that $12 \leq CAN(8,3) \leq 13$.  If
it is true that $\omega(UQI(12,3)) \leq 7$, it still may be true that
a $CA(12,8,3)$ exists, but it would not be balanced.  So, if a
$CA(12,8,3)$ exists and the above is an association scheme, then the
graph $UQI(12,3)$ would not be a core of $QI(12,3)$. This would mean
the answer to Question~\ref{quest:generalcores} is negative.

\begin{conj}
The set of graphs, $\{ G_i \; :\; i=2,\dots, 9\}$, defined above form
an association scheme on $\mathcal{U}^{12}_3$.
\end{conj}

\begin{conj}
A $CA(12,8,3)$ does not exist, so $CAN(8,3)=13$.
\end{conj}

\subsection{Partitions in $\mathcal{U}^{15}_3$}\label{sec:U15-3}
 
For $\mathcal{U}^{15}_3$ there are 13 non-isomorphic meet tables. Use
$P=$ 1 2 3 4 5 $\mid$ 6 7 8 9 10 $\mid$ 11 12 13 14 15 as the fixed
partition.  Below is a list of representatives of each non-isomorphic
meet table with respect to $P$.  In this case, the last three
representatives are qualitatively independent with the fixed
partition.

\[
\begin{tabular}{l@{\extracolsep{1cm}}l}
$1\; 2\; 3\; 4\; 5 \; \mid \; 6\; 7\; 8\; 9\; 10\; \mid \; 11\; 12\; 13\; 14\; 15$ &
$1\; 2\; 3\; 4\; 5 \; \mid \; 6\; 7\; 8\; 9\; 11\; \mid \; 10\; 12\; 13\; 14\; 15$ \\
$1\; 2\; 3\; 4\; 5 \; \mid \; 6\; 7\; 8\; 11\; 12\; \mid \; 9\; 10\; 13\; 14\; 15$ &
$1\; 2\; 3\; 4\; 6 \; \mid \; 5\; 7\; 8\; 9\; 11\; \mid \; 10\; 12\; 13\; 14\; 15$ \\
$1\; 2\; 3\; 4\; 6 \; \mid \; 5\; 7\; 8\; 11\; 12\; \mid \; 9\; 10\; 13\; 14\; 15$ &
$1\; 2\; 3\; 4\; 6 \; \mid \; 5\; 7\; 11\; 12\; 13\; \mid \; 8\; 9\; 10\; 14\; 15$ \\
$1\; 2\; 3\; 4\; 6 \; \mid \; 5\; 11\; 12\; 13\; 14\; \mid \; 7\; 8\; 9\; 10\; 15$ &
$1\; 2\; 3\; 6\; 7 \; \mid \; 4\; 5\; 8\; 11\; 12\; \mid \; 9\; 10\; 13\; 14\; 15$ \\
$1\; 2\; 3\; 6\; 7 \; \mid \; 4\; 5\; 11\; 12\; 13\; \mid \; 8\; 9\; 10\; 14\; 15$ & 
$1\; 2\; 3\; 6\; 11 \; \mid \; 4\; 5\; 7\; 8\; 12\; \mid \; 9\; 10\; 13\; 14\; 15$ \\
$1\; 2\; 3\; 6\; 11 \; \mid \; 4\; 7\; 8\; 9\; 12\; \mid \; 5\; 10\; 13\; 14\; 15$ &
$1\; 2\; 3\; 6\; 11 \; \mid \; 4\; 7\; 8\; 12\; 13\; \mid \; 5\; 9\; 10\; 14\; 15$ \\
$1\; 2\; 6\; 7\; 11 \; \mid \; 3\; 4\; 8\; 12\; 13\; \mid \; 5\; 9\; 10\; 14\; 15$ & 
\end{tabular}
\]

Similar to the case for $\mathcal{U}^{12}_3$, for each of the 13
non-isomorphic meet tables, it is possible to a build graph $G_i$, for
$i = 1,\dots,13$, with vertex set $\mathcal{U}^{15}_3$ and vertices
adjacent in $G_i$ if and only if their meet table is isomorphic to the
meet table corresponding to $G_i$.  The table below is the modified
matrix of eigenvalues of these graphs.

\begin{table}[h]
{\footnotesize
\begin{eqnarray*}
\left(
\begin{tabular}{c|ccc cccc cccc c}
2002  & 15 & -18 & 126 & 144 & -276 & -50 & 234 & 80 & -540 & 176 & 162 & -45 \\
32032 & 0 & -3 & -24 & 9 & 54 & -5 & 54 & -25 & -90 & 56 & -108 & 81 \\
90    & 47 & 132 & 450 & 1440 & 1240 & 110 & 600 & 320 & 720 & 160 & -1980 & -3240 \\
910   & 27 & 102 & 0 & 240 & -60 & 10 & 150 & 20 & -480 & -640 & -180 & 810 \\
7007  & 15 & -18 & 36 & -36 & 96 & 10 & 54 & -100 & 360 & -64 & -108 & -54 \\
25025 & -9 & 6 & 12 & 36 & -72 & 10 & -54 & 20 & 72 & 32 & -108 & 54 \\
350   & 36 & 66 & 252 & 186 & 228 & 55 & -984 & -340 & -1184 & 512 & 792 & 324 \\
1638  & 20 & 102 & -84 & -306 & -116 & -24 & 24 & 20 & 306 & 256 & 72 & -324 \\
1     & 75 & 300 & 1500 & 9000 & 6000 & 250 & 9000 & 2000 & 36000 & 8000 & 27000 &27000\\
11583 & 5 & -8 & 16 & -16 & 64 & -30 & -156 & 40 & 160 & -64 & -48 & 36 \\
22113 & 5 & -8 & -44 & 64 & -16 & 10 & -26 & 20 & -40 & -24 & 172 & -114 \\
1925  & 21 & -24 & 60 & -432 & 96 & 58 & 54 & 164 & -216 & 56 & -108 & 270 \\
21450 & -9 & 6 &30 & -72 & 36 & 2 & 54 & -16 & -36 & -64 & 162 & -90 
\end{tabular}
\right)
\end{eqnarray*}
}
\caption{Modified matrix of eigenvalues for the graphs $G_i$ (with vertex set $\mathcal{U}^{15}_3$)}
\end{table}

There are three representatives that are qualitatively independent
with the partition $P$, this means there are three graphs, namely $G_{11}, G_{12}$ and
$G_{13}$, in which partitions that are adjacent are qualitatively independent.
Thus the graph $QI(15,3)$ is not a single graph in this scheme.  So we
end with one last conjecture.

\begin{conj}
The graphs $G_i$, for $i=2,\dots , 13$,  described above give an association scheme for $\mathcal{U}^{15}_3$.
\end{conj}

%%%%%%%%%%%%%%%%%%%%%%%%%%%%%%%%%%%%%%%%%%%%%%%%%%%%
\chapter{Partition Systems}\label{partitionsystems}

In Chapter~\ref{EST}, we considered extremal problems for set
systems. These problems included Sperner's Theorem and the
Erd\H{o}s-Ko-Rado Theorem. There have been attempts to extend these
results to extremal problems in which the elements in the system are
families of sets, called {\em clouds}\index{cloud}, rather than sets.  Such
problems are called {\em higher order extremal problems}.  Ahlswede,
Cai and Zhang~\cite{MR1739309} give a good overview of such problems.
Most problems considered in~\cite{MR1739309} require that the clouds
be pairwise disjoint, that is, no set can occur in two distinct
clouds. With this restriction, the direct generalization of the
Erd\H{o}s-Ko-Rado Theorem for disjoint clouds proved to be
false~\cite{intersectingsystems}. Erd\H{o}s and
Sz\'ekely~\cite{MR2001f:05148} survey higher order Erd\H{o}s-Ko-Rado
Theorems where clouds are substituted by set systems with additional
structure and the disjointness requirement for pairs of set systems is
dropped.  They consider, among other cases, the particular case in
which each structure is a set partition. It is this type of problem
that is connected to covering arrays.

Sperner's Theorem and the Erd\H{o}s-Ko-Rado Theorem completely
determine the maximum cardinality of a set of qualitatively
independent sets. Since the rows of a binary covering array correspond
to subsets of an $n$-set, this determines the minimum size of a binary
covering array (Theorem~\ref{thm:katona73}).  The rows of a general
covering array, $CA(n,r,k)$, correspond to $k$-partitions of an $n$-set,
so it seems likely that results on extremal partition systems could
help determine the size of general covering arrays. In this chapter,
we give theorems for partition systems that are extensions of
Sperner's Theorem and the Erd\H{o}s-Ko-Rado Theorem to partition
systems.
%These are the first results in {\em extremal partition systems}

Where appropriate, the notation in this chapter will follow the notation
used in Chapter~\ref{EST}. In particular, the collection of all
$c$-sets of an $n$-set is denoted by ${[n] \choose c}$, the set of all
$k$-partitions of an $n$-set is denoted by $\mathcal{P}^n_k$ with
$S(n,k) = |\mathcal{P}^n_k|$ and the the set of all uniform
$k$-partitions of an $n$-set is denoted $\mathcal{U}^n_k$ with
$U(n,k) = |\mathcal{U}^n_k|$ .

The results in this chapter are contained in the
papers~\cite{EKR:partitions, karen4}.

\section{Sperner Partition Systems}\label{sec:spernerpartitions}

In this section, theorems similar to Sperner's Theorem
(Theorem~\ref{thm:sperner}) for partition systems are proven.

Recall from Section~\ref{sec:sperner} that sets $A$ and $B$ are {\em
incomparable} if they have the property that $A \not\subseteq B$ and
$B \not\subseteq A$.  From Definition~\ref{defn:sperner}, a Sperner
set system is a set system $\mathcal{A}$ with the property that any
distinct $A,B \in \mathcal{A}$ are incomparable.

We will extend this property to partitions; this definition
does not coincide with the extension of incomparability to
sequences given by Gargano, K\"orner and
Vaccaro~\cite{gargano:92,MR94e:05024,gargano:94}.

\begin{defn}[Sperner Property for Partitions]\index{Sperner property for partitions}
\label{def:spernerpartitions}
Let $n$ and $k$ be positive integers. 
Two $k$-partitions of an $n$-set, $P,Q \in \mathcal{P}^n_k$, with
$P=\{P_1,\dots,P_k\}$ and $Q=\{Q_1,\dots,Q_k\}$ have the {\em Sperner property} if
\[
P_i \not\subseteq Q_j \mbox{ and }  Q_i \not\subseteq P_j \mbox{ for all } i,j \in \{1,\dots,k\}.
\]
\end{defn}

\begin{defn}[Sperner Partition System]\index{Sperner partition system}
\label{def:spernerpartitionsystems}Let $n$ and $k$ be positive integers.
A partition system $\mathcal{P} \subseteq\mathcal{P}^n_k$ is a {\em
Sperner partition system} if all distinct $P,Q \in \mathcal{P}$ have
the Sperner property.
\end{defn}

Let $SP(n,k)$ denote the maximum cardinality of a Sperner
partition system in $\mathcal{P}^n_k$. 

If $\mathcal{P}$ is a partition system, then $\mathcal{P}$ is a
Sperner partition system if and only if all the partitions in
$\mathcal{P}$ are disjoint (no two partitions have a common class) and
the union of all the partitions in $\mathcal{P}$ is a Sperner set
system. In this sense, Sperner partition systems can be considered
as resolvable Sperner set systems. 

Further, any set of qualitatively independent partitions is a Sperner
partition system. Moreover, for $k=2$, the definition of the Sperner
property for partitions is equivalent to the definition of qualitative
independence.

%KAren check title
\subsection{Sperner's Theorem for Partition Systems in $\mathcal{P}^{ck}_k$}

In this section, we show that for positive integers $n,k,c$ with
$n=ck$, the largest Sperner $k$-partition system of an $n$-set is a
uniform partition system with cardinality ${n-1 \choose c-1}$.

%KAren move this
When $n=ck$, a 1-factor of the complete uniform hypergraph $K_n^{(c)}$
is equivalent to a uniform $k$-partition of an $n$-set, and a
1-factorization of $K_n^{(c)}$ corresponds to a Sperner partition
system.  From Theorem~\ref{thm:1factorization}, if $c$ divides $n$,
the hypergraph $K_n^{(c)}$ has a 1-factorization with ${n-1 \choose
c-1}$ factors.

\begin{lemma}\label{lem:maxSPS}
Let $n,k,c$ be positive integers with $n=ck$, then there exists a Sperner
partition system in $\mathcal{P}_k^n$ of cardinality ${n-1 \choose c-1}$.
\end{lemma}

The proof that this is the largest Sperner partition system uses a
result by Kleitman and Milner on Sperner set systems. For a set system
$\mathcal{A}$, define the {\em volume}\index{volume of a set system}
of $\mathcal{A}$ to be $t(\mathcal{A}) = \sum_{A \in \mathcal{A}}|A|$.

%KAREN check the terms here, sperner on a n set

\begin{thm}[Kleitman and Milner~\cite{MR0337638}]\label{thm:KM}
Let $\mathcal{A}$ be a Sperner set system on an $n$-set with $|\mathcal{A}| \geq
{n\choose c}$ and $c \leq \frac{n}{2}$, then
\begin{eqnarray*}
\frac{t( \mathcal{A} )}{|\mathcal{A}|} \geq c.
\end{eqnarray*}
This inequality is strict in all cases except when $\mathcal{A} = {[n] \choose c}$. 
\end{thm}

%KAREN
%sentence aobut maire?

\begin{thm}\label{thm:spernerpartition}
Let $n,k,c$ be positive integers with $n=kc$.  Then ${SP}(n,k) =
\frac{1}{k}{n \choose c}$.  Moreover, a Sperner partition system has
cardinality $\frac{1}{k}{n \choose c}$ only if it is a $c$-uniform
partition system.
\end{thm}

\begin{proof}
If $k=1$ then $\mathcal{P}^n_k$ has only one partition, namely
$\{ \{1,\dots,n\} \}$. So ${SP}(n,k) = 1$.

Let $\mathcal{P} \subseteq \mathcal{P}^n_k$ be a Sperner partition system.  Let
$\mathcal{A}$ be the Sperner set system formed by taking all classes from
all the partitions in $\mathcal{P}$. Thus $|\mathcal{A}| = k
|\mathcal{P}|$.

We show that $|\mathcal{A}| \leq {n \choose c}$.
Let $p_i$, $i \in \{1,\dots,n\}$ be the number of sets in  $\mathcal{A}$ with size $i$.

By the LYM Inequality (Inequality~(\ref{equationLYM}), Section~\ref{sec:sperner}), we have
\[
\sum_{i=1}^n \frac{p_i}{ {n \choose i} } \leq 1.
\]
Following the notation from~\cite{MR1448245}, define the function
$f(i) = {n \choose i}^{-1}$. With this, we get
\begin{equation}\label{eqn:LYMlike}
\sum_{i=1}^n \frac{p_i}{|\mathcal{A}|} f(i)  \leq \frac{1}{|\mathcal{A}|}.
\end{equation}

In~\cite{MR1448245} it is shown that the function $f(i)$ can be extended to a convex
function on the real numbers by
\[
f(i + u) = (1-u)f(i) + u f(i+1) \quad \mbox{ for } \quad 0 \leq u \leq 1. 
\]
%write the proof that it is convex?
Since the set system $\mathcal{A}$ is formed from a partition system with $n=ck$,
\begin{eqnarray}\label{eq:above}
\sum_{i=1}^n i p_i = \sum_{A\in \mathcal{A}} |A| = n|\mathcal{P}| = ck|\mathcal{P}| =c|\mathcal{A}|.
\end{eqnarray}
Using Equation~(\ref{eq:above}); the fact that $f$ is a convex function
and that $\sum_{i=1}^n \frac{p_i}{|\mathcal{A}|} = 1$; and Inequality~(\ref{eqn:LYMlike}),
\[
f(c) = f \left( \sum_{i=1}^n i \frac{p_i}{|\mathcal{A}|} \right) 
\leq \sum_{i=1}^n f(i) \frac{p_i}{|\mathcal{A}|} \leq  \frac{1}{|\mathcal{A}|}.
\]
Thus, $|\mathcal{A}| \leq {n \choose c}$ and 
\[
|\mathcal{P}| \leq \frac{1}{k}{n \choose c}.
\]

Next, we need to show that a Sperner partition system meets this bound,
then it is $c$-uniform. If $k=1$ this is clearly true since the only
partition in $\mathcal{P}^n_1$ is the $n$-uniform partition $\{
\{1,\dots,n\} \}$.  Assume that for $k \geq 2$, $|\mathcal{P}|
=\frac{1}{k}{n \choose c}$. Let $\mathcal{A}$ be the Sperner set
system formed from all the classes in $\mathcal{P}$, then
$|\mathcal{A}| = {n \choose c}$ and $c \leq \frac{n}{2}$.

Since $\frac{t(\mathcal{A})}{|\mathcal{A}|} = \frac{\Sigma_{A \in
\mathcal{A}} |A|}{|\mathcal{A}|}= c$, from Theorem~\ref{thm:KM}
it follows that $\mathcal{A}={[n] \choose c}$ and $\mathcal{P}$ is a uniform
partition system.

Finally, a Sperner partition system $\mathcal{P}$ with $|\mathcal{P}| =
\frac{1}{k}{n \choose c}$ exists for all $n,k,c$ with $n=ck$ from
Lemma~\ref{lem:maxSPS}.
\end{proof}

Theorem~\ref{thm:spernerpartition} is a natural extension of Sperner's
Theorem for sets (Theorem~\ref{thm:sperner}). Sperner's Theorem for set
systems says that the Sperner set system with maximum cardinality on
an $n$-set is the collection of all $\lfloor \frac{n}{2}
\rfloor$-sets.  Theorem~\ref{thm:spernerpartition} says that for
integers $n,k$ such that $k$ divides $n$, the Sperner $k$-partition
system on an $n$-set with the largest cardinality is the collection of
all $(\frac{n}{k})$-sets arranged in resolution classes.

A collection of qualitatively independent partitions is a Sperner
partition system, but the converse of this is not true in general. In
fact, for $k>2$, from Theorem~\ref{thm:spernerpartition}, the
cardinality of the largest Sperner $k$-partition system can be much
larger than the cardinality of a collection of qualitatively
independent $k$-partitions.  For example, for a positive integer $k$,
there can be at most $k+1$ qualitatively independent $k$-partitions of
a $k^2$-set (Inequality~(\ref{eq:maxMOLS}), Section~\ref{sec:MOLS} and
Theorem~\ref{thm:equivdesigns}), while ${SP}(k^2,k) = {k^2-1 \choose
k-1}$.

%When $k=2$, any Sperner $k$-partition system is also a qualitatively
%independent $k$-partition system.

\subsection{A Bound on the Cardinality of Sperner Partition Systems in $\mathcal{P}^n_k$}

For integers $k$ and $n$, if $k$ does not divide $n$, then we have
a bound on the cardinality of a Sperner partition system in $\mathcal{P}^n_k$.

\begin{thm}\label{thm:AUbound}
Let $n,k,c$ be positive integers with $n=ck+r$, where $0 \leq r < k$.
Then, 
\[
{SP}(n,k) \leq \frac{ {n \choose c} }{(k-r) + \frac{r(c+1)}{n-c}} .
\]
\end{thm}

\begin{proof}
Let $\mathcal{P} \subseteq \mathcal{P}^n_k$ be a Sperner partition system.
Let $\mathcal{A}$ be the Sperner set system formed by taking all
classes from all the partitions in $\mathcal{P}$. Thus 
$|\mathcal{A}| = k |\mathcal{P}|$.

Let $p_i$, $i \in \{1,\dots,n\}$ be the number of sets in
$\mathcal{A}$ with size $i$.

Again, using the function $f(i) = {n \choose i}^{-1}$ from~\cite{MR1448245}
and the LYM Inequality (Inequality~(\ref{equationLYM}) in Section~\ref{sec:sperner})
\[
\sum_{i=1}^n \frac{p_i}{|\mathcal{A}|} f(i)  \leq \frac{1}{|\mathcal{A}|}.
\]

Since the function $f(i)$ can be extended to a convex function on the real
numbers,
\[
f\left(\frac{n}{k}\right)    =  f \left( \sum_{A\in \mathcal{A}} \frac{|A|}{|\mathcal{A}|} \right)  
                  =   f \left( \sum_{i=1}^n i \frac{p_i}{|\mathcal{A}|} \right)  
                \leq  \sum_{i=1}^n f(i) \frac{p_i}{|\mathcal{A}|}  
                \leq  \frac{1}{|\mathcal{A}|}.
\]
By the definition of $f(i)$, 
\[
f\left(\frac{n}{k}\right) = f\left(\frac{ck+r}{k}\right)= f\left(c+\frac{r}{k}\right)
= \left(1-\frac{r}{k}\right){n \choose c}^{-1} + \frac{r}{k}{n \choose c+1}^{-1}.
\]

Thus, 
%\begin{eqnarray*}
\[
|\mathcal{A}| \leq \frac{1}{(1-\frac{r}{k}){n \choose c}^{-1} + \frac{r}{k}{n \choose c+1}^{-1}} 
              = {n \choose c} \frac{k}{(k-r) + \frac{r(c+1)}{n-c}} 
\]
%\end{eqnarray*}
 and 
$|\mathcal{P}| \leq \frac{ {n \choose c}}{(k-r) + \frac{r(c+1)}{n-c}} $.
\end{proof}

It would be better to know the exact cardinality and structure of the
largest Sperner partition system. We conjecture that the largest
Sperner partition system is an almost-uniform partition system (see
Section~\ref{setpartitions} for definition of almost-uniform partition systems).

\begin{conj}
Let $n,k$ be positive integers.  The largest Sperner partition system
in $\mathcal{P}^n_k$ is an almost-uniform partition system.
\end{conj}

\remove{For the case when $k$ does not divide $n$, the classes in a
$k$-partition of an $n$-set can not all have the same cardinality so it is
harder to see what a natural extension of the Sperner's Theorem for
set systems to partitions system would be.}

Similar to the case where $k$ divides $n$, the maximum cardinality of a system
of qualitatively independent $k$-partitions of an $n$-set is much
smaller than ${SP}(n,k)$. This is best seen by considering the
asymptotic growth of ${SP}(n,k)$ and comparing it to the
asymptotic growth of the maximum cardinality of a system of
qualitatively independent $k$-partitions of an $n$-set,
${N}(n,k)$.

\subsection{Asymptotic Growth of Maximal Sperner Partition Systems}\label{SpernerAsymptotics}

In this section, we give the asymptotic growth of Sperner partition
systems. This result is similar to the result for the asymptotic growth
of qualitatively independent partition systems cited in Section~\ref{asymptotics}.

\begin{lemma}\label{Spernerorder}
Let $n,k,c$ be positive integers with $n=ck+r$ and $0 \leq r < k$.
Then  ${SP}(n,k) \leq {SP}(n+1,k)$.
\end{lemma}
\begin{proof}
Let $\mathcal{P} \subset \mathcal{P}^n_k$ be a Sperner partition
system in $\mathcal{P}^n_k$ with $| \mathcal{P} | = SP(n,k)$. For
each partition in $\mathcal{P}$, add the element $n+1$ to the class in
the partition that contains the element 1.  This is a Sperner
partition system in $\mathcal{P}^{n+1}_k$ and the result follows.
\end{proof}

The following bounds on $SP(n,k)$ follow from
Theorem~\ref{thm:spernerpartition} and Lemma~\ref{Spernerorder}.

\begin{cor}\label{cor:bounds}
Let $n,k,c$ be positive integers with $k \geq 2$ and $n=ck+r$, where $0 \leq r < k$.
Then
\[
\frac{1}{k}{ck \choose c} \leq SP(n,k) \leq \frac{1}{k}{(c+1)k \choose c+1}.
\remove{\leq \frac{1}{k-r + \frac{r(c+1)}{n-c}}{n \choose c}.}
\]
\end{cor}

The above upper bound on $SP(n,k)$ is weaker than the upper bound in Theorem~\ref{thm:AUbound}.

The bounds in Corollary~\ref{cor:bounds} will be used to produce a
result on the asymptotic growth of $SP(n,k)$. Following the
result of Gargano, K\"orner and Vaccaro~\cite{MR94e:05024}
(see Section~\ref{asymptotics}), we consider the growth of
 $\lim_{n \rightarrow \infty}\frac{\log SP(n,k) }{n}$.

First we need a technical result.
\begin{lemma}[\cite{MR1397498}]\label{lem:lnn}
For all positive integers $n$,
\begin{eqnarray*}
\log_e{n!} = n \log_e{n} - n + \frac{\log_e{n}}{2} + O(1).
\end{eqnarray*}
\end{lemma}

\begin{thm}\label{thm:spernerasymp}
For positive integers $n,k,c$ with $k \geq 2$ and $n=kc+r$ where $0 \leq r < k$,
\[
\limsup_{n \rightarrow \infty} \frac{\log SP(n,k)}{n} = \log \frac{k}{k-1}+\frac{\log(k-1)}{k}.    
\]
\end{thm}

\begin{proof}
Since different logarithm bases differ by a constant factor, to prove
the theorem it is sufficient to prove that for a positive
integer $k$ with $k \geq 2$,
\[
\lim_{n \rightarrow \infty} \frac{\log_e SP(n,k)}{n} = \log_e k - \frac{k-1}{k}\log_e (k-1).
\]
Further, by Corollary~\ref{cor:bounds}, it is enough to show that
\begin{eqnarray}
\lim_{n \rightarrow \infty}\frac{ \log_e \frac{1}{k} {ck \choose c}}{n} 
                          &=& \log_e \frac{k}{k-1}+\frac{\log_e(k-1)}{k}, 
\label{eqn1}
\end{eqnarray}
and
\begin{eqnarray}
\lim_{n \rightarrow \infty}\frac{\log_e  \frac{1}{k} {(c+1)k \choose c+1} } {n} 
                                    &=&  \log_e\frac{k}{k-1} +\frac{\log_e(k-1)}{k},
                                    \label{eqn2}
\end{eqnarray}
where $c = c(n) = \lfloor \frac{n}{k}\rfloor$ and $r=r(n)$ such that $n=ck+r$ and $0\leq r(n) <k$.

First expand $\log_e{{ck \choose c}}$ using Lemma~\ref{lem:lnn}:
\begin{eqnarray*}
\log_e{{ck \choose c}} &=&\log_e{(ck)!} - \log_e{c!} - \log_e{(ck-c)!}\\
          &=& ck \log_e(ck) - ck + \frac{ \log_e (ck)}{2}- c \log_e c + c - \frac{ \log_e c}{2} \\
          &&  - (ck-c) \log_e(ck-c) + ck - c - \frac{ \log_e(ck-c)}{2} + O(1)     \\
%          &=& ck \log_e c + ck \log_e k - ck + \frac{\log_e c}{2}+ \frac{\log_e k}{2} - c \log_e c + c - \frac{\log_e c }{2} \\
%          &&  - ck\log_e(k-1)+c\log_e(k-1) - ck\log_e c + c\log_e c + ck - c - \frac{\log_e(k-1)}{2}-\frac{\log_e c}{2} \\
  &=&  ck\log_e k + \frac{\log_e k}{2} - ck\log_e(k-1) + c\log_e(k-1) -\frac{\log_e(k-1)}{2} \\
   && -\frac{\log_e c}{2}+O(1).
\end{eqnarray*}    
\orphan

{}From this and the fact that $n=ck+r$, we obtain Equation~(\ref{eqn1}):
\begin{eqnarray*}
\lim_{n \rightarrow \infty } \frac{\log_e \frac{1}{k} {ck \choose c} }{n}
  &=& \lim_{n \rightarrow \infty }\left( \frac{ck\log_e k}{n} + \frac{\log_e k}{2n} - \frac{ck\log_e(k-1)}{n} \right.\\ 
  && \ \ \ \ \left. + \frac{c\log_e(k-1)}{n} -\frac{\log_e(k-1)}{2n}-\frac{\log_e c}{2n} + \frac{O(1)}{n} - \frac{\log_e k}{n}\right) \\ 
  &=& \lim_{n \rightarrow \infty }\left( \frac{(n-r)\log_e k}{n} + \frac{\log_e k}{2n} - \frac{(n-r)\log_e(k-1)}{n} \right.\\ 
  && \ \ \ \left. + \frac{(n\!-\!r)\log_e(k\!-\!1)}{kn} -\frac{\log_e(k\!-\!1)}{2n}-\frac{\log_e \frac{n-r}{k}}{2n} + \frac{O(1)}{n} - \frac{\log_e k}{n}\right) \\ 
 &=&     \log_e k  - \log_e(k-1)  + \frac{1}{k} \log_e(k-1) 
\end{eqnarray*}

Similarly, to obtain Equation~(\ref{eqn2}), we first use Lemma~\ref{lem:lnn}:
\begin{eqnarray*}
\log_e{{(c+1)k \choose c+1}}
 &=&  (c+1)k\log_e k + \frac{\log_e k}{2} - (c+1)k\log_e(k-1)\\
  &&   + (c+1)\log_e(k-1) -\frac{\log_e(k-1)}{2}-\frac{\log_e (c+1)}{2}+O(1).
\end{eqnarray*}

Dividing by $n$, and taking the limit in the above equation, we obtain Equation~(\ref{eqn2}):
\begin{eqnarray*}
\lim_{n\rightarrow \infty} \frac{ \log_e\left(\frac{1}{k} {(c+1)k \choose c+1} \right) }{n}
&=& \lim_{n\rightarrow \infty} \left(
    \frac{(c+1)k\log_e k}{n} + \frac{\log_e k}{2n} - \frac{(c+1)k\log_e(k-1)}{n} \right. \\
  && \ \ \ \ + \frac{(c+1)\log_e(k-1)}{n} -\frac{\log_e(k-1)}{2n}-\frac{\log_e (c+1)}{2n} \\
   && \ \ \ \ \left.  +\frac{O(1)}{n} - \frac{\log_e k}{n} \right) \\
&=& \lim_{n\rightarrow \infty} \left(
    \frac{ (n-r+k)\log_e k}{n} + \frac{\log_e k}{2n} - \frac{(n-r+k)\log_e(k-1)}{n} \right. \\
  && \ \ \ \ + \frac{(n-r+k)\log_e(k-1)}{kn} -\frac{\log_e(k-1)}{2n}-\frac{\log_e \frac{n-r+k}{k}}{2n} \\
   && \ \ \ \ \left.  +\frac{O(1)}{n} - \frac{\log_e k}{n} \right) \\
&=&  \log_e k  - \log_e (k-1) + \frac{1}{k}\log_e (k-1).
\end{eqnarray*}
\end{proof}

Recall from Section~\ref{asymptotics} that for all $k \geq 2$,
\[
\limsup_{n \rightarrow \infty} \frac{\log_2 N(n,k)}{n} = \frac{2}{k},
\]
and from Theorem~\ref{thm:spernerasymp} above
\begin{eqnarray*}
\limsup_{n \rightarrow \infty} \frac{\log_2 SP(n,k)}{n}
&=& \limsup_{n \rightarrow \infty} \frac{1}{\log_e 2}\frac{ \log_e SP(n,k)}{n} \\
& =& \frac{1}{\log_e 2} \left( \log_e \frac{k}{k-1}+\frac{\log_e(k-1)}{k} \right).
\end{eqnarray*}
It is no surprise that for $k=2$, the asymptotic growth of $N(n,k)$ is
the same as the asymptotic growth for $SP(n,k)$.  For $k > 2$,
asymptotic growth of $SP(n,k)$ is larger than that of $N(k,n)$.  These
results on Sperner partition systems are not enough to give a bound on
the maximal collection of qualitatively independent partitions. In
the next section we extend the Erd\H{o}s-Ko-Rado Theorem to
partition systems.

\section{Intersecting Partitions}\label{intersectingpartitions}

Following Section~\ref{sect:EKRapps}, it would be interesting to find
a version of the Erd\H{o}s-Ko-Rado Theorem for partitions.  It is not
obvious how to extend the definition of intersection from sets to
partitions.  Two sets intersect if they have an element in common, and
since partitions are a collection of sets, a natural extension of
intersection is that two partitions intersect if they have a common
class. Further, two partitions {\em $t$-intersect} if they have at least $t$
classes in common.  This is the type of intersection considered in
this section and Section~\ref{matchings}, while in Section~\ref{PIP},
a different type of intersection, which we call {\em partial
intersection}, is considered.

\begin{defn}[$t$-Intersecting Partition System]\index{$t$-intersecting partition system}
\label{defn:intersectingpartitions}
A partition system $\mathcal{P} \subseteq \mathcal{P}^n_k$ is {\em
$t$-intersecting} if $|P \cap Q|\geq t$, for all $P,Q \in
\mathcal{P}$.  
\end{defn}

Recall from Section~\ref{sec:ISS}, for positive integers $t,k$ and
$n$, a $k$-uniform trivially $t$-intersecting set system of an
$n$-set is a set system formed by all $k$-sets of an $n$-set that
contain a given $t$-set. Similarly, a partition system $\mathcal{P}
\subseteq \mathcal{P}^n_k$ is a {\em trivially $t$-intersecting partition
system} if $\mathcal{P}$ is equal, up to a permutation on $\{1,\dots,n\}$, to
\[
\mathcal{Q}(n,k,t)=\left\{P \in \mathcal{P}^n_k: 
   \{\{1\},\{2\},\ldots, \{t\}\}\subseteq P\right\}. 
\]
The cardinality of a trivially $t$-intersecting partition system in
$\mathcal{P}^n_k$ is $S(n-t,k-t) = |\mathcal{P}^{n-t}_{k-t}|$ (recall
that $S(n,k)$ is the Stirling number of the second type, defined in
Section~\ref{setpartitions}).

For positive integers $n,k,c$ with $n=ck$, $\mathcal{P} \subseteq
\mathcal{U}^n_k$ is a {\em trivially $t$-intersecting uniform partition
system} \index{trivially $t$-intersecting uniform partition system}
if $\mathcal{P}$ is equal, up to a permutation on $[1,n]$, to
\begin{eqnarray*}
\mathcal{P}(n,k,t)&=&\left\{P \in \mathcal{U}^n_k:\right.\\ 
           && \left. \{[1,c],[c+1,2c],\ldots, [(t-1)c+1,tc]\}\subseteq P\right\}. 
\end{eqnarray*}
The cardinality of a trivially $t$-intersecting uniform partition
system in $\mathcal{U}^{ck}_k$ is  
\begin{eqnarray}\label{eq:sizeoftrivial}
U((k-t)c,k-t)=|\mathcal{U}^{(k-t)c}_{k-t}|.
\end{eqnarray}

Erd\H{o}s and Sz\'ekely observe that the following 
Erd\H{o}s-Ko-Rado type theorem for $t$-intersecting partition systems holds.

\begin{thm}[\cite{MR2001f:05148}]
Let $k \geq t \geq 1$ be positive integers. There exists a function
$n_0(k,t)$ such that if $n\geq n_0(k,t)$, and
$\mathcal{P}\subseteq \mathcal{P}^n_k$ is a $t$-intersecting partition
system, then $|\mathcal{P}| \leq S(n-t,k-t).$ This bound is
attained by a trivially $t$-intersecting partition system.
\end{thm}
We prove analogous theorems for uniform partition systems that
guarantee the uniqueness, up to isomorphism, of the maximal system.  
Our first theorem completely settles the case $t=1$.

\begin{thm}\label{two}
Let $n\geq k\geq 1$ be positive integers such that $k$ divides
$n$. Let $\mathcal{P}\subseteq \mathcal{U}^n_k$ be a $1$-intersecting uniform
partition system. Let $c=n/k$ be the size of a class in each
partition. Then, $|\mathcal{P}| \leq U(n-c,k-1).$ Moreover, this bound
holds with equality if and only if $\mathcal{P}$ is a trivially $1$-intersecting
uniform partition system.
\end{thm}

Our second theorem deals with general $t$ and determines the
cardinality and structure of maximal $t$-intersecting uniform
partition systems when $n$ is sufficiently large.  In this theorem,
$n$ can be sufficiently large with respect to $k$ and $t$. Alternately, if $c
\geq t+2$, then $n$ can be sufficiently large with respect to $c$ and $t$.

\begin{thm}\label{one}
Let $k \geq t \geq 1$. There exist functions $n_0(k,t)$ and $n_1(c,t)$
such that if ($n\geq n_0(k,t)$) or ($c\geq t+2$ and $n \geq n_1(c,t)$)
and $\mathcal{P}\subseteq \mathcal{U}^n_k$ is a $t$-intersecting
uniform partition system and $c=n/k$ is the size of a class in each
partition, then
\begin{enumerate}
\item $|\mathcal{P}| \leq U(n-tc,k-t)$;
\item this bound is tight if and only if $\mathcal{P}$
is a trivially $t$-intersecting uniform partition system.
\end{enumerate}
\end{thm}
In Section~\ref{proof1}, we give a straightforward lemma from which
we can easily prove Theorem~\ref{two} for all cases except $c=2$ and
Theorem~\ref{one}.  Indeed, the proof of Theorem~\ref{two} for $c=2$
is the only more involved case.  Since this proof applies to all $c$,
it is presented in this generality in Section~\ref{t1}.
In the proofs of Lemmas~\ref{lemma1}--\ref{lemma4} in the following 
sections, we apply a version of the {\em kernel method} introduced by 
Hajnal and Rothschild~\cite{MR0325398}.

\subsection{Erd\H{o}s-Ko-Rado Theorem for Partitions for $c \neq 2$}\label{proof1}

A {\em blocking set} $\mathcal{B} \subset {[n] \choose c}$ for a
uniform partition system $\mathcal{P}\subseteq \mathcal{U}^n_k$ is a
collection of $c$-sets, where $c=n/k$, such that $|\mathcal{B}\cap
P|\geq 1$, for all $P \in \mathcal{P}$. For an intersecting partition
system the set of classes in any partition in the system forms a blocking
set. In particular, if $\mathcal{P}\subseteq \mathcal{U}^n_k$, then
$\mathcal{P}$ has a blocking set with $k$ sets.

Let $\mathcal{P}\subseteq \mathcal{U}^n_k$, $c=n/k$ 
and let $A$ be a $c$-set of $[1,n]$; define $\mathcal{P}_A=\{ P \in \mathcal{P}: A \in P\}$.

\begin{lemma}\label{lemma1}
Let $n\geq k\geq t \geq 1$ be positive integers, and let $\mathcal{P}\subseteq \mathcal{U}^n_k$ be 
a $t$-intersecting partition system. Let $c=n/k$ be the size of
a class in each partition.
Assume that there does not exist a $c$-set that occurs as a class in every
partition in $\mathcal{P}$. Then,
\[
|\mathcal{P}| \leq k{k-2 \choose t} U(n-(t+1)c, k-(t+1) ).
\]
\end{lemma}

\begin{proof}
Let $A$ be a class from a partition in $\mathcal{P}.$ Since no single
class occurs in every partition in $\mathcal{P}$, there is a partition
$Q \in \mathcal{P}$ that does not contain $A$.  Every partition in
$\mathcal{P}_A$ must $t$-intersect $Q$.  There are at most $k-2$
classes in $Q$ that do not contain an element in $A$.  Each partition
in $\mathcal{P}_A$ must contain at least $t$ of these $k-2$ classes.
Thus, for any class $A$, $|\mathcal{P}_A| \leq {k-2 \choose
t}U(n-(t+1)c, k-(t+1) )$.

Let $R \in \mathcal{P}$. Then, $R$ is a blocking set of $\mathcal{P}$,
and $\mathcal{P}=\cup_{A\in R} \mathcal{P}_A$. Thus, since $|R|=k$, we 
get
\begin{eqnarray*}
|\mathcal{P}| \leq k{k-2 \choose t} U(n-(t+1)c, k-(t+1) ).
\end{eqnarray*}
\end{proof}

{\bf Proof of Theorem~\ref{two} for $c \neq 2$.} 
Let $n \geq k \geq 1$ and $c=n/k \neq 2$. Let $\mathcal{P}\subseteq
\mathcal{U}^n_k$ be a maximal 1-intersecting uniform partition system that is
not trivially 1-intersecting.  The theorem is clearly true when $c=1$
and when $k = 1$.  By Lemma~\ref{lemma1}
\begin{eqnarray*}
|\mathcal{P}| \leq k(k-2) U(n-2c, k-2).
\end{eqnarray*}

For $k\geq 2$ and $c\geq 3$,
\begin{eqnarray*} \nonumber
{kc-c \choose c} 
&\geq& {3k-3 \choose 3} 
> k(k-1)(k-2).
\end{eqnarray*}
Thus, we have $|\mathcal{P}| < \frac{1}{k-1}{n-c \choose c}U(n-2c,
k-2) = U(n-c,k-1)$ and any 1-intersecting uniform partition system
that is not trivially 1-intersecting has cardinality strictly less
than $U(n-c,k-1)$.

Finally, from Equation~(\ref{eq:sizeoftrivial}), a trivially
1-intersecting uniform partition system has cardinality $U(n-c,k-1)$
and the theorem holds for $c \neq 2$.
\qed

\vspace{.5cm}

{\bf Proof of Theorem~\ref{one}.}  Let $n \geq k \geq t\geq 1$. Let
$\mathcal{P}\subseteq \mathcal{U}^n_k$ be a maximal $t$-intersecting
uniform partition system that is not trivially $t$-intersecting.  Let
$c = n/k$ be the size of a class in each partition.  {}From
Equation~(\ref{eq:sizeoftrivial}), a trivially $t$-intersecting uniform
partition system has cardinality $U(n-tc,k-t)$, thus, it is enough to
show that for $n$ sufficiently large $|\mathcal{P}| < U(n-tc,k-t)$.

For $c=1$, there is only one partition and $|\mathcal{P}|=1$, so we
assume $c\geq 2$.  If $t = k$ or $k-1$, then two partitions are
$t$-intersecting if and only if they are identical. So we may also
assume that $t \leq k-2$.

Let $\mathcal{A}$ be the set of all $c$-sets that occur in every
partition in $\mathcal{P}$.  Let $s=|\mathcal{A}|$ and since
$\mathcal{P}$ is not trivially $t$-intersecting, we have $0 \leq s < t$.  
Consider the system
$\mathcal{P}' = \{ P\backslash \mathcal{A} : P \in \mathcal{P} \}$.  
The system $\mathcal{P}'$ is a $t'$-intersecting partition
system contained in $\mathcal{U}^{n'}_{k'}$, with $k'=k-s$, $t'=t-s$ and 
$n'=n-sc = c(k-s)$, and 
$|\mathcal{P}| = |\mathcal{P}'|$.  Furthermore, there exists no 
$c$-set in every partition in $\mathcal{P}'$, so by Lemma~\ref{lemma1},
\begin{eqnarray*}
|\mathcal{P}'| &\leq& k'{k'-2 \choose t'}U(n'-(t'+1)c, k'-(t'+1) ) \\
               &=& (k-s){k-s-2 \choose t-s}  U(n-(t+1)c, k-(t+1) ) \\
               &\leq& k{k-2 \choose t}  U(n-(t+1)c, k-(t+1) ). 
\end{eqnarray*}

Fix the value of $t$ and $k$. Then for some $n$ sufficiently large,
relative to $t$ and $k$ we have
\begin{eqnarray}\label{ineq:bound}
k{k-2 \choose t} &<& \frac{1}{k-t} {n-tc \choose c}.
\end{eqnarray}
Thus, there exists a function $n_0(k,t)$ such that for $n \geq
n_0(k,t)$ Inequality~(\ref{ineq:bound}) holds.

Since 
\[
\frac{1}{k-t} {n-tc \choose c} = {kc-tc-1 \choose c-1},
\]
if $c$ and $t$ are fixed with $1 \leq t < c-1$, for $k$ sufficiently
large Inequality~(\ref{ineq:bound}) holds.  Thus, if $1 \leq t < c-1$
then there exists a function $n_1(c,t)$ such that for $n \geq
n_1(c,t)$ Inequality~(\ref{ineq:bound}) holds.

Therefore, 
\begin{eqnarray*}
|\mathcal{P}'| &<& \frac{1}{k-t} {n-tc \choose c} U(n-(t+1)c, k-(t+1))\\
 &=& U(n-tc,k-t).
\end{eqnarray*}
\qed

\subsection{General Erd\H{o}s-Ko-Rado Theorem for Partitions}\label{t1}
It only remains to prove the case $c=2$ of Theorem~\ref{two}, but we 
give the proof for general $c$.	

Let $\mathcal{P}\subseteq \mathcal{U}^n_k$, with $n=ck$, and let $A$ be a
$c$-set of an $n$-set.  We denote $\mathcal{P}_A=\{ P \in \mathcal{P}:
A \in P\}$.  Further, if for some positive integer $\ell$, $\mathcal{A} =
\{A_1,A_2,\dots,A_\ell \}$ is a disjoint collection of $c$-sets of an $n$-set, then
$\mathcal{P}_{\mathcal{A}}=\{ P \in \mathcal{P}: A_i \in P \mbox{ for
all } i=1,\dots, \ell \}$.

\begin{lemma}\label{lemma2}
Let $n \geq k \geq 1$, and let $\mathcal{P}\subseteq \mathcal{U}^n_k$ be a
$1$-intersecting partition system that is not trivially $1$-intersecting. 
Let $c=n/k$ be the size of a class in each partition. Let $\ell$ be the size
of the smallest blocking set for $\mathcal{P}$.  Then, for any $1\leq i
< \ell$, any given set of $i$ classes of a partition can occur together in at most
\[
(k-i)(k-(i+1))\cdots(k-(\ell-1))U(c(k-\ell), k-\ell)
\]
partitions in $\mathcal{P}$.
\end{lemma}

\begin{proof}
First, since the set of classes of any partition from
$\mathcal{P}$ is a blocking set, we have $\ell \leq k$.  

Use induction on $\ell-i$. If $i = \ell-1$, consider a set of
$(\ell-1)$ disjoint $c$-sets $\mathcal{A} =
\{A_1,A_2,\dots,A_{\ell-1} \}$.  
Since $|\mathcal{A}| <\ell$, the set $\mathcal{A}$ is not a blocking
set for $\mathcal{P}$.
So, there exists a partition $Q \in \mathcal{P}$ that does not
contain any of the $A_j \in \mathcal{A}$.  
Since the $c$-sets $A_i$ are disjoint, $|\cup_{i=1}^{\ell-1} A_i| =
c(\ell-1)$. This means there are at least $\ell-1$ classes in $Q$ that
contain some element of $A_1 \cup A_2 \cup \cdots \cup A_{\ell-1}$. 
So, there are at most $k-(\ell-1)$ classes in $Q$ that
could appear in a partition in $\mathcal{P}_{\mathcal{A}}$.  Each
partition in $\mathcal{P}_{\mathcal{A}}$ must contain at least one of
these $k-(\ell-1)$ classes.  Thus,
\begin{eqnarray*}
|\mathcal{P}_{\mathcal{A}}| & \leq& ( k-(\ell-1) ) U(n-(\ell-1)c-c, k-(\ell-1)-1)\\ 
         &=& (k-\ell+1)U(n-\ell c, k-\ell).
\end{eqnarray*}
This completes the case $\ell=i+1$.

Now, for $\ell \geq i+1$, we assume that any set of $i$ disjoint $c$-sets can
occur together in at most
\[
 (k-i)(k-(i+1))\cdots(k-\ell+1)U(n-\ell c, k-\ell)
\]
partitions in $\mathcal{P}$.  Consider any set of $(i-1)$ disjoint
$c$-sets $\mathcal{A} = \{A_1,A_2,\dots,A_{i-1} \}$.  Since $i-1 <
\ell$, there exists a partition $Q \in \mathcal{P}$ that does not
contain any of the $A_j \in \mathcal{A}$.  There are at most $k-(i-1)$
classes in $Q$ that could appear in a partition in
$\mathcal{P}_{\mathcal{A}}$.  By the induction hypothesis, each of
these $k-(i-1)$ classes can occur together with all $A_j \in
\mathcal{A}$ in at most $(k-i)(k-i+1)\cdots(k-\ell+1)U(n-\ell c,
k-\ell)$ partitions. Thus,
\[
|\mathcal{P}_{ \mathcal{A} }| \leq  (k-(i-1)) (k-i)\cdots(k-\ell+1)U(n-\ell c,k-\ell).
\]
\end{proof}

A slightly stronger version of the previous lemma is needed for $i=1$.
\begin{lemma}\label{lemma3}
Let $n \geq k \geq 1$, and let $\mathcal{P}\subseteq \mathcal{U}^n_k$
be a $1$-intersecting system that is not trivially
$1$-intersecting. Let $c=n/k$ be the size of a class in each
partition. If the size of a smallest blocking set for $\mathcal{P}$ is
$\ell < k$, then any class can occur in at most
\[
(k-2) \left( \prod_{i=2}^{\ell-1}(k-i) \right) U(n-\ell c,k-\ell)
\]
partitions in $\mathcal{P}$.
\end{lemma}

\begin{proof}
Let $A$ be a class in a partition in $\mathcal{P}$.  Since the system
is not trivially 1-intersecting, there exists a partition $Q \in
\mathcal{P}$ which does not contain $A$.  Any partition in
$\mathcal{P}_A$ must intersect $Q$.  The elements from $A$ must be in
at least two separate classes in $Q$, thus there are at most $k-2$
classes in $Q$ which could be in this intersection. 

If $\ell=2$ then each of these
$k-2$ classes can occur in at most $(k-2)U(n-2c, k-2)$
partitions in $\mathcal{P}_A$.
If $\ell >2$ from Lemma~\ref{lemma2}, for the case $i=2$, we have that
any pair of classes can occur in at most
$(k-2)\cdots(k-(\ell-1))U(n-\ell c, k-\ell)$ partitions.
So each of these
$k-2$ classes can occur in at most $(k-2)\cdots(k-(\ell-1))U(n-\ell c, k-\ell)$
partitions in $\mathcal{P}_A$.  Thus, for all $\ell$
\[
|\mathcal{P}_A| \leq (k-2)  \left( \prod_{i=2}^{\ell-1}(k-i) \right) U(n-\ell c,k-\ell).
\]
\end{proof}

\begin{lemma}\label{lemma4}
Let $n \geq k \geq 4$, and let $\mathcal{P}\subseteq \mathcal{U}^n_k$ be a 
$1$-intersecting partition system that is not trivially $t$-intersecting. 
Let $\ell$ be the size of the
smallest blocking set for $\mathcal{P}$. If $\ell = k-1$ or $k$, then any
set of $i<k-2$ classes of a partition 
can occur in at most $(k-i)(k-i-1)(k-i-2) \cdots 3$ partitions in
$\mathcal{P}$.
\end{lemma}

\begin{proof}
Since $\ell \geq k-1$, for any set $\mathcal{A}$
of $k-2$ pairwise disjoint classes, there exists a partition
$Q \in \mathcal{P}$ that does not contain any of the classes 
in $\mathcal{A}$.
Any partition in $\mathcal{P}_{\mathcal{A}}$ must intersect
$Q$ and there are at most 2 classes in $Q$ which could be in this
intersection. Let $Q_0$ and $Q_1$ be these two classes.
For $P \in \mathcal{P}_{\mathcal{A}}$, if $Q_0 \in P$ then the first
$k-2$ classes of $P$ form $\mathcal{A}$, and another class is $Q_0$. Since
$Q_1$ is disjoint from all the sets in $\mathcal{A}$ and from $Q_0$,
the final class of $P$ must be $Q_1$. Thus there is only one partition in
$\mathcal{P}_{\mathcal{A}}$, and any set of $k-2$ classes can occur
in at most one partition in $\mathcal{P}$.

We will use induction on $k-i$. If $i = k-3$, consider a set $\mathcal{A}$ of
$k-3$ classes. Since $|\mathcal{A}|<\ell-1$, there is a partition $Q\in
\mathcal{P}$ that does not contain any of the classes in $\mathcal{A}$.
There are at most $k - (k-3) = 3$ classes in $Q$ that could appear
in a partition in $\mathcal{P}_\mathcal{A}$. Since no set of $(k-2)$ $c$-sets can
occur in more than one partition, $|\mathcal{P}_\mathcal{A}| \leq 3$.

Now, if $i \leq k-3$ we assume that any set of $i$ classes of a partition can
occur together in at most $(k-i)(k-i-1)(k-i-2) \cdots 3$ partitions in $\mathcal{P}$.
Consider any set $\mathcal{A}$ of $i-1$ classes. There exists a
partition $Q \in \mathcal{P}$ which does not contain any of the classes
in $\mathcal{A}$. There are at most $k-(i-1)$ classes in $Q$ which could
occur in a partition in $\mathcal{P}_\mathcal{A}$.
Thus, $|\mathcal{P}_\mathcal{A}| \leq  (k-(i-1))(k-i)(k-i-1)(k-i-2) \cdots 3$.
\end{proof}

Before giving the proof of Theorem~\ref{two}, for general $c$, we need two technical lemmas.
\begin{lemma}\label{lemma5}
For $k-2 \geq \ell \geq 2$,
\[
\ell(k-2) \prod_{i=2}^{\ell-1} (k-i) <  \prod_{i=1}^{\ell-1} (2(k-i)-1).
\]
\end{lemma}

\begin{proof}
We prove this by induction on $\ell$.
If $\ell=2$,
\[
2(k-2) < 2k-3.
\]

For $\ell \geq 2$, assume
\begin{eqnarray}
\ell(k-2) \prod_{i=2}^{\ell-1} (k-i) < \prod_{i=1}^{\ell-1}(2k-2i-1), \label{eqn:2}
\end{eqnarray}
then
\[
\begin{array}{ll}
(\ell+1)(k-2) \prod_{i=2}^{\ell} (k-i) =&\\
\ \ \ \ = \frac{(\ell+1)}{\ell}(k-\ell) \ell (k-2) \prod_{i=2}^{\ell-1} (k-i)  &  \\
\ \ \ \ < \frac{(\ell+1)}{\ell}(k-\ell)  \prod_{i=1}^{\ell-1}(2k-2i-1)    &\mbox{(by Inequality~(\ref{eqn:2}))} \\ 
\ \ \ \ \leq (2k-2\ell-1)  \prod_{i=1}^{\ell-1}(2k-2i-1)       \ \ \     &\mbox{(since $2 \leq \ell \leq k-2$)} \\
\ \ \ \ = \prod_{i=1}^{\ell}(2k-2i-1). &\\
\end{array}
\]
\end{proof}

\begin{lemma}\label{lemma6}
For integers $k>j \geq 2$, $c\geq 1$ and $n=ck$,
\begin{equation}
U(c(k-1), k-1) = \left( \prod_{i=1}^{j-1}{ ck-ic-1 \choose c-1} \right) U(c(k-j),k-j). \label{formula1}
\end{equation}
\end{lemma}
\begin{proof}
We prove this by induction on $j$.

If $j=2$, rewriting the size of a trivially $1$-intersecting system, we get
\begin{eqnarray*}
U(c(k-1), k-1) &=& \frac{1}{k-1}{ ck-c \choose c} U(c(k-2), k-2)  \\
               &=& { ck-c-1 \choose c-1} U(c(k-2), k-2).  \\
\end{eqnarray*}

For $k-1> j \geq 2$, assume Equation~(\ref{formula1}) holds,  
then
\begin{eqnarray*}
&&U(c(k-1), k-1)\\
&& = \left( \prod_{i=1}^{j-1}{ ck-ic-1 \choose c-1} \right) 
                        \frac{1}{k-j}{ ck-jc \choose c}   U(c(k-j)-c,k-j-1)\\
&& = \left( \prod_{i=1}^{j-1}{ ck-ic-1 \choose c-1} \right) 
                        { ck-jc-1 \choose c-1}   U(c(k-j)-c,k-j-1)\\
&& = \left( \prod_{i=1}^{j}{ ck-ic-1 \choose c-1} \right) U(c( k-(j+1) ), k-(j+1) ).
\end{eqnarray*}
\end{proof}

{\bf Proof of Theorem~\ref{two} for all $c$.}  Let
$\mathcal{P}\subseteq \mathcal{U}^n_k$ be a maximal $1$-intersecting
partition system that is not trivially $1$-intersecting.  It is enough
to show that $|\mathcal{P}| < U(n-c, k-1)$.  If $k=2$ then
every maximal $1$-intersecting partition system is trivially
$1$-intersecting. 

Consider the case when $k=3$. Assume that $P,Q,R \in \mathcal{U}^n_3$ with
$P=\{P_1,P_2,P_3\}$, $Q=\{Q_1,Q_2,Q_3\}$ and $R=\{R_1,R_2,R_3\}$.
Further assume $P,Q$ and $R$ are distinct intersecting partitions that
do not have a common class.  Since $P$, and $Q$ are distinct, $|P\cap
Q|=1$, so we may assume $P_1 =Q_1$. Since $R$ intersects $P$ and $Q$
but the partitions do not contain a common class, we may assume $R_1 =
P_2$ and $R_2 = Q_2$. Thus $P_2 \cap Q_2 = \emptyset$, but this means $P_2 = Q_3$
and $Q_2 = P_3$ contradicting that $P$ and $Q$ are distinct.  Thus $R$
can not intersect $P$ and $Q$ and for $k=3$ it is not possible to have
a non-trivial intersecting partition system.

We can assume that $k \geq 4$. For the same reason, we know $c \geq 2$.
Let $\ell$ be the size of a smallest blocking set for $\mathcal{P}$.
Since $\mathcal{P}$ is not trivially $1$-intersecting, we know that
$\ell>1$.

{\bf Case 1.} $2 \leq \ell \leq k-2$.

There exists a blocking set $\mathcal{B}$ for $\mathcal{P}$ with $|\mathcal{B}|=\ell$,
and from Lemma~\ref{lemma3} each class in $\mathcal{B}$ can be in at most
$(k-2)  \left( \prod_{i=2}^{\ell-1}(k-i) \right) U(c(k-\ell), k-\ell)$ partitions in $\mathcal{P}$.
Thus,
\begin{equation}
|\mathcal{P}| \leq \ell (k-2) \left( \prod_{i=2}^{\ell-1} (k-i) \right)  U(c(k-\ell), k-\ell). \label{e1}
\end{equation}

From Lemma~\ref{lemma5}, and the fact that $2(k-i)-2 < {c(k-i)-1 \choose c-1}$
for all $c \geq 2$, we get
\begin{equation}
\ell(k-2) \prod_{i=2}^{\ell-1} (k-i) <  \prod_{i=1}^{\ell-1} (2(k-i)-1)
                              \leq \prod_{i=1}^{\ell-1}{ c(k-i)-1 \choose c-1}.\label{e3}
\end{equation}

Therefore, 
\[
\begin{array}{lll}
|\mathcal{P}| &\leq \ell (k-2) \left( \prod_{i=2}^{\ell-1} (k-i) \right)  U(c(k-\ell), k-\ell) \
                           &\mbox{(by Inequality~(\ref{e1}) )} \\
              &<  \left( \prod_{i=1}^{\ell-1}{ c(k-i)-1 \choose c-1} \right)  U(c(k-\ell), k-\ell)  
                           &\mbox{(by Inequality~(\ref{e3}) )} \\           
              &= U(c(k-1), k-1)     &\mbox{(by Lemma~\ref{lemma6} with $j=\ell$)} \\
              &= U(n-c, k-1).     & \\
\end{array}
\]

{\bf Case 2.} $k \geq \ell \geq k-1$. 

By Lemma~\ref{lemma4}, any single class can occur
in at most $(k-1)(k-2) \cdots 3$ partitions in $\mathcal{P}$. 
Since there exists a blocking set of cardinality $k$,
\begin{eqnarray}
|\mathcal{P}| \leq k(k-1)(k-2)(k-3) \cdots 3 = \prod_{i=1}^{k-2} (k - i+1). \label{e4}
\end{eqnarray}

We have
\begin{eqnarray}
       &&k-i+1  < 2k-2i-1,\ \  \mbox{ for all $i \leq k-3$}     \label{e5}\\
\mbox{and}\ \      && 2(k-i)-1 \leq  { c(k-i)-1 \choose c-1}\ \   \mbox{ for all $c \geq 2$}. \label{e6}
\end{eqnarray}

Therefore,
\begin{eqnarray*}
|\mathcal{P}| &\leq& \prod_{i=1}^{k-2} (k - i  + 1) \ \ \ \ \ \ \ \ \mbox{(by Inequality~(\ref{e4}))}  \\
              &=& 3 \prod_{i=1}^{k-3} (k - i  + 1)     \\
              &<& 3\prod_{i=1}^{k-3}(2k-2i-1)   \ \ \       \mbox{(by Inequality~(\ref{e5}) 
                                                                   and $k\geq 4$)} \\
   &=& \prod_{i=1}^{k-2}(2k-2i-1)           \\
   &\leq& \prod_{i=1}^{k-2}{ c(k-i)-1 \choose c-1}\ \ \   \mbox{(by Inequality~(\ref{e6}) and $c\geq 2$)}  \\
   &=& U(c(k-1),k-1) \ \ \  \mbox{(by Lemma~\ref{lemma6} with $j=k-1$ and $U(c,1)=1$)} \\
   &=& U(n-c,k-1).
\end{eqnarray*}
\qed

%%%%%%%%%%%%%%%%%%%%%%%%%%%%%%%%%%%%%%%%%%%%%%%%%%
%% start of new stuff

\section{Intersecting Packing Systems}\label{matchings}

% change namemathcing?

These theorems can be seen as results on maximal families of 1-regular
$c$-uniform hypergraphs on $n$ vertices that intersect in at least $t$
edges. Alternatively, these hypergraphs can be thought of as perfect
matchings on $K_n^{(c)}$, the complete $c$-uniform hypergraph on $n$
vertices.  Thus, we can generalize our results for the case when $c$
does not divide $n$, by considering maximal matchings in place of
perfect ones.

Define an {\em $(n,c)$-packing}\index{$(n,c)$-packing} to be a set of
disjoint $c$-sets of an $n$-set.  Let $\mathcal{PC}_{n,c}$ denote the set of all
maximum $(n,c)$-packings, that is, all $(n,c)$-packings with $\lfloor
\frac{n}{c} \rfloor$ $c$-sets.  Set $P(n,c)=|\mathcal{PC}_{n,c}|$, then for $k
= \lfloor \frac{n}{c} \rfloor$, 
\[
P(n,c) = \frac{1}{k!} {n \choose c} {n-c \choose c} \cdots {n-(k-1)c \choose c}.
\] 
An $(n,c)$-packing system is a collection of $(n,c)$-packings.

An $(n,c)$-packing system $\mathcal{P} \subseteq \mathcal{PC}_{n,c}$ is 
{\em $t$-intersecting} if $|P \cap Q|\geq t$ for all $P,Q \in
\mathcal{P}$. It is straightforward to define a trivially 
intersecting $t$-intersecting $(n,c)$-packing system.

Generalizations of Theorems~\ref{two} and~\ref{one} are stated next
without proof.  The proofs for these are very similar to the ones
used for the original theorems.  Indeed, the only change to the
original proofs is that Lemma~\ref{lemma1} needs to have $k-1$ in
place of $k-2$ in the upper bound on $|\mathcal{P}|$.

\begin{lemma}\label{lemma7}
Let $n,c,k$ and $t$ be positive integers such that $n \geq c$ and $t \leq k = \lfloor n/c \rfloor$.  Let
$\mathcal{P}\subseteq \mathcal{P}_{n,c}$ be a $t$-intersecting $(n,c)$-packing
system. Assume that there does not exist a $c$-set that
occurs in every $(n,c)$-packing in $\mathcal{P}$.  Then,
\[
|\mathcal{P}| \leq k{k-1\choose t} P(n-(t+1)c,c).
\]
\end{lemma}

\begin{thm}\label{three}
Let  $n,c,k$ and $t$ be positive integers such that $n \geq c$ and $k = \lfloor \frac{n}{c} \rfloor$.
Let $\mathcal{P}\subseteq \mathcal{P}_{n,c}$ be a $1$-intersecting
$(n,c)$-packing system. Then
$|\mathcal{P}| \leq P(n-c,c).$
Moreover, this bound is tight 
if and only if $\mathcal{P}$ is a trivially $1$-intersecting 
$(n,c)$-packing system.
\end{thm}

\begin{thm}\label{four}
Let $n,c,k$ and $t$ be positive integers such that $n \geq c$ and $t
\leq k = \lfloor \frac{n}{c} \rfloor$.  Let $\mathcal{P}\subseteq
\mathcal{P}_{n,c}$ be a $t$-intersecting $(n,c)$-packing system. Then there exist functions
$n_0(k,t)$ and $n_1(c,t)$ such that if
($n\geq n_0(k,t)$) or ($c\geq t+2$ and $n \geq n_1(c,t)$) then,
\begin{enumerate}
\item $|\mathcal{P}| \leq P(n-tc,c)$;
\item moreover, this bound is tight if and only if $\mathcal{P}$
is a trivially $t$-intersecting $(n,c)$-packing system.
\end{enumerate}
\end{thm}

%%  end of new stuff
%%%%%%%%%%%%%%%%%%%%%%%%%%%%%%%%%%%%%%%%%%%%%%%%%%%%%%%%%%%%%%%%%%

\section[Complete Theorem for Intersecting Partition Systems]{Towards a Complete Theorem for $t$-Intersecting Partition
Systems}\label{sec:completetheorem}

Ahlswede and Khachatrian~\cite{MR97m:05251} have extended the 
{E}rd\H{o}s-{K}o-{R}ado Theorem for set systems by determining the cardinality
and structure of all maximal $t$-intersecting set systems  
$P\subseteq{[n]\choose k}$ for all possible $n \leq (k-t+1)(t+1)$ (see Theorem~\ref{thm:ahlswede}). 
This remarkable result went beyond proving a conjecture by Frankl~\cite{MR80c:05014} that stated 
a specific list of candidates for maximal set systems. 
Next, we state conjectures for uniform $t$-intersecting partition
systems, which parallel the conjecture of Frankl and the theorem
of Ahlswede and Khachatrian, respectively. 

For $0 \leq i \leq \lfloor(k-t)/2 \rfloor$, define the partition system
\begin{eqnarray*}
\mathcal{P}_i(n,k,t)&=&\{P \in \mathcal{U}^n_k:| P \cap \{[1,c],[c+1,2c],\\ 
  && \left. \ldots, [(t+2i-1)c+1, (t+2i)c]\} |   \geq t+i \right\}. 
\end{eqnarray*}

Note that $\mathcal{P}_0(n,k,t) = \mathcal{P}(n,k,t)$.
Theorem~\ref{one} says that for $n$ sufficiently large
$\mathcal{P}_0(n,k,t)$ is the unique (up to permutations on $[1,n]$)
largest $t$-intersecting uniform partition system in
$\mathcal{U}^{n}_k$.  We conjecture that for any $n$ the unique (up to
permutations on $[1,n]$) largest $t$-intersecting uniform partition
system in $\mathcal{U}^{n}_k$ is one of $\mathcal{P}_i(n,k,t)$ for $i
\in \{0, \dots \lfloor(k-t)/2
\rfloor \}$.

\begin{conj}
Let $n,k$ and $t$ be positive integers with $n \geq k \geq t \geq 1$,
and let $\mathcal{P}\subseteq \mathcal{U}^n_k$ be a $t$-intersecting
partition system. Then
\begin{eqnarray*}
|\mathcal{P}| \leq \max_{0 \leq i \leq  \frac{k-t}{2}}|\mathcal{P}_i(n,k,t)|.
\end{eqnarray*}
\end{conj}

\begin{conj}\label{conj:completeintersection}
Let $c,t$ and $i$ be positive integers.  There exists a function
$n_0(c,t,i)$, such that, for an integer $k$ with $k \geq t$ and
$n_0(c,t,i+1) < ck < n_0(c,t,i)$, if $\mathcal{P}\subseteq
\mathcal{U}^{ck}_k$ is a $t$-intersecting partition system, then
$|\mathcal{P}| \leq |\mathcal{P}_i(ck,k,t)|$.  Moreover, this bound is
tight if and only if $\mathcal{P}$ is equal (up to permutations on
$[1,ck]$) to $\mathcal{P}_i(ck,k,t)$.
\end{conj}

One could hope to be able to use the ideas in \cite{MR97m:05251} to
prove these conjectures; however, key techniques such as {\em left
compression}, which are used in their proofs, do not seem to have an
extension to partition systems.

We conclude with an infinite sequence of parameters $(n,k,t)$ for
which $|\mathcal{P}_1(n,k,t)| > |\mathcal{P}(n,k,t)|$. This is not a
counter example to Conjecture~\ref{conj:completeintersection} since,
in this example, we require that $t=k-3$, thus $n=c(t+3)$.  In fact,
this example gives a lower bound for the function $n_0(c,t,0)$,
specifically $n_0(c,t,0) \geq c(t+3)$.

\begin{prop}\label{prop:HMpartitions}
For positive integers $c,k,t$ with $k \geq 3$ and $t=k-3$, let
$n=ck$. For a function $k_0(c)$, if $k> k_0(c)$, then
$|\mathcal{P}_1(n,k,t)| > |\mathcal{P}(n,k,t)|$.
\end{prop}
\begin{proof}
From Equation~(\ref{eq:sizeoftrivial}),
\begin{eqnarray*}
|\mathcal{P}(n,k,k-3)| &=& U(n-(k-3)c, k-(k-3)) \\
                       &=& \frac{1}{3!}{3c \choose c}{2c \choose c} \\
                       &=& {3c-1 \choose c}{2c-1 \choose c-1}.
\end{eqnarray*}

Since 
\begin{eqnarray*}
\mathcal{P}_1(n,k,t)&=&\{P \in \mathcal{U}^n_k:| P \cap \{[1,c],[c+1,2c], \\ 
  && \left.      \ldots, [(t+1)c+1, (t+2)c]\} |   \geq t+1 \right \},
\end{eqnarray*}
for $t=k-3$,
\begin{eqnarray*}
|\mathcal{P}_1(n,k,t)| &=& {t+2 \choose t+1}U(n-(t+1)c, k-(t+1))  \\
                        && \ \ - {t+2 \choose t+2} U(n-(t+2)c, k-(t+2)) +1\\
                       &=& (t+2){2c-1 \choose c-1} - t - 1.
\end{eqnarray*}

For $k$ (and $t$, since $t=k-3$) sufficiently large
\[
{3c-1 \choose c-1}{2c-1 \choose c-1} < (t+2){2c-1 \choose c-1} -t-1.
\]
\end{proof}

\section{Partially Intersecting Partitions}\label{PIP}

Erd\H{o}s and Sz\'ekely~\cite{MR2001f:05148} define another
type of intersecting partitions, which we call here {\em partially $t$-intersecting}.
\index{partially $t$-intersecting partitions}
Two partitions $P, Q \in \mathcal{P}^n_k$ are said to be partially
$t$-intersecting if there exist classes $P_i \in P$ and $Q_j \in
Q$ such that $|P_i \cap Q_j|\geq t$.  A partition system $\mathcal{P}
\subseteq \mathcal{P}^n_k$ that is pairwise partially $t$-intersecting is called a 
{\em partially $t$-intersecting partition system}\index{partially
$t$-intersecting partition system}. 
 
A partition system is called a {\em trivially partially
$t$-intersecting partition system}\index{trivially partially
$t$-intersecting partition system} if it is equal, up to permutations on
$[1,n]$, to
\begin{eqnarray*}
\mathcal{R}(n,k,t)&=&\left\{P \in \mathcal{P}^n_k :  [1,t] \subseteq A, \mathrm{\ for\ some \ } A \in P \right\}.
\end{eqnarray*}

\begin{conj}{\bf(Czabarka's Conjecture, see~\cite{MR2001f:05148})}
Let $n\leq 2k-1$ and $\mathcal{P}\subseteq \mathcal{P}^n_k$ be a partially $2$-intersecting
partition system. Then, $|\mathcal{P}| \leq S(n-1,k)$. 
\end{conj}

This bound is attained by the system
\[ \{P \in \mathcal{P}^n_k: [1,2]\subseteq P_i,\mathrm{\ for\ some\ }P_i \in P\}.\]

We pose a similar conjecture for uniform partition systems.

\begin{conj} \label{conj:partialintersect}
Let $k,c,t$ be positive integers with $t\leq c$ and $n=ck$. Let
$\mathcal{P}\subseteq \mathcal{U}^n_k$ be a partially $t$-intersecting
uniform partition system. Then, $|\mathcal{P}| \leq {n-t \choose c-t}
U(n-c,k-1)$. Moreover, this bound is tight if and only if
$\mathcal{P}$ is equal, up to permutations of $[1,n]$, to
\[
\{P\in \mathcal{U}^n_k: [1,t] \subseteq P_i, \mathrm{\ for\ some\ }P_i \in P\}.
\]
\end{conj}
Note that Theorem~\ref{two} confirms Conjecture~\ref{conj:partialintersect} for $t=c$.

For some values of $t,k,n$, a partially $t$-intersecting $k$-partition
system of an $n$-set corresponds to an independent set in the graph
$UQI(n,k)$. 

\begin{prop}\label{prop:uniform}
For positive integers $c,k,t$ with $t \geq c-k+2$, a partially $t$-intersecting $c$-uniform
$k$-partition system is an independent set in the graph $UQI(ck,k)$.
\end{prop}
\begin{proof}
Assume $\mathcal{P}$ is a partially $t$-intersecting $c$-uniform
$k$-partition system in $\mathcal{P}^{ck}_k$ with  $t \geq c-k+2$.
Note that $\mathcal{P} \subseteq V(UQI(ck,k))$.

To show $\mathcal{P}$ is an independent set in $UQI(ck,k)$ we will
show that for any $P,Q \in V(UQI(ck,k))$, $P$ and $Q$ are not
qualitatively independent.  Since $P$ and $Q$ are partially
$t$-intersecting there are classes $P_1 \in P$ and $Q_1 \in Q$ such
that $|P_1 \cap Q_1| \geq t$. Since $P$ is a $c$-uniform partition and
$t \geq c-k+2$ there are at most $c - (c-k+2) = k-2$ elements in $P_1$
which are not in $Q_1$. Thus the class $P_1$ can intersect at most
$k-2$ classes in $Q$, other than $Q_1$. Since $Q$ is a $k$-partition,
$P_1$ can not intersect every class in $Q$. Hence, $P$ and $Q$ are not qualitatively independent.
\end{proof}

\begin{prop}\label{prop:notuniform}
Let $c,k,r$ be positive integers with $n=ck+r$ and $0 \leq r < k$.
For $t\geq c-k+3$, a partially $t$-intersecting almost-uniform $k$-partition
system on an $n$-set is an independent set in $AUQI(n,k)$.
\end{prop}
\begin{proof}
Assume $\mathcal{P}$ is a partially $t$-intersecting almost-uniform
$k$-partition system in $\mathcal{P}^{n}_k$ with  $t \geq c-k+3$ and $n=ck+r$.

To show $\mathcal{P}$ is an independent set in $AUQI(n,k)$ we will
show that for any $P,Q \in V(AUQI(n,k))$, $P$ and $Q$ are not
qualitatively independent.  Since $P$ and $Q$ are partially
$t$-intersecting there are classes $P_1 \in P$ and $Q_1 \in Q$ such
that $|P_1 \cap Q_1| \geq t$. Since $P$ is an almost uniform partition $|P_1| \leq c+1$.
Since, $t\geq c-k+3$ there are at most $(c+1) - (c-k+3) = k-2$ elements in $P_1$
which are not in $Q_1$. Thus the class $P_1$ can intersect at most
$k-2$ classes in $Q$, other than $Q_1$. Since $Q$ is a $k$-partition,
$P_1$ can not intersect every class in $Q$. Hence, $P$ and $Q$ are not qualitatively independent.
\end{proof}

We can prove the bound in Conjecture~\ref{conj:partialintersect} for
specific partition systems; to do this we need the following theorem.

\begin{thm}\label{PIPS:thm}
For positive integers $n,k,c$ and $t$ with $n=ck$ and $t \leq c$, if
there exists a resolvable $t$-$(n,c,1)$ design, then the size of the
largest partially $t$-intersecting uniform $k$-partition system from
$\mathcal{U}^n_k$ is bounded above by
\[
\frac{1}{(k-1)!} {n-t \choose c-t}{n-c \choose c} \cdots {c \choose c}
= {n-t \choose c-t} U(n-c, k-1).
\]
Further, this bound is attained by a trivially partially $t$-intersecting
$k$-partition system.
\end{thm}

\begin{proof}
Define a graph $PIP(n,k)$ whose vertex set is the set of all uniform
$k$-partitions of an $n$-set and two vertices are adjacent if and only
if they are partially $t$-intersecting.  The cardinality of the vertex set is
\[
U(n,k) = \frac{1}{k!} {n \choose c}{n-c \choose c} \cdots {c \choose c}.
\]
Further, the graph $PIP(n,k)$ is vertex transitive and a clique in $PIP(n,k)$
corresponds to a partially $t$-intersecting partition system.

An independent set in $PIP(n,k)$ is a set $\mathcal{Q}$ of partitions
such that no two distinct partitions in $\mathcal{Q}$ are partially
$t$-intersecting.  Consider a resolvable $t$-$(n,c,1)$ design. Each
resolution class in the resolvable $t$-$(n,c,1)$ design corresponds to
a $c$-uniform $k$-partition of an $n$-set. Further, for each $t$-set
there is exactly one resolution class in the resolvable
$t$-$(n,c,1)$ design which contains the $t$-set in one of its blocks.  Thus, a
resolvable $t$-$(n,c,1)$ design corresponds to an independent set in
$PIP(n,k)$.

The number of blocks in a resolvable $t$-$(n,c,1)$ design is
$\frac{n!(c-t)!}{c!(n-t)!}$, and the number of resolution classes is
$\frac{(n-1)!(c-t)!}{(c-1)!(n-t)!}$. The maximum independent set in
$PIP(n,k)$ has size at least $\frac{(n-1)!(c-t)!}{(c-1)!(n-t)!}$.
\remove{$\frac{n-1}{c-1}$.}

Since $PIP(n,k)$ is vertex transitive, with Inequality~(\ref{eq:vt2}) from Section~\ref{knesergraphs} we
have
\begin{eqnarray*}
\omega(PIP(n,k)) &\leq& \frac{ |V(PIP(n,k))|}{ \alpha(PIP(n,k))} \\ 
&\leq& \frac{ \frac{1}{k!} {n \choose c}{n-c \choose c} \cdots {c \choose c} }
     { \frac{(n-1)!(c-t)!}{(c-1)!(n-t)!}  } \\ 
&=& \frac{1}{(k-1)!} {n-t \choose c-t}    {n-c \choose c} \cdots {c \choose c}.
\end{eqnarray*}
\orphan
\end{proof}

\begin{cor}\label{PIPS:app1}
For $n = 3k$ and $k$ odd, the largest partially 2-intersecting uniform
$k$-partition system from
$\mathcal{U}^n_k$ has cardinality
\[
(3k-2)U(3k-3,k).
%\frac{1}{(k-1)!}{3k-3 \choose 3}{3k-6 \choose 3} \cdots {3 \choose 3}.
\]
This bound is attained by a trivially partially 2-intersecting $k$-partition system.
\end{cor}

\begin{proof}
If $k$ is odd, then $3k \equiv 3 \pmod{6}$ and a resolvable
$(3k,3,1)$-BIBD (equivalently, a resolvable 2-$(3k,3,1)$ design)
exists~\cite{colbourn:96}. The result follows from
Theorem~\ref{PIPS:thm}
\end{proof}

\begin{cor}\label{PIPS:app2}
For $n \equiv 4 \pmod{12}$ and $k= n/4$, the largest partially 2-intersecting uniform
$k$-partition system from
$\mathcal{U}^n_k$ has cardinality
\[
{4k-2 \choose 2} U(n-4,k).
%\frac{1}{(k-1)!}{n-4 \choose 4} \cdots {4 \choose 4}.
\]
This bound is attained by a trivially partially 2-intersecting $k$-partition system.
\end{cor}

\begin{proof}
If $n \equiv 4 \pmod{12}$, then a resolvable $(n,4,1)$-BIBD
(equivalently, a 2-$(n,4,1)$ design) exists~\cite{colbourn:96}. The
result follows from Theorem~\ref{PIPS:thm}
\end{proof}

\begin{cor}\label{PIPS:app3}
For $k$ a prime power and $n = k^2$, the largest partially
2-intersecting uniform $k$-partition system from $\mathcal{U}^n_k$ has
cardinality
\[
{k^2 -2 \choose k-2} U(k^2-k, k-1).
\]
This bound is attained by a trivially partially 2-intersecting $k$-partition system.
\end{cor}

\begin{proof}
For all $k$ a prime power there exists a resolvable $(k^2,k,1)$-BIBD (Section 8.3 \cite{street:87}).
\end{proof}

We have already seen a stronger version of this result for all $k$ in
another form.  For $n=k^2$, by Proposition~\ref{prop:uniform}, a
partially 2-intersecting uniform $k$-partition system of an $n$-set is
an independent set in $QI(k^2,k)$ and from
Lemma~\ref{rationprimepower}, $\alpha( QI(k^2,k) ) = {k^2-2 \choose
k-2} U(k^2-k, k-1)$.  Thus the bound from
Corollary~\ref{PIPS:app3} holds for all values of $k$.

\begin{thm}
For all integers $k$ and $n = k^2$, the largest partially 2-intersecting uniform
$k$-partition system from $\mathcal{U}^n_k$ has cardinality
\[
{k^2 -2 \choose k-2} U(k^2-k, k-1).
\]
This bound is attained by a trivially partially 2-intersecting $k$-partition system.
\end{thm}

It is interesting to note that
Lemma~\ref{rationprimepower} was proven with the ratio bound and the
eigenvalues of $QI(k^2,k)$. It may be possible to prove
Conjecture~\ref{conj:partialintersect} using the eigenvalues of the
graph $PIP(n,k)$ and the ratio bound.

\begin{quest}
Is it possible to find the eigenvalues of $PIP(n,k)$?  Does the ratio
bound (Theorem~\ref{thm:tauomegabound}) give a good bound on
$\alpha(PIP(n,k))$?
\end{quest}

\begin{cor}\label{PIPS:app4}
For $v$ a positive integer with $v \equiv 4,8 \pmod {12}$, except
possibly for $v \in \{220, 236, 292, 364, 460, 596, 676, 724,
1076, 1100, 1252, 1316, 1820, 2236, 2308,$ $2324,$ $2380,$ $2540, 2740,
2812, 3620, 3820, 6356\}$, the largest partially 3-intersecting uniform $v/4$-partition system has
cardinality
\[
(v-3)U(v-4, v/4-1).
\]
This bound is attained by a trivially partially 3-intersecting $v/4$-partition system.
\end{cor}

\begin{proof}
A resolvable 3-$(v,4,1)$ design exists for all such $v$~\cite{colbourn:96}.
\end{proof}

There is another approach that can be used to bound partially
$t$-intersecting uniform $k$-partition systems.  Recall from
Definition~\ref{defn:packing} that a relaxation of a design is a {\em
packing}. A $t$-$(n,k,\lambda)$ packing is a set system with the
property that every $t$-set occurs in at most $\lambda$ sets.  Thus,
the argument used in Theorem~\ref{PIPS:thm} follows exactly for
$t$-$(n,k,1)$ packings except with a different (possibly weaker)
bound. We state this theorem without proof.

\begin{thm}\label{PIPS:thm2}
For positive integers $k$ and $c$ with $n=ck$, if there exists a
resolvable $t$-$(n,c,1)$ packing with $B$ blocks, then the cardinality of the
largest partially $t$-intersecting uniform $k$-partition system is
bounded by
\[
 \frac{ \frac{1}{(k-1)!} {n \choose c}{n-c \choose c} \cdots {c \choose c} } {B} \\ 
\]
\end{thm}

Less is known about resolvable $t$-$(n,c,1)$ packings, but we do have a result for $c=3$.

\begin{cor}
Let $v=3k$ where $k$ is even and $k>4$. Then the largest partially 2-intersecting uniform
$k$-partition system has cardinality no more than
\[
(3k-1) U(3k-3,k).
%\frac{1}{(k-1)!}{3k-3 \choose 3}{3k-6 \choose 3} \cdots {3 \choose 3}.
\] 
\end{cor}
\begin{proof}
From~\cite{colbourn:96}, an optimal resolvable 2-$(v,3,1)$ packing
exists for all $v \equiv 0 \pmod 6$ except $v = 6,12$. This packing
has $\lfloor \frac{3k}{3} \lfloor \frac{3k-1}{2} \rfloor\rfloor$
blocks, which is the maximum possible.
\end{proof}

\chapter{Higher Order Problems}\label{higherparts}
\thispagestyle{empty}

%is EKR in here??

In the previous chapter, we consider ways to extend Sperner's
Theorem and the Erd\H{o}s-Ko-Rado Theorem to partition systems.  In
this chapter, these results are reformatted as part of a more general
scheme of extremal results on partition systems.

The scheme we use here is based on the framework used by Ahlswede, Cai
and Zhang~\cite{MR1349259, MR1739309} for disjoint {\em
clouds}. Recall from the introduction to
Chapter~\ref{partitionsystems} that a cloud is a collection of sets
and that a $c$-uniform cloud is a collection of $c$-sets. Ahlswede, Cai and
Zhang consider bounds on the cardinality of the largest system of pairwise
disjoint clouds for which certain kinds of binary relations hold.
They consider the following four binary relations between the sets in
the clouds: comparable, incomparable, disjoint, and intersecting.
Further, the binary relations can hold for clouds in four different
ways called {\em types}.  For example, for two clouds, a binary
relation holds with type $(\forall, \forall)$ if it holds between {\em
all} sets from the first cloud and {\em all} sets from the second
cloud. The three other types that they consider are $(\exists,
\forall)$, $(\forall, \exists)$ and $(\exists, \exists)$; each of
these types are defined in the following section.
 
Ahlswede, Cai and Zhang consider all four types of all four binary
relations for both disjoint $c$-uniform clouds and disjoint clouds.
For many of these 32 different problems they give
either an exact solution (with the system that attains the maximum
cardinality) or an asymptotic solution.  These problems include
Sperner's Theorem and the Erd\H{o}s-Ko-Rado Theorem. Other results for
specific problems in this scheme are also given
in~\cite{intersectingsystems} and~\cite{MR1368831}.

We will consider the same 32 problems for extremal set-partition
systems --- that is, both uniform $k$-partitions and $k$-partitions and
all types of the four binary relations.  Since we require that the
families of sets be partitions, rather than clouds, we require more
structure than a cloud. Moreover, we are looking for systems from
$\mathcal{P}^n_k$, not systems of subsets of the $n$-set, so we do not
require that the partitions be disjoint. Also, we require that
partitions in our systems be $k$-partitions, whereas
in~\cite{MR1739309} two clouds in a system could have a different
number of sets.  Most of the problems in this chapter are new; it is
not clear yet which may have applications.

Where appropriate, we will use the notation for partition systems and
uniform partition systems given in Chapter~\ref{EST}. 

\section{Types of Relations}

In this section, we detail the different problems that are
considered in this chapter.

Let $n$ be a positive integer and $A$ and $B$ be subsets of an $n$-set.
The four binary relations that we consider are the following:

\begin{enumerate}
\item{comparable ($A \subseteq B$ or $B \subseteq A)$,}
\item{incomparable ($A \not\subseteq B$ and $B \not\subseteq A$),}
\item{disjoint ($A \cap B = \emptyset$),}
\item{intersecting ($A \cap B \neq \emptyset$).}
\end{enumerate}

Unlike previous work on clouds~\cite{intersectingsystems, MR1349259,
MR1739309, MR1368831}, we do not require that the partitions in the
system be pairwise disjoint. Hence, the definition of comparable sets
includes sets that are equal and the definition of incomparable sets
requires that sets not be equal.

For each binary relation, there are four {\em types} of problems.  The
first type is {\em type $(\forall,\forall)$}.  The partition system
$\mathcal{P}$ is of type $(\forall, \forall)$ for a binary relation if
for any distinct $P,Q \in \mathcal{P}$ {\em for all} classes $P_{i}
\in P$ and {\em for all} classes $Q_{j} \in Q$, $P_{i}$ and $Q_{j}$
satisfy the binary relation.  A partition system $\mathcal{P}$ is of
type $(\exists, \forall)$ if for any distinct $P,Q \in \mathcal{P}$
{\em there exists} a class $P_{i} \in P$ such that {\em for all}
classes $Q_{j} \in Q$ the binary relation holds for $P_i$ and
$Q_j$. This is a weaker condition than type $(\forall,\forall)$.  For
a partition system $\mathcal{P}$ to be of type $(\forall,\exists)$,
for any distinct $P,Q \in \mathcal{P}$ it is required that {\em for
all} classes $P_{i} \in P$ {\em there exists} a class $Q_{j} \in Q$
such that the relation holds for $P_i$ and $Q_j$.  This is a weaker condition than both
type $(\exists,\forall)$ and $(\forall,\forall)$.  The final type of
problem is {\em type $(\exists, \exists)$} and it is the weakest
condition. All that is required for a partition system $\mathcal{P}$
is that for every distinct $P,Q \in \mathcal{P}$ {\em there exists} a
class $P_{i} \in P$ and that {\em there exists} a class $Q_{j} \in Q$
such that the relation holds for $P_i$ and $Q_j$.

Note that each of the relations, comparable, incomparable, disjoint
and intersecting, are symmetric, but not all the types of relations are
symmetric. In particular, only the types $(\forall, \forall)$ and
$(\exists, \exists)$ are symmetric for the relations we consider
here.

Following the notation in~\cite{MR1739309}, $\mathcal{A}_{n}(X,Y,k)$
will denote a $k$-partition system of an $n$-set of type $(X,Y)$ for
the relation $A$, and $A_{n}(X,Y,k)$ will denote the cardinality of the
largest such partition system. The parameter $A$ can
take the following values: $N$ for incomparable, $C$ for comparable,
$D$ for disjoint, and $I$ for intersecting. The values of $X$ and $Y$
can each be either $\forall$ or $\exists$.  Thus, $D_{n}(\forall,
\exists, k)$ denotes the cardinality of the largest partition system
$\mathcal{P} \subseteq \mathcal{P}^n_k$ with the property that for any
two distinct partitions $P,Q \in \mathcal{P}$, for all $P_i \in P$
there exists a $Q_j \in Q$ such that $P_i \cap Q_j = \emptyset$.

Similarly, we will denote a {\em uniform} $k$-partition system of an
$n$-set of type $(X,Y)$ for the relation $A$ by
$\mathcal{A}^{\ast}_{n}(X,Y,k)$ and the size of the the largest such
partition system by ${A}^{\ast}_{n}(X,Y,k)$.

If partitions $P$ and $Q$ satisfy a binary relation of type $(\exists,
\forall)$, then $Q$ and $P$ satisfy the same binary relation of type
$(\forall, \exists)$.  To see this, let $P_0$ be the class in $P$ such
that the relation holds between $P_0$ and all classes $Q_i \in
Q$. Then for every $Q_i \in Q$ the relation holds for some class in
$P$, namely $P_0$. Thus, for all binary symmetric relations $A$,
\begin{eqnarray}\label{order}
    A_{n}(\forall, \forall, k) \leq  A_{n}(\exists, \forall, k) \leq 
    A_{n}(\forall, \exists, k) \leq  A_{n}(\exists, \exists, k),
\end{eqnarray}
and
\begin{eqnarray}\label{uniformorder}
    A^{\ast}_{n}(\forall, \forall, k) \leq  A^{\ast}_{n}(\exists, \forall, k) \leq 
    A^{\ast}_{n}(\forall, \exists, k) \leq  A^{\ast}_{n}(\exists, \exists, k) .
\end{eqnarray}
Further, for every type $(X,Y)$ and every binary relation $A$, if $k$ divides $n$, then
$A^{\ast}_{n}(X, Y,k) \leq A_{n}(X, Y, k)$.

\section{Uniform Partition Systems}

First, we consider uniform $k$-partition systems.  Throughout this
section, $n,c,k$ will denote positive integers with $n=ck$. 

The requirement that the partitions in the system be uniform
$k$-partitions makes this problem more tractable than the general
case. In particular, many of the higher order extremal problems
considered here, with uniform partitions, can be translated to finding
maximum cliques and independent sets in a vertex-transitive graph.

\subsection{Uniform Partition System Graphs}\label{graph}

In the same way that we used the vertex-transitive graph $UQI(n,k)$
to find bounds on the cardinality of a set of qualitatively
independent uniform partitions (Section~\ref{sec:boundsVT}) and the
graph $PIP(n,k)$ to bound the cardinality of a set of partially
intersecting partitions (Section~\ref{PIP}), we can build
vertex-transitive graphs that can be used to find bounds on the
cardinality of uniform partition systems that satisfy a given type of
binary relation.

We define the graph $G_{\mathcal{A}^\ast_n(X,X,k)}$\index{$G_{\mathcal{A}^\ast_n(X,X,k)}$}
to be the graph whose vertex set is $\mathcal{U}^n_k$ and two
partitions are adjacent if and only if they satisfy the binary
relation $A$ with type $(X,X)$ (where $A \in \{ N,C,D,I\}$ and the
value of $X$ can be either $\forall$ or $\exists$). Since for the two
types $(\forall,\forall)$ and $(\exists,\exists)$ each of the binary
relations are symmetric, the graph $G_{\mathcal{A}^\ast_n(X,X,k)}$ is
well-defined.

In the graph $G_{\mathcal{I}^\ast_n(\forall,\forall, k)}$, two partitions $P, Q
\in \mathcal{U}^n_k$ are adjacent if and only if for all $P_{i} \in P$ and for all
$Q_j \in Q$, $P_{i}\cap Q_{j} \neq \emptyset$. This requirement is
exactly that $P$ and $Q$ be qualitatively independent, so
$G_{\mathcal{I}^\ast_n(\forall,\forall, k)}$ is exactly $UQI(n,k)$. The size of
the maximum clique in the graph is $I^\ast_n(\forall,\forall, k)$.
Moreover, two partitions are not adjacent in
$G_{\mathcal{I}^\ast_n(\forall,\forall, k)}$ if and only if they do not satisfy
this relation. The only way two partitions $P,Q$ could not satisfy
this relation is if there exists a class $P_{i} \in P$ and a class
$Q_j \in Q$ such that $P_{i}\cap Q_{j} = \emptyset$.  This means $P$
and $Q$ are not adjacent in $G_{\mathcal{I}^\ast_n(\forall,\forall, k)}$ if and
only if they satisfy the relation disjoint of type $(\exists,
\exists)$.  Thus, $D^\ast_n(\exists,\exists, k)$ is the size of the
maximum independent set in the graph
$G_{\mathcal{I}^\ast_n(\forall,\forall, k)}$.  Moreover,
$G_{\mathcal{I}^\ast_n(\forall,\forall, k)}$ and
$G_{\mathcal{D}^\ast_n(\exists,\exists, k)}$ are graph complements of
each other.

For each of these symmetric relations there is a {\em converse
relation}\index{converse relation}.  Below, each relation is matched
with its converse:

\begin{center}
\begin{tabular}{ccc}
$N^{\ast}_{n}(\forall, \forall, k)$ & $\longleftrightarrow$ & $C^{\ast}_{n}(\exists, \exists, k)$ \\
$C^{\ast}_{n}(\forall, \forall, k)$ & $\longleftrightarrow$ & $N^{\ast}_{n}(\exists, \exists, k)$ \\
$I^{\ast}_{n}(\forall, \forall, k)$ & $\longleftrightarrow$ & $D^{\ast}_{n}(\exists, \exists, k)$ \\
$D^{\ast}_{n}(\forall, \forall, k)$ & $\longleftrightarrow$ & $I^{\ast}_{n}(\exists, \exists, k)$ \\
\end{tabular}
\end{center}

All the graphs $G_{\mathcal{A}^\ast_n(X,X,k)}$, for $A \in \{ N,C,D,I\}$ and
$X$ either $\forall$ or $\exists$, are vertex transitive.
From Corollary~\ref{cor:vtbound}, if $G$ is a vertex-transitive graph,
$\omega(G) \leq \frac{|V(G)|}{\alpha(G)}$. Thus, for uniform
$k$-partition systems we have more bounds on the size of
$A^\ast_n(X,X,k)$.

\begin{prop}\label{vtbound}
Let $n,k$ be positive integers such that $k$ divides $n$. Then \hfill
\begin{enumerate}
\item{$N^{\ast}_{n}(\forall, \forall, k) \leq \frac{U(n,k)}{C^{\ast}_{n}(\exists, \exists, k)}$,}
\item{$C^{\ast}_{n}(\forall, \forall, k) \leq \frac{U(n,k)}{N^{\ast}_{n}(\exists, \exists, k)}$,}
\item{$I^{\ast}_{n}(\forall, \forall, k) \leq \frac{U(n,k)}{D^{\ast}_{n}(\exists, \exists, k)}$,}
\item{$D^{\ast}_{n}(\forall, \forall, k) \leq \frac{U(n,k)}{I^{\ast}_{n}(\exists, \exists, k)}$.}
\end{enumerate}
\end{prop}

We could define graphs $G_{\mathcal{A}^\ast_n(X,Y,k)}$ for all types of
relations, that is, for $A \in \{ N,C,D,I\}$ and the values of $X$
and $Y$ each either $\forall$ or $\exists$. But, if $X \neq Y$, then
these graphs may be directed and we would not have bounds equivalent
to the bounds in Proposition~\ref{vtbound}.

\subsection{Summary of Results}

The following chart organizes the results for all the different
types of problems for uniform partition systems. The entries either
give the exact value of $A^\ast_n(X,Y, k)$ or a bound on
$A^\ast_n(X,Y, k)$.

\begin{center}
\begin{tabular}{|l||c|c|c|c|}
\hline
$A^\ast \big \backslash (X,Y)$ & $(\forall, \forall)$ & $(\exists, \forall)$ & $(\forall, \exists)$ & 
 $(\exists, \exists)$ \\
\hline \hline
\!\small{Incomparable}\!  & ${n-1 \choose c-1}$ & $U(n,k)$ & $U(n,k)$ & $U(n,k)$\\ \hline\hline
\small{Comparable} & 1 & 1 & 1 & $U(n\!-\!c,k\!-\!1)$ \\ \hline\hline
            &   &   & $U(n,k)$ if $c<k$ & $U(n,k)$ if $c<k$\\ \cline{4-4} \cline{5-5}
\small{Disjoint}   & 1 & 1 & $\geq U(n\!-\!c,k\!-\!1)$ & $\!\geq\!{n-c+k-2 \choose k-2}U(n\!\!-\!\!c,k\!\!-\!\!1)$ \\
            &   &   &   if $c \geq k$   &   if $c\geq k$ \\ \hline\hline
            &  1 if $c<k$  & 1 if $c<k$ &  &  \\ \cline{2-2} \cline{3-3}
\small{Intersecting}& $\leq \frac{ {n-1 \choose c -1} }{ {n-(c-k+2) \choose k-2}}$ & $\leq {n-1 \choose c-1}$  &  $U(n,k)$   &  $U(n,k)$ \\
              &  if $c \geq k$  &   if $c \geq k$  &     &   \\
\hline
\end{tabular}
\end{center}

We will look at the different types of problems paired with their
converse problem: incomparable with comparable, and intersecting with disjoint.

\subsection{Incomparable and Comparable}

For uniform partitions, we can solve all types of problems where the
binary relation is either incomparable or comparable. This is due to the
fact that any two sets of the same cardinality will be comparable if
and only if they are identical.

\begin{thm}\label{thm:unifromcomparable} For all positive integers $c,k,n$ with $n=ck$, \hfill
\begin{enumerate}
\item{$N^\ast_{n}(\exists, \forall, k) = N^\ast_{n}(\forall, \exists, k) = N^\ast_{n}(\exists, \exists, k) = U(n,k)$,}
\label{unifomrcomparable1}
\item{$C^\ast_{n}(\forall,\exists, k) = C^\ast_{n}(\exists,\forall, k)= C^\ast_{n}(\forall,\forall, k)=1$,}\label{unifomrcomparable2}
\item{$N^\ast_{n}(\forall, \forall, k) = {n-1 \choose c-1}$,}\label{unifomrcomparable3}
\item{$C^\ast_{n}(\exists,\exists, k) = U(n-c,k-1)$.} \label{unifomrcomparable4}
\end{enumerate}
\end{thm}

{\bf Proof of  case~(\ref{unifomrcomparable1})}:
Let $P,Q \in \mathcal{U}^n_k$. Since $P$ and $Q$ are uniform, any classes $P_i
\in P$ and $Q_j \in Q$ will be incomparable if and only if they are not equal.
This means any two distinct partitions are incomparable of type
$(\exists, \forall)$. Thus $N^\ast_{n}(\exists, \forall, k) = U(n,k)$
and by Inequality~(\ref{uniformorder}) the other equations hold.

{\bf Proof of case~(\ref{unifomrcomparable2})}: Two partitions from
$\mathcal{U}^n_k$ will be comparable with type $(\forall, \exists)$ if
and only if they are equal. Thus $C^{\ast}_{n}(\forall, \exists, k) = 1$
and the other equations follow from Inequality~(\ref{uniformorder}).

{\bf Proof of cases~(\ref{unifomrcomparable3})
and~(\ref{unifomrcomparable4})}: In any $\mathcal{N}^\ast_n(\forall,
\forall, k)$, each $c$-set of the $n$-set can be a class in at most one
partition in the system. Since each partition in the system has
exactly $k$ $c$-sets, this gives the bound
\[
N_{n}^{\ast}(\forall,\forall,k) \leq \frac{1}{k} {n \choose c}.
\]

Moreover, a 1-factorization of the complete hypergraph $K^{(c)}_{n}$
corresponds to a partition system that is incomparable of type $(\forall, \forall)$.
{}From Theorem~\ref{thm:1factorization}, since $n=ck$, there is
a 1-factorization of the complete hypergraph $K^{(c)}_{n}$, so
$N_{n}^{\ast}(\forall,\forall,k) = \frac{1}{k}{n \choose c}$.

{}From the last equation and Proposition~\ref{vtbound}, we have
\[
C_{n}^{\ast}(\exists,\exists,k) \leq \frac{U(n,k)}{ \frac{1}{k} {n \choose c}} = U(n-c,k-1).
\] 
It is possible to construct a $C_n^\ast(\exists, \exists, k)$ of
cardinality $U(n-c,k-1)$.  Simply fix any $c$-set of the $n$-set, then
the system of all partitions from $\mathcal{U}^n_k$ that have the
given $c$-set as a class (from Section~\ref{intersectingpartitions}
this partition system is called a trivially 1-intersecting uniform
partition system) is a comparable system of type $(\exists, \exists)$
with cardinality $U(n-c,k-1)$.  Thus, $C_{n}^{\ast}(\exists,\exists,k)
= U(n-c,k-1)$.
\qed

The result that $C_{n}^{\ast}(\exists,\exists,k) = U(n-c,k-1)$ is not
new; it is a weaker version of Theorem~\ref{two}. Theorem~\ref{two}
further states that the only system that meets this bound is a trivially
1-intersecting uniform partition system. It is interesting that it is
so easy to get the bound in Theorem~\ref{two} using the graph
$G_{\mathcal{C}_{n}^{\ast}(\exists,\exists,k)}$.

Godsil and Newman~\cite{Godsil:Newman} (see Section~\ref{assscheme93})
use graph properties of $QI(9,3)$ to show that every maximum
independent set in $QI(9,3)$ is a trivially partially 2-intersecting
uniform partition system. Would it be possible to completely prove
Theorem~\ref{two} using properties of the graph
$G_{\mathcal{C}_{n}^{\ast}(\exists,\exists,k)}$? In particular, could
we prove that every maximum clique in
$G_{\mathcal{C}_{n}^{\ast}(\exists,\exists,k)}$ (or equivalently,
every maximum independent set in
$G_{\mathcal{N}_{n}^{\ast}(\forall,\forall,k)}$) is a trivially
1-intersecting uniform partition system using graph properties of
either $G_{\mathcal{C}_{n}^{\ast}(\exists,\exists,k)}$ or
$G_{\mathcal{N}_{n}^{\ast}(\forall,\forall,k)}$?

\remove{It is interesting to note for each of the cases in
Theorem~\ref{thm:unifromcomparable}, the bounds in
Proposition~\ref{vtbound} hold with equality.
\begin{quest}
Is there some property of the binary relations
comparable and incomparable for partition systems that makes the bound
in Proposition~\ref{vtbound} tight?
\end{quest}
}

\subsection{Disjoint and Intersecting} \label{sec:DIhigher}

The binary relations disjoint and intersecting are more complicated
than comparable and incomparable. We start by stating the trivial
cases without proof.

\begin{thm}\label{thm:uniformdisjoint} For positive integers $c,k,n$ with $n = ck$, \hfill
\begin{enumerate}
\item{$D^\ast_{n}(\exists, \forall, k) = 1$;}
\item{if $n<k^2$, then $D_n^\ast(\forall,\exists,k) = U(n,k)$;}
\item{$I^\ast_{n}(\forall,\exists, k)=U(n,k)$;}
\item{if $n<k^2$, then $I_n^\ast(\exists, \forall,k) =1$;}
\item{$D^\ast_{n}(\forall, \forall, k) = 1$ and $I^\ast_{n}(\exists,\exists, k) = U(n,k)$;}
\item{if $n<k^2$, then $D_n^\ast(\exists,\exists,k) = U(n,k)$ and $I_n^\ast(\forall, \forall,k) =1$.}
\end{enumerate}
\end{thm}

\subsubsection{Bounds for $D_n^\ast(\forall,\exists,k)$ for $k^2 \leq n$}
A trivially 1-intersecting uniform partition system is a 
$\mathcal{D}_n^\ast(\forall,\exists,k)$ with cardinality
$U(n-c,k-1)$. Thus
\[
D_n^\ast(\forall,\exists,k) \geq U(n-c,k-1).
\]

\begin{conj} For positive integers $c,k,n$ with $n=ck$ and $k \leq c$,
\[
 D_n^\ast(\forall,\exists,k) = U(n-c,k-1).
\]
\end{conj}

\subsubsection{Bounds for $I_n^\ast(\exists, \forall,k)$ for $k^2 \leq n$}

In any $\mathcal{I}_n^\ast(\exists, \forall,k)$ the partitions in the system must be
disjoint; as such, no class can be repeated. Thus,
\[
I_n^\ast(\exists, \forall,k) \leq {n-1 \choose c-1}.
\]

%\begin{quest}
%Is there a construction for a partition system that is intersecting of type $(\exists, \forall)$
%with cardinality ${n-1 \choose c-1}$?
%\end{quest}

\begin{conj} For positive integers $c,k,n$ with $n=ck$ and $k \leq c$,
\[
I_n^\ast(\exists, \forall,k) = {n-1 \choose c-1}.
\]
\end{conj}

Ahlswede, Cai and Zhang~\cite{MR1739309} conjecture that the
largest system of $c$-uniform clouds which is intersecting with type
$(\exists, \forall)$ has cardinality ${n-1 \choose
c-1}$. In~\cite{intersectingsystems} this conjecture is shown to be
false; they give a system of clouds that is larger than this bound for
all $c\geq 8$. The counter example given in~\cite{intersectingsystems}
requires that each cloud in the system contain intersecting sets.

\subsubsection{Bounds for $D_n^\ast(\exists,\exists,k)$ and $I_n^\ast(\forall,\forall,k)$ where $k\leq c$}

The next pair to consider is $I_n^\ast(\forall, \forall,k)$ and
$D_n^\ast(\exists, \exists, k)$.  Two partitions are intersecting of
type $(\forall, \forall)$ if and only if they are qualitatively
independent. Thus, the graph $G_{\mathcal{I}_n^\ast(\forall, \forall,k)}$ is
exactly the graph $UQI(n,k)$.  The value of $I_n^\ast(\forall,
\forall,k)$ is the size of the maximum clique in the uniform
qualitative independence graph.  This is not an easy question;
perhaps, trying to find the value of $D_n^\ast(\exists, \exists, k)$
might be easier or at least provide a new approach to tackle this question.

{}From the proof of Theorem~\ref{thm:generalvtbound}, with $n=ck$, we
have the following two bounds:
\begin{eqnarray*}
D_n^\ast(\exists, \exists, k) &\geq& {n-c+k-2 \choose k-2} U(n-c,k-1),
\\ I_n^\ast(\forall, \forall, k) &\leq& \frac{\frac{1}{ k!}{n
\choose c}} {\frac{1}{(k-1)!}{ n-(c-(k-2)) \choose k-2 }} = \frac{
{n-1 \choose c -1} }{ {n-(c-k+2) \choose k-2} }.
\end{eqnarray*}

For $k=2$, from Theorem~\ref{thm:katona73},
$I_n^\ast(\forall,\forall,2) = {n-1 \choose \lfloor n/2 \rfloor -1}$
and $D_n^\ast(\exists,\exists,2) = 1$.

In Section~\ref{sec:blockrecursive}, for $k$ a prime power, the
block-size recursive construction can be used to construct a uniform
covering array $CA( (i+1)k^2-ik, k^i(k+1), k)$. So, for $k$ a prime
power,
\[
\frac{\log_2 k}{k(k-1)} \leq \limsup_{n \rightarrow \infty} \frac{\log_2(I_n^\ast(\forall,\forall,k))}{n}. 
\]

It is not clear if this is the exact asymptotic growth of
$I_n^\ast(\forall,\forall,k)$.  From Equation~(\ref{eq:asymptotic}),
Section~\ref{asymptotics}, we also have the following upper bound from
a cardinality of sets of qualitatively independent partitions (not necessarily uniform)
\[
\limsup_{n \rightarrow \infty} \frac{\log_2(I_n^\ast(\forall,\forall,k))}{n} \leq \frac{2}{k}.
\]

The exact value of $\limsup_{n \rightarrow \infty}
\frac{\log_2(I_n^\ast(\forall,\forall,k))}{n}$ is an open question.
If it is strictly less then $\frac{2}{k}$ then the size of a maximum
clique in $UQI(n,k)$ would be strictly smaller than the size of a
maximum clique in $QI(n,k)$. In particular $UQI(n,k)$ would not be a
core of $QI(n,k)$.  This gives us a new way to consider
Question~\ref{quest:generalcores}. In particular, a negative answer to
the question below would yield a negative answer to
Question~\ref{quest:generalcores}.
\begin{quest}
Is it true that 
\[
\limsup_{n \rightarrow \infty} \frac{\log_2(I_n^\ast(\forall,\forall,k))}{n}  = \frac{2}{k} \; \; ? 
\]
\end{quest}

\section{Non-Uniform Partitions}

Throughout this section, $n,k$ and $c$ will be positive integers with
$n=ck+r$ for some integer $r$ with $0 \leq r <k$.  We will consider all partition systems
from $\mathcal{P}^k_n$ and not just the uniform ones. We summarize
our results in the table below. The entry ``??'' denotes that we have
no bounds and ``asymptotics'' indicates that we know the asymptotic growth of
$\limsup_{n \rightarrow \infty} \frac{\log_2(A_n(X,Y,k))}{n}$.
 
\begin{center}
\begin{tabular}{|l||c|c|c|c|}
\hline
$A \big\backslash (X,Y)$ & $(\forall, \forall)$ & $(\exists, \forall)$ &
$(\forall, \exists)$ & $(\exists, \exists)$ \\
\hline \hline
\small{Incomparable} & ${n-1 \choose c-1}$ if $n=ck$ & ??  & $S(n,k)$ & $S(n,k)$\\ \cline{2-2}
              &  asymptotics & & & \\  \hline\hline
\small{Comparable}   & 1 & 1 & ?? & $\geq S(n-1,k-1)$\\ \hline\hline
\small{Disjoint}     & 1 & 1 & $\geq S(n-1,k-1)$ & ??\\ \hline\hline
\small{Intersecting} &  $\frac{1}{2}{n \choose \lfloor n/2 \rfloor}$ for $k=2$ & ?? & $S(n,k)$ & $S(n,k)$\\ \cline{2-2}
              & asymptotics &  & & \\
\hline
\end{tabular}
\end{center}

We do not have bounds analogous to the ones in
Proposition~\ref{vtbound} for the non-uniform case. For the symmetric
relations, $(\forall,\forall)$ and $(\exists, \exists)$, graphs
$G_{\mathcal{A}_n(X,X, k)}$ can be defined as before, but since they are 
not vertex transitive, we do not get bounds from the converse relations.

\remove{ it is possible to define graphs
$G_{\mathcal{A}_n(X,X, k)}$ that have vertex set $\mathcal{P}^n_k$ and
two partitions are adjacent if and only if they satisfy the relation
$A$ with type $(X,X)$ (where $A \in \{ N,C,D,I\}$ and the value of $X$
can be either $\forall$ or $\exists$).  The value of $A_n(X,X,k)$
would be equal to the size of the maximum clique in
$G_{\mathcal{A}_n(X,X, k)}$, but, since this graph is not
vertex transitive we can not use Corollary~\ref{cor:vtbound} to get
any bounds.
}

\subsection{Incomparable}
It is trivial to see that $N_n(\forall,\exists, k) = S(n,k)$ and
$N_n(\exists,\exists, k) = S(n,k)$.

Two partitions are incomparable of type $(\forall,\forall)$ if and
only if they have the Sperner property for partitions
(Section~\ref{sec:spernerpartitions}). Thus,
Theorem~\ref{thm:spernerpartition} can be restated as follows.
\begin{thm}
If $k$ divides $n$, then $N_n(\forall,\forall, k) = {n-1 \choose c-1}$.
\end{thm}
From Theorem~\ref{thm:AUbound} we have a bound for general $N_n(\forall,\forall, k)$.
\begin{thm}
If $n = ck+r$ with $0 \leq r< k$, then $N_n(\forall,\forall, k) \leq
\frac{1}{(k-r)+\frac{r(c+1)}{n-c}}{n \choose c}$.
\end{thm}
Finally, from Theorem~\ref{thm:spernerasymp} the asymptotic growth of
$N_n(\forall,\forall, k)$ is:
\[
\limsup_{n \rightarrow \infty} \frac{\log N_n(\forall,\forall, k) }{n} = \log\left( \frac{k}{k-1}\right)+\frac{\log(k-1)}{k}.
\]

%$N_n(\exists,\forall, k)\geq U(n-1,k-1)$ 

\subsection{Comparable}
For the relation comparable we have only the trivial results:

\begin{thm}\label{thm:comparable} For positive integers $k,n$, \hfill
\begin{enumerate}
\item{$C_n(\forall, \forall, k) = 1$,}
\item{$C_n(\exists, \forall, k) = 1$,}
\item{$C_n(\exists, \exists, k) \geq S(n-1,k-1)$}.
\end{enumerate}
\end{thm}

\subsection{Disjoint}

We start by stating the trivial results for disjoint partition systems without proof.

\begin{thm}\label{thm:disjoint} For positive integers $k,n$,\hfill
\begin{enumerate}
\item{$D_n(\forall, \forall, k) = 1$,}
\item{$D_n(\exists, \forall, k) = 1$,}
\item{if $n < k^2$ then $D_n(\exists, \exists, k) = S(n,k)$.}
\end{enumerate}
\end{thm}

The system of all partitions that contain the set $\{1\}$ as a class
is a disjoint system of type $(\forall, \exists)$. Thus,
\[
D_n(\forall, \exists, k) \geq S(n-1,k-1).
\]

Finally, for $k^2 \leq n$ the system of all partitions from
$\mathcal{P}^n_k$ with at least one class of size strictly less than
$k$ is a system of disjoint type $(\exists, \exists)$.  
It is not clear how to count the number of partitions in such a system.

\subsection{Intersecting}

It is straightforward that $I_n(\forall, \exists, k) = S(n,k)$ and
$I_n(\exists, \exists, k) = S(n,k)$.

Recall from Section~\ref{CA}, that for positive integers $n$ and $k$,
$N(n,k)$ is the largest $r$ for which a $CA(n,r,k)$ exists. So, by
definition, $I_n(\forall,\forall,k) = N(n,k)$. Finding an exact
solution for $I_n(\forall,\forall,k)$ is very hard, since it is
equivalent to finding the maximum set of qualitatively independent
partitions, which is also equivalent to determining the minimum size of
covering arrays.
 
We know the exact result for $k=2$, from Theorem~\ref{thm:katona73},
$I_n(\forall,\forall,2) = {n-1 \choose \lfloor n/2 \rfloor-1}$.
Poljak and Tuza~\cite{poljak:89} prove (using Bollab\'as's Theorem)
that $N(n,k) \leq \frac{1}{2}{ \lfloor{2n/k}\rfloor \choose
\lfloor{n/k}\rfloor}$ (see Section~\ref{asympbollobas}).
Thus, 
\[
I_n(\forall,\forall,k) \leq  \frac{1}{2}{ \lfloor{2n/k}\rfloor \choose \lfloor{n/k}\rfloor}.
\]
For general $k$, Gargano, K\"orner and Vaccaro give an asymptotic
result (Equation~(\ref{eq:asymptotic}), Section~\ref{asymptotics})
\[
\limsup_{n \rightarrow \infty} \frac{log_2(I_n(\forall,\forall,k))}{n} = \frac{2}{k}.
\]

Finally, there are many more constructions (see
Section~\ref{CAconstructions}) and results from heuristic searches
that give bounds on $N(n,k)$, and hence $I_n(\forall,\forall,k)$, for
specific values of $n$ and $k$~\cite{Colbourn3, Colbourn2, Colbourn1,
brett:97a, brett:97b}.

%\chapter{Conclusion}?

%%%%%%%%%%%%%%%%%%%%%%%%%%%%%%%%%%%%%%%%%%%%%%%%%%%%%%%%%%%% Open Questions

\chapter{Conclusion}\label{Questions}
\thispagestyle{empty}

In this thesis, we introduce the qualitative independence graphs and
the uniform and almost-uniform qualitative independence graphs. We
also introduce extremal partition theory, in particular we prove
versions of Sperner's Theorem and the Erd\H{o}s-Ko-Rado Theorem for
partitions. These are all interrelated and have applications to both
covering arrays and covering arrays on graphs.

The qualitative independence graphs are particularly useful for
studying covering arrays.  For a graph $G$ and integers $n,k$, a
$CA(n,G,k)$ exists if and only if there is a homomorphism from $G$ to
$QI(n,k)$.  In addition, an $r$-clique in the qualitative independence
graph $QI(n,k)$ corresponds to a covering array $CA(n,r,k)$ and an
independent set in $QI(n,k)$ is related to a partially intersecting
partition system in $\mathcal{P}^n_k$.

We give several results on binary covering arrays on graphs.  In
particular, we prove that a core of $QI(n,2)$ is $UQI(n,2)$ if $n$ is
even, and $AUQI(n,2)$ if $n$ is odd. This raises the question: does
there exist a homomorphism from $QI(n,k)$ to $UQI(n,k)$ or to
$AUQI(n,k)$ for all integers $n$ and $k$? If there is such a
homomorphism, bounds on the size of the maximum clique in $UQI(n,k)$
and $AUQI(n,k)$ would provide many new bounds on the size of covering
arrays. Moreover, we would know that there exists an optimal covering
array whose rows correspond to uniform or almost-uniform partitions
rather than to general partitions.  Since we know that the asymptotic
growth of the maximum cardinality of a set of qualitatively
independent $k$-partitions is $\frac{2}{k}$ (see
Section~\ref{asympbollobas}), it is interesting to consider the same
asymptotic growth for a set of qualitatively independent uniform
$k$-partitions.  In particular, if the limit $\limsup_{n \rightarrow
\infty}
\frac{\log_2(I_n^\ast(\forall,\forall,k))}{n}$ is strictly smaller
than $\frac{2}{k}$, then we would know that there is no homomorphism
from $QI(n,k)$ to $UQI(n,k)$.

%spectral stuff
The uniform qualitative independence graphs $UQI(ck,k)$, where $c,k$
are positive integers, are particularly interesting.  We exhibit an
equitable partition on the vertices of $UQI(ck,k)$ that can be used to
find the eigenvalues of these graphs; we give the spectra for many
small uniform qualitative independence graphs.  Mathon and
Rosa~\cite{MR86e:05068} give an association scheme that has $QI(9,3)$
as one of its graphs.  In Section~\ref{sec:otherschemes}, we describe
sets of graphs that are an extension of this scheme. It seems likely
that these sets of graphs form either an association scheme or some
generalization of an association scheme (perhaps an asymmetric
association scheme). This is very interesting, since association
schemes are difficult to construct and can give more information about
their graphs. These results open a new and exciting direction for
covering array research.
 
\remove{
We show that the technique of translating designs into graphs is useful.
We have results for 
(uniform) qualitative independent sets of partition $QI(n,k)$ $(UQI(n,k))$
intersecting partitions systems from $\mathcal{P}^n_k$  $IP(n,k)$
$G_{\mathcal{A}^\ast(X,X,k)}$
Moreover, in a surprising number of cases the
bound from Inequality~(\ref{eq:vt2}) in Section~\ref{knesergraphs} holds
with equality. Is there a class of designs for which this bound holds?
}

%extremal partition systems
Sperner's Theorem for set systems is essential in our proof that a
core of $QI(n,2)$ is $UQI(n,2)$ if $n$ is even, and $AUQI(n,2)$ if
$n$ is odd.  Motivated by this, we give an extension of Sperner's
Theorem to partition systems.  We prove that, for integers $n$ and $k$
where $k$ divides $n$, the largest Sperner partition system in
$\mathcal{P}^{n}_k$ is a uniform partition system.  We also
conjecture, for all values of $n$ and $k$, that the largest Sperner
partition system in $\mathcal{P}^{n}_k$ is an almost-uniform partition
system.  Unlike the case for $k=2$, this extension of Sperner's
Theorem to partition systems is not enough to prove the existence of a
homomorphism from $QI(n,k)$ to $UQI(n,k)$. It would be interesting to
know if the largest qualitatively independent partition system is a
uniform (or an almost-uniform) partition system. This would be true if
there is a homomorphism from $QI(n,k)$ to $UQI(n,k)$ (or to
$AUQI(n,k)$).

There are two ways to extend the Erd\H{o}s-Ko-Rado Theorem to
partition systems. For the first extension, we define two partitions
to be intersecting if they have a class in common. We prove that the
largest uniform partition system with this type of intersection is a
trivially intersecting system. For the second extension, we consider a
different type of intersection called partial intersection. This type
of intersection is related to independent sets in qualitative
independence graphs. We conjecture that the largest partially
$t$-intersecting partition system is a trivially partially
$t$-intersecting partition system.  We can prove this conjecture for
partially $t$-intersecting partition systems from $\mathcal{U}^n_k$
for many specific values of $t,n$ and $k$. For example, a partially
2-intersecting partition system in $\mathcal{U}^{k^2}_k$ is equivalent
to an independent set in the graph $QI(k^2,k)$. Using the eigenvalues
of $QI(k^2,k)$, we prove that a trivially partially 2-intersecting
partition system is the largest partially 2-intersecting partition
system in $\mathcal{U}^{k^2}_k$.

There are many open problems related to the extension of Sperner's Theorem and
the Erd\H{o}s-Ko-Rado Theorem to partition systems.  First, it would be nice
to find a tighter connection between Sperner's Theorem for partition
systems, the Erd\H{o}s-Ko-Rado Theorem for partition systems and
covering arrays (see Question~\ref{quest1}).  Second, it would be
interesting to know the exact bound for $n$ in the Erd\H{o}s-Ko-Rado Theorem
for partition systems (see Theorem~\ref{one}) and to find an extension of the
complete Erd\H{o}s-Ko-Rado Theorem for partition systems (see
Conjecture~\ref{conj:completeintersection}).  Wilson's
proof~\cite{MR86f:05007} of the exact bound for $n$ in the
Erd\H{o}s-Ko-Rado Theorem used the ratio bound,
Lemma~\ref{lem:smallevbound}, with the eigenvalues of the Johnson
scheme (Example~\ref{exmp:johsonscheme},
Section~\ref{sec:introassscheme}). Is it possible to prove a complete
Erd\H{o}s-Ko-Rado Theorem for intersecting partition systems (or
partially intersecting partitions) in a similar manner?  Finally, it
would be better to have a version of the Erd\H{o}s-Ko-Rado Theorem for
all partially intersecting partition systems, rather than the
collection of specific cases (see Corollaries~\ref{PIPS:app1},
\ref{PIPS:app2}, and \ref{PIPS:app3}). 

%higher order extremal partition systems
These extensions of Sperner's Theorem and the Erd\H{o}s-Ko-Rado
Theorem to partition systems are a starting point for the study of
extremal partition theory. Throughout the thesis, we give motivation for why these problems
are interesting and in Chapter~\ref{higherparts}, we present a framework in which to consider a
variety of extremal problems for partition systems.  For several of
these problems, we give exact results and for others, we give a bound on the
cardinality of the maximum system. There are still many
problems to solve and other possible extensions to be considered.
We conclude with a list of questions followed by several conjectures.

\section{Questions}\label{sec:quest}

\begin{quest}
In Section~\ref{sec:groupconst}, the group construction for covering
arrays is described. This construction builds a $CA( r(k-1) + 1,r,k)$
for many values of $k$ and $n$, and for $k+1 \leq r \leq 2k$ this
construction often gives a good upper bound for $CAN(r,k)$. Can this
construction be extended or generalized to give a good upper bound for
$CAN(r,k)$ for $r \geq 2k$?  This construction uses a ``starter
vector'' which is currently found either by an exhaustive search or a
heuristic search. Is it possible to construct the starter vectors
directly?
\end{quest}

\begin{quest}
Gargano, K\"orner and Vaccaro~\cite{gargano:94} give the asymptotic
growth of the maximum size of a set of qualitative independent
partitions.  Is it possible to find the asymptotic growth of the
maximum size of a set of qualitative independent uniform partitions?
In particular, if $N'(n,k)$ denotes the largest integer $r$ such that
a balanced $CA(n,r,k)$ exists, then what is the value of
\[
\limsup_{n \rightarrow \infty} \frac{\log_2 N'(n,k)}{n}?
\]
\end{quest}

%qi graphs
\begin{quest}
Let $n$ and $k$ be positive integers.
What are the cores of qualitative independence graphs $QI(n,k)$ for
$k>2$? A more daring version of this question is: for
positive integers $c$ and $k$, is the graph $UQI(ck,k)$ a core of
$QI(ck,k)$? And, for positive integers $n$ and $k$, is the graph
$AUQI(n,k)$ a core of $QI(n,k)$?
\end{quest}

\begin{quest}
What are $\alpha(QI(n,k))$, $\chi(QI(n,k))$ and $\chi^\ast(QI(n,k))$
for $k>2$?
\end{quest}

\begin{quest}
The natural extension of strength-2 covering arrays on graphs to
higher strengths are strength-$t$ covering arrays on $t$-uniform
hypergraphs. Many of the basic results from Sections~\ref{sec:homomor}
and \ref{sec:QIgraphs} can be extended to strength-$t$ covering arrays
on $t$-uniform hypergraphs. What can be proven for higher strength
covering arrays on hypergraphs with a binary alphabet?
\end{quest}

\begin{quest}
We can define a graph, similar to the qualitative independence graphs,
that correspond to orthogonal arrays. Let $n,k,\lambda$ be positive
integers with $n = \lambda k^2$. Define a graph $O(n,k)$ to have
vertex set $\mathcal{U}^n_k$. Vertices $P,Q \in \mathcal{U}^n_k$
are adjacent in $O(n,k)$ if and only if for all $P_i \in P$ and $Q_j
\in Q$, where $i,j \in \{1,\dots, k\}$, $P_i \cap Q_j = \lambda$.  Can
we find results for the graph $O(n,k)$ similar to the ones for the
graph $UQI(n,k)$ given in Chapter~\ref{chp:ass}? Moreover, can we use
the graph $O(n,k)$ to find bounds on the size of $r$ in $OA(n,r,k,2)$?
\end{quest}

%agt
\begin{quest}
In Example~\ref{exmp:kneserevalues}, Section~\ref{sec:ep}, the
eigenvalues of the Kneser graphs are given. These eigenvalues are
found by partitioning the vertices of the Kneser graph into an
equitable partition with a singleton class.  The adjacency matrix for
the quotient graph of the Kneser graph, with respect to this equitable
partition, has the same eigenvalues as the Kneser graph.  In
Section~\ref{sec:betterep}, an equitable partition on the vertices of
$UQI(n,k)$ with a singleton class is given. This partition is used to
reduce the calculations to find the eigenvalues of some uniform
qualitative independence graphs. Can this equitable partition be used
to find a formula for all the eigenvalues for all the uniform
qualitative independence graphs?
\end{quest}

\begin{quest}
What are the eigenvalues of the graph $QI(k^2+1,k)$?  What bound on
$\omega(QI(k^2+1,k))$ do we get from the ratio bound for maximum
cliques?
\end{quest}

\begin{quest}
Can the equitable partition from Section~\ref{sec:betterep} be
generalized to the vertices of the graphs $AUQI(n,k)$ for any positive
integers $n$ and $k$? What are the eigenvalues for $AUQI(n,k)$?
\end{quest}

\begin{quest}
Do the graphs in Section~\ref{sec:U12-3} and \ref{sec:U15-3} describe
an association scheme for $\mathcal{U}^n_3$ where $n=12,15$?  Do the
graphs given in Section~\ref{sec:otherschemes} describe an asymmetric
association scheme on $\mathcal{U}^n_k$ for all $n$ and $k$?
\end{quest}

\begin{quest}\label{quest:IP}
In the proof of Theorem~\ref{PIPS:thm}, for positive integers $n,k$,
where $k$ divides $n$, a graph $PIP(n,k)$ is defined. The vertex set
of this graph is $\mathcal{U}^n_k$ and two partitions are adjacent if
and only if they are partially 2-intersecting.  Similarly, define a
(vertex-transitive) graph $IP(n,k)$ with vertex set all uniform
$k$-partitions of an $n$-set and two partitions are adjacent if and
only if they are intersecting. The graph $IP(n,k)$ is also the graph
$G_{\mathcal{I}_n^\ast(\forall,\forall,k)}$ described in
Section~\ref{sec:DIhigher}.  A clique in $IP(n,k)$ is an intersecting
uniform $k$-partition system.  The bound from Theorem~\ref{two} can be
found with the formula $\omega(IP(n,k)) \leq \frac{|V(IP(n,k))|}
{\alpha(IP(n,k))}$ and the fact that $\alpha(IP(n,k)) = {n-1
\choose \frac{n}{k}-1}$.  Theorem~\ref{two} also states that only
trivially intersecting partition systems meet this bound. Is there a
way to use the graph $IP(n,k)$ to prove this fact?
\end{quest}

\begin{quest}
Let $n,k$ be positive integers. For $i=0,1,\dots, k$, define a graph
$G_i$ on vertex set $\mathcal{U}^n_k$ with partitions $P,Q \in
\mathcal{U}^n_k$ adjacent if and only if $P$ and $Q$ have exactly $i$
classes in common. What are the eigenvalues of this graph?
\end{quest}

\begin{quest}
Why is the ratio bound for independent sets (Lemma~\ref{lem:smallevbound}) tight for the graph
$QI(k^2,k)$ for all positive integers $k$?  What other types of
designs can be converted into independent sets or cliques in graphs?
For what type of design is such a graph regular and when is the ratio
bound for independent sets tight for such graphs?
\end{quest}

\begin{quest}\label{quest1}
If a partition system is qualitatively independent, then it is a
Sperner partition system.  Does this relation have a converse?  For
example, if a 2-partition system is both an intersecting and a Sperner
partition system, then it is qualitatively independent.  What extra
conditions would a general Sperner partition system need to meet to
also be a qualitatively independent partition system?
\end{quest}

\begin{quest}
Are there other applications of Sperner partition systems?
\end{quest}

\begin{quest}
In Theorem~\ref{one}, we prove that a $t$-intersecting uniform
partition system has cardinality no more than $U(n-tc,k-t)$ for $n$
sufficiently large. What is the exact lower bound for $n$?  In
particular, can we set up a graph whose vertex set is $\mathcal{U}^n_k$
and vertices are adjacent if and only if the partitions are
$t$-intersecting, and use eigenvalues of this graph to find the exact
lower bound?
\end{quest}

\begin{quest}
Recall the graphs $G_{\mathcal{A}^\ast_n(X,X,k)}$ (where $A \in \{ N,C,D,I\}$
and $X$ can be either $\forall$ or $\exists$) defined in
Section~\ref{graph}.  Can we use an equitable partition, similar to
the one defined in Section~\ref{sec:betterep}, to find the eigenvalues
of these graphs?
\end{quest}

\section{Conjectures}
\setcounter{thm}{0}
\begin{conj}
For all $n \geq k^2$, $CAN(QI(n,k),k) = n$. This is equivalent to there being no
homomorphism
\[
QI(n,k) \rightarrow QI(n-1,k).
\]
This was given in Conjecture~\ref{conj:canQIn} and Conjecture~\ref{conj:nohomo}.
\end{conj}

\begin{conj}
From Theorem~\ref{biggchrom}, we have for all integers $k$, $\chi(
QI(k^2,k) ) \leq {k+1 \choose 2}$; we conjecture that for all integers $k$,
$\chi( QI(k^2,k) ) = {k+1 \choose 2}$.
\end{conj}

\remove{\begin{conj}
For all $k$, for any $n \geq k+2$, the graph $AUQI(n,k)$ is not
arc transitive.
\end{conj}}

\begin{conj}\label{conj:alphadim}
Let $k$ be a positive integer. Consider the graphs $QI(k^2,k)$. For
all distinct $i,j \in \{1,\dots,k^2\}$, let $S_{ \{i,j\} }$ denote the
set of all partitions in $V(QI(k^2,k))$ with $i$ and $j$ in the same
class. Then every maximum independent set in $QI(k^2,k)$ is a set $S_{
\{i,j\} }$ for some distinct $i,j \in \{1,\dots,k^2\}$.

In addition, for all $n$, the set of characteristic vectors of $S_{ \{i,j\} }$ in
$UQI(n,k)$ for all distinct $i,j \in \{1,\dots,n\}$ spans a vector
subspace of $\mathbb{R}^{U(n, k)}$ with dimension ${n \choose 2}-{n
\choose 1}$.  Note that in~\cite{mikesthesis} it is shown that this is
true for $n=k^2$.
\end{conj}

\begin{conj}
For all $k$, the graph $QI(k^2,k)$ is a core. We further conjecture
that this can be proved using Conjecture~\ref{conj:alphadim} above.
\end{conj}

\begin{conj}
Let $n,k,c,r$ be positive integers with $n=ck+r$ and $0 < r <k$.  The
largest Sperner partition system in $\mathcal{P}^n_k$ is an
almost-uniform Sperner partition system.
\end{conj}

\begin{conj}
For positive integers $n,k,c,t$ with $n=ck$ and $t \leq c$, let $\mathcal{P}$ be a
partially $t$-intersecting uniform $k$-partition system on an
$n$-set. Then
\[
|\mathcal{P}| \leq {n-t \choose c-t}\frac{1}{(k-1)!} {n-c \choose
c}{n-2c \choose c} \cdots{c \choose c}.
\]
Moreover, equality holds if and only if $\mathcal{P}$ is a trivially
partially $t$-intersecting partition system.

Note that Corollaries~\ref{PIPS:app1}, \ref{PIPS:app2},
\ref{PIPS:app3} and \ref{PIPS:app4} prove the bound in this
conjecture for several specific values of $t,k$ and $n$.
\end{conj}

\appendix

\chapter{Tables of Bounds for $CAN(\MakeLowercase{r},\MakeLowercase{k})$}
\label{appendix}
\thispagestyle{empty}

Tables~\ref{results3:tab} and \ref{results4:tab} gives a list of starter vectors that
improve the previously best known upper bounds for covering arrays. 

\begin{table}
\caption{New upper bounds and corresponding starter vectors for $k \leq 12$. \label{results3:tab}}
{\small
\begin{tabular} {clcc} \hline
$k=6$ & starter vector & new bound & old bound\\  \hline
$r = 9$ & 0 1 1 2 1 1 3 5 3 & 46 & 48\\ 
$r = 10$ & 0 1 1 1 1 2 4 3 1 2 & 51 & 52\\
\hline
$k=7$ & \\  \hline
$r = 10$ & 0 1 1 1 3 4 1 3 2 6 & 61 & 63 \\
$r = 11$ & 0 1 1 1 1 2 1 4 6 5 3 & 67 & 73\\
$r = 12$ & 0 1 1 1 1 1 2 1 4 6 5 3 & 73 & 76\\
\hline
$k=8$ & \\  \hline
$r = 11$ & 0 1 1 2 2 4 2 5 6 3 2         & 78 & 80\\
$r = 12$ & 0 1 1 1 2 6 2 6 1 1 6 5       & 85 & 99\\
$r = 13$ & 0 1 1 1 1 2 1 3 7 5 1 3 4     & 92 & 102\\
$r = 14$ & 0 1 1 1 1 1 2 1 3 7 5 1 3 4   & 99 & 104\\
%$r = 15$ & 0 1 1 1 1 1 1 2 1 3 7 5 1 3 4 & 106 & 107\\
\hline
$k=9$ \\  \hline
$r = 13$ & 0 1 1 1 3 2 1 6 2 5 5 3 4           & 105 & 120\\
$r = 14$ & 0 1 1 1 1 2 8 5 6 2 1 3 6 7         & 113 & 131\\
$r = 15$ & 0 1 1 1 1 1 2 3 2 7 1 5 4 2 5       & 121 & 135\\
$r = 16$ & 0 1 1 1 1 1 1 2 3 2 7 1 5 4 2 5     & 129 & 145\\
$r = 17$ & 0 1 1 1 1 1 1 1 2 1 2 6 8 5 3 6 7   & 137 & 148\\
$r = 18$ & 0 1 1 1 1 1 1 1 1 2 1 2 6 8 5 3 6 7 & 145 & 151\\
\hline
$k=10$ \\  \hline
$r = 15$ & 0 1 1 1 1 4 2 8 1 9 4 5 6 8 4             & 136 & 166\\
$r = 16$ & 0 1 1 1 1 1 2 9 5 7 1 5 5 3 2 8           & 145 &177\\
$r = 17$ & 0 1 1 1 1 1 1 2 4 9 5 8 7 2 1 8 5         & 154 & 180\\
$r = 18$ & 0 1 1 1 1 1 1 1 2 4 9 5 8 7 2 1 8 5       & 163 & 180\\
$r = 19$ & 0 1 1 1 1 1 1 1 1 2 1 4 9 7 8 2 6 8 5     & 172 & 180\\
%$r = 20$ & 0 1 1 1 1 1 1 1 1 1 2 1 4 9 7 8 2 6 8 5   & 181 & 180\\
$r = 21$ & 0 1 1 1 1 1 1 1 1 1 1 1 2 5 6 8 4 8 6 5 2 & 190 & 202\\
$r = 22$ & 0 1 1 1 1 1 1 1 1 1 1 1 1 2 5 6 8 4 8 6 5 2  & 199 & 202\\
\hline
$k=11$ \\  \hline
$r = 17$ & 0 1 1 1 1 1 3 10 6 9 4 10 9 7 1 4 5 & 171 & 221 \\ 
$r = 18$ & 0 1 1 1 1 1 1 2 5 8 3 1 3 10 6 5 8 2 & 181 & 225 \\
$r = 19$ & 0 1 1 1 1 1 1 1 2 5 8 3 1 3 10 6 5 8 2   & 191 & 231 \\
$r = 20$ & 0 1 1 1 1 1 1 1 1 2 5 8 3 1 3 10 6 5 8 2  & 201 & 231 \\
$r = 21$ & 0 1 1 1 1 1 1 1 1 1 2 5 8 3 1 3 10 6 5 8 2   & 211 & 231 \\
$r = 22$ & 0 1 1 1 1 1 1 1 1 1 1 2 5 8 3 1 3 10 6 5 8 2  & 221 & 231 \\
%$r = 23$ & 0 1 1 1 1 1 1 1 1 1 1 1 2 5 8 3 1 3 10 6 5 8 2  & 231 & 231 \\
\hline
$k=12$ \\  \hline
$r =18 $ & 0 6 5 10 7 3 9 1 3 2 3 7 6 6 3 1 4 8  & 199 & 255\\ 
$r =19 $ & 0 8 3 3 6 7 10 9 7 10 5 9 5 7 4 2 10 9 3   & 210 & 276\\ 
$r = 20$ &0 7 9 2 10 10 2 9 8 9 4 9 2 3 1 4 5 6 1 9 & 221 & 276\\
$r = 21$ &0 5 1 9 1 4 2 7 9 4 7 6 5 5 2 3 7 7 2 7 9 & 232 & 276\\
$r = 22$ & 0 1 1 3 10 6 9 4 10 9 1 1 1 4 8 6 3 2 7 8 4 11     & 243 & 288\\
$r = 23$ & 0 1 1 1 3 10 6 9 4 10 9 1 1 1 1 9 2 11 5 11 1 6 5  & 254 & 288\\
$r = 24$ & 0 1 1 1 1 3 10 6 9 4 10 9 1 1 1 1 6 5 2 11 1 5 6 2 & 265 & 288\\
$r = 25$ & 0 7 10 7 8 2 2 1 1 10 1 5 5 5 1 7 10 1 6 1 6 4 9 9 10 & 276 & 288\\
%$r = 26$ &0 5 9 8 3 1 5 7 1 4 1 4 5 9 1 6 2 3 1 1 2 10 9 2 6 1 & 287 & 288\\
%\\ \hline
\end{tabular}
%\caption{Starter vectors improving known bounds \label{results2:tab} }
}
\end{table}

\begin{table}
\caption{New upper bounds and corresponding starter vectors for $13 \leq k \leq 18$. \label{results4:tab}}
{\small
\begin{tabular} {clcc} \hline
% & starter vector & new & old \\  
% &  & bound & bound\\  \hline
 $k=13$ \\ \hline
$r = 21$ & 0 6 9 5 7 6 2 8 8 9 6 2 3 10 5 3 5 9 8 1 8 & 253 & 325\\ 
$r = 22$ &0 3 9 11 4 1 8 7 3 7 1 3 3 4 1 1 4 2 6 11 6 4  & 265  & 325\\ 
$r = 23$ &0 6 7 10 2 5 4 11 8 2 4 1 9 3 5 3 3 1 10 10 11 4 6  & 277  & 325\\ 
$r = 24$ & 0 7 3 7 7 3 10 2 8 3 5 4 2 4 5 8 9 2 11 7 8 7 3 10 & 289  & 325\\ 
$r = 25$ & 0 1 2 4 10 11 5 1 10 2 7 6 4 3 6 8 6 7 10 3 10 1 2 8 8  & 301  & 325\\ 
$r = 26$ & 0 7 4 9 2 2 3 7 10 6 6 10 10 7 2 7 9 4 5 4 2 4 11 10 4 4  & 313  & 325\\ 
%$r = 27$ & 0 10 1 4 5 1 5 7 1 8 9 10 2 3 3 10 8 1 3 1 1 10 5 10 9 9 1  & 325 & 325\\ 
\hline
 $ k=14$ \\  \hline
$r = 23$ & 0 1 5 4 1 9 4 6 1 1 6 1 5 10 7 8 8 1 4 5 1 12 6 & 300 & 437\\ 
$r = 24$ &0 5 4 12 9 2 7 10 10 10 2 1 6 2 8 12 9 7 8 8 1 8 10 7  &313  & 437\\ 
$r = 25$ & 0 8 8 2 12 1 12 10 6 1 6 6 10 4 7 4 3 6 12 1 3 4 9 11 7 & 326 & 437\\ 
$r = 26$ & 0 7 2 9 10 7 6 3 11 1 10 2 6 8 3 5 10 5 5 6 12 5 3 2 6 5  & 339 & 437\\ 
$r = 27$ & 0 8 4 5 10 12 4 8 3 11 4 3 8 9 2 1 5 12 5 8 7 5 2 2 3 3 12 & 352 & 437\\ 
$r = 28$ &0 5 3 11 6 5 5 10 3 10 6 1 4 11 9 4 12 6 8 1 8 1 10 1 11 12 4 4  & 365 & 437\\ 
$r = 29$ &0 2 4 1 1 2 5 10 3 1 2 9 8 1 1 1 1 9 3 11 10 9 10 1 1 4 12 12 8  & 378 & 437\\ 
 \hline
$k=15$ \\   \hline
$r = 27$ & 0 11 8 9 11 1 7 12 4 4 11 9 4 2 6 12 10 6 7 6 6 5 1 8 11 4 12  & 379 & 450\\ 
$r = 28$ & 0 12 7 1 1 5 13 11 12 8 3 4 5 10 11 2 13 1 4 12 4 3 3 10 4 3 4 3 & 393 & 450\\ 
$r = 29$ & 0 2 13 8 12 8 1 13 1 8 11 3 8 2 9 8 13 5 7 7 13 10 11 6 11 6 5 1 12  & 407 &450 \\ 
$r = 30$ & 0 1 7 13 11 2 2 2 12 8 3 5 8 9 6 9 8 8 11 11 5 4 4 2 13 9 2 13 3 6 & 421  &450 \\ 
$r = 31$ &0 9 8 1 7 1 5 2 5 3 1 4 3 12 12 1 2 6 1 2 1 13 1 6 1 11 7 7 10 3 10  & 435 & 450\\ 
 \hline
$ k=16$ \\   \hline
$r = 29$ & 0 5 2 6 2 5 11 4 13 1 3 14 9 6 4 8 9 10 12 3 8 3 14 10 9 9 2 9 6 & 436  & 496\\ 
$r = 30$ &0 11 8 9 13 7 3 8 2 10 4 14 6 3 2 12 4 8 6 4 6 13 3 6 6 12 12 4 2 1  & 451 &496 \\ 
$r = 31$ & 0 9 5 5 13 2 10 14 8 1 12 4 9 14 6 1 4 6 10 6 6 7 1 14 4 7 13 7 6 3 5 & 466 & 496\\ 
$r = 32$ & {\footnotesize 0 7 7 3 7 5 1 7 12 3 5 10 11 11 8 10 1 13 6 2 10 11 5 12 5 8 7 13 14 11 9 4} & 481  & 496\\ 
%$r = 33$ & 0 5 1 5 4 1 8 13 10 14 2 3 9 3 5 1 4 10 5 7 8 2 14 5 14 10 11 1 14 11 8 1 1 & 496  & 496\\
 \hline
$ k=17$ \\  \hline
$r = 33$ & {\scriptsize 0 13 4 10 9 5 8 9 4 14 8 13 2 11 6 14 13 10 14 2 4 9 13 5 5 5 7 10 4 7 13 11 4} & 529 & 561\\ 
$r =34 $ & {\scriptsize 0 8 4 5 7 10 1 7 6 4 14 6 13 1 5 10 10 11 8 5 2 15 3 6 2 13 8 11 8 3 12 2 8 12} & 545 & 561\\ 
%$r=35$&0 13 4 4 3 6 14 11 8 4 13 11 11 11 14 15 7 7 10 4 13 15 9 15 13 7 12 12 7 11 13 12 3 2 10&561& 561\\ 
 \hline
$ k=18$ \\  \hline
$r=35$& {\scriptsize 0 2 13 12 8 4 5 12 9 12 4 8 8 15 10 10 16 5 5 4 6 13 15 13 1 11 8 13 14 5 13 10 7 10 8} & 596 & 666\\
$r=36$& {\scriptsize 0 3 16 12 4 5 4 1 15 15 16 11 13 6 9 16 3 2 10 8 14 1 14 8 13 15 15 14 10 4 1 10 2 6 9 4} & 613 & 666 \\
$r=37$& {\scriptsize 0 6 14 3 5 8 12 15 4 15 5 2 11 7 2 3 2 7 15 5 12 12 10 8 4 1 5 5 7 7 10 16 10 16 9 2 4} & 630 & 666\\ 
%$k=19$ \\  
%$r = 39$ & 0 13 8 13 8 1 13 17 13 2 7 15 9 11 4 15 16 7 11 11 7 5 2 4 6 12 15 17 8 17 10 17 16 4 13 14 2 11 3 & 703 & 703\\
\end{tabular}
}
\end{table}

\addcontentsline{toc}{chapter}{Bibliography}

\bibliographystyle{plain}
\bibliography{bibliography}%{publishers,journals,bibliography}

\addcontentsline{toc}{chapter}{Index}
\printindex

\end{document}